\documentclass[12pt,twoside]{amsbook}
\usepackage{etex}

\usepackage[dvips]{graphics}
\usepackage{exscale}
\usepackage{amssymb}
\usepackage{tikz}
\usepackage{epsfig}
\usepackage{pgf}
\usepackage{pgfplots}
\usepackage{wrapfig}
\usepackage{xcolor}
\usepackage{lipsum}
\usepackage{calrsfs}%para escrever os letras caligráficas mais arredondadas
\usepackage{colortbl}
\usepackage{bm}

\usepackage{stackrel}

\usepackage{filecontents}

\usepackage[latin1]{inputenc}
\usepackage[T1]{fontenc}
\usepackage{pstricks,pst-3dplot}
\usepackage{exscale}
\usepackage{psfrag,epsfig,wrapfig}
\usepackage{infix-RPN}
\usepackage{amssymb,latexsym,amsthm,amsfonts,framed,color,fancyhdr}
\usepackage{pst-plot,pst-infixplot,pst-node,pstricks,graphicx}
\usepackage{mathrsfs}
\usepackage{enumitem}
\usepackage{extarrows}
\usepackage{multicol}
\usetikzlibrary{shapes.geometric,fit}

\newrgbcolor{verde}{0.0 0.5 0.0}
\newrgbcolor{morado}{.5 0 .5}
\newrgbcolor{azul}{0 0 .5}
\newrgbcolor{rojo}{0.5 0 0}
\newrgbcolor{brown}{0.55 0.14 0.14}
\newrgbcolor{orange}{1.0 0.3 0.0}
\newrgbcolor{seagreen}{0.33 1.00 0.62}
\newrgbcolor{darkgreen}{0.00 0.39 0.00}
\newrgbcolor{azulnovo}{0.0 0.75 1.0}%azul turquesa

\definecolor{shadecolor}{rgb}{0.9,1,1} % RECOMENDADO
\numberwithin{equation}{chapter}

\newtheorem{theorem}{Theorem}[chapter]

\newtheorem{corollary}[theorem]{Corollary}
\newtheorem{proposition}[theorem]{Proposition}

\newtheorem{lemma}[theorem]{Lemma}
\newtheorem{remark}{Remark}[chapter]
\newtheorem{definition}{Definition}[chapter]

\usepackage{amsmath} % for \genfrac

\newmuskip\pFqskip
\mathchardef\pFcomma=\mathcode`, % keep a copy of the comma

\newcommand*\pFq[5]{%
  \begingroup
  \begingroup\lccode`~=`,
    \lowercase{\endgroup\def~}{\pFcomma\mkern\pFqskip}%
  \mathcode`,=\string"8000
  {}_{#1}F_{#2}\left(\genfrac..{0pt}{}{#3}{#4};#5\right)%
  \endgroup
}

\begin{document}

%%%%%%%%%%%%%%%% colocar % a partir daqui (para vir só o capitulo pretendido, dado pelo includeonly) %%%%%%%%%%%%%%%%%%

\thispagestyle{empty}

\pagenumbering{roman} \vspace*{5em}

\begin{center}

%{\LARGE{\bf Texts of}} \bigskip

%{\Large{\bf \blue ORTHOGONAL POLYNOMIALS} \bigskip
%
%{\bf \blue AND SPECIAL FUNCTIONS}}\vspace{2em}

%Este é bom {\blue\bf {\large ORTHOGONAL POLYNOMIALS} $\&$ {\large SPECIAL FUNCTIONS}}\vspace{2em}

{\blue\bf {\Large Orthogonal Polynomials and Special Functions}}\vspace{2em}

%{\LARGE{\bf  Interuniversitary PhD Program in Mathematics Universities of Coimbra and Porto}}
{\large{\bf UC|UP Joint PhD Program in Mathematics}}
\end{center}
\vspace{5em}

\begin{center}
{\bf Jos\'e Carlos Soares Petronilho}
\bigskip

%({\sl Professor Associado})
\end{center}

\vspace{10em}

\begin{center}
\begin{pspicture}(0,0)(0,0)
%\rput(5.5,-6)
\rput(0,0)
{\includegraphics[width=10em,scale=0.14]{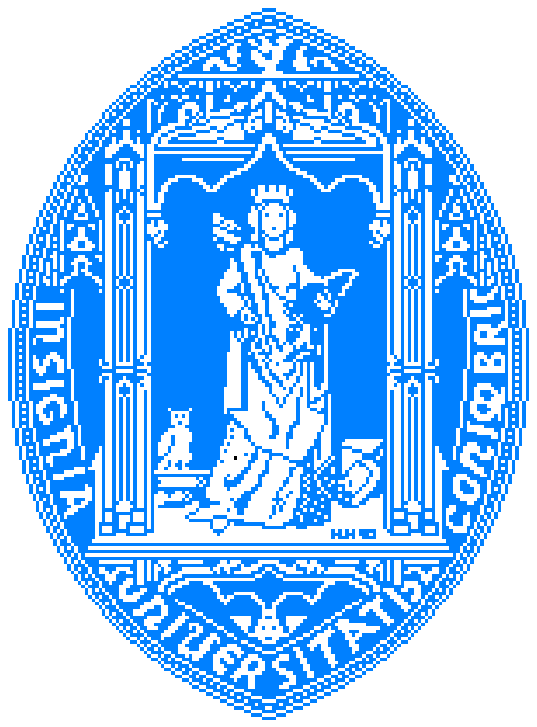}}
\end{pspicture}
\end{center}

\vfill
\begin{center}
{\bf Department of Mathematics }

{\bf Faculty of Sciences and Technology}

{\bf University of Coimbra}

{\bf 2016--2017}
\end{center}

%%%%%%%%%%%%% colocar % até aqui %%%%%%%%%%%%%%%%%%%%%%%%%%%%%%%%%%%%%%%%%%%%%%%%%%%%%%%%%%%%%%%%%%

%\newpage
%
%\vspace{10em} {\bf Nota pr\'evia}\vspace{1em}
%
%Textos de apoio \`as aulas da disciplina
%{\sl An\'alise Real e Funcional} no ano lectivo 2014--2015.
%
%\vspace{10em} {\bf Nota pr\'evia}\vspace{1em}
%
%O texto que se apresenta serviu de apoio \`as aulas da disciplina
%de {\sl An\'alise Funcional Aplicada} leccionada
%pelo autor nos tr\^es \'ultimos anos. Deve ter-se em mente que o
%texto n\~ao abarca a totalidade dos t\'opicos estudados, assim
%como inclui outros que n\~ao foram abordados, e que a ordem de
%apresenta\c c\~ao n\~ao corresponde, em muitas situa\c c\~oes,
%\`aquela em que a mat\'eria foi, efectivamente, leccionada.
%Trata-se, pois, apenas de umas notas ainda em constru\c c\~ao.

\newpage

\pagestyle{empty}
\vspace*{5em}

%The present text contains the
These notes contain part of the lectures of an introductory course on
orthogonal polynomials and special functions that I gave
in the joint PhD Program in Mathematics UC|UP in the
academic years 2015-2016 (at University of Porto) and 2016-2017 (at University of Coimbra).
\medskip

The notes were written for students who have never contacted with the above topics.
%Thus detailed proofs of most results presented were included.
Most results presented here can be found in the available bibliography at the end of each text/chapter,
although in general more detailed proofs have been included %(``the devil is in the details''),
(a few of them different from the ones presented in the source references),
hoping this helps the beginner student.
Besides the topics contained in the notes, several other subjects were covered in the course, including an introduction to
discrete orthogonal polynomials, orthogonal polynomials on the unit circle, spectral theory of Jacobi operators,
%and monotonicity properties of their zeros,
%the statement and proof of Ap\'ery's theorem on the irrationality of $\zeta(3)$,
%and Rivoal's wonderful theorem about the existence of infinitely many irrational numbers
%among the values of the Riemann zeta function at the odd integer numbers.
%This text must be regarded as simple class notes which are still under construction...
and the study of the arithmetic nature of the values of the Riemann zeta function at the integer numbers,
including Ap\'ery's theorem %(proved in the framework of the orthogonal polynomials theory)
 and Ball and Rivoal's results about the existence of infinitely many irrational numbers
among the values of the Riemann zeta function at the odd integer numbers.
%Several open questions were also mentioned along the course.

\medskip

%%%%%%%%%%%%%%%%%%%%%%%%%%%%%%%%%%%%%%%%%%%%%%%%%%%%%%%%%%%%%%%%%%%%%%%%%%%%%
%I would like to thanks to the students who chose the course
%(specially for the patience they had with my poor English...),
%and to some of them for pointing out several misprints,
%including Ali Moghanni, Dieudonn\'e Mbouna, Lili Song, Nikolaus Tsopanidis,
%Peter Lombaers, R\'uben Sousa, and Willian Ribeiro.
%Their questions and comments helped me to improve the present version of the text.
%\medskip
%
%Finally, I thanks deeply my Colleagues Renato \'Alvarez-Nodarse, Kenier Castillo, and M\'arcio do Nascimento de Jesus
%for their interest and comments that helped me to improve the present version of the text,
%or simply for supporting me since ever.
%%%%%%%%%%%%%%%%%%%%%%%%%%%%%%%%%%%%%%%%%%%%%%%%%%%%%%%%%%%%%%%%%%%%%%%%%%%%%
I would like to thank all the students who chose the course, %attended the course,
 Ali Moghanni, Dieudonn\'e Mbouna, Lili Song, Nikolaus Tsopanidis,
Peter Lombaers, R\'uben Sousa, and Willian Ribeiro.
Their questions and comments helped me to improve %the present version of these notes.
earlier versions of these notes.
\medskip

%Finally, I thanks deeply to Renato \'Alvarez-Nodarse, Kenier Castillo,
%and M\'arcio do Nascimento de Jesus for their interest and comments that helped me to
%improve the text, or simply for supporting me since ever.
%%or simply for always supporting me.
%\bigskip

\hfill
J. Petronilho (July 2017).\footnote{\,This version: revised on April 2018.}

%A instrução seguinte retira a palavra capítulo no título de cada capítulo (e deixa apenas o número do capítulo)

\def\chaptername{\bf\huge \,}

\tableofcontents

%%%%%%%%%%%%%%%%%%%%%%%%%%%%%%%%%%%%%%%%%%%%%%%%%%%%%%%%%%%%%%%%%%%%%%
%%%%%%%% Nota: o comando "\def\contentsname{\'Indice}"      %%%%%%%%%%
%%%%%%%% tem que ser colocado imediatamente antes           %%%%%%%%%%
%%%%%%%% do comando "\tableofcontents"                      %%%%%%%%%%
%%%%%%%% ( Sen\~ao sai "Conte\'udo" em vez de "\'Indice" )  %%%%%%%%%%
%%%%%%%%%%%%%%%%%%%%%%%%%%%%%%%%%%%%%%%%%%%%%%%%%%%%%%%%%%%%%%%%%%%%%%

%\pagenumbering{roman}
%\tableofcontents

\pagenumbering{arabic}
\setcounter{page}{1}

\setcounter{page}{-1} %colocar esta instrucao se for para vir o texto todo
%\setcounter{chapter}{-1}
%\include{opsf_appx1_2015-16}
%\include{opsf_tables_test2_2016-17}
%%%%%%%%%%%%%%%%%%%%%%%%%%%%%%%%%%%%%%%%%%%%%%%%%%%%%%%%%%%%%
%\setcounter{page}{1}
%\include{opsf_t1_2016-17}
%\setcounter{page}{1}
%\include{opsf_t2_2016-17}
%\setcounter{page}{1}
%\include{opsf_t3_2016-17}
%\setcounter{page}{1}
%\include{opsf_t4_2016-17}
%\setcounter{page}{1}
%\include{opsf_t5_2016-17}
%\setcounter{page}{1}
%\include{opsf_t6_2016-17}
%\setcounter{page}{1}
%\include{opsf_t7_2016-17}
%\setcounter{page}{1}
%\include{opsf_t8_2016-17}
%\setcounter{page}{1}
%\include{opsf_t9_2016-17}
%\setcounter{page}{1}
%%\include{opsf_t9a_2016-17}
%%\setcounter{page}{1}
%%\include{opsf_t10_2016-17}
%%\setcounter{page}{1}
%%\include{opsf_t11_2016-17}
%%\setcounter{page}{1}
%%\include{opsf_t12_2016-17}
%%\setcounter{page}{1}
%%\include{opsf_t13_2016-17}
%%\setcounter{page}{1}
%%\include{opsf_texer_2016-17}
%%%%%%%%%%%%%%%%%%%%%%%%%%%%%%%%%%%%%%%%%%%
%\appendix
%\setcounter{page}{1}
%\include{opsf_appx1_2016-17}
%\setcounter{page}{1}
%\include{opsf_appx2_2016-17}
%%\setcounter{page}{1}
%%\include{opsf_tables_2016-17}
%%\setcounter{page}{1}
%%\include{opsf_appx3_2016-17}
%%\include{ref}
%%%%%%%%%%%%%%%%%%%%%%%%%%%%%%%%%%%%%%%%%%%%%%%%%%%%%%%%%%%%

%%%%%%%%%%%%%%%%%%%%%%%%%%%%%%%%%%%%%%%%%%%%%%%%%%%%%%%%%%%%
\setcounter{page}{1}

%\setcounter{page}{1}
%\chapter{Preliminaries on Functional Analysis}
%\chapter{A crash course on locally convex spaces}
%\chapter{Topological preliminaries}

\chapter{Orthogonal polynomials: foundations}\label{OP-foundations}

\pagestyle{myheadings}\markright{Orthogonal polynomials: foundations}
\pagestyle{myheadings}\markleft{J. Petronilho}

\section{The spaces $\mathcal{P}$ and $\mathcal{P}'$}

Orthogonal polynomials (OP) can be studied from several different points of view.
%As remarked by the French mathematician Pascal Maroni,
Following the French mathematician Pascal Maroni,
from an algebraic viewpoint (meaning that orthogonality will be considered with respect to a moment linear functional, not necessarily represented by a weight function or a positive Borel measure),
it is very useful to consider OP as test functions living in an appropriate locally convex space (LCS), which we will denote by $\mathcal{P}$. This LCS is the set of all polynomials (with real or complex coefficients) endowed with a strict inductive limit topology, so that
\begin{snugshade}
\begin{equation}\label{P=indlimPn}
\mathcal{P}=\bigcup_{n=0}^\infty\mathcal{P}_n=\mbox{\rm ind\,lim}_n\,\mathcal{P}_n\; ,
\end{equation}
\end{snugshade}
\noindent
where $\mathcal{P}_n$ is the space of all polynomials of degree at most $n$.\footnote{\,For the sake of simplicity, we do not distinguish between {\it polynomial} and {\it polynomial function}.}
$\mathcal{P}_n$ being a finite dimensional vector space, all its norms are equivalent,
so there is no need to specify any one in particular.
%Being a finite dimensional vector space, all the norms in $\mathcal{P}_n$ are equivalent,
%so there is no need to specify any one in particular.
%Thus we may consider the following one, which is useful for what cames next:
%\begin{snugshade}
%$$
%\|p\|_n:=\sum_{j=0}^n|a_j|\; ,\quad p\equiv a_0+a_1x+\cdots+a_nx^n\in\mathcal{P}_n
%$$
%\end{snugshade}
%\noindent
%$(a_1,\cdots,a_n\in\mathbb{C})$.
For the development of the theory to be presented here it is not important to know much about the above topology (the definition and basic properties of LCS, including inductive limit topologies, can be found, e.g., in the book by
%Michael Reed and Barry Simon,
M. Reed and B. Simon, Chapter V---see also Appendix \ref{Appx1}), but the reader should keep in mind that the reason why such topology is introduced is because it implies the following fundamental property: the topological and the algebraic dual spaces of $\mathcal{P}$ coincide.

\begin{snugshade}
%\begin{teo}[Maroni]\label{expDB1}
\begin{theorem}\label{P*=Plinha}
%Let $\mathcal{P}$ be endowed with the strict inductive limit topology $(\ref{P=indlimPn})$,
Let $\mathcal{P}:=\mbox{\rm ind\,lim}_n\,\mathcal{P}_n$, as in $(\ref{P=indlimPn})$,
and let $\mathcal{P}^*$ and $\mathcal{P}'$ be the algebraic and the topological duals of $\mathcal{P}$,
respectively. Then %the set equality
\begin{snugshade}\vspace*{-1em}
\begin{equation}\label{P*equalsPprime}
\mathcal{P}'=\mathcal{P}^*\,.
\end{equation}
%\noindent holds.
\end{snugshade}
\end{theorem}
\end{snugshade}

{\it Proof.}
Obviously, $\mathcal{P}'\subseteq\mathcal{P}^*$.
To prove that $\mathcal{P}^*\subseteq\mathcal{P}'$, take ${\bf u}\in\mathcal{P}^*$.
%From the basic properties of the inductive limit topologies,
Taking into account Theorem \ref{Teo-Tcontin-LimInd}, to prove that ${\bf u}\in\mathcal{P}'$
it suffices to show that the restriction ${\bf u}|\mathcal{P}_n$ is continuous for every $n$.
But this is a trivial assertion, since ${\bf u}|\mathcal{P}_n$ is a linear functional
defined on a finite dimensional normed space.
\qed
\medskip

Equality (\ref{P*equalsPprime}) means that every linear functional defined in $\mathcal{P}$ is continuous (for the strict inductive limit topology in $\mathcal{P}$). This is a curious property, since we know that, $X$ being a normed vector space, it is true that $X'=X^*$ if ${\rm dim}\,X<\infty$, but, for the contrary, $X'\neq X^*$ whenever ${\rm dim}\,X=\infty$.
This last fact can be easily stated by using Zorn's lemma. (Exercise \ref{Ex-cp1-1}) Of course there is no contradiction between (\ref{P*equalsPprime}) and the fact that $\mathcal{P}$ is an infinite dimensional vector space, because $\mathcal{P}$ (carried with the inductive limit topology) is not a normed space. Indeed, being a strict inductive limit of the spaces $\mathcal{P}_n$, and taking into account that each $\mathcal{P}_n$ is a proper closed subspace of $\mathcal{P}_{n+1}$ (so that $\mathcal{P}$ is indeed an hyper strict inductive limit of the spaces $\mathcal{P}_n$), the general theory of LCS (cf. Theorem \ref{TeoContinutyLimInd}) ensures that $\mathcal{P}$ cannot be a metrizable space, and so \emph{a fortiori} it is not a normed space---or, to be more precise, it is not possible to provide $\mathcal{P}$ with a norm that generates in it the above inductive limit topology.

%In the dual space
In $\mathcal{P}'$ we consider the weak dual topology, which, by definition, is generated by the family of semi-norms
$s_p:\mathcal{P}'\to[0,+\infty[$, $p\in\mathcal{P}$, defined by
\begin{snugshade}
\begin{equation}\label{semip1}
s_p\big({\bf u}\big):=|\langle {\bf u},p\rangle|\; ,\quad {\bf u}\in\mathcal{P}'\; .
\end{equation}
\end{snugshade}
\noindent
It turns out that this family of semi-norms $s_p$ is equivalent to the family of semi-norms
$|\cdot|_n:\mathcal{P}'\to[0,+\infty[$, $n\in\mathbb{N}_0$, defined by
\begin{snugshade}
\begin{equation}\label{sharp1}
|{\bf u}|_n:=\max_{0\leq k\leq n}|\langle {\bf u},x^k\rangle|\; ,\quad {\bf u}\in\mathcal{P}'\; .
\end{equation}
\end{snugshade}
\noindent
Indeed, the following proposition holds.
\begin{snugshade}
%\begin{teo}[Maroni]\label{expDB1}
\begin{theorem}\label{Snumeravel}
$\mathcal{S}:=\{s_p:p\in\mathcal{P}\}$ and $\mathcal{S}_\sharp:=\{|\cdot|_n:n\in\mathbb{N}_0\}$,
with $s_p$ and $|\cdot|_n$ given by $(\ref{semip1})$--$(\ref{sharp1})$,
are equivalent families of seminorms in  $\mathcal{P}'$, provided $\mathcal{P}:=\mbox{\rm ind\,lim}_n\,\mathcal{P}_n$.
\end{theorem}
\end{snugshade}

{\it Proof.}
Given $p\in\mathcal{P}$, putting $p(x)=\sum_{j=0}^na_jx^j$ and $C(p):=\sum_{j=0}^n|a_j|$, we have
$$
s_{p}({\bf u})=|\langle {\bf u},p\rangle|=\Big| \sum_{j=0}^na_j\langle {\bf u},x^j\rangle \Big|
\leq C(p)|{\bf u}|_n\,,\quad \forall{\bf u}\in\mathcal{P}'\,.
$$
On the other hand, given $n\in\mathbb{N}_0$, setting $p_j(x):=x^j$ ($j=0,1,\ldots,n$), we have
$$
|{\bf u}|_n=\max_{0\leq j\leq n}|\langle {\bf u},x^j\rangle|
\leq\sum_{j=0}^n|\langle {\bf u},x^j\rangle|=\sum_{j=0}^n s_{p_j}({\bf u})
\,,\quad \forall{\bf u}\in\mathcal{P}'\,.
$$
Therefore, by Proposition \ref{teo-LCS-equiv-semin},
$\mathcal{S}$ and $\mathcal{S}_\sharp$ are equivalent families of semi-norms. %in $\mathcal{P}'$.
%(see \ref{ReedSimon1972i}, p.126).
\qed
\smallskip

Since $\mathcal{S}_\sharp$ is a countable family of seminorms, then by Theorem \ref{Snumeravel},
%and the theory of LCS (cf. Theorem \ref{LCS-metrizable} and Remark \ref{LCS-metrizable-metric})
together with Theorem \ref{LCS-metrizable} and Remark \ref{LCS-metrizable-metric},
%allow us to we conclude that $\mathcal{P}'$ is a metrizable space,
$\mathcal{P}'$ is a metrizable space, a metric being given by
\begin{snugshade}
\begin{equation}\label{metricPprime}
\varrho({\bf u},{\bf v}):=\sum_{n=0}^\infty\frac{1}{2^n}\frac{|{\bf u}-{\bf v}|_n}{1+|{\bf u}-{\bf v}|_n}\; ,
\quad {\bf u},{\bf v}\in\mathcal{P}'\; .
\end{equation}
\end{snugshade}
\noindent
Moreover, $\mathcal{P}'$ is a Fr\'echet space. (Exercise \ref{Ex-cp1-3})

\section{Dual basis in $\mathcal{P}^*$}

Since we will work in the dual space $\mathcal{P}^*$, it would be useful to explicitly building bases in $\mathcal{P}^*$.
This makes sense, since (\ref{P*equalsPprime}) allow us writing expansions (finite or infinite sums) of the elements of $\mathcal{P}^*$ in terms of the elements of a given basis, in the sense of the weak dual topology in $\mathcal{P}'$. Such basis in $\mathcal{P}^*$ may be achieved in a natural way, using simple sets of polynomials. A {\sl simple set} in $\mathcal{P}$ is a sequence of polynomials, $\{R_n\}_{n\geq0}$, such that ${\rm deg}\,R_n=n$ for every $n\in\mathbb{N}_0$ (where $R_0\equiv\mbox{\rm const.}\neq0$). To any simple set in $\mathcal{P}$, $\{R_n\}_{n\geq0}$, we may associate a {\sl dual basis}, which, by definition, is a sequence of linear functionals $\{{\bf a}_n\}_{n\geq0}$, being ${\bf a}_n:\mathcal{P}\to\mathbb{C}$, such that
\begin{snugshade}
$$
\langle {\bf a}_n,R_k\rangle:=\delta_{n,k}\quad (n,k=0,1,2,\cdots)\; ,
$$
\end{snugshade}
\noindent
where $\delta_{n,k}$ represents the Kronecker symbol ($\delta_{n,k}=1$ if $n=k$; $\delta_{n,k}=0$ if $n\neq k$).
The following is a fundamental result.
Together with equality (\ref{P*equalsPprime}) it is on the foundations of the (algebraic) theory of OP.

\begin{snugshade}
%\begin{teo}[Maroni]\label{expDB1}
\begin{theorem}\label{expDB1}
Let $\{R_n\}_{n\geq0}$ be a simple set in $\mathcal{P}$ and $\{{\bf a}_n\}_{n\geq0}$ the associated dual basis.
Let ${\bf u}\in\mathcal{P}^*$. Then
\begin{equation}\label{expDB2}
{\bf u}=\sum_{n=0}^\infty\langle{\bf u},R_n\rangle\,{\bf a}_n\, , %\quad \lambda_n:=\langle{\bf u},R_n\rangle \,,
\end{equation}
in the sense of the weak dual topology in $\mathcal{P}'$.
\end{theorem}
\end{snugshade}

{\it Proof.}
Notice first that the assertion makes sense, according with (\ref{P*equalsPprime}). To prove it,
fix $N\in\mathbb{N}$ and let $${\bf s}_N:=\sum_{n=0}^{N-1}\lambda_n{\bf a}_n\quad(\lambda_n:=\langle{\bf u},R_n\rangle)$$
be the partial sum of order $N$ of the series appearing in (\ref{expDB2}).
We need to show that
$$
\lim_{N\to\infty}\langle {\bf s}_N-{\bf u},p\rangle=0\;,\quad\forall p\in\mathcal{P}\;.
$$
Clearly, it suffices to prove that this equality holds for $p\in\{R_0,R_1,R_2,\cdots\}$.
Indeed, fix $k\in\mathbb{N}_0$. Then, for $N>k$,
$$
\langle {\bf s}_N-{\bf u},R_k\rangle
=\sum_{n=0}^{N-1}\langle{\bf u},R_n\rangle\langle{\bf a}_n,R_k\rangle -\langle{\bf u},R_k\rangle
=0\; ,
$$
hence $\;\lim_{N\to\infty}\langle {\bf s}_N-{\bf u},R_k\rangle=0$.
\qed
\bigskip

%\begin{example}
%\rm
%The dual basis corresponding to the simple set $\{x^n\}_{n\geq0}$ is given by
%$$
%{\bf a}_n:=\frac{(-1)^n}{n!}\bm{\delta}^{(n)}\;,\quad n\in\mathbb{N}_0\;,
%$$
%where $\bm{\delta}:\mathcal{P}\to\mathbb{C}$ is the Dirac functional,
%defined by $\langle\bm{\delta},f\rangle=f(0)$, $f\in\mathcal{P}$,
%and $\bm{\delta}^{(n)}$ is the (distributional) derivative of $\bm{\delta}$.
%\end{example}

\section{Basic operations in $\mathcal{P}$ and in $\mathcal{P}'$}

In this section we introduce some fundamental operations in the framework of the algebraic theory of OP.
Given a functional ${\bf u}\in\mathcal{P}'$, we will denote by
\vspace*{-0.5em}
\begin{snugshade}
$$
u_n:=\langle{\bf u},x^n\rangle\;,\quad n\in\mathbb{N}_0\; ,
$$
\end{snugshade}
\noindent
the {\sl moment} of order $n$ of ${\bf u}$.
Clearly, if ${\bf u}$ and ${\bf v}$ are two functionals in $\mathcal{P}'$
such that the corresponding sequences of moments satisfy
$u_n=v_n$ for all $n\in\mathbb{N}_0$, then ${\bf u}={\bf v}$.
Therefore, {\it  each functional ${\bf u}\in\mathcal{P}'$ is uniquely determined
by its sequence of moments.} % $\{u_n\}_{n\geq0}$}.
%uniquely determines}.

Define operators $M_\phi$, $T$, and $\theta_c$, from $\mathcal{P}$ into $\mathcal{P}$, by
\vspace*{-0.5em}
\begin{snugshade}\vspace*{-0.5em}
\begin{equation}\label{thetac}%$$
M_\phi p(x):= \phi(x)p(x)\;,\quad Tp(x):= -p'(x)\;,\quad \theta_cp(x):=\frac{p(x)-p(c)}{x-c}\quad (p\in\mathcal{P})\;,
\end{equation}%$$
\end{snugshade}
\noindent
where $\phi\in\mathcal{P}$ (fixed), $'$ denotes
derivative with respect to $x$, and $c\in\mathbb{C}$.
%and $\theta_c:\mathcal{P}\to\mathcal{P}$ is defined by
%\begin{snugshade}
%$$
%\theta_cp(x):=\frac{p(x)-p(c)}{x-c}
%$$
%\end{snugshade}
%\noindent
Note that $\theta_cp(x)$ is defined as above if $x\neq c$,
with the obvious definition $\theta_cp(c):=p'(c)$ if $x=c$
(so that, indeed, $\theta_cp\in\mathcal{P}$). %for each $p\in\mathcal{P}$).
By Theorem \ref{Tdual-contin}, the dual operators (cf. Appendix \ref{Appx1})
%corresponding the above three operators $M_\phi$, $T$, and $\theta_c$,
$M_\phi'$, $T'$, and $\theta_c'$
belong to $\mathcal{L}(\mathcal{P}',\mathcal{P}')$.
For each ${\bf u}\in\mathcal{P}'$, the images $M_\phi'{\bf u}$, $T'{\bf u}$, and $\theta_c'{\bf u}$,
%by each of those dual operators over an element ${\bf u}\in\mathcal{P}'$
%give rise to three new
are elements (functionals) in $\mathcal{P}'$,
hereafter denoted by $\phi\,{\bf u}$, $D{\bf u}$, and $(x-c)^{-1}{\bf u}$. %, respectively.

%are the mappings $M_q'$, $T'$, and $\theta_c'$,
%from $\mathcal{P}'$ into $\mathcal{P}'$, given for each ${\bf u}\in\mathcal{P}'$ and $p\in\mathcal{P}$ by
%$$
%\langle M_q'{\bf u}, p\rangle:=\langle{\bf u}, qp\rangle\;,\quad
%\langle T'{\bf u}, p\rangle:=-\langle{\bf u}, p'\rangle\;,\quad
%\langle \theta_c'{\bf u}, p\rangle:=\langle{\bf u}, \theta_cp\rangle\;.
%$$

\begin{snugshade}
\begin{definition}\label{def-left-mult}
Let ${\bf u}\in\mathcal{P}'$, $\phi\in\mathcal{P}$, and $c\in\mathbb{C}$.
\begin{enumerate}
\item[{\rm (i)}]
the {\sl left multiplication} of ${\bf u}$ by $\phi$, denoted by $\phi{\bf u}$,
is the functional defined by
$$
\langle \phi {\bf u}, p\rangle:=\langle {\bf u}, \phi p\rangle\;,\quad p\in\mathcal{P}\;;
$$
\item[{\rm (ii)}]
the {\sl derivative} of ${\bf u}$, denoted by $D{\bf u}$,
is the functional in $\mathcal{P}'$ defined by
$$
\langle D{\bf u}, p\rangle:=-\langle{\bf u}, p'\rangle\;,\quad p\in\mathcal{P}\;;
$$
\item[{\rm (iii)}]
the {\sl division} of ${\bf u}$ by $x-c$, denoted by $(x-c)^{-1}{\bf u}$,
is the functional defined by
$$
\langle (x-c)^{-1}{\bf u}, p\rangle:=\langle{\bf u}, \theta_cp\rangle
=\Big\langle{\bf u},\frac{p(x)-p(c)}{x-c} \Big\rangle\;,\quad p\in\mathcal{P}\;.
$$
\end{enumerate}
\end{definition}
\end{snugshade}

Note that these definitions, introduced here by duality with respect to the
operators defined in (\ref{thetac}), are in accordance with those usually given in
Theory of Distributions (this explains the minus sign appearing in the second definition).

\begin{snugshade}
%\begin{teo}[Maroni]\label{expDB1}
\begin{theorem}%\label{expDB1}
Let ${\bf u}\in\mathcal{P}'$ and $\phi\in\mathcal{P}$.
Then
$$
D(\phi{\bf u})=\phi'{\bf u}+\phi\,D{\bf u}\;.
$$
\end{theorem}
\end{snugshade}

{\it Proof.}
Indeed, for each $p\in\mathcal{P}$, we have
$$
\begin{array}{rcl}
\langle D(\phi{\bf u}),p \rangle &=&
-\langle \phi{\bf u},p' \rangle \,=\, -\langle {\bf u},\phi p' \rangle
\,=\,-\langle{\bf u},-\phi'p+(\phi p)'\rangle \\ [.25em]
&=& \langle {\bf u},\phi' p\rangle-\langle {\bf u},(\phi p)'\rangle
\,=\,\langle \phi'{\bf u}, p\rangle+\langle D{\bf u},\phi p\rangle\\ [.25em]
&=& \langle \phi'{\bf u}+\phi\,D{\bf u},p \rangle\;,
\end{array}
$$
hence the desired equality holds.
\qed

\medskip

By definition, the left multiplication of a functional in $\mathcal{P}'$ by a polynomial
is another functional in $\mathcal{P}'$. We may also define a
right multiplication of a linear functional by a polynomial.
The result is a polynomial.

\begin{snugshade}
\begin{definition}\label{def-right-mult}
Let ${\bf u}\in\mathcal{P}'$ and $\psi\in\mathcal{P}$.
The {\sl right multiplication} of ${\bf u}$ by $\psi$, denoted by ${\bf u}\psi$,
is the polynomial defined by
\vspace*{-0.5em}
$$
{\bf u}\psi(x):=\langle {\bf u}_\xi, \theta_\xi(x\psi)\rangle
=\Big\langle{\bf u}_\xi,\frac{x\psi(x)-\xi\psi(\xi)}{x-\xi} \Big\rangle\;,
$$
where the subscript $\xi$ in ${\bf u}_\xi$ means that ${\bf u}$ acts in polynomials of the variable $\xi$.
\end{definition}
\end{snugshade}

Setting $\psi(x):=\sum_{i=0}^na_ix^i$, the polynomial ${\bf u}\psi$ is explicitly given by
\begin{snugshade}
\begin{equation}\label{def-right-mult-explic}
{\bf u}\psi(x)=\sum_{i=0}^n\Big(\sum_{j=i}^na_ju_{j-i}\Big)x^i\;,
\end{equation}
\end{snugshade}
\noindent
and it also admits the following useful matrix representation:
\begin{snugshade}
\begin{equation}\label{def-right-mult-mz}
{\bf u}\psi(x)=[a_0\;a_1\cdots a_n]
\left[
\begin{array}{ccccc}
u_0 & 0 & 0 & \cdots & 0 \\
u_1 & u_0  & 0 & \cdots & 0 \\
u_2 & u_1 & u_0  & \cdots & 0 \\
\vdots & \vdots & \vdots & \ddots & \vdots \\
u_n & u_{n-1} & u_{n-2}  & \cdots & u_0
\end{array}
\right]
\left[
\begin{array}{c}
1 \\ x \\ x^2 \\ \vdots \\ x^n
\end{array}
\right]\;.
\end{equation}
\end{snugshade}

The right multiplication of a functional by a polynomial enable us
to introduce a product in $\mathcal{P}'$, by duality.
Indeed, fix ${\bf u}\in\mathcal{P}'$,
and let $T_{\bf u}:\mathcal{P}\to\mathcal{P}$ be defined by
$$
T_{\bf u}p:={\bf u}p\; ,\quad p\in\mathcal{P}\,.
$$
The dual operator, $T_{\bf u}^\prime:\mathcal{P'}\to\mathcal{P'}$
(${\bf v}\mapsto T_{\bf u}^\prime{\bf v}$), is given by
$$
\langle T_{\bf u}^\prime{\bf v},p\rangle:=\langle {\bf v},T_{\bf u}p\rangle
=\langle {\bf v},{\bf u}p\rangle\;,\quad p\in\mathcal{P}\;.
$$
Thus, we may introduce a product in $\mathcal{P}'$,
by duality with respect to the right multiplication of a functional by a polynomial.

\begin{snugshade}
\begin{definition}\label{def-prod-in-Pprime}
Let ${\bf u},{\bf v}\in\mathcal{P}'$.
The {\sl product} ${\bf uv}$ is the functional in $\mathcal{P}'$ given by
$$
\langle {\bf uv}, p\rangle:=\langle{\bf v},{\bf u}p\rangle\;,\quad p\in\mathcal{P}\,.
$$
\end{definition}
\end{snugshade}

This product is commutative. This fact may be seen easily by noticing that
the moments of the functionals ${\bf uv}$ and ${\bf vu}$ coincide:
\vspace*{-0.5em}
\begin{snugshade}
\begin{equation}\label{uv-comutes}
\langle{\bf uv},x^n\rangle=\sum_{i+j=n}u_iv_j=\langle{\bf vu},x^n\rangle\,,\quad n\in\mathbb{N}_0\;.
\end{equation}
\end{snugshade}
%\noindent
Further, there exists unit element in $\mathcal{P}'$, namely the Dirac functional at the origin,
$\bm{\delta}\equiv\bm{\delta}_0$. Indeed, using (\ref{uv-comutes}), it is easy to prove that
\vspace*{-0.5em}
\begin{snugshade}
$${\bf u}\bm{\delta}={\bf u}\,,\quad {\bf u}\in\mathcal{P}'\;.$$
\end{snugshade}
\noindent
Recall that the {\sl Dirac functional} at a point $c\in\mathbb{C}$,
$\bm{\delta}_c:\mathcal{P}\to\mathbb{C}$, is defined by
\vspace*{-0.5em}
\begin{snugshade}
$$\langle\bm{\delta}_c,p\rangle:=p(c)\,,\quad p\in\mathcal{P}\;.$$
\end{snugshade}

The next proposition lists some basic properties concerning the above operations.
The proof is left to the reader. (Exercise \ref{Ex-cp1-6})

\begin{snugshade}
\begin{proposition}\label{prop-PandPprime}
Let ${\bf u},{\bf v},{\bf w}\in\mathcal{P}'$, $p\in\mathcal{P}$, and $c\in\mathbb{C}$.
Then:\vspace*{-0.5em}
\begin{multicols}{2}
\begin{enumerate}
\item[{\rm 1.}] $\bm{\delta}p=p$
\item[{\rm 2.}] ${\bf v}({\bf u}p)=({\bf v}{\bf u})p$
\item[{\rm 3.}] $({\bf u}+{\bf v}){\bf w}={\bf u}{\bf w}+{\bf v}{\bf w}$
\item[{\rm 4.}] $({\bf u}{\bf v}){\bf w}={\bf u}({\bf v}{\bf w})$
\item[{\rm 5.}] ${\bf u}$ {\rm has an inverse iff} $u_0\neq0$
%\item[{\rm 6.}] $p({\bf u}{\bf v})=(p{\bf u}){\bf v}+x({\bf u}\theta_0p){\bf v}$
\item[{\rm 6.}\hspace*{3em}] \hspace*{-3em} $p({\bf u}{\bf v})=(p{\bf u}){\bf v}+x({\bf u}\theta_0p){\bf v}$
\item[{\rm 7.}\hspace*{3em}] \hspace*{-3em} $D({\bf u}p)=(D{\bf u})p+{\bf u}p'+{\bf u}\theta_0p$
\item[{\rm 8.}\hspace*{3em}] \hspace*{-3em} $D({\bf u}{\bf v})=(D{\bf u}){\bf v}+{\bf u}D{\bf v}+x^{-1}({\bf u}{\bf v})$
\item[{\rm 9.}\hspace*{3em}] \hspace*{-3em} $(x-c)\big((x-c)^{-1}{\bf u}\big)={\bf u}$
\item[{\rm 10.}\hspace*{3em}] \hspace*{-3em} $(x-c)^{-1}\big((x-c){\bf u}\big)={\bf u}-u_0\bm{\delta}_c$\,.
%\item[{\rm 11.}\hspace*{3em}] \hspace*{-3em} $(x-c){\bf u}=(x-c){\bf v}\;$ {\rm iff} $\;{\bf u}={\bf v}+(u_0-v_0)\delta_c$\,.
\end{enumerate}
\end{multicols}
\end{proposition}
\end{snugshade}

%{\it Proof.}
%The proof is left to the reader. (Exercise \ref{Ex-cp1-6})
%%%%%%%%%%%%%%%%%%%%%%%%%%%%%%%%%%%%%%%%%%%%%%%%%%%%%%%%%%%%%%%%%%%%%%%%%%%%%%%%%
%\footnote{%\begin{remark}\em
%\;On the course of the proof of several properties listed in Proposition \ref{prop-PandPprime}, it may be useful to use the following identities, valid for any array $\;\{\alpha_{i,j}\}_{0\leq j\leq i\leq n}\;$ of $\frac{(n+1)(n+2)}{2}$ complex numbers:
%\begin{snugshade}
%$$
%\sum_{i=0}^n\sum_{j=0}^i\alpha_{i,j}=\sum_{i=0}^n\sum_{j=i}^n\alpha_{j,j-i}=\sum_{j=0}^n\sum_{i=j}^n\alpha_{i,j}\;.
%$$
%\end{snugshade}%}
%\noindent
%Indeed, dispose the array elements to form a right triangle:
%$$
%\begin{array}{cccccc}
%\alpha_{0,0} &&&&& \\
%\alpha_{1,0} &\alpha_{1,1}&&&& \\
%\alpha_{2,0} &\alpha_{2,1}&\alpha_{2,2}&&& \\
%\vdots&\vdots&\vdots&\ddots&& \\
%%\alpha_{n-1,0} &\alpha_{n-1,1}&\alpha_{n-1,2}&\ldots&\alpha_{n-1,n-1}& \\
%%\alpha_{n,0} &\alpha_{n,1}&\alpha_{n,2}&\ldots&\alpha_{n,n-1}& \alpha_{n,n}
%\alpha_{n,0} &\alpha_{n,1}&\alpha_{n,2}&\ldots&\alpha_{n,n}&
%\end{array}
%$$
%Then we only need to notice that the first (double) sum corresponds to adding the elements by horizontal lines (from top to %bottom), the second sum corresponds to adding the elements by diagonal lines (starting from the ``hypotenuse''), and the third one %corresponds to adding the elements by vertical lines (from left to right).}
%%%%%%%%%%%%%%%%%%%%%%%%%%%%%%%%%%%%%%%%%%%%%%%%%%%%%%%%%%%%%%%%%%%%%%%%%%%%%%%%%
%\qed
%\smallskip

Notice that from property 10 one also obtains
\begin{snugshade}
$$%\begin{equation}\label{xcu-xcv}
(x-c){\bf u}=(x-c){\bf v}\quad \mbox{\rm iff}\quad{\bf u}={\bf v}+(u_0-v_0)\bm\delta_c\;.
$$%\end{equation}
\end{snugshade}

We conclude this section by stating the following

\begin{snugshade}
\begin{proposition}\label{prop-uvzero}
Let ${\bf u}\in\mathcal{P}'$ and $p,q\in\mathcal{P}\setminus\{0\}$, and denote by $Z_p$ and $Z_q$ the zeros of $p$ and $q$, respectively.
Then the following property holds:
\begin{equation}\label{pu0qu}
Z_p\cap Z_q=\emptyset\quad\wedge\quad p{\bf u}=q{\bf u}={\bf 0}\qquad\Rightarrow\qquad{\bf u}={\bf 0}\;.
\end{equation}
\end{proposition}
\end{snugshade}

{\it Proof.}
Denote the degrees of $p$ and $q$ by $m$ and $n$, respectively.
%Without loss of generality, we may assume that $p$ and $q$ are monic.
Without loss of generality, we may assume that $q$ is monic (i.e., the coefficient of $x^n$ is equal to $1$).
Moreover, if $m=0$ or $n=0$ (i.e., if $p$ or $q$ is a nonzero constant) then trivially ${\bf u}={\bf 0}$,
so we may assume as well that $m,n\geq1$.
%Accordingly, write $$p(x)=\sum_{j=0}^ma_jx^j\;,$$
%where $a_m=1$ and $a_j\in\mathbb{C}$ for each $j\in\{1,\ldots,m-1\}$.
We will prove that ${\bf u}={\bf 0}$ by induction on the degree $n$ of $q$ (and keeping $p$ fixed).

Suppose that $n=1$, so that $q(x)=x-b$, with $b\in\mathbb{C}$ and $p(b)\neq0$.
Then, from $q{\bf u}={\bf 0}$ we have $x{\bf u}=b{\bf u}$,
hence $x^j{\bf u}=b^j{\bf u}$ for each $j\in\mathbb{N}_0$.
%Therefore, writing $p(x)=\sum_{j=0}^ma_jx^j$ (with $a_j\in\mathbb{C}$ for each $j\in\{1,\ldots,m\}$),
%we have $p{\bf u}=\sum_{j=0}^ma_jx^j{\bf u}=p(b){\bf u}$,
Therefore, $p{\bf u}=p(b){\bf u}$,
and so equation $p{\bf u}={\bf 0}$ is equivalent to $p(b){\bf u}={\bf 0}$.
Thus ${\bf u}={\bf 0}$.

Suppose now (induction hypothesis) that property (\ref{pu0qu}) holds for each polynomial $q$ of degree $n$
(which does not share zeros with $p$ and fulfills $q{\bf u}={\bf 0}$).
Let $\widetilde{q}$ be a polynomial of degree $n+1$
which does not share zeros with $p$ and fulfills $\widetilde{q}{\bf u}={\bf 0}$.
Then we may write $\widetilde{q}=(x-b)q$ where $q$ is a polynomial of degree $n$
which does not share zeros with $p$ and $p(b)\neq0$. Let ${\bf v}:=q{\bf u}$.
Then $(x-b){\bf v}=\widetilde{q}{\bf u}={\bf 0}$ and $p{\bf v}=q(p{\bf u})={\bf 0}$,
and so, by the case $n=1$ already proved, we conclude that ${\bf v}={\bf 0}$, i.e., $q{\bf u}={\bf 0}$.
Thus, one has $p{\bf u}=q{\bf u}={\bf 0}$, and since $q$ is a polynomial of degree $n$
which does not share zeros with $p$, it follows by the induction hypothesis that ${\bf u}={\bf 0}$.
\qed
\medskip

\section{The formal Stieltjes series}

Let $\Delta'$ %$\Delta'\equiv\mathbb{C}[[z]]$
the vector space of the formal series in the variable $z$
with coefficients in $\mathbb{C}$:\vspace*{-0.5em}
\begin{snugshade}
$$
\Delta':=\Big\{\sum_{n=0}^\infty c_nz^n\;\Big|\; c_n\in\mathbb{C}\;\;\mbox{for all $n\in\mathbb{N}_0$}\Big\}\;.
$$
\end{snugshade}
\noindent
In $\Delta'$ the operations of addition, multiplication and scalar multiplication
are defined in the usual way (see e.g. Tr\`eves's book \ref{Treves1967i}).
Endowing $\Delta'$ with the family of seminorms $\rho_n:\Delta'\to[0,+\infty)$, $n\in\mathbb{N}_0$, where
\vspace*{-0.5em}
\begin{snugshade}
$$
\rho_n\Big(\sum_{n=0}^\infty c_nz^n\Big):=\max_{0\leq j\leq n}|c_j|\; ,
$$
\end{snugshade}
\noindent
$\Delta'$ becomes a metrizable LCS. This space can be identified with $\mathcal{P}'$.

\begin{snugshade}
%\begin{teo}[Maroni]\label{expDB1}
\begin{theorem}\label{operF-Delta}
The operator $F:\mathcal{P}'\to\Delta'$ given by
$$
{\bf u}\in\mathcal{P}'\quad\mapsto\quad F({\bf u}):=\sum_{n=0}^\infty u_nz^n
$$
is a topological isomorphism ($\mathcal{P}'$ being endowed with the weak dual topology).
\end{theorem}
\end{snugshade}

{\it Proof.}
%It is easily seen that $F$ is a linear and bijective operator.
Clearly, $F$ is linear and bijective.
Moreover, for each ${\bf u}\in\mathcal{P}'$ and $n\in\mathbb{N}_0$,
$$
\rho_n\big(F({\bf u})\big)=\max_{0\leq j\leq n}|u_j|=|{\bf u}|_n\;.
$$
Therefore, since the family of seminorms $\{|\cdot|_n:n\in\mathbb{N}_0\}$
generates the topology in $\mathcal{P}'$ (cf. Theorem \ref{Snumeravel}),
we deduce, for each sequence $\{{\bf u}_j\}_{j\geq0}$ in $\mathcal{P}'$,
$$
\begin{array}{rcl}
F({\bf u}_j)\to0 \;\; \mbox{\rm (in $\Delta'$)}
&\mbox{\rm iff}& \rho_n\big(F({\bf u}_j)\big)\to0\;\;\mbox{for each $n\in\mathbb{N}_0$} \\ [0.5em]
&\mbox{\rm iff}& |{\bf u}_j|_n\to0 \;\;\mbox{for each $n\in\mathbb{N}_0$} \\ [0.5em]
&\mbox{\rm iff}& {\bf u}_j\to\textbf{0} \;\; \mbox{\rm (in $\mathcal{P}'$)} \;.
\end{array}
$$
Thus, $F$ is bicontinuous, which proves the theorem.
\qed
\smallskip

Note that the isomorphism $F$ allow us to transfer the algebraic structure
from $\Delta'$ into $\mathcal{P}'$.
This fact is accomplished through the formal Stieltjes series.
\begin{snugshade}
\begin{definition}\label{def-serie-Stieltjes}
Let ${\bf u}\in\mathcal{P}'$.
The {\sl formal Stieltjes series} associated with ${\bf u}$ is
\vspace*{-0.5em}
$$
S_{\bf u}(z):=-\sum_{n=0}^\infty\frac{u_n}{z^{n+1}}\equiv
-\frac{\;1\;}{z}\,F\left(\frac{\;1\;}{z}\right)\;.
$$
\end{definition}
\end{snugshade}
%\noindent
Note that $S_{\bf u}(z)$ gives a representation for the sequence of moments, $\{u_n\}_{n\geq0}$, of ${\bf u}$.
The formal Stieltjes series is an important tool in the theory of OP,
allowing us to state characterizations theorems concerning certain important
classes of OP, e.g., classical OP, semiclassical OP, and Laguerre-Hahn OP.
$S_{\bf u}(z)$ and its formal derivative,
\begin{snugshade}
\vspace*{-0.1em}
$$
S_{\bf u}'(z):=\sum_{n=0}^\infty\frac{(n+1)u_n}{z^{n+2}}\;,
$$
\end{snugshade}
\noindent
become tools of major importance in the study of these classes of OP.
\medskip

%$S_{\bf u}(z)$ is a representation for the moments $u_n$ of ${\bf u}$.
%Thus, since the sequence $\{u_n\}_{n\geq0}$ characterizes ${\bf u}$,
%then so does $S_{\bf u}(z)$.
%Formally, $S_{\bf u}(z)$ admits the representation
%%($\textbf{u}_x$ means that $\textbf{u}$ acts on functions of the variable $x$):
%$$S_{\bf u}(z)=\Big\langle\textbf{u}_x,\dfrac{1}{x-z}\Big\rangle\;,$$
%where $\textbf{u}_x$ means that $\textbf{u}$ acts on functions of the variable $x$.

%\section{Remarks and historical notes}
%
%Most of the material presented in this chapter is contained in the work by Pascal Maroni, specially in references \ref{},
%although some aspects of the presentation reflect the author's viewpoint.

\section*{Exercises}
%\bigskip

{\small
%\noindent
\begin{enumerate}[label=\emph{\bf \arabic*.},leftmargin=*]
\item\label{Ex-cp1-1}
Show that in any infinite dimensional normed space there are linear functionals which are not continuous---hence, an equality like (\ref{P*equalsPprime}) cannot holds on an infinite dimensional normed space. As a consequence, being $X$ a normed space, there holds:
$$
\left\{
\begin{array}{rcl}
X^\prime= X^* & \mbox{\rm if} & {\rm dim}\,X<\infty \; ; \\ [0.25em]
X^\prime\subsetneqq X^* & \mbox{\rm if} & {\rm dim}\,X=\infty \; .
\end{array}
\right.
$$
%\smallskip

\noindent
({\sl Hint:} If $\mbox{\rm dim}\,X=\infty$, there exists a denumerable subset $E=\{e_n\,|\,n\in\mathbb{N}\}\subset X$ whose elements are linearly independent unit vectors in $X$. Set $Y:=\langle E\rangle$. Then, Zorn's Lemma ensures that $Y$ has a complementary subspace in $X$, say, $Z$, so that $X=Y+Z$, with $Y\cap Z=\{0\}$, and each $x\in X$ admits a unique representation as $x=y+z$, with $y\in Y$ and $z\in Z$ [see e.g. Lax's book, p. 14, Lemma 9]. Denoting by $\mathbb{K}$ ($=\mathbb{R}$ or $\mathbb{C}$) the field of scalars associated with the vector space $X$, define $\ell:X\to\mathbb{K}$ by $\ell(x):=\ell_0(y)$, being $\ell_0:Y\to\mathbb{K}$ the linear functional given by
$$
\ell_0\Big(\sum_{k=1}^N\alpha_ke_k\Big)=\sum_{k=1}^Nk\alpha_k\;\quad
\big(\alpha_1,\cdots,\alpha_N\in\mathbb{K}\,;\;\;N\in\mathbb{N}\big)\; .
$$
Show that $\ell\in X^*\setminus X^\prime$.)
\medskip

\item
Let $\mathbb{P}:=\mathbb{K}[x]$ be the set of all polynomials (regarded as polynomial functions) with coefficients in
$\mathbb{K}$ ($=\mathbb{R}$ or $\mathbb{C}$). Prove that the mapping $\|\cdot\|:\mathbb{P}\to[0,+\infty)$ defined by
$$
\|f\|:=\max_{0\leq k\leq N}|a_k|\; ,\quad f(x)\equiv \sum_{k=0}^Na_kx^k\in\mathbb{K}[x]\quad (N:=\deg f)\,,
% (N\in\mathbb{N}_0)\; ,
$$
is a norm in $\mathbb{P}$, but with this norm $\mathbb{P}$ is not a complete (Banach) space.
\smallskip

\noindent
({\sl Hint:} To prove noncompleteness use the Banach-Steinhauss theorem.)
\medskip

\item\label{Ex-cp1-3}
Prove that $\mathcal{P}'$ is a Fr\'echet space.
\smallskip

({\sl Hint.} The weak dual topology in $\mathcal{P}'$ is generated by the countable family of seminorms $\mathcal{S}_\sharp:=\{|\cdot|_n:n\in\mathbb{N}_0\}$, hence a given sequence $\{{\bf u}_n\}_{n\geq0}$ in $\mathcal{P}'$ is Cauchy if and only if
$$
\forall\epsilon>0\;\forall k\in\mathbb{N}_0\;\exists n_0=n_0(\epsilon,k)\in\mathbb{N}\;:\;\forall n,m\in\mathbb{N}\;,\;
n,m\geq n_0\;\Rightarrow\; |{\bf u}_n-{\bf u}_m|_k<\epsilon\;.
$$
To prove that $\mathcal{P}'$ is Fr\'echet we have to show that each Cauchy sequence $\{{\bf u}_n\}_{n\geq0}$ in $\mathcal{P}'$ converges, i.e., there exists ${\bf u}\in\mathcal{P}'$ such that $|{\bf u}_n-{\bf u}|_k\to 0$ as $n\to\infty$, for every $k\in\mathbb{N}_0$.)
\medskip

\item
Let ${\bf u}\in\mathcal{P}'$, $\phi\in\mathcal{P}$, and $n\in\mathbb{N}_0$.
Prove Leibniz formula
$$
D^n(\phi{\bf u})=\sum_{k=0}^n{n\choose k}\phi^{(k)}\,D^{n-k}{\bf u}\;.
$$
\medskip

%\item
%Find the dual basis associated with the simple set of polynomials $\{x^n\}_{n\geq0}$.
%\medskip

\item
Show that the dual basis $\{{\bf a}_n\}_{n\geq0}$ corresponding to the simple set $\{x^n\}_{n\geq0}$ is given by
$$
{\bf a}_n:=\frac{(-1)^n}{n!}\bm{\delta}^{(n)}\;,  %\quad n\in\mathbb{N}_0\;,
$$
%where $\bm{\delta}:\mathcal{P}\to\mathbb{C}$ is the Dirac functional,
%defined by $\langle\bm{\delta},f\rangle=f(0)$, $f\in\mathcal{P}$,
where $\bm{\delta}^{(n)}$ is the (distributional) derivative of order $n$
of the Dirac functional $\bm{\delta}\equiv\bm{\delta}_0$.
Conclude that each functional ${\bf u}\in\mathcal{P}'$ admits the representation
$$
{\bf u}=\sum_{n=0}^\infty (-1)^n \frac{u_n}{n!}\bm{\delta}^{(n)}\;,
$$
in the sense of the weak dual topology in $\mathcal{P}'$.
\medskip

\item\label{Ex-cp1-6}
Prove the properties listed in Proposition \ref{prop-PandPprime}.\footnote{\,On the course of the proof of several properties listed in Proposition \ref{prop-PandPprime}, it may be useful to use the following identities, valid for any array $\;\{\alpha_{i,j}\}_{0\leq j\leq i\leq n}\;$ of $\frac{(n+1)(n+2)}{2}$ complex numbers:
\begin{snugshade}
$$
\sum_{i=0}^n\sum_{j=0}^i\alpha_{i,j}=\sum_{i=0}^n\sum_{j=i}^n\alpha_{j,j-i}=\sum_{j=0}^n\sum_{i=j}^n\alpha_{i,j}\;.
$$
\end{snugshade}%}
\noindent
Indeed, dispose the array elements to form a right triangle:
$$
\begin{array}{cccccc}
\alpha_{0,0} &&&&& \\
\alpha_{1,0} &\alpha_{1,1}&&&& \\
\alpha_{2,0} &\alpha_{2,1}&\alpha_{2,2}&&& \\
\vdots&\vdots&\vdots&\ddots&& \\
%\alpha_{n-1,0} &\alpha_{n-1,1}&\alpha_{n-1,2}&\ldots&\alpha_{n-1,n-1}& \\
%\alpha_{n,0} &\alpha_{n,1}&\alpha_{n,2}&\ldots&\alpha_{n,n-1}& \alpha_{n,n}
\alpha_{n,0} &\alpha_{n,1}&\alpha_{n,2}&\ldots&\alpha_{n,n}&
\end{array}
$$
Then we only need to notice that the first double sum corresponds to adding the elements by horizontal lines, from top to bottom, the second sum corresponds to adding the elements by diagonal lines, starting from the ``hypotenuse'', and the third one corresponds to adding the elements by vertical lines, from left to right.}
%)
%\medskip

%\item
%Prove that, being $F$ the operator introduced in Theorem \ref{operF-Delta}, the following property holds:
%$$F({\bf uv})=F({\bf u})F({\bf v})\;,\quad\forall {\bf u},{\bf v}\in\mathcal{P}'\,.$$
%\medskip
%
%\item
%Colocar um exercicio com v.p. Cauchy
%\medskip

%\item
%Prove that the following defines an equivalence relation in $\mathcal{P}'\,$:
%$$
%{\bf u}\sim{\bf v}\qquad\mbox{iff}\qquad
%\exists (a,b)\in\mathbb{C}\setminus\{0\}\times\mathbb{C}\;:\quad{\bf u}=\big(h_{a^{-1}}\circ\tau_{-b}\big){\bf v}\;.
%$$
%\medskip

%\item
%Let ${\bf u}^{(\alpha,\beta)}:\mathcal{P}\to\mathbb{C}$, with $\alpha,\beta>-1$, be defined by
%$$
%\langle {\bf u}^{(\alpha,\beta)},p\rangle:=\int_{-1}^1p(x)(1-x)^\alpha(1+x)^\beta\,{\rm d}x\; ,\quad p\in\mathcal{P}\;.
%$$
%\begin{enumerate}
%\item Show that there exist polynomials $\phi\in\mathcal{P}_2$ and $\psi\in\mathcal{P}_1$ such that ${\bf u}^{(\alpha,\beta)}$ fulfils the (distributional) differential equation
%    $$ D(\phi{\bf u}^{(\alpha,\beta)})=\psi{\bf u}^{(\alpha,\beta)}\;, $$
%    and determine explicitly $\phi$ and $\psi$.
%\item Let ${\bf u}$ and ${\bf v}$ be defined by
%$$
%{\bf u}:={\bf u}^{(-1/2,3/2)}-\pi\frac{1+2a}{1+a}\bm{\delta}_1\;,\quad
%{\bf v}:={\bf u}^{(1/2,-1/2)}-\frac{\pi}{2}\frac{1+2a}{1+a}\bm{\delta}_{-1}\;.
%$$
%Prove that $(1-x){\bf u}=(1+x^2){\bf v}$.
%\end{enumerate}
%\medskip

\end{enumerate}
%\bigskip
}
\medskip

\section*{Final remarks}

The theoretical foundations on the (so called) algebraic theory of OP are contained in the articles 
\ref{Maroni1985i}, \ref{Maroni1988i}, and \ref{Maroni1991i} by Pascal Maroni,
which were our source references for writing this text/chapter.
The basic facts about LCS needed to understanding this text are contained e.g. in the book \ref{ReedSimon1972i} by Michael Reed and Barry Simon, that we have used also for writing Appendix \ref{Appx1} on LCS.
An alternative/complementar reference is the book \label{Treves1967i} by Fran\c cois Tr\`eves \ref{Treves1967i}.
Exercises {\bf 1} and {\bf 2} can be found in several books on Functional Analysis.
Exercises {\bf 3} up to {\bf 6} were elaborated using as source, essentially, the articles by Maroni mentioned above
(where some proofs of the results presented therein were not included).
\medskip

%\section*{Notas finais}
%Os fundamentos te\'oricos sobre a (chamada) teoria alg\'ebrica dos PO est\~ao contidos
%nos artigos \ref{Maroni1985i}, \ref{Maroni1988i} e \ref{Maroni1991i} de Pascal Maroni,
%os quais serviram de base \`a elabora\c c\~ao deste texto/cap\'\i tulo.
% %(reflectindo tamb\'em, naturalmente, o ponto de vista do autor deste relat\'orio).
%Os pr\'e-requisitos sobre ELC necess\'arios \`a compreens\~ao do texto est\~ao
%contidos, e.g., no livro \ref{ReedSimon1972i} de Michael Reed e Barry Simon, no qual
%base\'amos o Ap\^endice sobre ELC contido nas notas de curso disponibilizadas aos estudantes.
%Um texto alternativo/complementar \'e o livro \label{Treves1967i} de Fran\c cois Tr\`eves \ref{Treves1967i}.
%Os exerc\'\i cios {\bf 1} e {\bf 2} podem encontrar-se em livros de An\'alise Funcional.
%Os exerc\'\i cios {\bf 3} a {\bf 6} foram elaborados, essencialmente, a partir dos artigos de Maroni referidos
%(nos quais muitas das provas dos resultados apresentados n\~ao s\~ao dadas).

\section*{Bibliography}
%\medskip

{\small
\begin{enumerate}[label=\emph{\rm [\arabic*]},leftmargin=*]
\item\label{Lax2002i} Peter D. Lax, {\sl Functional Analysis}, John Wiley $\&$ Sons (2002).
\item\label{Maroni1985i} P. Maroni, {\it Sur quelques espaces de distributions qui sont des formes lin\'eaires sur l'espace
vectoriel des polyn\^omes}, In C. Brezinski et al. Eds., Simposium Laguerre, Bar-le-Duc, Lecture
Notes in Math. {\bf 1171}, Springer-Verlag (1985) 184--194.
%\item\label{Maroni1987i} P. Maroni, {\it Prol\'egom\`enes \`a l'\'etude des polyn\^omes orthogonaux semiclassiques}, Ann. Mat. Pura Appl. {\bf 149} (4) (1987) 165--184.
\item\label{Maroni1988i} P. Maroni, {\it Le calcul des formes lin\'eaires et les polyn\^omes orthogonaux semiclassiques}, In M. Alfaro et al. Eds., Orthogonal Polynomials and Their Applications, Lecture Notes in Math.
{\bf 1329}, Springer-Verlag (1988) 279--290.
\item\label{Maroni1991i} P. Maroni, {\it Une th\'eorie alg\'ebrique des polyn\^omes orthogonaux. Applications aux polyn\^omes
orthogonaux semiclassiques}, In C. Brezinski et al. Eds., Orthogonal Polynomials and Their Applications, Proc. Erice 1990, IMACS, Ann. Comp. App. Math. {\bf 9} (1991) 95--130.
\item\label{ReedSimon1972i} M. Reed and B. Simon, {\sl Methods of Modern Mathematical Physics I: Functional Analysis},
Academic Press (1972).
\item\label{Treves1967i} F. Tr\`eves, {\sl Topological Vector Spaces, Distributions and Kernels}, Academic Press (1967).
\end{enumerate}
}

\chapter{Basic theory of orthogonal polynomials}

\pagestyle{myheadings}\markright{Basic theory of orthogonal polynomials}
\pagestyle{myheadings}\markleft{J. Petronilho}

%Neste texto introduzem-se as no\c c\~oes fundamentais e as propriedades alg\'ebricas b\'asicas
%da teoria das sucess\~oes de polin\'omios ortogonais.

%In this text we introduce the fundamental concepts and basic algebraic properties 
%of the theory of orthogonal polynomial sequences.

\section{Orthogonal polynomial sequences}

\begin{snugshade}
\begin{definition}\label{def-OPS}
Let ${\bf u}\in\mathcal{P}'$ and $\{P_n\}_{n\geq0}$ a sequence in $\mathcal{P}$.
\begin{enumerate}
\item[{\rm (i)}]
$\{P_n\}_{n\geq0}$ is called an {\sl orthogonal polynomial sequence (OPS)} with respect to ${\bf u}$
if $\{P_n\}_{n\geq0}$ is a simple set (so that $\deg P_n=n$ for all $n$)
and there exists a sequence $\{h_n\}_{n\geq0}$, with $h_n\in\mathbb{C}\setminus\{0\}$, such that
$$
%\begin{equation}\label{knOPS}
\langle {\bf u},P_mP_n \rangle=h_n\delta_{m,n}\;,\quad m,n=0,1,2,\ldots\;;
%\end{equation}
$$
\item[{\rm (ii)}]
${\bf u}$ is called {\sl regular} (or {\sl quasi-definite})
if there exists an OPS with respect to ${\bf u}$.
\end{enumerate}
\end{definition}
\end{snugshade}

\begin{remark}
\em
Whenever ${\bf u}$ is regular and $\{P_n\}_{n\geq0}$ is an OPS with respect to ${\bf u}$,
we will use such phrases as
``$\{P_n\}_{n\geq0}$ is an OPS for ${\bf u}$'', or
``$\{P_n\}_{n\geq0}$ is an OPS associated with ${\bf u}$'', or
``${\bf u}$ is regular and $\{P_n\}_{n\geq0}$ an associated OPS'', etc.
%or ``${\bf u}$ is regular and $\{P_n\}_{n\geq0}$ the corresponding monic OPS'', etc.
\end{remark}

Next we notice that if ${\bf u}$ is regular and $\{P_n\}_{n\geq0}$ an associated OPS,
then every polynomial admits a Fourier-type expansion in terms of a finite subset of $\{P_n\}_{n\geq0}$.

\begin{snugshade}
\begin{theorem}\label{Fourier-coeff-OPS}
Let ${\bf u}\in\mathcal{P}'$, regular, and $\{P_n\}_{n\geq0}$ an OPS with  respect to ${\bf u}$.
Let $\pi_k$ be a polynomial of degree $k$.
Then,
$$
\pi_k(x)=\sum_{j=0}^kc_{k,j}P_j(x)\;,\quad c_{k,j}
:=\frac{\langle{\bf u},\pi_kP_j\rangle}{\langle{\bf u},P_j^2\rangle}\; .
$$
\end{theorem}
\end{snugshade}

{\it Proof.}
Since $\{P_j\}_{j\geq0}$ is a simple set in $\mathcal{P}$,
then $\{P_j\}_{j=0}^k$ is an algebraic basis in $\mathcal{P}_k$.
Therefore, since $\pi_k\in\mathcal{P}_k$, there exist complex numbers $c_{k,j}$ such that
$$
\pi_k(x)=\sum_{j=0}^kc_{k,j}P_j(x)\;.
$$
Multiplying both sides of this equality by $P_\ell$, being $\ell$ fixed, $0\leq\ell\leq k$, and then
taking the action of the functional ${\bf u}$ in both sides of the resulting equality, we deduce
$$
\langle{\bf u},\pi_kP_\ell\rangle
=\sum_{j=0}^kc_{k,j}\langle{\bf u},P_jP_\ell\rangle=c_{k,\ell}\langle{\bf u},P_\ell^2\rangle\;,
$$
hence the desired result follows.
\qed

%The following result gives some useful characterizations of OPS.

\begin{snugshade}
\begin{theorem}\label{charcOPS}
Let ${\bf u}\in\mathcal{P}'$ and let $\{P_n\}_{n\geq0}$ be a simple set in $\mathcal{P}$.
Then the following are equivalent:
\begin{enumerate}
\item[{\rm (i)}] $\{P_n\}_{n\geq0}$ is an OPS with respect to ${\bf u}$;
\item[{\rm (ii)}] for each $n\in\mathbb{N}_0$ and $\pi\in\mathcal{P}_n\setminus\{0\}$,
there is $h_n=h_n(\pi)\in\mathbb{C}\setminus\{0\}$, such that
$$\langle {\bf u},\pi P_n \rangle=h_n\delta_{m,n}\;,\quad m:=\deg\pi\;.$$
%$$
%\langle {\bf u},\pi P_n \rangle
%\left\{
%\begin{array}{lcl}
%=0 & \mbox{\rm if} &  \deg\pi<n\; , \\ [0.25em]
%\neq0 & \mbox{\rm if} &  \deg\pi=n\; ;
%\end{array}
%\right.
%$$
\item[{\rm (iii)}] for each $n\in\mathbb{N}_0$, there exists $h_n\in\mathbb{C}\setminus\{0\}$ such that
$$\langle {\bf u},x^m P_n \rangle=h_n\delta_{m,n}\;,\quad m=0,1,\ldots,n\;.$$
\end{enumerate}
\end{theorem}
\end{snugshade}

{\it Proof.}
Assume that (i) holds. Fix $n\in\mathbb{N}_0$ and let $\pi\in\mathcal{P}_n$.
Setting $m:=\deg\pi$, from Theorem \ref{Fourier-coeff-OPS}
we know that there exists complex numbers $c_{m,j}$ such that
$$
\pi(x)=\sum_{j=0}^mc_{m,j}P_j(x)\; .
$$
Clearly, $c_{m,m}\neq0$ (since $\pi\not\equiv0$, $\deg\pi=m$,
and $\{P_j\}_{j\geq0}$ is a simple set in $\mathcal{P}$). Thus,
$$
\langle {\bf u},\pi P_n \rangle
=\sum_{j=0}^mc_{m,j}\langle {\bf u},P_j P_n \rangle
=\left\{
\begin{array}{lcl}
0 & \mbox{\rm if} &  m<n\; , \\ [0.25em]
c_{n,n}\langle {\bf u},P_n^2 \rangle\neq0 & \mbox{\rm if} &  m=n\; ,
\end{array}
\right.
$$
hence (i)$\Rightarrow$(ii), being $h_n:=c_{n,n}\langle {\bf u},P_n^2\rangle$.
Taking $\pi(x):=x^m$ in (ii), it is clear that (ii)$\Rightarrow$(iii).
Finally, assume that (iii) holds.
Fix $j,n\in\mathbb{N}_0$ and, without loss of generality, assume that $j\leq n$.
Since $\{P_k\}_{k\geq0}$ is a simple set, there exists complex numbers $c_{j,m}$, with $c_{j,j}\neq0$, such that
$\,P_j(x)=\sum_{m=0}^jc_{j,m}x^m\,$. Therefore, we deduce
$$
\langle {\bf u},P_j P_n \rangle
=\sum_{m=0}^jc_{j,m} \langle {\bf u},x^m P_n\rangle
=\sum_{m=0}^jc_{j,m} h_n\delta_{m,n}
=\widetilde{h}_n\delta_{j,n}\;,
$$
where $\widetilde{h}_n:=c_{n,n}h_n\neq0$.
Thus (iii)$\Rightarrow$(i), which completes the proof.
\qed
\medskip

The next proposition states that, up to normalization, there exists only one OPS
associated with a given regular functional.

\begin{snugshade}
\begin{theorem}\label{uniqueOPS}
Let ${\bf u}\in\mathcal{P}'$ be regular, and let
$\{P_n\}_{n\geq0}$ and $\{Q_n\}_{n\geq0}$ be two OPS with  respect to ${\bf u}$.
Then, there exists a sequence $\{c_n\}_{n\geq0}$, with $c_n\in\mathbb{C}\setminus\{0\}$, such that
$$
Q_n(x)=c_nP_n(x)\;,\quad n=0,1,2,\cdots\;.
$$
\end{theorem}
\end{snugshade}

{\it Proof.}
Fix $k\in\mathbb{N}$. Since $\{Q_n\}_{n\geq0}$ is an OPS with respect to ${\bf u}$, then
$$
\langle {\bf u},Q_kP_j \rangle=0\quad\mbox{\rm if}\quad j<k\;.
$$
Thus, by Theorem \ref{Fourier-coeff-OPS}, taking $\pi_k(x)=Q_k(x)$, we obtain
$Q_k(x)=c_{k,k}P_k(x)$, being $c_{k,k}=\langle {\bf u},P_kQ_k \rangle/\langle {\bf u},P_k^2 \rangle\neq0$,
which concludes the proof, by taking $c_k:=c_{k,k}$.
\qed
\medskip

Theorem \ref{uniqueOPS} implies that an OPS is uniquely determined if it satisfies
a condition that fixes the leading coefficient of each $P_n$ (i.e., the coefficient of $x^n$).
In particular, if $\{P_n\}_{n\geq0}$ is an OPS and the leading coefficient of each $P_n$ is $1$,
we say that $\{P_n\}_{n\geq0}$ is the {\sl monic} OPS (with respect to ${\bf u}$).
In general, being $\{P_n\}_{n\geq0}$ an OPS with respect to ${\bf u}$ (not necessarily monic),
and being $k_n$ the leading coefficient of $P_n$, so that
$$P_n(x)=k_nx^n+\mbox{\rm lower degree terms}\;,$$
the corresponding monic OPS is $\{\widehat{P}_n\}_{n\geq0}$, where\vspace*{-0.5em}
\begin{snugshade}
$$
\widehat{P}_n(x):=k_n^{-1}P_n(x)\; .
$$
\end{snugshade}
\noindent
On the other hand, if $\{P_n\}_{n\geq0}$ is an OPS with respect to ${\bf u}$ and
$$
\langle {\bf u},P_n^2 \rangle=1\quad (n=0,1,2,\ldots)\;,
$$
we say that $\{P_n\}_{n\geq0}$ is an {\sl orthonormal} polynomial sequence (with respect to ${\bf u}$).
In general, being $\{P_n\}_{n\geq0}$ an OPS with respect to ${\bf u}$ (not necessarily orthonormal),
the sequence $\{p_n\}_{n\geq0}$, where\vspace*{-0.5em}
\begin{snugshade}
$$
p_n(x):=\langle{\bf u},P_n^2\rangle^{-1/2}P_n(x)\; ,
$$
\end{snugshade}
\noindent
is orthonormal with respect to ${\bf u}$.
Here the square root needs not be real, but, as noticed above, $p_n(x)$ may be uniquely determined by requiring
an additional condition on its leading coefficient (e.g., that its leading coefficient be positive).

Finally, we notice the following obvious fact:
{\it if $\{P_n\}_{n\geq0}$ is an OPS with respect to the functional ${\bf u}\in\mathcal{P}'$, then
$\{P_n\}_{n\geq0}$ is also an OPS with respect to the functional ${\bf v}:=c{\bf u}$,
for every constant $c\in\mathbb{C}\setminus\{0\}$}.

\section{Existence of OPS}

In this section we analyze the question of wether a given functional ${\bf u}\in\mathcal{P}'$ is regular,
i.e., we ask for necessary and sufficient conditions
that guarantee the existence of an OPS with respect to ${\bf u}$.
To answer this question, we introduce the so called Hankel determinants.
Denoting, as usual, by $u_j:=\langle{\bf u},x^j\rangle$, $j\in\mathbb{N}_0$, the moments of ${\bf u}$,
we define the associated {\sl Hankel determinant} $H_n\equiv H_n({\bf u})$ as
%\vspace*{-0.5em}
\begin{snugshade}
\begin{equation}\label{def-Hn}
H_n:=\det\big\{[u_{i+j}]_{i,j=0}^n\big\}=
\left|
\begin{array}{ccccc}
u_0 & u_1 & \cdots & u_{n-1} & u_n \\ [0.25em]
u_1 & u_2 & \cdots & u_{n} & u_{n+1} \\ [0.25em]
\vdots&\vdots&\ddots&\vdots&\vdots \\ [0.25em]
u_{n-1} & u_n & \cdots & u_{2n-2} & u_{2n-1} \\ [0.25em]
u_n & u_{n+1} & \cdots & u_{2n-1} & u_{2n}
\end{array}
\right|\;,\quad n\in\mathbb{N}_0\;.
\end{equation}
\end{snugshade}
\noindent
Notice that $H_n$ is a determinant of order $n+1$. It is also useful to set
\vspace*{-0.5em}
\begin{snugshade}
\begin{equation}\label{def-H-1}
H_{-1}:=1\;.
\end{equation}
\end{snugshade}

\begin{snugshade}
\begin{theorem}\label{exist-OPS}
Let ${\bf u}\in\mathcal{P}'$. Then,
${\bf u}$ is regular if and only if
\begin{equation}\label{Hn-not-zero}
H_n\neq0\;,\quad \forall n\in\mathbb{N}_0\,.
\end{equation}
Under such conditions, the monic OPS $\{P_n\}_{n\geq0}$ with respect to ${\bf u}$ is given by $P_0(x)=1$ and
\begin{equation}\label{Hn-Pn}
%P_0(x)=1\;,\quad
P_n(x)=\frac{1}{H_{n-1}}\left|
\begin{array}{ccccc}
u_0 & u_1 & \cdots & u_{n-1} & u_n \\ [0.25em]
u_1 & u_2 & \cdots & u_{n} & u_{n+1} \\ [0.25em]
\vdots&\vdots&\ddots&\vdots&\vdots \\ [0.25em]
u_{n-1} & u_n & \cdots & u_{2n-2} & u_{2n-1} \\ [0.25em]
1 & x & \cdots & x^{n-1} & x^n
\end{array}
\right|\;,\quad n\in\mathbb{N}\;.
\end{equation}
\end{theorem}
\end{snugshade}

{\it Proof.}
Suppose that ${\bf u}$ is regular.
Let $\{P_n\}_{n\geq0}$ be an OPS with respect to ${\bf u}$.
Fix $n\in\mathbb{N}_0$.
Then, there exists $c_{n,0},c_{n,1},\ldots,c_{n,n}\in\mathbb{C}$ such that
\begin{equation}\label{Hn-Pn0a}
P_n(x)=\sum_{k=0}^nc_{n,k}x^k\;.
\end{equation}
By Theorem \ref{charcOPS}, there exists $h_n\in\mathbb{C}\setminus\{0\}$ such that
\begin{equation}\label{Hn-Pn0aA}
h_n\delta_{m,n}=\langle {\bf u},x^m P_n \rangle
=\sum_{k=0}^nc_{n,k}u_{k+m}\;,\quad m=0,1,\ldots,n\;.
\end{equation}
This may be written in matrix form as
\begin{equation}\label{kn-system}
\left[
\begin{array}{ccccc}
u_0 & u_1 & \cdots & u_{n-1} & u_n  \\
u_1 & u_2 & \cdots & u_{n} & u_{n+1}  \\
\vdots&\vdots&\ddots&\vdots& \vdots \\
u_{n-1} & u_{n} & \cdots & u_{2n-2} & u_{2n-1} \\
u_{n} & u_{n+1} & \cdots & u_{2n-1} & u_{2n}
\end{array}
\right]
\left[
\begin{array}{c}
c_{n,0} \\ c_{n,1} \\ \vdots \\ c_{n,n-1} \\ c_{n,n}
\end{array}
\right]=
\left[
\begin{array}{c}
0 \\ 0 \\ \vdots \\ 0 \\ h_n
\end{array}
\right]\;.
\end{equation}
Since the sequence $\{h_n\}_{n\geq0}$ in (\ref{Hn-Pn0aA}) uniquely determines
the OPS $\{P_n\}_{n\geq0}$\footnote{\,Indeed, if $\{Q_n\}_{n\geq0}$ is also an OPS with respect to ${\bf u}$ such
that (\ref{Hn-Pn0aA}) holds with $Q_n$ instead of $P_n$, with the same $h_n$, then,
since, by Theorem \ref{uniqueOPS}, $Q_n(x)=c_nP_n(x)$ for some complex number $c_n\neq0$, one would have
$h_n=\langle {\bf u},x^n Q_n \rangle=c_n\langle{\bf u},x^n P_n \rangle=c_nh_n$, hence $c_n=1$, so that $Q_n=P_n$.},
then this system (\ref{kn-system}), where the coefficients $c_{n,0}, c_{n,1},\ldots,c_{n,n}$ of $P_n$ are the unknowns,
has a unique solution.
Hence $H_n\neq0$, because $H_n$ is the determinant of such system.

Conversely, suppose that $H_n\neq0$ for all $n\in\mathbb{N}_0$.
Then, for any fixed $n\in\mathbb{N}_0$, to each constant $h_n\in\mathbb{C}\setminus\{0\}$
corresponds a unique vector $(c_{n,0}, c_{n,1},\ldots,c_{n,n})$, solution of the system (\ref{kn-system}).
Using the components of this vector, we may define a polynomial $P_n(x)$ by expression (\ref{Hn-Pn0a}).
This polynomial fulfils (\ref{Hn-Pn0aA}), since (\ref{Hn-Pn0aA}) and (\ref{kn-system}) are equivalent.
To conclude that $\{P_n\}_{n\geq0}$ is an OPS with respect to ${\bf u}$, it remains to prove that
it is a simple set, i.e., $\deg P_n=n$ for all $n$.
Indeed, solving (\ref{kn-system}) for $c_{n,n}$ by Crammer's rule,
and taking into account the hypothesis $H_n\neq0$, we obtain
\begin{snugshade}
\begin{equation}\label{kn-hn-Hn}
c_{n,n}=\frac{h_n H_{n-1}}{H_{n}}\; ,\quad n\in\mathbb{N}_0\;,
\end{equation}
\end{snugshade}
\noindent
hence $c_{n,n}\neq0$, which proves that, indeed, $\deg P_n=n$ for all $n$.

It remains to prove (\ref{Hn-Pn}). We will present two proofs.
The first one is a constructive proof. The second one is much more concise.
%\smallskip

\textit{First proof of} (\ref{Hn-Pn}).
Let $\{P_n\}_{n\geq0}$ be the monic OPS with respect to ${\bf u}$.
For each fixed $n\in\mathbb{N}$, $P_n(x)$ may be written as in (\ref{Hn-Pn0a}),
%\begin{equation}\label{Hn-Pn0}
%P_n(x)=\sum_{k=0}^nc_kx^k\;,
%\end{equation}
being $c_{n,0},c_{n,1},\ldots,c_{n,n-1}\in\mathbb{C}$, and $c_{n,n}=1$.
As above, for each $m\in\{0,1,\cdots,n-1\}$, we deduce
$0=\langle{\bf u}, x^mP_n\rangle=\sum_{k=0}^nc_{n,k}u_{k+m}$.
From this, and taking into account that $c_{n,n}=1$,
we obtain the following system of $n$ equations in the $n$
unknowns $c_{n,0},c_{n,1},\ldots,c_{n,n-1}$:
%$$
%\sum_{k=0}^{n-1}u_{k+m}c_{n,k}=-u_{n+m}\,,\quad m=0,1,\cdots,n-1\;,
%$$
%i.e., in matrix notation,
$$
\left[
\begin{array}{cccc}
u_0 & u_1 & \cdots & u_{n-1}  \\
u_1 & u_2 & \cdots & u_{n}  \\
\vdots&\vdots&\ddots&\vdots \\
u_{n-1} & u_n & \cdots & u_{2n-2}
\end{array}
\right]
\left[
\begin{array}{c}
c_{n,0} \\ c_{n,1} \\ \vdots \\ c_{n,n-1}
\end{array}
\right]=
\left[
\begin{array}{c}
-u_n \\ -u_{n+1} \\ \vdots \\ -u_{2n-1}
\end{array}
\right]\;.
$$
The determinant of this system is $H_{n-1}\neq0$.
Solving by Crammer's rule, we obtain
$$
c_{n,k}=\frac{1}{H_{n-1}}\,
\left|
\begin{array}{ccccccc}
u_0 & \cdots & u_{k-1} & -u_n & u_{k+1} & \cdots & u_{n-1} \\ [0.25em]
u_1 & \cdots & u_{k} & -u_{n+1} & u_{k+2} & \cdots & u_{n} \\ [0.25em]
\vdots&\ddots&\vdots&\vdots&\vdots&\ddots&\vdots \\ [0.25em]
u_{n-1} & \cdots & u_{n+k-2} & -u_{2n-1} & u_{n+k} & \cdots & u_{2n-2} \\ [0.25em]
\end{array}
\right| %\;,\quad  k=0,1,\cdots,n-1\;.
$$
for each $k=0,1,\cdots,n-1$.
Performing elementary operations on the columns of this determinant, by
moving successively the $(k+1)$th column to its right
(so that $n-k-1$ permutations on columns must be done), we deduce %the last determinant may be written as
\begin{equation}\label{Hn-Pn1}
c_{n,k}=\frac{(-1)^{n-k}}{H_{n-1}}\,
\left|
\begin{array}{ccccccc}
u_0 & \cdots & u_{k-1} & u_{k+1} & \cdots & u_{n-1} & u_n \\ [0.25em]
u_1 & \cdots & u_{k} & u_{k+2} & \cdots & u_{n} & u_{n+1} \\ [0.25em]
\vdots&\ddots&\vdots&\vdots&\ddots&\vdots&\vdots \\ [0.25em]
u_{n-1} & \cdots & u_{n+k-2} & u_{n+k} & \cdots & u_{2n-2} & u_{2n-1} \\ [0.25em]
\end{array}
\right| %\;,\quad  k=0,1,\cdots,n-1\;.
\end{equation}
for each $k=0,1,\cdots,n-1$.
Clearly, (\ref{Hn-Pn1}) is also true for $k=n$, since in that case
the right-hand side of (\ref{Hn-Pn1}) reduces to $1$.
Therefore, substituting (\ref{Hn-Pn1}) into (\ref{Hn-Pn0a}), we obtain
$$
P_n(x)=\frac{1}{H_{n-1}}\sum_{k=0}^n(-1)^{n-k}
\left|
\begin{array}{ccccccc}
u_0 & \cdots & u_{k-1} & u_{k+1} & \cdots & u_{n-1} & u_n \\ [0.25em]
u_1 & \cdots & u_{k} & u_{k+2} & \cdots & u_{n} & u_{n+1} \\ [0.25em]
\vdots&\ddots&\vdots&\vdots&\ddots&\vdots&\vdots \\ [0.25em]
u_{n-1} & \cdots & u_{n+k-2} & u_{n+k} & \cdots & u_{2n-2} & u_{2n-1} \\ [0.25em]
\end{array}
\right|\,x^k\;,
$$
hence formula (\ref{Hn-Pn}) follows by Laplace's Theorem,
developing the determinant in the right-hand side of (\ref{Hn-Pn}) along its last row.
%\smallskip

\textit{Second proof of} (\ref{Hn-Pn}).
%Formula (\ref{Hn-Pn}) can be easily proved by noticing that
Let $Q_n(x)$ be the (monic) polynomial of degree $n$ defined by the right-hand side of (\ref{Hn-Pn}).
If $m<n$ then, clearly, $\langle{\bf u},x^mQ_n\rangle=0$
(since $\langle{\bf u},x^mQ_n\rangle$ becomes
a determinant whose $m+1$ row and $n+1$ row are equal).
If $m=n$, then we simply notice that
$\langle{\bf u},x^nQ_n\rangle=H_n/H_{n-1}\neq0$.
Thus, by Theorem \ref{charcOPS}, $\{Q_n\}_{n\geq0}$ is an OPS with respect to ${\bf u}$,
and since each $Q_n$ is a monic polynomial,
one should conclude that (\ref{Hn-Pn}) holds.
%$P_n(x)=Q_n(x)$ for each $n\in\mathbb{N}_0$.
%This completes the proof.
\qed

\begin{snugshade}
\begin{corollary}\label{leading-coeff-Hn}
Let $\{P_n\}_{n\geq0}$ be an OPS with respect to ${\bf u}\in\mathcal{P}'$, and let
$\pi_n$ be a polynomial of degree $n$.
Denote by $k_n$ and $a_n$ the leading coefficients of $P_n$ and $\pi_n$, respectively, so that
$$
\begin{array}{l}
P_n(x)=k_nx^n+\mbox{\rm lower degree terms}\,, \\ [0.25em]
\pi_n(x)=a_nx^n+\mbox{\rm lower degree terms}
\end{array}
$$
for each $n\in\mathbb{N}_0$. Then
\begin{equation}\label{Hn-not-zero}
\langle{\bf u},\pi_nP_n\rangle=a_n\langle{\bf u},x^nP_n\rangle
=\frac{a_nk_nH_n}{H_{n-1}}\,,\quad n\in\mathbb{N}_0\;.
\end{equation}
\end{corollary}
\end{snugshade}

{\it Proof.}
Writing $\pi_n(x)=a_nx^n+\pi_{n-1}(x)$, with $\pi_{n-1}\in\mathcal{P}_{n-1}$,
and taking into account Theorem \ref{charcOPS}, we deduce, for each $n\in\mathbb{N}_0$,
$$
\langle{\bf u},\pi_nP_n\rangle=a_n\langle{\bf u},x^nP_n\rangle+\langle{\bf u},\pi_{n-1}P_n\rangle
=a_n\langle{\bf u},x^nP_n\rangle=a_nh_n=\frac{a_nk_nH_n}{H_{n-1}}\;,
$$
where the last equality follows from (\ref{kn-hn-Hn}),
noticing that $c_{n,n}=k_n$.
\qed

%(\ref{Hn-Pn0aA})
%c_{n,n}=\frac{h_n H_{n-1}}{H_{n}}

\section{OPS in the positive-definite sense}

In many important occurrences of OP, the functional ${\bf u}\in\mathcal{P}^*$
with respect to which the polynomials are orthogonal admits an integral representation
involving a weight function, or, in the most general situation,
a positive Borel measure, $\mu$, whose support is an infinite subset of $\mathbb{R}$,
and with finite moments of all orders, so that
\begin{equation}\label{up+1}
\langle{\bf u},p\rangle=\int_\mathbb{R}p(x)\,{\rm d}\mu(x)\;,\quad p\in\mathcal{P}\;.
\end{equation}
One easily verifies that, under such conditions, the property
\begin{equation}\label{up+2}
\langle{\bf u},p\rangle>0
\end{equation}
holds for each polynomial $p\in\mathcal{P}$ which is nonzero (i.e., it doesn't vanishes identically)
and nonnegative for all $x\in\mathbb{R}$.
It turns out that this property characterizes functionals
${\bf u}\in\mathcal{P}'$ such that an integral representation as (\ref{up+1}) holds,
under the conditions described above.
This ``equivalence'' between (\ref{up+1}) and (\ref{up+2}) is a nontrivial fact,
and it will be proved latter.
We start the study of such functionals by introducing the following definition.

\begin{snugshade}
\begin{definition}\label{def-OPS-posit}
A functional ${\bf u}\in\mathcal{P}'$ is called {\sl positive-definite} if the condition
\begin{equation}\label{u-def+}
\langle{\bf u},p\rangle>0
\end{equation}
holds for each polynomial $p$ which is nonzero  and nonnegative for all real $x$.
\end{definition}
\end{snugshade}

Next we state some basic properties of positive-definite linear functionals in $\mathcal{P}$.

\begin{snugshade}
\begin{theorem}\label{def-OPS-positive}
Let ${\bf u}\in\mathcal{P}'$ be positive-definite. Then, the moments
$u_n:=\langle{\bf u},x^n\rangle$ are real numbers.
More precisely, the following holds:
\begin{equation}\label{un+}
u_{2n}>0\;,\quad u_{2n+1}\in\mathbb{R}\;,\quad \forall n\in\mathbb{N}_0\;.
\end{equation}
\end{theorem}
\end{snugshade}

{\it Proof.}
On the first hand, since ${\bf u}$ is positive-definite and $x^{2n}\geq0$ for $x\in\mathbb{R}$, then
$$
u_{2n}=\langle{\bf u},x^{2n}\rangle>0\;.
$$
On the other hand, using again the positive-definiteness of ${\bf u}$ and Newton's
binomial formula, we may write
$$
0<\langle{\bf u},(1+x)^{2n}\rangle=\sum_{k=0}^{2n}\binom{2n}{k}u_{k}\;,
$$
hence it follows by induction that $u_{2n+1}$ is a real number.
\qed
\smallskip

Given a positive-definite functional defined in $\mathcal{P}$, a step-by-step method of constructing a corresponding
orthonormal polynomial sequence can be described, known as {\it Gram-Schmidt process}. This method produces
real orthonormal polynomials. %, described as follows.

\begin{snugshade}
\begin{theorem}[Gram-Schmidt process]\label{GramShmidt}
Let ${\bf u}\in\mathcal{P}'$ be positive-definite.
Define a sequence of polynomials $\{p_n\}_{n\geq0}$, constructed step-by-step, as follows:
\begin{equation}\label{GShm1}
%p_0(x):=u_0^{-1/2}\; ,\quad %\\ [0.5em]
p_{n}(x):=\langle{\bf u},P_{n}^2\rangle^{-1/2}P_{n}(x)\;,\quad n\in\mathbb{N}_0\;,
\end{equation}
where $\{P_n\}_{n\geq0}$ is a simple set of monic polynomials, constructed step-by-step as
%$P_0(x)=1$ and %, for each $n\in\mathbb{N}_0$ and $1\leq k\leq n$,
\begin{equation}\label{GShm2}
P_{n}(x):=x^{n}-\sum_{k=0}^{n-1} \langle{\bf u},x^{n}p_k\rangle\, p_k(x)\;,\quad n\in\mathbb{N}_0\,.
%P_{n+1}(x)=x^{n+1}-\sum_{k=0}^n \langle{\bf u},x^{n+1}p_k\rangle\, p_k(x)\;,\quad n\in\mathbb{N}_0\,.
\end{equation}
Then, $\{p_n\}_{n\geq0}$ is orthonormal with respect to ${\bf u}$,
being each $p_n(x)$ a real polynomial (i.e., with real coefficients).
Moreover, $\{P_n\}_{n\geq0}$ is the corresponding monic OPS,
being also each $P_n(x)$ a real polynomial.
%As a consequence, every positive-definite functional is regular.
\end{theorem}
\end{snugshade}

{\it Proof.}
To state the theorem we prove that each $P_n(x)$ is a real polynomial
(this implies $\langle{\bf u},P_{n}^2\rangle>0$, since ${\bf u}$ is positive-definite,
so $p_n(x)$ is also a real polynomial), and
\begin{equation}\label{GSproof}
\langle{\bf u},p_{n}^2\rangle=1\;,\quad
\langle{\bf u},p_jp_{n+1}\rangle=0\;,\quad n\in\mathbb{N}_0\;,\;\; 0\leq j\leq n\;.
\end{equation}
This will be proved by induction over $n$.
For $n=0$, we have
$$P_0(x)=1\;,\quad p_0(x)=u_0^{-1/2}\,,$$
hence $p_0(x)$ is real---notice that, by Theorem \ref{def-OPS-positive}, conditions (\ref{un+}) hold---and
$$
\langle{\bf u},p_{0}^2\rangle=u_0^{-1}\langle{\bf u},1\rangle=u_0^{-1}u_0=1\; .
$$
Now, we compute
$$P_1(x)=x-\langle{\bf u},xp_0\rangle p_0(x)=x-u_1/u_0\;,$$
hence $P_1(x)$ is a real polynomial, and since ${\bf u}$ is positive-definite,
we have $\langle{\bf u},P_1^2\rangle>0$.
Thus, $p_1(x):=\langle{\bf u},P_1^2\rangle^{-1/2}P_1(x)$ is well defined,
it is a real polynomial, and
$$
\langle{\bf u},p_0p_1\rangle =u_0^{-1/2}\langle{\bf u},P_1^2\rangle^{-1/2}\langle{\bf u},x-u_1/u_0\rangle=0\;,
$$
and we conclude that (\ref{GSproof}) holds for $n=0$.
Assume now (induction hypothesis) that, for some $m\in\mathbb{N}_0$,
the polynomials $P_1(x),\ldots,P_{m+1}(x)$ are real, and (\ref{GSproof}) holds
for all positive integers $n\leq m$ .
We need to prove that $P_{m+2}(x)$ is also a real polynomial and
(\ref{GSproof}) remains true if $n$ is replaced by $m+1$.
Indeed, since $P_{m+1}(x)$ is real and ${\bf u}$ is positive-definite,
then $\langle{\bf u},P_{m+1}^2\rangle>0$, and so
$$
\langle{\bf u},p_{m+1}^2\rangle
=\big\langle{\bf u},\langle{\bf u},P_{m+1}^2\rangle^{-1} P_{m+1}^2\big\rangle=1\;.
$$
Moreover, since, by the induction hypothesis, $P_1(x),\ldots,P_{m+1}(x)$ are real,
then so are $p_0(x),p_1(x),\ldots,p_{m+1}(x)$, hence so is $P_{m+2}(x)$.
Then, $\langle{\bf u},P_{m+2}^2\rangle>0$, and so
$p_{m+2}(x):=\langle{\bf u},P_{m+2}^2\rangle^{-1/2}P_{m+2}(x)$ is well defined.
Thus, for each $j\in\{0,1,\ldots,m+1\}$,
$$
\langle{\bf u},p_jp_{m+2}\rangle=\langle{\bf u},P_{m+2}^2\rangle^{-1/2}
\Big(\langle{\bf u},p_jx^{m+2}\rangle-\sum_{k=0}^{m+1}\langle{\bf u},x^{m+2}p_k\rangle\langle{\bf u}, p_jp_k \rangle\Big)\, .
$$
%$$
%\begin{array}{rcl}
%\langle{\bf u},p_jp_{n+2}\rangle
%&=& \displaystyle\big\langle{\bf u},P_{n+2}^2\big\rangle^{-1/2}
%\Big\langle{\bf u},p_jx^{n+2}-\sum_{k=0}^{n+1}\langle{\bf u},x^{n+2}p_k\rangle\, p_jp_k \Big\rangle \\ [1em]
%&=& \displaystyle\langle{\bf u},P_{n+2}^2\rangle^{-1/2}
%\Big(\langle{\bf u},p_jx^{n+2}\rangle-\sum_{k=0}^{n+1}\langle{\bf u},x^{n+2}p_k\rangle\langle{\bf u}, p_jp_k \rangle\Big)\, .
%\end{array}
%$$
Since, by the induction hypothesis,
$\langle{\bf u},p_jp_k\rangle=\delta_{j,k}$ if $j,k\in\{0,1,\ldots,m+1\}$, we deduce
$$
\langle{\bf u},p_jp_{m+2}\rangle=\langle{\bf u},P_{m+2}^2\rangle^{-1/2}
\big(\langle{\bf u},p_jx^{m+2}\rangle-\langle{\bf u},x^{m+2}p_j\rangle\big)=0\;,\quad
0\leq j\leq m+1\,.
$$
This completes the proof.
\qed

%The following proposition is an immediate consequence of Theorem \ref{GramShmidt}.

\begin{snugshade}
\begin{corollary}\label{corGS}
Let ${\bf u}\in\mathcal{P}'$ be positive-definite.
Then, ${\bf u}$ is regular.
\end{corollary}
\end{snugshade}

Next we state the connection between positive-definite functionals defined in $\mathcal{P}$
and the Hankel determinants introduced in (\ref{def-Hn}).
We will need the following classical result characterizing non-negative polynomials.

\begin{snugshade}
\begin{lemma}\label{poly+}
Let $\pi(x)$ be a polynomial that is non-negative for all real $x$.
Then, there are real polynomials $P(x)$ and $Q(x)$ such that
\begin{equation}\label{poly+1}
\pi(x)=P^2(x)+Q^2(x)\;.\end{equation}
\end{lemma}
\end{snugshade}

{\it Proof.}
Since $\pi(x)\geq0$ for $x\in\mathbb{R}$, then $\pi$ is a real polynomial
(i.e., its coefficients are all real numbers) such that its real zeros have even multiplicity
and its non-real zeros occur in conjugate pairs. Thus, we can write
$$
\pi(x)=R^2(x)\prod_{k=1}^m(x-a_k+ib_k)(x-a_k-ib_k)\;,
%=R^2(x)\prod_{k=1}^m(x-a_k+ib_k)\overline{\prod_{k=1}^m(x-a_k+ib_k)}\;,
$$
where $R$ is a real polynomial and $a_k,b_k$ real numbers.
Therefore, since we may write
$$
\prod_{k=1}^m(x-a_k+ib_k)=A(x)+iB(x)\;,
$$
where $A$ and $B$ are real polynomials, we deduce
$$
\pi(x)=R^2(x)[A(x)+iB(x)]\overline{[A(x)+iB(x)]}=R^2(x)\big[A^2(x)+B^2(x)\big]\;,
$$
hence the desired result follows by taking $P:=RA$ and $Q:=RB$.
\qed

\begin{snugshade}
\begin{theorem}\label{exist-OPS+}
Let ${\bf u}\in\mathcal{P}'$. Then, ${\bf u}$ is positive-definite
if and only if the following two conditions hold:
\begin{enumerate}
\item[{\rm (i)}] the moments $u_n:=\langle {\bf u},x^n\rangle$ are real for each $n\in\mathbb{N}_0$;
\item[{\rm (ii)}] the Hankel determinants $(\ref{def-Hn})$ are all positive:
\begin{equation}\label{Hn-positive}
H_n>0\;,\quad \forall n\in\mathbb{N}_0\,.
\end{equation}
\end{enumerate}
\end{theorem}
\end{snugshade}

{\it Proof.}
Suppose that ${\bf u}$ is positive-definite.
Then by Theorem \ref{def-OPS-positive}
all the moments $u_n$ are real.
Moreover, by Theorem \ref{GramShmidt}, a monic OPS $\{P_n\}_{n\geq0}$ with respect to ${\bf u}$ exists,
with each $P_n(x)$ a real polynomial, and so $\langle{\bf u},P_n^2\rangle>0$ for each $n\in\mathbb{N}_0$
(since ${\bf u}$ is positive-definite). Then, and taking into account Corollary \ref{leading-coeff-Hn},
we have
$$
0<\langle{\bf u},P_n^2\rangle=\frac{H_n}{H_{n-1}}\,,\quad n\in\mathbb{N}_0\;.
$$
Therefore, since $H_{-1}=1$, it follows by induction that $H_n>0$ for all $n\in\mathbb{N}_0$.

Conversely, suppose that conditions (i) and (ii) hold.
(ii) and Theorem \ref{exist-OPS} ensure that ${\bf u}$ is regular,
hence there exists a monic OPS $\{P_n\}_{n\geq0}$ with respect to ${\bf u}$.
Since $P_n(x)$ admits the representation (\ref{Hn-Pn}), it follows from (i) and (ii)
that each $P_n(x)$ is a real polynomial.
%Writing each $P_n(x)$ as in (\ref{Hn-Pn0a}), we see that the coefficients of $P_n(x)$
%are solutions of the matrix system (\ref{kn-system}), hence hypothesis (i) guarantees that
%such coefficients are real, i.e., each $P_n(x)$ is a real polynomial.
Also, again by Corollary \ref{leading-coeff-Hn} and by (ii), we have
$$
\langle{\bf u},P_n^2\rangle=\frac{H_n}{H_{n-1}}>0\,,\quad n\in\mathbb{N}_0\;.
$$
Let $Q(x)$ be a nonzero real polynomial of degree $m$.
Since each $P_n(x)$ is real, we may write
$\,Q(x)=\sum_{j=0}^ma_{j}P_j(x)\,$, where $a_j\in\mathbb{R}$ for all $j$, with $a_{m}\neq0$.
Therefore,
$$
\langle{\bf u},Q^2\rangle=\sum_{j=0}^ma_j^2\langle{\bf u},P_j^2\rangle>0\;.
$$
Thus, it follows from Lemma \ref{poly+} that ${\bf u}$ is positive-definite.
\qed

\begin{snugshade}
\begin{corollary}\label{Cor-exist-OPS+}
Let ${\bf u}\in\mathcal{P}'$. Suppose that ${\bf u}$ is regular
and let $\{P_n\}_{n\geq0}$ be the associated monic OPS.
Assume further that $P_n(x)$ is real for each $n\in\mathbb{N}_0$, and
\begin{equation}\label{eqPn+}%$$
\langle {\bf u},P_n^2\rangle>0\;,\quad \forall n\in\mathbb{N}_0\;.
\end{equation}%$$
Then, ${\bf u}$ is positive-definite.
\end{corollary}
\end{snugshade}

{\it Proof.}
The hypothesis allow us to proceed as in the last part of the proof of Theorem \ref{exist-OPS+},
in order to obtain $\langle{\bf u},Q^2\rangle>0$ for every
nonzero real polynomial $Q$, so that, by Lemma \ref{poly+}, ${\bf u}$ is positive-definite.
\qed
%The proof is immediate if we use the ideas in the last part of the proof of Theorem \ref{exist-OPS+},
%noticing that the hypothesis here allow us to obtain $\langle{\bf u},Q^2\rangle>0$ for every
%nonzero real polynomial $Q$, so that, by Lemma \ref{poly+}, ${\bf u}$ is positive-definite.
%Since ${\bf u}$ is regular, then $H_n\neq0$ for all $n\in\mathbb{N}_0$.
%By Corollary \ref{leading-coeff-Hn} and the hypothesis (\ref{eqPn+}), we have
%$$
%\frac{H_n}{H_{n-1}}=\langle{\bf u},P_n^2\rangle>0\,,\quad n\in\mathbb{N}_0\;.
%$$
%Therefore, since $H_{-1}=1$, it follows by induction that $H_n>0$ for all $n\in\mathbb{N}_0$.

%\section{Favard's Theorem and Christoffel-Darboux identities}
\section{Favard's Theorem}

One of the most important characterizations of OPS is the fact that any
three consecutive polynomials are connected by a very simple relation.
This is the content of Favard's Theorem.
We begin by stating the following proposition.

\begin{snugshade}
\begin{theorem}\label{thm-TTRR}
Let ${\bf u}\in\mathcal{P}'$, regular, and
$\{P_n\}_{n\geq0}$ the corresponding monic OPS.
Then, $\{P_n\}_{n\geq0}$ satisfies the three-term recurrence relation
\begin{equation}\label{TTRR}
P_{n+1}(x)=(x-\beta_n)P_n(x)-\gamma_nP_{n-1}(x)\;,\quad n=0,1,2,\cdots
\end{equation}
with initial conditions %$P_{-1}(x)=0$ and $P_{0}(x)=1$,
\begin{equation}\label{TTRRinitial}
P_{-1}(x)=0\;,\quad P_{0}(x)=1\;,
\end{equation}
where $\{\beta_n\}_{n\geq0}$ and $\{\gamma_n\}_{n\geq1}$ are
sequences of complex numbers such that
%two sequences of complex numbers fulfilling
\begin{equation}\label{TTRRbgn}
\gamma_n\neq0\;,\quad n=1,2,3,\ldots\;.
\end{equation}
Moreover, if ${\bf u}$ is positive-definite, then
\begin{equation}\label{TTRRbgn+}
\beta_{n-1}\in\mathbb{R}\;,\quad\gamma_n>0\;,\quad n=1,2,3,\ldots\;.
\end{equation}
\end{theorem}
\end{snugshade}

{\it Proof.}
Since $xP_n(x)$ is a polynomial of degree $n+1$ then, by Theorem \ref{Fourier-coeff-OPS},
$$
xP_n(x)=\sum_{j=0}^{n+1}c_{n,j}P_j(x)\;,\quad c_{n,j}
:=\frac{\langle{\bf u},xP_nP_j\rangle}{\langle{\bf u},P_j^2\rangle}\quad(0\leq j\leq n+1)\; .
$$
Clearly, $\langle{\bf u},xP_nP_j\rangle=0$ if $0\leq j\leq n-2$
(since $\{P_n\}_{n\geq0}$ is an OPS with respect to ${\bf u}$),
and $c_{n,n+1}=1$ (since each $P_j$ is a monic polynomial).
Hence, %(\ref{}) reduces to
$$
xP_n(x)=P_{n+1}(x)+c_{n,n}P_{n}(x)+c_{n,n-1}P_{n-1}(x)\;,\quad n=0,1,2,\cdots\;.
$$
Therefore, we obtain (\ref{TTRR}), with
$$%\begin{equation}\label{TTRRb}
\beta_n:=c_{n,n}=\frac{\langle{\bf u},xP_n^2\rangle}{\langle{\bf u},P_n^2\rangle}\quad (n=0,1,2,\cdots)
$$%\end{equation}
and
$$%\begin{equation}\label{TTRRg}
\gamma_n:=c_{n,n-1}=\frac{\langle{\bf u},xP_nP_{n-1}\rangle}{\langle{\bf u},P_{n-1}^2\rangle}
=\frac{\langle{\bf u},P_n^2\rangle}{\langle{\bf u},P_{n-1}^2\rangle}\in\mathbb{C}\setminus\{0\}\quad (n=1,2,\cdots)\;.
$$%\end{equation}
If ${\bf u}$ is positive-definite, then, by Theorem \ref{GramShmidt}, each $P_n(x)$ is a real polynomial.
Hence it follows from the previous expressions for $\beta_n$ and $\gamma_n$ %(\ref{TTRRb}) and (\ref{TTRRg})
that conditions (\ref{TTRRbgn+}) hold.
\qed

\begin{remark}\em
Since $P_{-1}(x)=0$, then it doesn't matter how to define $\gamma_0$.
Often we will make the useful choice $\gamma_0:=u_0$. % or $\gamma_0:=1$.
\end{remark}

\begin{snugshade}
\begin{corollary}\label{cor-TTRR}
Under the conditions of Theorem \ref{thm-TTRR}, the following holds:
\begin{enumerate}
\item[{\rm (i)}] the $\beta-$parameters are given by
\begin{equation}\label{TTRRb}
\beta_n=\frac{\langle{\bf u},xP_n^2\rangle}{\langle{\bf u},P_n^2\rangle}\;,\quad n=0,1,2,\ldots\;;
\end{equation}
\item[{\rm (ii)}] the $\gamma-$parameters are given by
\begin{equation}\label{TTRRg}
\gamma_n=\frac{\langle{\bf u},P_n^2\rangle}{\langle{\bf u},P_{n-1}^2\rangle}
=\frac{H_{n-2}H_{n}}{H_{n-1}^2}\;,\quad n=1,2,\ldots\;;
\end{equation}
\item[{\rm (iii)}] setting
\begin{equation}\label{Pn=fngn}
P_n(x)=x^n+f_nx^{n-1}+g_nx^{n-2}+\mbox{\rm lower degree terms}\;,
\end{equation}
the coefficients $f_n$ and $g_n$ are given in terms of the $\beta$ and $\gamma-$parameters by
\begin{equation}\label{coefxn-1}
f_n=-\sum_{j=0}^{n-1}\beta_j\;,\quad n=1,2,\cdots
\end{equation}
and
\begin{equation}\label{coefxn-2}
g_n=\sum_{0\leq i<j\leq n-1}\beta_i\beta_j-\sum_{k=1}^{n-1}\gamma_k\;,\quad n=2,3,\cdots
\end{equation}
\end{enumerate}
\end{corollary}
\end{snugshade}

{\it Proof.}
(i) and (ii) follow from the proof of Theorem \ref{thm-TTRR} and
taking into account (\ref{Hn-not-zero}) in Corollary \ref{leading-coeff-Hn}, so that
$\langle{\bf u},P_{n}^2\rangle=H_n/H_{n-1}$ for all $n\in\mathbb{N}_0$.
%$$\langle{\bf u},P_{n}^2\rangle=\frac{H_n}{H_{n-1}}\;,\quad n=0,1,2\ldots\;.$$
%The proof of (\ref{coefxn-1})--(\ref{coefxn-2}) is left to the reader. (Exercise \ref{Ex-cp2-4})
To prove (iii), substitute
$P_{n}(x)=x^{n}+f_{n}x^{n-1}+g_{n}x^{n-2}+\cdots$ and the corresponding expressions
for $P_{n+1}(x)$ and $P_{n-1}(x)$ in the recurrence relation (\ref{TTRR}), so that
$$
\begin{array}{l}
x^{n+1}+f_{n+1}x^n+g_{n+1}x^{n-1}+\cdots \\  [0.5em]
\qquad\qquad =(x-\beta_n)\big(x^{n}+f_{n}x^{n-1}+g_{n}x^{n-2}+\cdots\big)
-\gamma_n\big(x^{n-1}+f_{n-1}x^{n-2}+\cdots\big) \\  [0.5em]
\qquad\qquad =x^{n+1}+(f_n-\beta_n)x^{n}+(g_n-\beta_nf_n-\gamma_n)x^{n-1}+\cdots\;.
\end{array}
$$
Therefore, by comparing coefficients, and defining $f_0=g_1=0$, we obtain
%\begin{equation}\label{fn+1}
%%f_0=0\,,\quad
%f_{n+1}=f_n-\beta_n\;,\quad n\geq0
%\end{equation}
%\begin{equation}\label{gn+1}
%%g_1=0\,,\quad
%g_{n+1}=g_n-\beta_nf_n-\gamma_n\;,\quad n\geq1\;.
%\end{equation}
%Henceforth, (\ref{coefxn-1}) and (\ref{coefxn-2}) follow easily by induction
%(or by the telescoping method) from (\ref{fn+1}) and (\ref{gn+1}).
$$
\begin{array}{c}
f_{n+1}=f_n-\beta_n\;,\quad n\geq0\;; \\ [0.25em]
g_{n+1}=g_n-\beta_nf_n-\gamma_n\;,\quad n\geq1\;.
\end{array}
$$
Hence, (\ref{coefxn-1}) and (\ref{coefxn-2}) follow easily by induction
(or by applying the telescoping property for sums). % from (\ref{fn+1}) and (\ref{gn+1}).
\qed

\begin{remark}\label{Remark-uPn2}\em
Regarding to Corollary \ref{cor-TTRR}, notice also the relations
(with the convention that empty product equals one)
\begin{snugshade}
\begin{equation}\label{usefull-R}
\langle{\bf u},P_{n}^2\rangle=\frac{H_n}{H_{n-1}}
=u_0\prod_{j=1}^n\gamma_j\;,\quad n=0,1,2\ldots\;.
\end{equation}
\end{snugshade}
\end{remark}

\begin{snugshade}
\begin{theorem}[Favard]\label{Favard}
Let $\{\beta_n\}_{n\geq0}$ and $\{\gamma_n\}_{n\geq0}$ be two arbitrary
sequences of complex numbers, and let $\{P_n\}_{n\geq0}$ be a sequence of (monic) polynomials
defined by the three-term recurrence relation
\begin{equation}\label{TTRRFavard}
P_{n+1}(x)=(x-\beta_n)P_n(x)-\gamma_nP_{n-1}(x)\;,\quad n=0,1,2,\cdots
\end{equation}
with initial conditions %$P_{-1}(x)=0$ and $P_{0}(x)=1$,
\begin{equation}\label{TTRRinitialFavard}
P_{-1}(x)=0\;,\quad P_{0}(x)=1\;.
\end{equation}
Then there exists a unique functional ${\bf u}\in\mathcal{P}'$ such that
\begin{equation}\label{Fv1}
\langle{\bf u},1\rangle=u_0:=\gamma_0\;,\qquad
\langle{\bf u},P_nP_m\rangle=0\quad\mbox{\rm if}\quad n\neq m\quad (n,m\in\mathbb{N}_0)\;.
\end{equation}
Moreover, ${\bf u}$ is regular and $\{P_n\}_{n\geq0}$ is the corresponding monic OPS
if and only if $\gamma_n\neq0$ for each $n\in\mathbb{N}_0$,
while ${\bf u}$ is positive-definite and $\{P_n\}_{n\geq0}$ is the corresponding monic OPS
if and only if $\beta_{n}\in\mathbb{R}$ and  $\gamma_n>0$ for each $n\in\mathbb{N}_0$.
\end{theorem}
\end{snugshade}

{\it Proof.}
Since $\{P_n\}_{n\geq0}$ defined by (\ref{TTRRFavard}) is clearly a simple set in $\mathcal{P}$
(so that it is an algebraic basis in $\mathcal{P}$),
we may define a functional ${\bf u}\in\mathcal{P}'$ by
\begin{equation}\label{TTRRF1a}
\langle{\bf u},P_0\rangle=\langle{\bf u},1\rangle:=\gamma_0\;,\qquad
\langle{\bf u},P_n\rangle=0\,,\quad n\geq1\;.
\end{equation}
Rewrite (\ref{TTRRFavard}) as
\begin{equation}\label{TTRRF1b}
xP_n(x)=P_{n+1}(x)+\beta_nP_n(x)+\gamma_nP_{n-1}(x)\;,\quad n\geq0\,.
\end{equation}
Then, $\langle{\bf u},xP_n\rangle
=\langle{\bf u},P_{n+1}\rangle+\beta_n\langle{\bf u},P_n\rangle+\gamma_n\langle{\bf u},P_{n-1}\rangle$
for each $n\geq0$, hence, by (\ref{TTRRF1a}),
\begin{equation}\label{TTRRF1c}
\langle{\bf u},xP_n\rangle=0\,,\quad n\geq2\;.
\end{equation}
Multiplying both sides of (\ref{TTRRF1b}) by $x$ and using (\ref{TTRRF1c}), we find
$$%\begin{equation}\label{TTRRF1d}
\langle{\bf u},x^2P_n\rangle=0\,,\quad n\geq3\;.
$$%\end{equation}
Continuing in this manner, we deduce
\begin{equation}\label{TTRRF1e}
\langle{\bf u},x^kP_n\rangle=0\,,\quad 0\leq k<n\;,\quad n\in\mathbb{N}\;.
\end{equation}
Therefore, if $m\neq n$, say $m<n$, then writing $P_m(x)=\sum_{k=0}^ma_{m,k}x^k$, we obtain
$$
\langle{\bf u},P_mP_n\rangle=\sum_{k=0}^ma_{m,k}\langle{\bf u},x^kP_n\rangle=0\;.
$$
This proves (\ref{Fv1}).
Next, for each $n\in\mathbb{N}$, multiplying both sides of (\ref{TTRRF1b}) by $x^{n-1}$, we find
$\langle{\bf u},x^nP_n\rangle%=\langle{\bf u},x^{n-1}xP_n\rangle
=\langle{\bf u},x^{n-1}P_{n+1}\rangle+
\beta_n\langle{\bf u},x^{n-1}P_n\rangle+\gamma_n\langle{\bf u},x^{n-1}P_{n-1}\rangle$,
hence, using (\ref{TTRRF1e}),
$$
\langle{\bf u},x^nP_n\rangle=\gamma_n\langle{\bf u},x^{n-1}P_{n-1}\rangle\;,\quad n\in\mathbb{N}\;.
$$
Applying successively this equality, we find
\begin{equation}\label{TTRRF1f}
\langle{\bf u},P_n^2\rangle=\langle{\bf u},x^nP_n\rangle=\gamma_0\gamma_1\cdots\gamma_n\;,\quad n\in\mathbb{N}_0\;.
\end{equation}
This holds for $n=0$ since $\langle{\bf u},1\rangle:=\gamma_0$.
Notice also that the first equality in (\ref{TTRRF1f}) holds taking into accout (\ref{TTRRF1e}), after
writing  $P_n(x)=x^n+\sum_{k=0}^{n-1}a_{n,k}x^k$. It follows from (\ref{Fv1}) and (\ref{TTRRF1f}) that
${\bf u}$ is regular and $\{P_n\}_{n\geq0}$ is the corresponding monic OPS
if and only if $\gamma_n\neq0$ for each $n\in\mathbb{N}_0$.

In addition, if ${\bf u}$ is positive-definite and $\{P_n\}_{n\geq0}$ is the corresponding monic OPS,
%Then, it follows from (\ref{TTRRF1f}) that
%$\gamma_0\gamma_1\cdots\gamma_n=\langle{\bf u},P_n^2\rangle>0$ for each $n\in\mathbb{N}_0$, hence,
%since $\gamma_0=\langle{\bf u},1\rangle>0$, it follows by induction that $\gamma_n>0$ for all $n\in\mathbb{N}_0$.
then $\gamma_0=\langle{\bf u},1\rangle>0$ and so, by Theorem \ref{thm-TTRR}, we may conclude that
$\beta_n\in\mathbb{R}$ and $\gamma_n>0$ for each $n\in\mathbb{N}_0$.
Conversely, assume that $\beta_n\in\mathbb{R}$ and $\gamma_n>0$ for each $n\in\mathbb{N}_0$.
Then, by (\ref{TTRRFavard}) we see that $P_n(x)$ is real (i.e., it  has real coefficients) for each $n\in\mathbb{N}_0$.
Moreover, from (\ref{TTRRF1f}), we have $\langle{\bf u},P_n^2\rangle>0$ for each $n\in\mathbb{N}_0$.
This, together with (\ref{Fv1}), proves that $\{P_n\}_{n\geq0}$ is the monic OPS with respect to ${\bf u}$.
By Corollary \ref{Cor-exist-OPS+}, ${\bf u}$ is positive-definite.
\qed

\begin{remark}\em
%It is worth mentioning that,
Since $\{P_n\}_{n\geq0}$ in Theorem
\ref{Favard} is independent of $\gamma_0$, and $u_0:=\gamma_0$, then the functional ${\bf u}$ is
unique up to the (given) choice of $\gamma_0$, i.e.,
up to the choice of its first moment $u_0:=\langle{\bf u},1\rangle$.
\end{remark}

\begin{remark}\em
The original theorem of Favard concerned only the positive-definite case and
the functional ${\bf u}$ was represented by a Stieltjes integral.
The corresponding result for regular functionals was subsequently observed by Shohat.
\end{remark}

\section{The Christoffel-Darboux identities}

In this section we state other important consequences of the three-term
recurrence relation characterizing a given OPS.

\begin{snugshade}
\begin{theorem}[Christoffel-Darboux identities]\label{CD-id}
Let $\{P_n\}_{n\geq0}$ be a monic OPS fulfilling the three-term recurrence relation
$(\ref{TTRR})$--$(\ref{TTRRinitial})$.
Then, for each $n\in\mathbb{N}_0$,
\begin{equation}\label{CD1}
\sum_{j=0}^n\frac{P_j(x)P_j(y)}{\gamma_1\gamma_2\cdots\gamma_j}=
\frac{1}{\gamma_1\gamma_2\cdots\gamma_n}\frac{P_{n+1}(x)P_n(y)-P_{n}(x)P_{n+1}(y)}{x-y}\quad\mbox{\rm if}\quad x\neq y
\end{equation}
(with the convention that empty product equals one), and
\begin{equation}\label{CD2}
\sum_{j=0}^n\frac{P_j^2(x)}{\gamma_1\gamma_2\cdots\gamma_j}=
\frac{P_{n+1}'(x)P_n(x)-P_{n}'(x)P_{n+1}(x)}{\gamma_1\gamma_2\cdots\gamma_n}\;.
\end{equation}
\end{theorem}
\end{snugshade}

{\it Proof.}
Since (\ref{CD2}) follows from (\ref{CD1}) by taking the limit $y\to x$, we only need to prove (\ref{CD1}).
From (\ref{TTRR})--(\ref{TTRRinitial}) we have, for each $n\in\mathbb{N}_0$,
$$
\begin{array}{l}
xP_n(x)P_{n}(y)=P_{n+1}(x)P_n(y)+\beta_nP_n(x)P_n(y)+\gamma_nP_{n-1}(x)P_n(y)\;, \\ [0.25em]
yP_n(y)P_{n}(x)=P_{n+1}(y)P_n(x)+\beta_nP_n(y)P_n(x)+\gamma_nP_{n-1}(y)P_n(x)\;.
\end{array}
$$
Subtracting the second equation from the first one yields
\begin{equation}\label{eqCDa}
(x-y)P_n(x)P_{n}(y)=G_{n+1}(x,y)-\gamma_nG_n(x,y)\;,\quad n\in\mathbb{N}_0\,,
\end{equation}
where
$$
G_n(x,y):=P_n(x)P_{n-1}(y)-P_n(y)P_{n-1}(x)\;.
$$
Dividing both sides of (\ref{eqCDa}) by $\gamma_1\gamma_2\cdots\gamma_n(x-y)$,
and then in the resulting equality changing $n$ into $j$, we obtain
$$
\frac{P_j(x)P_{j}(y)}{\gamma_1\gamma_2\cdots\gamma_j}
=\frac{G_{j+1}(x,y)}{\gamma_1\cdots\gamma_j(x-y)}
-\frac{G_{j}(x,y)}{\gamma_1\cdots\gamma_{j-1}(x-y)}\,,\quad j\in\mathbb{N}_0\;.
$$
Summing from $j=0$ to $j=n$, the right-hand side becomes a telescoping sum,
hence, taking into account that $G_0(x,y)=0$, we deduce (\ref{CD1}).
\qed

\section*{Exercises}
%\bigskip

{\small
%\noindent
\begin{enumerate}[label=\emph{\bf \arabic*.},leftmargin=*]

\item\label{Ex-cp2-0}
Is the simple set $\{x^n\}_{n\geq0}$ an OPS with respect to some ${\bf u}\in\mathcal{P}'\;$?
\medskip

\item\label{Ex-cp2-1} % Szego, p.22 ; tese FOUPOUAGNIGNI, p.18
Let ${\bf u}\in\mathcal{P}'$, regular, and $\{P_n\}_{n\geq0}$ the corresponding monic OPS. Show that
$$
P_n(x)=\frac{1}{\Delta_{n-1}}
\left|
\begin{array}{ccccc}
\langle{\bf u},R_0R_0\rangle & \langle{\bf u},R_0R_1\rangle & \cdots &
\langle{\bf u},R_0R_{n-1}\rangle & \langle{\bf u},R_0R_n\rangle \\ [0.25em]
\langle{\bf u},R_1R_0\rangle & \langle{\bf u},R_1R_1\rangle & \cdots &
\langle{\bf u},R_1R_{n-1}\rangle & \langle{\bf u},R_1R_n\rangle \\ [0.25em]
\vdots&\vdots&\ddots&\vdots&\vdots \\ [0.25em]
\langle{\bf u},R_{n-1}R_0\rangle & \langle{\bf u},R_{n-1}R_1\rangle & \cdots &
\langle{\bf u},R_{n-1}R_{n-1}\rangle & \langle{\bf u},R_{n-1}R_n\rangle \\ [0.25em]
R_0(x) & R_1(x) & \cdots & R_{n-1}(x) & R_n(x)
\end{array}
\right|
$$
where $\{R_n\}_{n\geq0}$ is any simple set of monic polynomials, and
$$
\Delta_{-1}:=1\;,\quad
\Delta_n:=\det\big\{\big[\langle{\bf u},R_iR_j\rangle\big]_{i,j=0}^{n}\big\}\;,\quad n\geq0\;.
$$

%\item\label{Ex-cp2-2}
%%Let $E=\{x_1,\ldots,x_N\}$ be a set with $N$ distinct real numbers, and let $h_1,\ldots,h_N>0$.
%Let $x_1,\ldots,x_N$ be any $N$ distinct real numbers, and let $h_1,\ldots,h_N>0$.
%Define ${\bf u}\in\mathcal{P}'$ by
%$$
%\langle {\bf u},x^n\rangle:=\sum_{j=1}^N h_jx_j^n\,,\quad n\in\mathbb{N}_0\;.
%$$
%Prove that:
%\begin{enumerate}
%\item ${\bf u}$ is positive-definite on $E:=\{x_1,\ldots,x_N\}$;
%\item ${\bf u}$ is not positive-definite on any set $S\subseteq\mathbb{R}$ such that $E$ is a proper subset of $S$.
%\end{enumerate}
%\medskip

%\label{Ex-cp2-1}

%\item
%Let ${\bf u}\in\mathcal{P}'$ and $H_n$ the associated Hankel determinant of order $n+1$, defined as in (\ref{def-Hn}).
%Show that if $H_n>0$ for all $n\in\mathbb{N}$ then all the moments $u_n:=\langle{\bf u},x^n\rangle$ are real numbers.
%More precisely, for each $n\in\mathbb{N}_0$ the following holds:
%$$
%u_{2n}>0\;,\quad u_{2n+1}\in\mathbb{R}\;.
%$$
%
%\item\label{Ex-cp2-4}
%Prove (\ref{coefxn-1}) and (\ref{coefxn-2}).

\item\label{Ex-cp2-3}
Let $\{T_n\}_{n\geq0}$ be the sequence of the Chebyshev polynomials of the first kind, defined by
$$
T_n(x)=\cos(n\theta)\;,\quad x=\cos\theta\quad(0\leq\theta\leq\pi\;;\;-1\leq x\leq 1).
$$
\begin{enumerate}
\item%[{\rm (a)}]
Prove that $\{T_n\}_{n\geq0}$ fulfills the three-term recurrence relation
$$
2xT_{n}(x)=T_{n+1}(x)+T_{n-1}(x)\;,\quad n=1,2,3,\ldots
$$
with initial conditions $T_0(x)=1$ and $T_{1}(x)=x$.
(Note that this shows that $T_n$ is a polynomial of degree $n$ for each $n\in\mathbb{N}$.)
%(Hence $U_n$ is a polynomial of degree $n$ and, by Favard's theorem, $\{U_n\}_{n\geq0}$ is an OPS w.r.t. a positive-definite functional ${\bf u}\in\mathcal{P}'$.)
%%\item%[{\rm (b)}]
%%Find an explicit expression for $U_n(x)$ without involving trigonometric functions.
%%%%%%%%%%%%
\item
Set $p_0(x):=\frac{1}{\sqrt{\pi}}\,T_0(x)$ and $p_n(x):=\sqrt{\frac{2}{\pi}}\,T_n(x)$ if $n\geq1$. Show that $\{p_n\}_{n\geq0}$ is orthonormal with respect to ${\bf u}\in\mathcal{P}'$ given by
$$
\langle{\bf u},p\rangle:=\int_{-1}^1\frac{p(x)}{\sqrt{1-x^2\,}}\;{\rm d}x\;.
$$
\item
Prove that $T_n(x)$ admits the explicit expression
$$
T_n(x)=\sum_{k=0}^{\lfloor n/2\rfloor}(-1)^k \binom{n}{2k} x^{n-2k}(1-x^2)^k\;,
$$
where $\lfloor s\rfloor$ denotes the greatest integer less than or equal to the real number $s$.
\end{enumerate}
\bigskip

\item\label{Ex-cp2-4}
Let $\{U_n\}_{n\geq0}$ be the sequence of the Chebyshev polynomials of the second kind, defined by
$$
U_n(x)=\frac{\sin(n+1)\theta}{\sin\theta}\;,\quad x=\cos\theta\quad(0\leq\theta\leq\pi\;;\;-1\leq x\leq 1).
$$
(It is assumed that $U_n(x)$ is defined by continuity whenever $\sin\theta=0$.)\smallskip
\begin{enumerate}
\item%[{\rm (a)}]
Prove that $\{U_n\}_{n\geq0}$ fulfills the three-term recurrence relation
$$
2xU_{n}(x)=U_{n+1}(x)+U_{n-1}(x)\;,\quad n=1,2,3,\ldots
$$
with initial conditions $U_0(x)=1$ and $U_{1}(x)=2x$.
(Hence $U_n$ is a polynomial of degree $n$ for each $n\in\mathbb{N}$.)
%(Hence $U_n$ is a polynomial of degree $n$ and, by Favard's theorem, $\{U_n\}_{n\geq0}$ is an OPS w.r.t. a positive-definite functional ${\bf u}\in\mathcal{P}'$.)
%%\item%[{\rm (b)}]
%%Find an explicit expression for $U_n(x)$ without involving trigonometric functions.
%%%%%%%%%%%%
\item
Set $p_n(x):=\sqrt{\frac{2}{\pi}}\,U_n(x)$. Show that $\{p_n\}_{n\geq0}$ is orthonormal with respect to ${\bf u}\in\mathcal{P}'$ given by
$$
\langle{\bf u},p\rangle:=\int_{-1}^1p(x)\,\sqrt{1-x^2\,}\;{\rm d}x\;.
$$
\item
Prove that $U_n(x)$ admits the explicit representation
$$
U_n(x)=\sum_{k=0}^{\lfloor(n+1)/2\rfloor}(-1)^k \binom{n+1}{2k+1} x^{n-2k}(1-x^2)^k\;.
$$
%where $\lfloor s\rfloor$ represents the greatest integer that is less than or equal to the real number $s$.\newpage
%$\;\displaystyle U_n(x)=\sum_{k=0}^{\lfloor(n+1)/2\rfloor}(-1)^k \binom{n+1}{2k+1} x^{n-2k}(1-x^2)^k\;$.
%%%%%%%%%%%%%%%%%%%%%%%
%\item%[{\rm (c)}]
%%Let $a,b,c\in\mathbb{R}$, with $bc>0$, and let $A_n$ be the tridiagonal Toeplitz matrix of order $n$
%%\medskip
%Let $a,b,c\in\mathbb{R}$, with $bc>0$, and set
%$$P_n(x):=(bc)^{n/2}U_n\Big(\frac{x-a}{2\sqrt{bc}}\Big)\;,\quad n=0,1,2,\ldots\;.$$
%%for each $n\in\mathbb{N}_0$.
%\begin{enumerate}
%\item
%Show that $\{P_n\}_{n\geq0}$ is a monic OPS with respect ro some
%positive-definite functional on $\mathcal{P}'$.
%\item
%Consider the tridiagonal Toeplitz matrix of order $n$
%$$
%A_n=\left(
%\begin{array}{cccccc}
%a & b &  &  &  &   \\
%c & a & b &  & &  \\
% & c & a & b &  &   \\
% & & \ddots & \ddots & \ddots &  \\
% & & & c & a & b  \\
% & & & & c & a \\
%\end{array}
%\right) \, .
%$$
%%Use the Chebyshev polynomials $\{U_n\}_{n\geq0}$
%%where $U_n(x):=\sin(n+1)\theta/\sin\theta$ $(x=\cos\theta)$,
%Prove that the eigenvalues of $A_n$ are
%$$
%\lambda_j:=a+2\sqrt{bc}\,\cos\frac{j\pi}{n+1}\quad (j=1,2,\cdots, n)\, ,
%$$
%and find explicitly corresponding eigenvectors.
%\smallskip
%
%\noindent%\hspace*{-1em}
%({\sl Hint.} Define $Q_n(x):=b^{-n}P_n(x)$, and write the TTRR for $\{Q_n\}_{n\geq0}$ in matrix form.) %$P_n(x):=(bc)^{n/2}U_n\big(\frac{x-a}{2\sqrt{bc}}\big)$.)
%\end{enumerate}
\end{enumerate}
\bigskip

\item\label{Ex-cp2-5}
Let $\{P_n\}_{n\geq0}$ be the sequence of the Legendre polynomials, defined by
$$
P_n(x):=\frac{1}{2^n n!}\frac{{\rm d}^n}{{\rm d}x^n}\big\{(x^2-1)^n\big\}\;.
$$
Notice that  the leading coefficient of $P_n$ is $2^{-n}{2n\choose n}$, hence it is not a monic polynomial.
For each $n\in\mathbb{N}_0$, set
$$
p_n(x):=\sqrt{\frac{2n+1}{2}}\,P_n(x)\,.
$$
Show that $\{p_n\}_{n\geq0}$ is orthonormal with respect to ${\bf u}\in\mathcal{P}'$ given by
$$
\langle{\bf u},p\rangle:=\int_{-1}^1 p(x)\;{\rm d}x\;.
$$
\end{enumerate}
%\bigskip
}
\medskip

\section*{Final remarks}

The main sources on the basis of this text are the books by 
Theodore S. Chihara \ref{Chihara1978-C2} (1978) and
Mourad E. H. Ismail \ref{Ismail2004-C2} (2005),
where the student may found most of the results presented here.
The notion of OPS introduced in Definition \ref{def-OPS} reflects our option to adopt the concept of
{\it formal orthogonality}, also called {\it regular orthogonality}.
Many researchers/authors prefer to adopt a definition of orthogonality that corresponds to the positive-definite case.
The book by Gabor Szeg\"o \ref{Szegoo1975-C2} (whose 1$^{\rm\scriptsize st}$ edition goes back to 1939)
is considered the first important book entirely dedicated to the theory of OP.
Other recommended references containing the general theory presented here are the books by
Geza Freud \ref{Freud1971-C2} (1976) and Walter Gautschi \ref{Gautschi2004-C2} (2004).
Exercise {\bf 2} may be found e.g. in Szeg\" o's book.
Exercises {\bf 3} up to {\bf 5} involve three families of OP that the students probably already meet on previous courses (Numerical Analysis, Linear Algebra, or Functional Analysis, among others) and they can be found in the books included on the bibliography (appearing therein as exercises or not).
\medskip

%\section*{Notas finais}
%\section*{Coment\'arios finais}
%
%As refer\^encias principais usadas na elabora\c c\~ao deste texto s\~ao os livros de
%Theodore S. Chihara \ref{Chihara1978-C2} (1978) e de
%Mourad E. H. Ismail\footnote{\,M. E. H. Ismail \'e {\it Fellow} da AMS (2015) ``for contributions to classical analysis and special function theory''.} \ref{Ismail2004-C2} (2005),
%nos quais se encontra a maioria dos resultados apresentados.
%A no\c c\~ao de SPO que introduzimos (Defini\c c\~ao \ref{def-OPS}) reflecte a op\c c\~ao de adoptar o conceito
%de {\it ortogonalidade formal}, tamb\'em chamada {\it ortogonalidade regular}.
%Muitos investigadores/autores preferem adoptar uma defini\c c\~ao de ortogonalidade
%que corresponde ao que design\'amos por caso definido-positivo.
%O livro de Gabor Szeg\"o \ref{Szegoo1975-C2} (cuja 1$^{\rm\scriptsize a}$ edi\c c\~ao data de 1939)
%\'e considerado o primeiro livro importante inteiramente dedicado \`a teoria dos PO.
%Outros livros recomendados contendo a teoria geral s\~ao os de
%Geza Freud \ref{Freud1971-C2} (1976) e de Walter Gautschi \ref{Gautschi2004-C2} (2004).
%O exerc\'\i cio {\bf 2} pode encontrar-se e.g. no livro de Szeg\" o.
%Os exerc\'\i cios {\bf 3} a {\bf 5} incidem sobre tr\^es fam\'\i lias de PO
%que os estudantes provavelmente j\'a encontraram noutras ocasi\~oes
%(em An\'alise Num\'erica, ou \'Algebra Linear, ou An\'alise Funcional)
%e podem ser encontrados nos livros indicados na bibliografia
%(sob a forma de exerc\'\i cio ou n\~ao).

\section*{Bibliography}
\medskip

{\small
\begin{enumerate}[label=\emph{\rm [\arabic*]},leftmargin=*]
\item\label{Chihara1978-C2} T. S. Chihara, {\sl An introduction to orthogonal polynomials}, Gordon and Breach (1978).
%\item\label{Deift1997-C2} P. Deift, {\sl Orthogonal Polynomials and Random Matrices: A Riemann-Hilbert Approach}, AMS Courant Lecture Notes {\bf 3} (2000).
\item\label{Freud1971-C2} G. Freud, {\sl Orthogonal polynomials}, Pergamon Press, Oxford (1971).
\item\label{Gautschi2004-C2} W. Gautschi, {\sl Orthogonal polynomials. Computation and Approximation}, Oxford University Press, Oxford (2004).
\item\label{Ismail2004-C2} M. E. H. Ismail, {\sl Classical and Quantum Orthogonal Polynomials in One Variable},
          Cambridge University Press (2005) [paperback edition: 2009].
%\item\label{Maroni1991-C2} P. Maroni, {\it Une th\'eorie alg\'ebrique des polyn\^omes orthogonaux. Applications aux polyn\^omes
%orthogonaux semiclassiques}, In C. Brezinski et al. Eds., Orthogonal Polynomials and Their Applications, Proc. Erice 1990, IMACS, Ann. Comp. App. Math. {\bf 9} (1991) 95--130.
%\item\label{Lax2002i} Peter D. Lax, {\sl Functional Analysis}, John Wiley $\&$ Sons (2002).
%\item\label{ReedSimon1972i} M. Reed and B. Simon, {\sl Methods of Modern Mathematical Physics I: Functional Analysis},
%Academic Press (1972).
\item\label{Szegoo1975-C2} G. Szeg\"o, {\sl Orthogonal Polynomials}, AMS Colloq. Publ. {\bf 230} (1975), 4th ed.
\end{enumerate}
}

\chapter{Zeros of orthogonal polynomials and quadrature formulas}

\pagestyle{myheadings}\markright{Zeros of orthogonal polynomials and quadrature formulas}
\pagestyle{myheadings}\markleft{J. Petronilho}

\section{Zeros of OPS}

When ${\bf u}\in\mathcal{P}'$ is %(regular and)
positive-definite, then the zeros of the corresponding OPS exhibit a certain regularity in their behavior.
In order to discuss this behavior  we need to make an extension of the concept of
positive-definiteness as introduced in Definition \ref{def-OPS-posit}.
To make it clear %To avoid misunderstandings  %To avoid inaccurate statements
we emphasize that a polynomial $p$ is said to be nonzero on a set $E$
(written $p\not\equiv0$ on $E$) if it does not vanish identically on $E$.

\begin{snugshade}
\begin{definition}\label{def-OPS-positive-E}
Let ${\bf u}\in\mathcal{P}'$ and $E\subseteq\mathbb{R}$.
\begin{enumerate}
\item[{\rm (i)}]
${\bf u}$ is said {\sl positive-definite on} $E$  if the condition
\begin{equation}\label{u-def+E}
\langle{\bf u},p\rangle>0
\end{equation}
holds for each real polynomial $p$ which is nonzero and nonnegative on $E$;
\item[{\rm (ii)}]
if ${\bf u}$ is positive-definite on $E$, then $E$ is called a {\sl supporting set} for ${\bf u}$.
\end{enumerate}
\end{definition}
\end{snugshade}

\begin{remark}\em
Notice that if $E=\mathbb{R}$ then positive-definiteness on $\mathbb{R}$ is the same as positive-definiteness as
introduced in Definition \ref{def-OPS-posit}.
\end{remark}

\begin{snugshade}
\begin{theorem}\label{OPS-zeros-Thm1}
Let $E\subseteq\mathbb{R}$, with ${\small\#} E=\infty$.
Let  ${\bf u}\in\mathcal{P}'$ and suppose that ${\bf u}$ is positive-definite on $E$.
Then, the following holds:
\begin{enumerate}
\item[{\rm (i)}]
if $E\subseteq S$, then ${\bf u}$ is positive-definite on $S$;
\item[{\rm (ii)}]
if $E\supseteq S$ and $\overline{S}=E$, then ${\bf u}$ is positive-definite on $S$.
\end{enumerate}
\end{theorem}
\end{snugshade}

{\it Proof.}
(i) Let $p$ be a real polynomial and suppose that $p(x)\geq0$ for all $x\in S$, and that $p\not\equiv0$ on $S$.
Since (by hypothesis) $E\subseteq S$, then also $p(x)\geq0$ for all $x\in E$ and, moreover,
$p\not\equiv0$ on $E$ (since $p\not\equiv0$ on $S$---hence $p$ does not vanishes identically
on $\mathbb{R}$---and ${\small\#} E=\infty$).
Therefore, since (by hypothesis) ${\bf u}$ is positive-definite on $E$, we deduce $\langle{\bf u},p\rangle>0$.
Thus, ${\bf u}$ is positive-definite on $S$.

(ii) Take a real polynomial $p$ such that $p(x)\geq0$ for all $x\in S$ and $p\not\equiv0$ on $S$.
Then
\begin{equation}\label{eqSE0}
p(x)\geq 0\;,\quad\forall x\in E\;.
\end{equation}
Indeed, suppose that there is $x_0\in E$ with $p(x_0)<0$.
Since $p$ is continuous on $\mathbb{R}$, then
\begin{equation}\label{eqSE1}
\exists\delta>0\;:\;\forall x\in\mathbb{R}\;,\; |x-x_0|<\delta\Rightarrow p(x)<0\;.
\end{equation}
Now, since $x_0\in E$ and $\delta>0$, taking into account that $\overline{S}=E$, we may ensure that
\begin{equation}\label{eqSE2}
\exists s\in S\;:\; |s-x_0|<\delta\;.
\end{equation}
From (\ref{eqSE1}) and (\ref{eqSE2}) we conclude that $p(s)<0$, in contradiction with the choice of $p$.
Henceforth, (\ref{eqSE0}) holds. %we have $p(x)\geq0$ for all $x\in E$, and $p\not\equiv0$ on $E$
Moreover, $p\not\equiv0$ on $E$ (since $p\not\equiv0$ on $S$ and $S\subseteq E$).
Thus, since (by hypothesis) ${\bf u}$ is positive-definite on $E$,
we conclude that $\langle{\bf u},p\rangle>0$,
hence ${\bf u}$ is positive-definite on $S$.
\qed
%\smallskip

\begin{remark}\em
Statement (ii) in Theorem \ref{OPS-zeros-Thm1} holds trivially if ${\small\#} E<\infty$, since in that case $S=E$.
On the contrary, statement (i) does not holds if ${\small\#} E<\infty$.
For instance, %\item\label{Ex-hwk2-2}
if $x_1,\ldots,x_N$ are any $N$ distinct real numbers, and $h_1,\ldots,h_N>0$,
then, being ${\bf u}\in\mathcal{P}'$ the functional defined by
$$
\langle {\bf u},x^n\rangle:=\sum_{j=1}^N h_jx_j^n\quad (n\in\mathbb{N}_0)\,,
$$
${\bf u}$ is positive-definite on $E:=\{x_1,\ldots,x_N\}$, but it is not positive-definite
on any set $S\subseteq\mathbb{R}$ such that $E$ is a proper subset of $S$.\footnote{
\,This can be seen immediately by noticing that ${\bf u}$ may be represented as a
Riemann-Stieltjes integral with respect to the right continuous step function
supported on $E$ with jump $h_j$ at the point $x_j$, for each $j\in\{1,\ldots,N\}$.}
%$\psi:\mathbb{R}\to\mathbb{R}$ such that
%$$
%\psi(-\infty)=0\;,\;
%\psi(x_j+0)-\psi(x_j-0)=h_j\;(1\leq j\leq N)\;,\;
%\psi(+\infty)=h_1+\ldots+h_N\,,
%$$
%and being $E$ the support of $\psi$. (Exercise \ref{Ex-cp3-2})
\end{remark}

%\begin{remark}\em
%Statement (ii) in Theorem \ref{OPS-zeros-Thm1} shows that, in general,
%there is no ``smallest'' infinite supporting set for a positive-definite functional
%${\bf u}\in\mathcal{P}'$.
%\end{remark}

\begin{snugshade}
\begin{theorem}\label{OPS-zeros-Thm2}
Let ${\bf u}\in\mathcal{P}'$ be positive-definite, and let $\{P_n\}_{n\geq0}$ be the monic OPS with respect to ${\bf u}$.
Let $I$ be an interval which is a supporting set for ${\bf u}$.
Then, for each $n\in\mathbb{N}$, the zeros of $P_n$ are all real, simple, and they are located in the interior of $I$.
\end{theorem}
\end{snugshade}

{\it Proof.}
Fix $n\in\mathbb{N}$.
Since ${\bf u}$ is positive-definite then $P_n$ is a real polynomial, i.e.,
its coefficients are real numbers (by Theorem \ref{GramShmidt}).
Moreover, since $\langle{\bf u},P_n\rangle=0$ then $P_n(x)$ must change sign at least once in the interior of the interval $I$.
[Indeed, if $P_n(x)\geq0$ for all $x\in I$ then, since $I$ is a supporting set for ${\bf u}$, we would have $\langle{\bf u},P_n\rangle>0$, a contradiction with $\langle{\bf u},P_n\rangle=0$, hence, there is at least one point $r_1\in I$ such that $P_n(r_1)<0$. Similarly, if $P_n(x)\leq0$ for all $x\in I$ then $-P_n(x)\geq0$ for all $x\in I$, so we would have $\langle{\bf u},P_n\rangle=-\langle{\bf u},-P_n\rangle<0$, again a contradiction with $\langle{\bf u},P_n\rangle=0$, hence, there is at least one point $r_2\in I$ such that $P_n(r_2)>0$. Therefore, $P_n(r_1)P_n(r_2)<0$, so $P_n(x)$ change sign at least once in the interval $(r_1,r_2)\subset I$.]
Therefore, $P_n(x)$ has at least one zero of odd multiplicity located in the interior of $I$.
Let $x_1,\ldots,x_k$ denote the distinct zeros of odd multiplicity of $P_n(x)$ which are located in the interior of $I$. Set
$$
\pi_k(x):=(x-x_1)(x-x_2)\cdots(x-x_k)\;.
$$
Then the polynomial $\pi_kP_n$ has no zeros of odd multiplicity
in the interior of $I$, hence $\pi_k(x)P_n(x)\geq0$ for each $x\in I$.
Therefore, since ${\bf u}$ is positive-definite on $I$,
\begin{equation}\label{Z1a}
\langle{\bf u},\pi_kP_n\rangle>0\,.
\end{equation}
On the other hand, since $\{P_n\}_{n\geq0}$ is an OPS with respect to ${\bf u}$, we must have
\begin{equation}\label{Z2a}
\langle{\bf u},\pi_kP_n\rangle\left\{
\begin{array}{rcl}
=0 & \mbox{\rm if} & k<n\,, \\ [0.5em]
\neq0 & \mbox{\rm if} & k=n\,.
\end{array}
\right.
\end{equation}
From (\ref{Z1a}) and (\ref{Z2a}) we deduce that $k=n$.
This means that $P_n(x)$ has $n$ distinct zeros of odd multiplicity in the interior of $I$,
and since ${\rm deg}\,P_n=n$, we may conclude that $P_n(x)$
has $n$ real and simple zeros, all in the interior of $I$.
\qed
\bigskip

Let $\{P_n\}_{n\geq0}$ be a monic OPS with respect to a positive-definite functional ${\bf u}\in\mathcal{P}'$.
According to Theorem \ref{OPS-zeros-Thm2},
the zeros $x_{n,1},\ldots,x_{n,n}$ of each $P_n(x)$ may be ordered by increasing size, so that
\begin{snugshade}
\begin{equation}\label{OPS-zeros-order1}
x_{n,1}<x_{n,2}<\ldots<x_{n,n}\,,\quad n\geq2\,.
\end{equation}
\end{snugshade}
\noindent
Since $P_n(x)$ as positive leading coefficient ($=1$), it follows that for each $n\geq1$,
\begin{snugshade}
\begin{equation}\label{OPS-zeros-order2}
P_n(x)>0\quad\mbox{\rm if}\quad x>x_{n,n}\,;
\end{equation}
\begin{equation}\label{OPS-zeros-order3}
\mbox{\rm sgn}\, P_n(x)=(-1)^n\quad\mbox{\rm if}\quad x<x_{n,1}\,,
\end{equation}
\end{snugshade}
\noindent
where $\mbox{\rm sgn}$ is the {\it signum} function, defined by
\begin{snugshade}
$$
\mbox{\rm sgn}\,(x):=\left\{
\begin{array}{rcl}
-1 & \mbox{\rm if} & x<0 \\
0 & \mbox{\rm if} & x=0 \\
1 & \mbox{\rm if} & x>0\,.
\end{array}
\right.
$$
\end{snugshade}

\begin{snugshade}
\begin{theorem}\label{OPS-zeros-Thm3}
Let ${\bf u}\in\mathcal{P}'$ be positive-definite, and $\{P_n\}_{n\geq0}$ the corresponding monic OPS.
Suppose (without loss of generality) that the zeros of $P_n$ fulfill $(\ref{OPS-zeros-order1})$ for each $n\geq2$.
Then, the following holds:
\begin{enumerate}
\item[{\rm (i)}] $P_{n}'$ has exactly one zero in each open interval $(x_{n,j},x_{n,j+1})$, $1\leq j\leq n-1$.
Moreover:
\begin{equation}\label{OPS-zeros-order4}
\mbox{\rm sgn}\, P_n'(x_{n,j})=(-1)^{n-j}\,,\quad j=1,2,\ldots,n\,;
\end{equation}
\item[{\rm (ii)}] the zeros of $P_{n}$ and $P_{n+1}$ fulfill the separating (or interlacing) property:.
\begin{equation}\label{OPS-zeros-order5}
x_{n+1,j}<x_{n,j}<x_{n+1,j+1}\,,\quad j=1,2,\ldots,n\,;
\end{equation}
\item[{\rm (iii)}] for each $j\in\mathbb{N}$, $\{x_{n,j}\}_{n\geq j}$ is a decreasing sequence, while
$\{x_{n,n-j+1}\}_{n\geq j}$ is an increasing sequence;
\item[{\rm (iv)}] for each $j\in\mathbb{N}$, the limits
\begin{equation}\label{OPS-zeros-order6}
\xi_j:=\lim_{n\to\infty} x_{n,j}\;,\quad
\eta_j:=\lim_{n\to\infty} x_{n,n-j+1}
\end{equation}
all exist (at least in the extended real number system).
\end{enumerate}
\end{theorem}
\end{snugshade}

{\it Proof.}
(i) Since $P_n(x)$ has $n$ real and distinct zeros $x_{n,1},\ldots,x_{n,n}$,
then by the Cauchy-Bolzano theorem the derivative $P_n'(x)$ has $n-1$ real and distinct zeros,
one zero in between each pair of consecutive zeros of $P_n(x)$. Henceforth, in each interval
$(x_{n,j},x_{n,j+1})$, $1\leq j\leq n-1$, there is exactly one zero of $P_n'(x)$.
Moreover, we see that for each $j\in\{1,2,\ldots,n-1\}$, $P_n'(x_{n,j})$ alternates in sign as $j$ varies from $0$ to $n$,
and since $P'_n(x)$ has positive leading coefficient ($=n$), we conclude that (\ref{OPS-zeros-order4}) holds.
%\smallskip

(ii) By the confluent Cristoffel-Darboux formula (\ref{CD2}), with $x=x_{n+1,j}$, we have
\begin{equation}\label{ZZpn1}
P_{n+1}'(x_{n+1,j})P_n(x_{n+1,j})>0\,,\quad n\in\mathbb{N}_0\;,\;\; 1\leq j\leq n+1\;.
\end{equation}
On the other hand, by (\ref{OPS-zeros-order4}) with $n$ replaced by $n+1$, we also have
\begin{equation}\label{ZZpn2}
\mbox{\rm sgn}\, P_{n+1}'(x_{n+1,j})=(-1)^{n+1-j}\,,\quad j=1,2,\ldots,n+1\,.
\end{equation}
It follows from (\ref{ZZpn1}) and (\ref{ZZpn2}) that
$$
\mbox{\rm sgn}\, P_n(x_{n+1,j})=(-1)^{n+1-j}\,,\quad j=1,2,\ldots,n+1\,.
$$
Therefore, $P_n(x)$ has at least one zero, and hence exactly one zero, on each of the intervals
$(x_{n+1,j},x_{n+1,j+1})$, $1\leq j\leq n$, which proves (\ref{OPS-zeros-order5}).

(iii) It is an immediate consequence of (ii).

(iv) It is an immediate consequence of (iii).
\qed
\medskip

Property (iv) in Theorem \ref{OPS-zeros-Thm3} motivates a very important definition:

\begin{snugshade}
\begin{definition}\label{trueIO}
Let ${\bf u}\in\mathcal{P}'$ be positive-definite,
and $\{P_n\}_{n\geq0}$ the corresponding monic OPS.
The closed interval $[\xi,\eta]$, where
\begin{equation}\label{trueIO1}
\xi:=\lim_{n\to\infty} x_{n,1}\;,\quad
\eta:=\lim_{n\to\infty} x_{n,n}
\end{equation}
is called the {\sl true interval of orthogonality} of $\{P_n\}_{n\geq0}$.
\end{definition}
\end{snugshade}

\begin{remark}\label{TIZ}\em
The true interval of orthogonality is the smallest closed interval that contains all the zeros of all the polynomials in the sequence $\{P_n\}_{n\geq0}$. Moreover, it can be shown that
%Moreover, it can be shown---cf. (\ref{SupCoSup})---that
{\it the true interval of orthogonality is the smallest closed interval that is a supporting set for ${\bf u}$}.
\end{remark}

\begin{remark}\em
The three-term recurrence relation for a given monic OPS (not necessarily with respect to a positive-definite moment linear functional)
$$
xP_{n-1}(x)=P_n(x)+\beta_{n-1}P_{n-1}(x)+\gamma_{n-1}P_{n-2}(x)\;,\quad n\geq1\,,
$$
with initial conditions $P_{-1}(x):=0$ and $P_0(x)=1$, may be written in matrix form as
\begin{snugshade}
\begin{equation}\label{JnPn}
x\left(\begin{array}{c} P_0(x) \\ P_1(x) \\ \vdots \\ P_{n-2}(x) \\ P_{n-1}(x) \end{array}\right)=
J_n\,\left(\begin{array}{c} P_0(x) \\ P_1(x) \\ \vdots \\ P_{n-2}(x) \\ P_{n-1}(x) \end{array}\right)+
P_n(x)\,\left(\begin{array}{c} 0 \\ 0 \\ \vdots \\ 0 \\ 1 \end{array}\right)\,,
\end{equation}
\end{snugshade}\noindent
where $J_n$ is a tridiagonal matrix of order $n$ given by
\begin{snugshade}
\begin{equation}\label{Jn}
J_n:=\left(
\begin{array}{cccccc}
\beta_0 & 1 &  &  &  &   \\
\gamma_1 & \beta_1 & 1 &  & &  \\
 & \gamma_2 & \beta_2 & 1 &  &   \\
 & & \ddots & \ddots & \ddots &  \\
 & & &  & \beta_{n-2} & 1  \\
 & & & & \gamma_{n-1} & \beta_{n-1} \\
\end{array}
\right) \, .
\end{equation}
\end{snugshade}\noindent
Clearly, for each $n\in\mathbb{N}$, the following holds:
\begin{enumerate}
\item[{\rm (i)}]
the eigenvalues of $J_n$ are the zeros of $P_n$, hence the spectrum of $J_n$ is
\begin{snugshade}
$$
\sigma(J_n)=\{x_{n,j}\,:\,j=1,\ldots,n\}\,;
$$
\end{snugshade}
\item[{\rm (ii)}]
an eigenvector $v_{n,j}$ corresponding to the eigenvalue $x_{n,j}$ is
\begin{snugshade}
$$
v_{n,j}:=\left(\begin{array}{c} P_0(x_{n,j}) \\ P_1(x_{n,j}) \\ \vdots \\ P_{n-2}(x_{n,j}) \\ P_{n-1}(x_{n,j}) \end{array}\right)\;,\quad j=1,\ldots,n\,.
$$
\end{snugshade}
\end{enumerate}\noindent
This establishes a connection between Orthogonal Polynomials and Linear Algebra.
\end{remark}

\begin{remark}\em
Often we will refer to $J_n$ as the {\it Jacobi matrix} associated with $P_n$,
although in the framework of Linear Algebra the name ``Jacobi'' is usually attached to symmetric tridiagonal matrices.
\end{remark}

\begin{remark}\em
As a consequence of the connection just mentioned, $P_n$ is the (monic) characteristic
polynomial associated with the matrix $J_n$, so that
\begin{snugshade}
$$P_n(x)=\det\big(xI_n-J_n\big)\;,$$
\end{snugshade}\noindent
where $I_n$ is the identity matrix of order $n$. Henceforth, $P_n(x)$ may be represented
as a determinant involving only the sequences of the $\beta$ and $\gamma-$parameters:
\begin{snugshade}
\begin{equation}\label{PndetJn}
P_n(x)=\left|
\begin{array}{cccccc}
x-\beta_{0} & 1  &  0& \dots & 0 & 0  \\
\gamma_{1} & x-\beta_{1}  &  1 &\dots & 0 & 0 \\
0 & \gamma_{2} & x-\beta_{2}&   \dots & 0 & 0 \\
\vdots & \vdots & \vdots   & \ddots & \vdots & \vdots \\
0 &0 & 0 & \cdots& x-\beta_{n-2} & 1 \\
0 &0 & 0 &  \cdots & \gamma_{n-1} & x-\beta_{n-1}
\end{array}\right|\;,\quad n\in\mathbb{N}\;.
\end{equation}
\end{snugshade}
\end{remark}

\section{Gauss-Jacobi-Christoffel quadrature formula}

Fix $n$ points $(t_j,y_j)\in\mathbb{R}^2$, $1\leq j\leq n$ ($n\in\mathbb{N}$).
Assume that $t_i\neq t_j$ if $i\neq j$.
It is well known that the only solution for the problem ---known as {\it Lagrange problem}--- of
constructing a polynomial of degree at most $n-1$ whose graph passes through all the points $(t_j,y_j)$ is the so--called {\it Lagrange interpolation polynomial}, $L_n$, defined by
\begin{snugshade}
\begin{equation}\label{LagrangeInterpPoly}
L_n(x):=\sum_{j=1}^n y_j\ell_j(x)\;,
\end{equation}
\end{snugshade}\noindent
where
\begin{snugshade}
\begin{equation}\label{Lagrange-ellj}
\ell_j(x):=\frac{F(x)}{(x-t_j)F'(t_j)}\;,\quad F(x):=\prod_{i=1}^n(x-t_i)\;.
\end{equation}
\end{snugshade}\noindent
Clearly, $\ell_j$ is a polynomial of degree $n-1$ for each $j=1,\ldots,n$, which fulfils
\begin{snugshade}
\begin{equation}\label{elljtk}
\ell_j(t_k)=\delta_{j,k}\;,\quad j,k=1,2,\ldots,n\;.
\end{equation}
\end{snugshade}
Moreover, the interpolation property implies that $L_n$ satisfies the property
\begin{snugshade}
\begin{equation}\label{Lntj=yj}
L_n(t_j)=y_j\;,\quad j=1,2,\ldots,n\;.
\end{equation}
\end{snugshade}

We will use the Lagrange interpolation polynomial to obtain the Gauss-Jacobi-Christoffel quadrature formula.
%This is a very useful tool in Numerical Analysis, e.g. for the approximation of integrals by numerical quadrature.\footnote{\,Numerical quadrature consists of approximating the integral of a given integrable function $f$,
%$I[f]:=\int_\mathbb{R}f(x)\,{\rm d}\mu(x)\,$, with respect some positive Borel measure $\mu$,
%by a finite sum which uses only the values of $f$ at $n$ points $t_j$ (called nodes),
%$I_n[f]:=\sum_{j=1}^nf(t_j)A_j\,$,
%where the coefficients $A_j$ (which may depend on $n$, as well as the notes $t_j$)
%have to be chosen properly so that the quadrature formula is
%correct,---i.e., the equality $I[f]=I_n[f]$ holds---, for as many functions $f$ as possible.}

\begin{snugshade}
\begin{theorem}[Gauss-Jacobi-Christoffel quadrature formula]\label{GaussJC}
Let ${\bf u}\in\mathcal{P}'$ be positive-definite, and $\{P_n\}_{n\geq0}$ the corresponding monic OPS.
For each $n\in\mathbb{N}$, denote by $x_{n,1},\ldots,x_{n,n}$ the zeros of $P_n$.
Then
\begin{equation}\label{GJC1}
\forall n\in\mathbb{N}\;,\;\exists A_{n,1},\ldots,A_{n,n}>0\;:\;\forall p\in\mathcal{P}_{2n-1}\,,\quad
\langle{\bf u},p\rangle=\sum_{j=1}^nA_{n,j}p(x_{n,j})\;.
\end{equation}
Moreover,\vspace*{-0.5em}
\begin{equation}\label{GJC2}
\sum_{j=1}^nA_{n,j}=u_0:=\langle{\bf u},1\rangle\;.
\end{equation}
\end{theorem}
\end{snugshade}

{\it Proof.}
Let $p\in\mathcal{P}_{2n-1}$.
Consider the Lagrange interpolation polynomial $L_n$ that passes through the points
$(t_j,y_j)\equiv(x_{n,j},p(x_{n,j}))$, $1\leq j\leq n$, i.e.,
$$
L_n(x):=\sum_{j=1}^n p(x_{n,j})\ell_{j,n}(x)\;,\quad
\ell_{j,n}(x):=\frac{P_n(x)}{(x-x_{n,j})P_n'(x_{n,j})}\;.
$$
Let $Q(x):=p(x)-L_n(x)$. Then, $Q\in\mathcal{P}_{2n-1}$
and $Q(x_{n,j})=p(x_{n,j})-L_n(x_{n,j})=y_j-y_j=0$ for each $j=1,\ldots,n$,
hence $Q(x)$ vanishes at the zeros of $P_n(x)$. Therefore,
$$
\exists R\in\mathcal{P}_{n-1}\;:\;\; Q(x)=R(x)P_n(x)\;.
$$
Since $\{P_n\}_{n\geq0}$ is an OPS with respect to ${\bf u}$, we deduce
$$
\langle{\bf u},p\rangle=\langle{\bf u},Q+L_n\rangle=
\langle{\bf u},RP_n\rangle+\langle{\bf u},L_n\rangle=\langle{\bf u},L_n\rangle=
\sum_{j=1}^n p(x_{n,j})\,\langle{\bf u},\ell_{j,n}\rangle\;.
$$
Thus, setting
\begin{snugshade}
\begin{equation}\label{GJC-Anj}
A_{n,j}:=\langle{\bf u},\ell_{j,n}\rangle\;,\quad j=1,2,\ldots,n\;,
\end{equation}
\end{snugshade}\noindent
we obtain
\begin{equation}\label{GJC-up}
\langle{\bf u},p\rangle=\sum_{j=1}^nA_{n,j}p(x_{n,j})\;.
\end{equation}
Therefore, (\ref{GJC1}) will become proved provided we can show that the $A_{n,j}$'s
defined by (\ref{GJC-Anj}) are all positive numbers.
Indeed, taking $p(x)\equiv \ell_{j,n}^2(x)$ in (\ref{GJC-up})---notice that each $\ell_{j,n}$, $1\leq j\leq n$,
is a polynomial of degree $n-1$, hence $\ell_{j,n}^2\in\mathcal{P}_{2n-1}$---, and taking into account that ${\bf u}$ is positive-definite, we have
$$
0<\langle{\bf u},\ell_{j,n}^2\rangle=\sum_{k=1}^nA_{n,k}\ell_{j,n}^2(x_{n,k})
=\sum_{k=1}^nA_{n,k}\delta_{j,k}=A_{n,j}
$$
for each $j=1,\ldots,n$.
Notice also that the $A_{n,j}$'s defined by (\ref{GJC-Anj}) do not depend on $p$.
Hence, (\ref{GJC1}) is proved.
Finally, choosing $p(x)\equiv1$ in (\ref{GJC1}), we obtain (\ref{GJC2}).
\qed

\begin{remark}\em
Quadrature formulas are very useful tools in Numerical Analysis, e.g. for computing integrals by approximation.
Indeed, numerical quadrature consists of approximating the integral of a given integrable function $f$,
$$I[f]:=\int_\mathbb{R}f(x)\,{\rm d}\mu(x)\,,$$ with respect some positive Borel measure $\mu$,
by a finite sum which uses only the values of $f$ at $n$ points $t_j$ (called nodes),
$$I_n[f]:=\sum_{j=1}^nf(t_j)A_j\,,$$
where the coefficients $A_j$ (which may depend on $n$, as well as the notes $t_j$)
have to be chosen properly so that the quadrature formula is
correct,---i.e., the equality $I[f]=I_n[f]$ holds---, for as many functions $f$ as possible.
\end{remark}

\section*{Exercises}
%\bigskip

{\small
%\noindent
\begin{enumerate}[label=\emph{\bf \arabic*.},leftmargin=*]

\item\label{Ex-cp3-1}
%Let $E=\{x_1,\ldots,x_N\}$ be a set with $N$ distinct real numbers, and let $h_1,\ldots,h_N>0$.
Let $x_1,\ldots,x_N$ be any $N$ distinct real numbers, and let $h_1,\ldots,h_N>0$.
Define ${\bf u}\in\mathcal{P}'$ by
$$
\langle {\bf u},x^n\rangle:=\sum_{j=1}^N h_jx_j^n\,,\quad n\in\mathbb{N}_0\;.
$$
Prove that:
\begin{enumerate}
\item ${\bf u}$ is positive-definite on $E:=\{x_1,\ldots,x_N\}$;
\item ${\bf u}$ is not positive-definite on any set $S\subseteq\mathbb{R}$ such that $E$ is a proper subset of $S$.
\end{enumerate}
\smallskip

%\item\label{Ex-cp3-1a}
%Let $x_1,\ldots,x_N$ be any $N$ distinct real numbers, and let $h_1,\ldots,h_N>0$.
%Define ${\bf u}\in\mathcal{P}'$ by
%$$
%\langle {\bf u},x^n\rangle:=\sum_{j=1}^N h_jx_j^n\,,\quad n\in\mathbb{N}_0\;.
%$$
%Prove that:
%\begin{enumerate}
%\item ${\bf u}$ is positive-definite on $E:=\{x_1,\ldots,x_N\}$;
%\item ${\bf u}$ is not positive-definite on any set $S\subseteq\mathbb{R}$ such that $E$ is a proper subset of $S$.
%\end{enumerate}
%\medskip

%\noindent
%({\sl Hint.} Represent ${\bf u}$ as a Riemann-Stieltjes integral with respect to a certain step function.)
%\medskip

%\item\label{Ex-cp3-2a}
%Let $a,b,c\in\mathbb{R}$, with $bc>0$, and let $A_n$ be the tridiagonal Toeplitz matrix of order $n$
%\medskip
%$$
%A_n=\left(
%\begin{array}{cccccc}
%a & b &  &  &  &   \\
%c & a & b &  & &  \\
% & c & a & b &  &   \\
% & & \ddots & \ddots & \ddots &  \\
% & & & c & a & b  \\
% & & & & c & a \\
%\end{array}
%\right) \, .
%$$
%Use the Chebyshev polynomials of the second kind, $\{U_n\}_{n\geq0}$,
%%where $U_n(x):=\sin(n+1)\theta/\sin\theta$ $(x=\cos\theta)$,
%to prove that the eigenvalues of $A_n$ are explicitly given by
%$$
%\lambda_j:=a+\sqrt{bc}\,\cos\frac{j\pi}{n+1}\quad (j=1,2,\cdots, n)\, ,
%$$
%with corresponding eigenvectors
%$$
%v_j:=\big(\big)^T\;,\quad j=1,2,\cdots, n\,.
%$$
%\medskip

\item\label{Ex-cp3-2}
%Let $a,b,c\in\mathbb{R}$, with $bc>0$, and let $A_n$ be the tridiagonal Toeplitz matrix of order $n$
%\medskip
Let $a,b,c\in\mathbb{R}$, with $bc>0$. For each $n\in\mathbb{N}_0$, set
$$
P_n(x):=(bc)^{n/2}U_n\Big(\frac{x-a}{2\sqrt{bc}}\Big) \;,
$$
where $\{U_n\}_{n\geq0}$ is the sequence of the Chebyshev polynomials of the second kind.
\begin{enumerate}
\item
Show that $\{P_n\}_{n\geq0}$ is a monic OPS w.r.t. a %with respect to a
positive-definite functional ${\bf u}\in\mathcal{P}'$.
\item
Consider the tridiagonal Toeplitz matrix of order $n$
$$
A_n=\left(
\begin{array}{cccccc}
a & b &  &  &  &   \\
c & a & b &  & &  \\
 & c & a & b &  &   \\
 & & \ddots & \ddots & \ddots &  \\
 & & & c & a & b  \\
 & & & & c & a \\
\end{array}
\right) \, .
$$
Prove that the eigenvalues of $A_n$ are
$$
\lambda_j:=a+2\sqrt{bc}\,\cos\frac{j\pi}{n+1}\quad (j=1,2,\cdots, n)\, ,
$$
with corresponding eigenvectors
$$
%v_j:=\big(\big)^T\;,\quad j=1,2,\cdots, n\,.
v_j:=\frac{1}{\sin\frac{j\pi}{n+1}}
\left( \begin{array}{c}
\sin\frac{j\pi}{n+1} \\ [0.2em]
(c/b)^{1/2}\sin\frac{2j\pi}{n+1} \\ \vdots \\ [0.2em]
(c/b)^{(n-1)/2}\sin\frac{nj\pi}{n+1}
\end{array} \right)\;, \quad j=1,2,\cdots, n\,.
$$
\smallskip

%\noindent\hspace*{-1em}
({\sl Hint.} Define $Q_n(x):=b^{-n}P_n(x)$, and write the TTRR for $\{Q_n\}_{n\geq0}$ in matrix form.) %$P_n(x):=(bc)^{n/2}U_n\big(\frac{x-a}{2\sqrt{bc}}\big)$.)
\end{enumerate}
\medskip

\item\label{Ex-cp3-2}
Let ${\bf u}\in\mathcal{P}'$ be positive-definite. Let $\{P_n\}_{n\geq0}$ be the corresponding monic OPS
and $\{p_n\}_{n\geq0}$ an associated orthonormal sequence.
Denote by $x_{n1},\ldots,x_{nn}$ the zeros of $P_n(x)$ and let $\{\gamma_n\}_{n\geq1}$ be the sequence of $\gamma-$parameters appearing in the TTRR fulfilled by $\{P_n\}_{n\geq0}$.
Show that the ``weights'' $A_{nk}$ in the associated Gauss quadrature formula admit the following representations:
$$
A_{nk}=-\frac{u_0\gamma_1\gamma_2\cdots\gamma_n}{P_{n+1}(x_{nk})P_{n}'(x_{nk})}
=\Big\{\sum_{j=0}^{n}p_j^2(x_{nk})\Big\}^{-1}\quad (1\leq k\leq n\;;\; n\in\mathbb{N})\;.
$$
\smallskip

\noindent
({\sl Hint.} Use the Christoffel-Darboux identities.)
\medskip

\item\label{Ex-cp3-3}
Let ${\bf u}\in\mathcal{P}'$ be defined as
$$
\langle{\bf u},p\rangle:=\int_\mathbb{R}p(x)\,{\rm d}\mu(x)\;,\quad p\in\mathcal{P}\,,
$$
where $\mu$ is a positive Borel measure with infinite support\footnote{\,The {\sl support} of $\mu$ is the set $\text{supp}(\mu):=\big\{x\in\mathbb{R}\,:\, \mu\big((x-\epsilon,x+\epsilon)\big)>0\,,\,
\forall\epsilon>0\big\}\,$.}
and finite moments of all orders.
\begin{enumerate}
\item
Prove that ${\bf u}$ is positive-definite.
\item
Let $\{P_n\}_{n\geq0}$ be the monic OPS with respect to ${\bf u}$.
Prove that the maximum of the ratio
$$
\frac{\displaystyle\int_{\mathbb{R}}xQ_n^2(x)\,{\rm d}\mu(x)}{\displaystyle\int_\mathbb{R}Q_n^2(x)\,{\rm d}\mu(x)}
%\int_{\mathbb{R}}xQ_n^2(x)\,{\rm d}\mu(x)\Big/\int_\mathbb{R}Q_n^2(x)\,{\rm d}\mu(x)
$$
taken over all real polynomials $Q_n$ of degree at most $n$ is equal to the largest zero ${x_{n+1,n+1}}$ of the polynomial $P_{n+1}$, and the minimum is equal to the smallest zero ${x_{n+1,1}}$ of $P_{n+1}$.
%$$
%\min_{q\in\mathcal{P}_n(\mathbb{R})}\frac{\displaystyle\int_{\mathbb{R}}xq_n^2(x)\,{\rm d}\mu(x)}{\displaystyle\int_\mathbb{R}q_n^2(x)\,{\rm d}\mu(x)}=x_{n+1,1}
%%\int_{\mathbb{R}}xQ_n^2(x)\,{\rm d}\mu(x)\Big/\int_\mathbb{R}Q_n^2(x)\,{\rm d}\mu(x)
%$$
\item
Determine polynomials $Q_n$ where these maximum and minimum ratios are attained.
\end{enumerate}
%\medskip

%\item\label{Ex-cp3-4}
%Let $\{P_n\}_{n\geq0}$ be a monic OPS with respect to a positive Borel measure $\mu$.
%%with an infinite support and whose moments of all orders exist.
%Prove that each $P_n(x)$ admits the representation
%$$
%P_n(x)=\frac{1}{n!H_{n-1}} \int_\mathbb{R}\cdots\int_\mathbb{R}\; %\int_{\mathbb{R}^n}
%\prod_{i=1}^n(x-x_i)\prod_{1\leq j<k\leq n}^n(x_j-x_k)^2\,{\rm d}\mu(x_1)\cdots{\rm d}\mu(x_n)\;,
%$$
%where $H_{n-1}$ is the Hankel determinant of order $n$.
%\medskip
%
%{\sl Remark.} This formula goes back at least as far as Heine (from Heine-Borel), 1878.
%Nowadays it has important applications in \emph{Random Matrix Theory}.

\end{enumerate}
%\bigskip
}
\medskip

\section*{Final remarks}

The presentation of the topics considered in this text follows Chihara's book \ref{Chihara1978-C2}.
These topics may be found also in most books containing chapters on the general theory of OP.
In particular, they are treated (with more or less detail) in the books included in the bibliography.
Exercises {\bf 1} up to {\bf 3} appear in Chihara's book,
being the results contained therein proved in several textbooks appearing in the bibliography.
The result expressed by exercise {\bf 2} appears in useful applications of OP.
Concerning exercise {\bf 4}, see e.g. the article \ref{VanAssche1997-C3a} by W. Van Assche.
\medskip

%\section*{Notas finais}
%\section*{Coment\'arios finais}
%
%
%A apresenta\c c\~ao dos t\'opicos abordados neste texto segue o livro de Chihara \ref{Chihara1978-C2}.
%Trata-se, contudo, de assuntos que se encontram na maioria dos livros
%com cap\'\i tulos sobre a teoria geral dos PO. Em particular, s\~ao abordados, com maior ou menor detalhe,
%na bibliografia indicada no final do texto.
%Os exerc\'\i cios {\bf 1} e {\bf 3} aparecem no livro de Chihara,
%sendo os resultados neles contidos provados em v\'arios dos textos da bibliografia.
%O resultado expresso pelo exerc\'\i cios {\bf 2} encontra-se em muitas aplica\c c\~oes dos PO.
%No tocante ao exerc\'\i cios {\bf 4}, veja-se, e.g., o artigo \ref{VanAssche1997-C3a} de W. Van Assche.

\section*{Bibliography}
\medskip

{\small
\begin{enumerate}[label=\emph{\rm [\arabic*]},leftmargin=*]
\item\label{Chihara1978-C3a} T. S. Chihara, {\sl An introduction to orthogonal polynomials}, Gordon and Breach (1978).
\item\label{Freud1971-C3a} G. Freud, {\sl Orthogonal polynomials}, Pergamon Press, Oxford (1971).
\item\label{Gautschi2004-C3a} W. Gautschi, {\sl Orthogonal polynomials. Computation and approximation}, Oxford University Press, Oxford (2004).
\item\label{Ismail2004-C3a} M. E. H. Ismail, {\sl Classical and quantum orthogonal polynomials in one variable},
          Cambridge University Press (2005) [paperback edition: 2009].
\item\label{Szegoo1975-C3a} G. Szeg\"o, {\sl Orthogonal polynomials}, AMS Colloq. Publ. {\bf 230} (1975), 4th ed.
\item\label{VanAssche1997-C3a} W. Van Assche, {\it Orthogonal polynomials in the complex plane and on the real line}, {\sl Special functions, $q-$series and related topics (Toronto, ON, 1995)}, Fields Inst. Commun. {\bf 14}, AMS (1997) 211-245.
%M. E. H. Ismail, D. R. Masson, and M. Rahman Eds., Fields Institute Communications {\bf 14} (1997) 211-245.
\end{enumerate}
}

\chapter{The spectral theorem for orthogonal polynomials}

\pagestyle{myheadings}\markright{The spectral theorem for orthogonal polynomials}
\pagestyle{myheadings}\markleft{J. Petronilho}

%\section{A brief review on the Riemann-Stieltjes integral}

%\section{A brief review on positive Borel measures}

Here we still concentrates in the study of OPS with respect to positive-definite
moment linear functionals ${\bf u}\in\mathcal{P}'$.
Our aim is to prove that any such functional admits an integral representation involving
a positive Borel measure $\mu$ on $\mathbb{R}$ (which needs not to be unique)
with infinite support and such that all its moments exist. %(i.e., they are finite).

\section{Helly's theorems}

In this section we state some preliminary results needed for the proof of the
representation theorem to be stated in the next section.

\begin{snugshade}
\begin{lemma}\label{LammaCountableContinuous}
Let $I\subseteq\mathbb{R}$ be an interval and let $f:I\to\mathbb{R}$ be a monotone function.
Then, $f$ has at most countably many discontinuity points.
\end{lemma}
\end{snugshade}

{\it Proof.}
This is a well known result in Real Analysis.
A nice proof can be found e.g. in the book \ref{Leoni2009-C3} by G. Leoni.
\qed

\begin{snugshade}
\begin{lemma}\label{LammaHelly}
Let $\{f_n\}_{n\geq1}$ be a sequence of real functions defined on a countable set $E$.
Suppose that $\{f_n(x)\}_{n\geq1}$ is a bounded sequence for each $x\in E$.
Then $\{f_n\}_{n\geq1}$ contains a subsequence $\{f_{n_j}\}_{j\geq1}$ that converges everywhere on $E$,
i.e., the (sub)sequence $\{f_{n_j}(x)\}_{j\geq1}$ converges for each $x\in E$.
\end{lemma}
\end{snugshade}

{\it Proof.}
Set $E:=\{x_1,x_2,x_3,\ldots\}$ and write $f_n^{(0)}\equiv f_n$.
Since $\{f_n^{(0)}(x_1)\}_{n\geq1}$ is a bounded sequence of real numbers, it contains a convergent subsequence,
i.e., there exists a subsequence $\{f_n^{(1)}\}_{n\geq1}$ of $\{f_n^{(0)}\}_{n\geq1}$ such that
$\{f_n^{(1)}(x)\}_{n\geq1}$ converges for $x=x_1$.
Now, since $\{f_n^{(1)}(x_2)\}_{n\geq1}$ is a bounded sequence, it contains a convergent subsequence, hence,
there exists a subsequence $\{f_n^{(2)}\}_{n\geq1}$ of $\{f_n^{(1)}\}_{n\geq1}$ such that
$\{f_n^{(2)}(x)\}_{n\geq1}$ converges for $x=x_2$.
Proceeding in this way, we obtain sequences
$$
\{f_n^{(0)}\}_{n\geq1}\,,\quad \{f_n^{(1)}\}_{n\geq1}\,,\quad \{f_n^{(2)}\}_{n\geq1}\,,\quad\ldots,\quad
\{f_n^{(k)}\}_{n\geq1}\;,\quad\ldots
$$
such that:
\smallskip

(i) $\{f_n^{(k)}\}_{n\geq1}$ is a subsequence of $\{f_n^{(k-1)}\}_{n\geq1}$, for each $k=1,2,3,\ldots$;
\smallskip

(ii) $\{f_n^{(k)}(x)\}_{n\geq1}$ converges for each $x\in E_k:=\{x_1,x_2,\ldots,x_k\}$.
\smallskip

\noindent
It follows from (i)---with a little care (passing to a subsequence if necessary)
to ensure that the relative order of terms is preserved---that the diagonal
sequence, $\{f_n^{(n)}\}_{n\geq1}$, is also a subsequence of $\{f_n\}_{n\geq1}$.
Since, for each $k\in\mathbb{N}$, except for the first $k-1$ terms,
$\{f_n^{(n)}\}_{n\geq1}$ is also a subsequence of $\{f_n^{(k)}\}_{n\geq1}$,
it follows from (ii) that
$$
\{f_n^{(n)}(x)\}_{n\geq1}\;\;\mbox{\rm converges for each}\;\; x\in\cup_{k=1}^\infty E_k=E\;.
$$
Therefore, since, for each $x\in E$, $\{f_n^{(n)}(x)\}_{n\geq1}$ is a subsequence of $\{f_n(x)\}_{n\geq1}$,
the proof is concluded.
\qed

\begin{snugshade}
\begin{theorem}[Helly's selection principle]\label{Helly1}
Let  $\{\phi_n\}_{n\geq1}$ be a uniformly bounded sequence of nondecreasing functions defined on $\mathbb{R}$.
Then, $\{\phi_n\}_{n\geq1}$ has a subsequence which converges on $\mathbb{R}$ to a bounded and nondecreasing function.
\end{theorem}
\end{snugshade}

{\it Proof.}
Consider the set of rational numbers, $\mathbb{Q}$.
According to Lemma \ref{LammaHelly}, there is a subsequence $\{\phi_{n_k}\}_{k\geq1}$
which converges everywhere on $\mathbb{Q}$.
Henceforth, we may define a function $\Phi_1:\mathbb{Q}\to\mathbb{R}$ as
\begin{equation}\label{phiHy1}
\Phi_1(r):=\lim_{k\to+\infty}\phi_{n_k}(r)\;,\quad r\in\mathbb{Q}\;.
\end{equation}
It follows from the hypothesis on $\{\phi_n\}_{n\geq1}$ that $\Phi_1$ is bounded and nondecreasing on $\mathbb{Q}$.
We now extend the domain of $\Phi_1$ to $\mathbb{R}$ by defining $\Phi_2:\mathbb{R}\to\mathbb{R}$ as
\begin{equation}\label{phiHy2}
\Phi_2(x):=\left\{
\begin{array}{ccl}
\Phi_1(x) & \mbox{\rm if} & x\in\mathbb{Q}\; , \\ [0.5em]
\displaystyle\sup_{\substack{r\in\mathbb{Q} \\ r<x}}\Phi_1(r) & \mbox{\rm if} & x\in\mathbb{R}\setminus\mathbb{Q}\;.
\end{array}\right.
\end{equation}
This function $\Phi_2$ is clearly bounded and nondecreasing on $\mathbb{R}$
(since the same properties are fulfilled by $\Phi_1$).
According to (\ref{phiHy1}), $\{\phi_{n_k}\}_{k\geq1}$
converges to $\Phi_2(x)$ at each point $x\in\mathbb{Q}$.
Next we show that $\{\phi_{n_k}\}_{k\geq1}$ also
converges to $\Phi_2(x)$ at each point $x$ where $\Phi_2$ is continuous.
To this end, suppose that $\Phi_2$ is continuous at the point $x\in\mathbb{R}\setminus\mathbb{Q}$.
Since $\mathbb{Q}$ is a dense subset of $\mathbb{R}$ and $\Phi_2$ is continuous at $x$, then
\begin{equation}\label{phiHy3}
\forall\epsilon>0\;,\;\exists x_2\in\mathbb{Q}\;:\quad
x<x_2\;\;\wedge\;\; \Phi_2(x_2)<\Phi_2(x)+\epsilon\;.
\end{equation}
Fix arbitrarily $x_1\in\mathbb{Q}$, with $x_1<x$.
Then, since (by hypothesis) $\phi_{n_k}$ is a nondecreasing function on $\mathbb{R}$, we have
$$%\begin{equation}\label{phiHy4}
\phi_{n_k}(x_1)\leq\phi_{n_k}(x)\leq\phi_{n_k}(x_2)\;.
$$%\end{equation}
Therefore, we deduce
$$
\begin{array}{rcl}
\Phi_2(x_1)&=&\displaystyle\Phi_1(x_1)=\lim_{k\to+\infty}\phi_{n_k}(x_1)=\liminf_{k\to+\infty}\phi_{n_k}(x_1)\leq
\liminf_{k\to+\infty}\phi_{n_k}(x) \\ [0.75em]
&\leq&\displaystyle\limsup_{k\to+\infty}\phi_{n_k}(x)\leq\displaystyle\limsup_{k\to+\infty}\phi_{n_k}(x_2)
=\lim_{k\to+\infty}\phi_{n_k}(x_2)=\Phi_1(x_2)=\Phi_2(x_2) \\ [1em]
&<&\displaystyle\Phi_2(x)+\epsilon\;.
\end{array}
$$
Summarizing, we proved that, if $\Phi_2$ is continuous at a point $x\in\mathbb{R}\setminus\mathbb{Q}$, then
%$$
%\forall x_1\in\mathbb{Q}\;,\quad x_1<x\;\;\Rightarrow\;\;
%\Phi_2(x_1)\leq\Phi_2(x_2)\;.
%$$
%Summarizing, we proved that for each $x_1\in\mathbb{Q}$ such that $x_1<x$,
$$%\begin{equation}\label{phiHy5}
\forall x_1\in\mathbb{Q}\;,\quad x_1<x\;\;\Rightarrow\;\;
\Phi_2(x_1)\leq\liminf_{k\to+\infty}\phi_{n_k}(x)\leq\limsup_{k\to+\infty}\phi_{n_k}(x)<\Phi_2(x)+\epsilon\;.
$$%\end{equation}
%Since (\ref{phiHy5}) holds for all $x_1\in\mathbb{Q}$ such that $x_1<x$, and
Therefore, %since $\Phi_2$ is continuous at $x$, we may write
$$
\Phi_2(x)=\lim_{\substack{x_1\to x^-\\  x_1\in\mathbb{Q}\;\;}}\Phi_2(x_1)\leq\liminf_{k\to+\infty}\phi_{n_k}(x)\leq
\limsup_{k\to+\infty}\phi_{n_k}(x)<\Phi_2(x)+\epsilon\;,
$$
%$$
%\begin{array}{rcl}
%\Phi_2(x)&=&\displaystyle\lim_{x_1\to x^-}\Phi_2(x_1)\leq\lim_{x_1\to x^-}\Big(
%\liminf_{k\to+\infty}\phi_{n_k}(x)\Big)=\liminf_{k\to+\infty}\phi_{n_k}(x) \\ [0.75em]
%&\leq&\displaystyle\limsup_{k\to+\infty}\phi_{n_k}(x)\leq\Phi_2(x_2)\leq\Phi_2(x)+\epsilon\;,
%\end{array}
%$$
hence, since $\epsilon>0$ is arbitrary, we deduce (letting $\epsilon\to0^+$),
$$
\Phi_2(x)\leq\liminf_{k\to+\infty}\phi_{n_k}(x)\leq
\limsup_{k\to+\infty}\phi_{n_k}(x)\leq\Phi_2(x)\;.
$$
Since the left-hand side and the right-hand side coincide, these inequalities are indeed equalities, hence
$$
\lim_{k\to+\infty}\phi_{n_k}(x)=\Phi_2(x)\;.
$$
Thus indeed $\{\phi_{n_k}(x)\}_{k\geq1}$ converges to $\Phi_2(x)$
at each point $x\in\mathbb{R}$ of continuity of $\Phi_2$.

Now, observe that $\Phi_2$ is a nondecreasing function, so (by Lemma \ref{LammaCountableContinuous})
the set of its points of discontinuity form an at most countable set.
Denote by $D$ the set of points of discontinuity of $\Phi_2$ which does not belong to $\mathbb{Q}$.
Applying Lemma \ref{LammaHelly} to $\{\phi_{n_k}\}_{k\geq1}$ and $D$, we deduce that there is a subsequence $\{\phi_{n_{k_j}}\}_{j\geq1}$ of $\{\phi_{n_k}\}_{k\geq1}$ which converges everywhere on $D$ to a limit function $\Phi_3:D\to\mathbb{R}$. Finally, define $\phi:\mathbb{R}\to\mathbb{R}$ by
$$
\phi(x):=\left\{
\begin{array}{ccl}
\Phi_2(x) & \mbox{\rm if} & x\in\mathbb{R}\setminus D\; , \\ [0.25em]
\Phi_3(x) & \mbox{\rm if} & x\in D\;.
\end{array}\right.
$$
It is clear that
$$
\lim_{j\to+\infty}\phi_{n_{k_j}}(x)=\phi(x)\;,\quad\forall x\in\mathbb{R}\;.
$$
Moreover, the hypothesis on $\{\phi_{n}\}_{n\geq1}$ ensure that $\phi$ is bounded and nondecreasing.
In conclusion: $\{\phi_{n}\}_{n\geq1}$ has a convergent subsequence $\{\phi_{n_{k_j}}\}_{j\geq1}$ which
converges to a bounded and nondecreasing function $\phi$ on $\mathbb{R}$.
The proof is complete.
\qed
\medskip

%$f(x) \overset{\text{def}}{=} x \ln(1+x)$

The next theorem involves the Riemann-Stieltjes integral.
%Appendix \ref{Appx2} contains a review on this integral,
%including most properties needed on this course.
The needed facts concerning this integral can be found e.g.
in the book \ref{KolmogorovFomini1968-A2} by Kolmogorov and Fomini.

\begin{snugshade}
\begin{theorem}[Helly's convergence theorem]\label{Helly2}
Let $\{\phi_n\}_{n\geq1}$ be a uniformly bounded sequence of nondecreasing functions defined on
a compact interval $[a,b]$, and suppose that this sequence converges on $[a,b]$ to a limit function $\phi$, so that
\begin{equation}\label{Hyp-phi1}
\phi(x):=\lim_{n\to\infty}\phi_n(x)\,,\quad x\in[a,b]\;.
\end{equation}
Then, for each continuous function $f:[a,b]\to\mathbb{R}$, the following holds:
\begin{equation}\label{Hyp-phi2}
\lim_{n\to\infty}\int_a^bf(x)\,{\rm d}\phi_n(x)=\int_a^bf(x)\,{\rm d}\phi(x)\;.
\end{equation}
\end{theorem}
\end{snugshade}

{\it Proof.}
Since $\{\phi_n\}_{n\geq1}$ is uniformly bounded and each $\phi_n$ is nondecreasing on $[a,b]$,
$$
\exists M>0\;:\;\forall n\in\mathbb{N}\;,\quad 0\leq\phi_n(b)-\phi_n(a)\leq M\; ,
$$
hence, by hypothesis (\ref{Hyp-phi1}), also
\begin{equation}\label{phi1}
0\leq\phi(b)-\phi(a)\leq M\; .
\end{equation}
Fix $\epsilon>0$.
Since, by hypothesis, $f$ is a continuous function on the compact set $[a,b]$,
then $f$ is uniformly continuous on $[a,b]$, hence there is a partition $P_\epsilon$ of $[a,b]$,
$$
%P_\epsilon:=\{x_0,x_1,\cdots,x_m\}
P_\epsilon\quad:\quad a=x_0<x_1<\cdots<x_{i-1}<x_i<\cdots<x_{m-1}<x_m=b\;,
$$
such that
\begin{equation}\label{Pexz}
x,z\in[x_{i-1},x_i]\quad\Rightarrow\quad\big|f(x)-f(z)\big|<\frac{\epsilon}{2M}\,,\quad \forall i\in\{1,2,\ldots,m\}\;.
\end{equation}
[Indeed, being $f$ uniformly continuous on $[a,b]$, this means that
$$
\forall\epsilon>0\;,\;\exists\delta=\delta(\epsilon)>0\;:\;\forall x,z\in[a,b]\;,\quad
|x-z|<\delta\;\;\Rightarrow\;\;\big|f(x)-f(z)\big|<\epsilon/(2M)\;;
$$
hence we choose the partition $P_\epsilon$
so that $|x_i-x_{i-1}|<\delta$ for each $i\in\{1,2,\ldots,m\}$.]
Now, for each $i\in\{1,2,\ldots,m\}$, choose an ``intermediate point''
$\xi_i\in[x_{i-1},x_i]$, and set
$$
\Delta_i\phi:=\phi(x_i)-\phi(x_{i-1})\;,\quad\Delta_i\phi_n:=\phi_n(x_i)-\phi_n(x_{i-1})\;.
$$
By the Mean Value theorem for the Riemann-Stieltjes integral, %(Theorem \ref{propret-RSintegral}),
for each $i\in\{1,2,\ldots,m\}$, there exists $\xi_i'\in[x_{i-1},x_i]$ such that
$$
\int_{x_{i-1}}^{x_i}f(x)\,{\rm d}\phi(x)=f(\xi_i')\big(\phi(x_i)-\phi(x_{i-1})\big)
=\big(f(\xi_i')-f(\xi_i)+f(\xi_i)\big)\Delta_i\phi\;,
$$
hence
$$
\int_{x_{i-1}}^{x_i}f(x)\,{\rm d}\phi(x)-f(\xi_i)\Delta_i\phi=
\big(f(\xi_i')-f(\xi_i)\big)\Delta_i\phi\;,\quad i\in\{1,2,\ldots,m\}\;.
$$
Summing over $i$ and taking into account that
$\sum_{i=1}^m\int_{x_{i-1}}^{x_i}f(x)\,{\rm d}\phi(x)=\int_a^bf(x)\,{\rm d}\phi(x)$,
and then applying the triangular inequality, we obtain
%$$
%\int_a^bf(x)\,{\rm d}\phi(x)-\sum_{i=1}^mf(\xi_i)\Delta_i\phi
%=\sum_{i=1}^m\big(f(\xi_i')-f(\xi_i)\big)\Delta_i\phi\;,
%$$
%hence, taking into account (\ref{Pexz}),
$$
\Big|\int_a^bf(x)\,{\rm d}\phi(x)-\sum_{i=1}^mf(\xi_i)\Delta_i\phi\Big|
\leq\sum_{i=1}^m\big|f(\xi_i')-f(\xi_i)\big|\Delta_i\phi\;,
$$
hence, since $\xi_i,\xi_i'\in[x_{i-1},x_i]$, so that we may apply (\ref{Pexz}), we deduce
\begin{equation}\label{ineqphi1}
\Big|\int_a^bf(x)\,{\rm d}\phi(x)-\sum_{i=1}^mf(\xi_i)\Delta_i\phi\Big|
<\frac{\epsilon}{2M}\sum_{i=1}^m\Delta_i\phi=\frac{\epsilon}{2M}\big(\phi(b)-\phi(a)\big)
\leq\frac{\epsilon}{2}\;,
\end{equation}
where the last inequality holds by (\ref{phi1}).
In the same way, replacing $\phi$ by $\phi_n$ in the previous reasoning, we deduce
\begin{equation}\label{ineqphi2}
\Big|\int_a^bf(x)\,{\rm d}\phi_n(x)-\sum_{i=1}^mf(\xi_i)\Delta_i\phi_n\Big|<\frac{\epsilon}{2}\;.
\end{equation}
Next, observe that
$$
\begin{array}{l}
\displaystyle\Big|\int_a^bf(x)\,{\rm d}\phi(x)-\int_a^bf(x)\,{\rm d}\phi_n(x)\Big|\leq
\Big|\int_a^bf(x)\,{\rm d}\phi(x)-\sum_{i=1}^mf(\xi_i)\Delta_i\phi\Big| \\ [1em]
\qquad\displaystyle+\Big|\sum_{i=1}^mf(\xi_i)\Delta_i\phi-\sum_{i=1}^mf(\xi_i)\Delta_i\phi_n\Big|
+\Big|\int_a^bf(x)\,{\rm d}\phi_n(x)-\sum_{i=1}^mf(\xi_i)\Delta_i\phi_n\Big|\;.
\end{array}
$$
Therefore, taking into account (\ref{ineqphi1}) and (\ref{ineqphi2}), and noticing that
$$
\Big|\sum_{i=1}^mf(\xi_i)\Delta_i\phi-\sum_{i=1}^mf(\xi_i)\Delta_i\phi_n\Big|
\leq\sum_{i=1}^m|f(\xi_i)|\,\big|\Delta_i\phi-\Delta_i\phi_n\big|\;,
%\leq\sum_{i=1}^m|f(\xi_i)|\,\big|\Delta_i\phi-\Delta_i\phi_n\big|
$$
we obtain
\begin{equation}\label{ineqphi3}
\Big|\int_a^bf(x)\,{\rm d}\phi(x)-\int_a^bf(x)\,{\rm d}\phi_n(x)\Big|
<\epsilon+\sum_{i=1}^m|f(\xi_i)|\,\big|\Delta_i\phi-\Delta_i\phi_n\big|\;.
\end{equation}
Keeping $P_\epsilon$ fixed, we have
$$
\lim_{n\to+\infty}\big(\Delta_i\phi-\Delta_i\phi_n\big)=
\lim_{n\to+\infty}\big\{\big(\phi(x_i)-\phi_n(x_i)\big)-\big(\phi(x_{i-1})-\phi_n(x_{i-1})\big)\big\}=0\;,
$$
where the last equality follows from (\ref{Hyp-phi1}), hence we conclude from (\ref{ineqphi3}) that
\begin{equation}\label{ineqphi4}
\limsup_{n\to+\infty}\Big|\int_a^bf(x)\,{\rm d}\phi(x)-\int_a^bf(x)\,{\rm d}\phi_n(x)\Big|
\leq\epsilon\;.
\end{equation}
Since $\epsilon$ is positive and arbitrary, the $\limsup$ in (\ref{ineqphi4}) must be equal to zero.
It turns out that in (\ref{ineqphi4}) we may replace the $\limsup$ by the limit,
and since this limit is equal to zero, we conclude that (\ref{Hyp-phi2}) holds.
\qed

\section{The representation theorem}

We are ready to state the important representation theorem
for a positive-definite functional ${\bf u}\in\mathcal{P}'$,
showing that such a functional admits an integral representation
as a Riemann-Stieltjes integral with respect to a real bounded nondecreasing function on $\mathbb{R}$
fulfilling some natural conditions (namely, finite moments of all orders, and infinite spectrum).
We begin by introducing some useful concepts.

%Let ${\bf u}\in\mathcal{P}'$ be positive-definite, and let
%$\{P_n\}_{n\geq0}$ be the monic OPS with respect to ${\bf u}$.
%We know that each $P_N(x)$ is a real polynomial having $N$ real simple zeros, so that
%$$
%x_{N,1}<x_{N,2}<\cdots<x_{N,N}\;.
%$$
%For each $N\in\mathbb{N}$, we introduce a distribution function $\psi_N:\mathbb{R}\to\mathbb{R}$,
%characterized as being a right continuous step function with support $\{x_{N,1},x_{N,2},\ldots,x_{N,N}\}$,
%%$\big\{x_{N,j}:1\leq j\leq N\big\}$ (the set of zeros of $P_N(x)$),
%having a jump $\rho(x_{N,j})$ at the $j$th zero $x_{N,j}$ of $P_N(x)$, so that
%\begin{snugshade}
%\begin{equation}\label{dist-psiN}
%\begin{array}{l}
%\psi_N(-\infty):=0\;,\quad \psi_N(+\infty):=u_0\;, \\ [0.5em]
%\psi_N(x_{N,j})-\psi_N(x_{N,j}-0)=\rho(x_{N,j})\;,\quad 1\leq j\leq N\,.
%\end{array}
%\end{equation}
%%\begin{equation}\label{dist-psiN}
%%\psi_N(-\infty):=0\;,\quad \psi_N(x_{N,j})-\psi_N(x_{N,j}-0)=\rho(x_{N,j})\;,\quad
%%\psi_N(+\infty):=u_0\;,\quad 1\leq j\leq N\,.
%%\end{equation}
%\end{snugshade}
%
%\begin{snugshade}
%\begin{lemma}\label{LemSpecThmOP}
%Let ${\bf u}\in\mathcal{P}'$ be positive-definite.
%Fix $j,N\in\mathbb{N}$, with $1\leq j\leq 2N-2$,
%and let $\psi_N$ be the distribution function $(\ref{dist-psiN})$.
%Then, the moment of order $j$ of $\psi_N$,
%\begin{equation}\label{mom1}
%\int_\mathbb{R}x^j\,{\rm d}\psi_N(x)\;,
%\end{equation}
%only depends on $h_0,h_1,\ldots,h_{\lfloor(j+1)/2\rfloor}$, being
%$$
%h_n:=\prod_{k=0}^n\gamma_k\,,\quad n\in\mathbb{N}_0\;,
%$$
%with the convention $\gamma_0:=u_0$.
%\end{lemma}
%\end{snugshade}

\begin{snugshade}
\begin{definition}\label{dist-func}
A function $\psi:\mathbb{R}\to\mathbb{R}$ is called a {\sl distribution function} if it is
bounded, nondecreasing and all its moments
$$
\int_{-\infty}^{+\infty}x^n\,{\rm d}\psi(x)\;,\quad n=0,1,2,\ldots
$$
are finite. The {\sl spectrum} of a distribution function $\psi$ is the set
$$
\sigma(\psi):=\big\{ x\in\mathbb{R}\,:\,\psi(x+\delta)-\psi(x-\delta)>0\;,\;\forall\delta>0 \big\}\;.
$$
\end{definition}
\end{snugshade}

\begin{remark}\em
Often, being $\psi$ a distribution function,
a point in $\sigma(\psi)$ is called a {\sl spectral point}, or an {\sl increasing point} of $\psi$.
\end{remark}

\begin{snugshade}
\begin{theorem}\label{sigma-psi-closed}
%If $\psi$ is a distribution function, then $\sigma(\psi)$ is a closed set in $\mathbb{R}$.
The spectrum $\sigma(\psi)$ of a distribution function $\psi$ is closed in $\mathbb{R}$.
\end{theorem}
\end{snugshade}

{\it Proof.}
We will prove that $\mathbb{R}\setminus\sigma(\psi)$ is an open set.
Let $x_0\in\mathbb{R}\setminus\sigma(\psi)$. Then
$$
\exists\delta>0\;:\;\; \psi(x_0+\delta)-\psi(x_0-\delta)\leq0\;.
$$
Since $\psi$ is nondecreasing, $\psi$ must be constant ($=C$)
on the interval $(x_0-\delta,x_0+\delta)$.
Therefore, we see that if $x\in\big(x_0-\delta/2,x_0+\delta/2\big)$,
then there is $\delta'>0$ (choose $\delta'$ such that $0<\delta'<\delta/2$) such that
$\psi(x+\delta')-\psi(x-\delta')=C-C=0$,
hence $x\not\in\sigma(\psi)$. Thus,
$$
\Big(x_0-\mbox{$\frac{\delta}{2}$},x_0+\mbox{$\frac{\delta}{2}$}\Big)
\subseteq\mathbb{R}\setminus\sigma(\psi)\;,
$$
so that $x_0$ is an interior point of $\mathbb{R}\setminus\sigma(\psi)$.
Since $x_0$ was arbitrarily fixed on the set $\mathbb{R}\setminus\sigma(\psi)$,
we conclude that this set is open in $\mathbb{R}$.
\qed
\medskip

Let ${\bf u}\in\mathcal{P}'$ be positive-definite and
$\{P_n\}_{n\geq0}$ the monic OPS with respect to ${\bf u}$.
%We know
By Theorems \ref{GramShmidt} and \ref{OPS-zeros-Thm2},
each $P_n(x)$ is a real polynomial having $n$ real simple zeros:
$$
x_{n,1}<x_{n,2}<\cdots<x_{n,n}\;.
$$
For each $n\in\mathbb{N}$, introduce a distribution function $\psi_n:\mathbb{R}\to\mathbb{R}$,
characterized as being a bounded and right continuous step function with spectrum $\sigma(\psi_n)=\{x_{n,1},\ldots,x_{n,n}\}$, and
having jump $A_{n,j}>0$ at the $j$th zero $x_{n,j}$ of $P_n(x)$, where the $A_{n,j}$'s
($1\leq j\leq n$) are the weights appearing in Gauss quadrature formula, so that
\begin{snugshade}
\begin{equation}\label{dist-psiN}
\begin{array}{l}
\psi_n(-\infty):=0\;,\quad \psi_n(+\infty):=u_0\;, \\ [0.5em]
\psi_n(x_{n,j})-\psi_n(x_{n,j}-0)=A_{n,j}\;,\quad 1\leq j\leq n\,.
\end{array}
\end{equation}
\end{snugshade}\noindent
Explicitly, we may write
\begin{snugshade}
\begin{equation}\label{dist-psiN-bis}
\psi_n(x):=\left\{\begin{array}{ccl}
0 & \mbox{\rm if} & x<x_{n,1}\;; \\ [0.5em]
A_{n,1}+\cdots+A_{n,j} & \mbox{\rm if} & x_{n,j}\leq x<x_{n,j+1} \;,\quad 1\leq j\leq n-1\,;\\ [0.5em]
u_0 & \mbox{\rm if} & x\geq x_{n,n}\;.
\end{array}\right.
\end{equation}
\end{snugshade}\noindent
Then, for each fixed $k\in\{0,1,\ldots,2n-1\}$, by using the Gauss-Jacobi-Christoffel quadrature formula
(Theorem \ref{GaussJC}) applied to the polynomial $p(x):=x^k$,
 %and then taking into account Theorem \ref{propret-RSintegral}-(v),
one sees that the moment $u_k:=\langle {\bf u},x^k\rangle$
may be represented as a Riemann-Stieltjes integral with respect to $\psi_n$ as:
\begin{snugshade}
\begin{equation}\label{GaussUk1}
u_k=\sum_{j=1}^nA_{n,j}x_{n,j}^k=\int_{-\infty}^{+\infty}x^k\,{\rm d}\psi_n(x)\;,\quad k=0,1,\ldots,2n-1\,.
\end{equation}
\end{snugshade}\noindent
Now, by Helly's selection principle, there is a subsequence $\{\psi_{n_j}\}_{j\geq0}$ of $\{\psi_n\}_{n\geq0}$
which converges on $\mathbb{R}$ to a bounded and nondecreasing function $\psi$:
%$\psi:\mathbb{R}\to\mathbb{R}$,
\begin{snugshade}
\begin{equation}\label{Psi=limPsinj}
\psi(x):=\lim_{j\to\infty}\psi_{n_j}(x)\;,\quad x\in\mathbb{R}\,.
\end{equation}
\end{snugshade}
%\noindent

\begin{snugshade}
\begin{definition}\label{Psi=NaturalRep}
Let ${\bf u}\in\mathcal{P}'$ be positive-definite.
A function $\psi:\mathbb{R}\to\mathbb{R}$ defined as in $(\ref{Psi=limPsinj})$ ---limit of a subsequence
of the step functions $(\ref{dist-psiN})$--- is called a {\sl natural representative} for ${\bf u}$.
\end{definition}
\end{snugshade}\noindent

%Notice that a natural representative for ${\bf u}$ is a bounded,
%nondecreasing and right continuous function which approaches zero as $x\to-\infty$
%(so that it is a distribution function).

\begin{remark}\em
A natural representative $\psi$ for ${\bf u}$ is a distribution function.
Indeed, as noted above, $\psi$ is bounded and nondecreasing.
Moreover, all the moments
$$
\int_{-\infty}^{+\infty}x^n\,{\rm d}\psi(x)\;,\quad n=0,1,2,\ldots
$$
are finite, as follows by the representation Theorem \ref{RepThmOP},
to be proved next.
\end{remark}

%A natural representative $\psi$ provides ${\bf u}$ with an integral representation
%expressed in terms of the Riemann-Stieltjes integral with respect to $\psi$.

\begin{snugshade}
\begin{theorem}[representation theorem for positive-definite functionals on $\mathcal{P}$]\label{RepThmOP}
Let ${\bf u}\in\mathcal{P}'$ be positive-definite. %, and let $\psi_n$ be defined by $(\ref{})$ for each $n\geq0$.
Then, there is a natural representative of ${\bf u}$, $\psi:\mathbb{R}\to\mathbb{R}$, whose spectrum is an infinite set, such that
%Then, there exists a subsequence of $\{\psi_n\}_{n\geq1}$ which converges on $\mathbb{R}$ to a distribution function $\psi$ (bounded, nondecreasing and right-continuous), whose spectrum is an infinite set, fulfilling
\begin{equation}\label{u-rep-Thm}
\langle{\bf u},p\rangle=\int_{-\infty}^{+\infty}p(x)\,{\rm d}\psi(x)\;,\quad p\in\mathcal{P}\;.
\end{equation}
\end{theorem}
\end{snugshade}

{\it Proof.}
We consider two cases.
\smallskip

{\sl Case 1}. Assume that the true interval of orthogonality $[\xi,\eta]$ is bounded (compact).
Then (cf. Remark \ref{TIZ}) from (\ref{dist-psiN}) and (\ref{Psi=limPsinj}) we see that $\psi(x)=0$ if $x<\xi$, and
$\psi(x)=u_0>0$ if $x>\eta$. Therefore, for each (fixed) $k\in\mathbb{N}_0$ we may write
\begin{equation}\label{GaussUk2}
\int_{-\infty}^{+\infty}x^k\,{\rm d}\psi(x)=
\int_{\xi}^{\eta}x^k\,{\rm d}\psi(x)=
\lim_{j\to\infty}\int_{\xi}^{\eta}x^k\,{\rm d}\psi_{n_j}(x)\;,
\end{equation}
where the last equality holds by Helly's convergence theorem.
Keeping $k$ fixed, and since $k\leq 2 n_j-1$ for $j$ sufficiently large,
we deduce from (\ref{GaussUk1}) that the limit in (\ref{GaussUk2}) equals $u_k$, and so
\begin{equation}\label{GaussUk3}
\int_{-\infty}^{+\infty}x^k\,{\rm d}\psi(x)=u_k=\langle{\bf u},x^k\rangle
\;,\quad k\in\mathbb{N}_0\,.
\end{equation}
Thus, (\ref{u-rep-Thm}) follows whenever $[\xi,\eta]$ is bounded.
\smallskip

{\sl Case 2}. Assume now that $[\xi,\eta]$ is unbounded.\footnote{\,In this case
Helly's convergence theorem cannot be applied (Exercise \ref{Ex-cp4-1a}).}
By Helly's selection theorem, there exists a subsequence
$\{\psi_{n_i}\}_{i\geq0}$ of $\{\psi_n\}_{n\geq0}$ which converges
on $\mathbb{R}$ to a bounded and nondecreasing function $\psi$.
Setting $\phi_i:=\psi_{n_i}$, according with (\ref{GaussUk1}) we deduce
\begin{equation}\label{GaussUk4}
\int_{-\infty}^{+\infty}x^k \,{\rm d}\phi_i(x)=u_k \quad\mbox{\rm if}\quad
n_i\geq\frac{k+1}{2}\; ,\quad k=0,1,2,\cdots\; .
\end{equation}
Fix $k\in\mathbb{N}_0$. For any compact interval $[\alpha,\beta]$,
by Helly's convergence theorem we have
\begin{equation}\label{IntLim}
\lim_{j\rightarrow +\infty}\int_\alpha^\beta
x^k \,{\rm d}\phi_j(x)=\int_\alpha^\beta x^k \,{\rm d}\psi(x)\; .
\end{equation}
Therefore, choosing $-\infty<\alpha < 0 < \beta<+\infty$ and $j$ such that $n_j > k+1$, we deduce
\begin{equation}\label{GaussUk5}
\begin{array}{l}
\displaystyle \Big|u_k-
\int_\alpha^\beta x^k \,{\rm d}\psi(x)\Big| \; =\; \Big|\int_{-\infty}^{+\infty}
x^k \,{\rm d}\phi_j(x) - \int_\alpha^\beta x^k \,{\rm d}\psi(x)\Big| \\ [1.5em]
\qquad \displaystyle =\Big|\int_{-\infty}^\alpha x^k
\,{\rm d}\phi_j(x)+\int_\alpha^\beta x^k \,{\rm d}\phi_j(x) +
\int_\beta^{+\infty}x^k \,{\rm d}\phi_j(x) - \int_\alpha^\beta x^k
\,{\rm d}\psi(x)\Big|\\ [1.5em]
\qquad \displaystyle \leq\Big|\int_{-\infty}^\alpha x^k
\,{\rm d}\phi_j(x)\Big|+\Big|\int_\beta^{+\infty}x^k
\,{\rm d}\phi_j(x)\Big|+\Big|\int_\alpha^\beta x^k \,{\rm d}\phi_j(x) -
\int_\alpha^\beta x^k \,{\rm d}\psi(x)\Big|\; .
\end{array}
\end{equation}
But,
$$
\Big|\int_{-\infty}^{\alpha}x^k\,{\rm d}\phi_j(x)\Big|
=\Big|\int_{-\infty}^{\alpha}\frac{x^{2k+2}}
{x^{k+2}}\,{\rm d}\phi_j(x)\Big| \leq \frac{1}{|\alpha|^{k+2}}\int_{-\infty}^{+\infty}x^{2k+2}
\,{\rm d}\phi_j(x) =\frac{u_{2k+2}}{|\alpha|^{k+2}}\; ,
$$
where the last equality follows from (\ref{GaussUk4}), since $n_j>k+1$.
Similarly,
$$
\Big|\int_{\beta}^{+\infty} x^k \,{\rm d}\phi_j(x)\Big| \leq
\frac{u_{2k+2}}{\beta^{k+2}}\,.
$$
Therefore,
sending $j\rightarrow +\infty$ in (\ref{GaussUk5}) and taking into account (\ref{IntLim}), we find
$$\Big|u_k - \int_{\alpha}^{\beta}x^k \,{\rm d}\psi(x)\Big| \leq
u_{2k+2}\left(\frac{1}{|\alpha|^{k+2}} + \frac{1}{\beta^{k+2}}\right)\,,\quad k\in\mathbb{N}_0\,.$$
Thus, taking the limits $\alpha\to -\infty$ and $\beta \to +\infty$,
we conclude that (\ref{GaussUk3}) holds also whenever the true interval of orthogonality is unbounded.

It remains to prove that the spectrum of any natural representative $\psi$
(of ${\bf u}$) fulfilling (\ref{u-rep-Thm}) is an infinite set.
Indeed, if $\sigma(\psi)=\{x_1,\cdots,x_N\}$ (a finite subset of $\mathbb{R}$),
define
%then considering the polynomial
$$
p(x):=(x-x_1)(x-x_2)\cdots(x-x_N)\;,
$$
and let $h_j:=\psi(x_j+0)-\psi(x_j-0)$ be the jump of $\psi$ at the point $x_j$.
Then, from (\ref{u-rep-Thm}) we would have
$$
\langle{\bf u},p^2\rangle=\int_{-\infty}^{+\infty}p^2(x)\,{\rm d}\psi(x)=
\sum_{j=1}^Np^2(x_j)h_j=0\;,
$$
in contradiction with the positive definiteness of ${\bf u}$.
\qed

\begin{remark}\em
We have remarked before that a natural representative $\psi$ for
a positive-definite functional ${\bf u}\in\mathcal{P}'$ is a distribution function.
Moreover, being a nondecreasing function, the set of points of discontinuity of $\psi$ is finite or denumerable.
 %Thus, since---by Proposition \ref{propret-RSintegral}-(iv)---changing the values of $\psi$ at its points of discontinuity
Thus, since changing the values of $\psi$ at its points of discontinuity
does not change the value of the Riemann-Stieltjes integral with respect to $\psi$
for continuous integrand functions (and so in particular for polynomials),
it follows that {\it there is a representative of ${\bf u}$,
in the sense of $(\ref{u-rep-Thm})$,
which is a bounded nondecreasing right-continuous function
with infinite spectrum and finite moments of all orders}.
\end{remark}

\begin{remark}\em
If $\psi$ is a distribution function which represents
a positive-definite functional ${\bf u}\in\mathcal{P}'$ in the sense of (\ref{u-rep-Thm}),
then so is any function obtained by adding a constant to $\psi$.
Such distribution functions are called {\sl essentially equal}.
The discussion about the existence of different distribution functions
(not essentially equal) which represent a given functional will be made later.
\end{remark}

\section{The spectral theorem}

In this section we present an alternative statement of the representation theorem
(Theorem \ref{RepThmOP}), called the {\it spectral theorem for orthogonal polynomials}.
It is worth mentioning that both the representation theorem and
the spectral theorem are equivalent versions of Favard's Theorem in the positive-definite case.

The spectral theorem asserts that any positive-definite
functional ${\bf u}\in\mathcal{P}'$ admits an integral representation involving
a positive Borel measure $\mu$ on $\mathbb{R}$ (which needs not to be unique)
with infinite support and such that all its moments exist (i.e., they are finite).
Recall that the {\sl support} of $\mu$ is the set
\vspace*{-0.5em}
\begin{snugshade}\begin{equation}\label{supp-mu}
\text{supp}(\mu):=\big\{x\in\mathbb{R}\,:\, \mu\big((x-\epsilon,x+\epsilon)\big)>0\;,\;
\forall\epsilon>0\big\}\;,
\end{equation}\end{snugshade}
\noindent
while saying that all the moments of $\mu$ exist (are finite) means that
\vspace*{-0.5em}
\begin{snugshade}\begin{equation}\label{mom-finite}
\int_\mathbb{R}|x|^n\,{\rm d}\mu<\infty\;,\quad\forall n\in\mathbb{N}_0\;.
\end{equation}\end{snugshade}
\noindent
Given a finite positive Borel measure $\mu$ on $\mathbb{R}$,
the function $F_\mu:\mathbb{R}\to\mathbb{R}$ defined by
\vspace*{-0.5em}
\begin{snugshade}\begin{equation}\label{dist-f1}
F_\mu(x):=\mu\big((-\infty,x]\big)
\end{equation}\end{snugshade}
\noindent
is called the {\sl distribution function} of $\mu$.
This function $F_\mu$ is bounded, nondecreasing, right-continuous, nonnegative, and it fulfills
$$
\lim_{x\to-\infty}F_\mu(x)=0\;.
$$
Conversely, any function $F:\mathbb{R}\to\mathbb{R}$ satisfying these five properties
is a distribution function for a finite positive Borel measure $\mu$, so that $F\equiv F_\mu$, and
\vspace*{-0.5em}
\begin{snugshade}\begin{equation}\label{intRS=intL}
\int_\mathbb{R}f(x)\,{\rm d}F_\mu(x)=\int_\mathbb{R}f(x)\,{\rm d}\mu(x)
\end{equation}\end{snugshade}
\noindent
for each continuous and $\mu-$integrable function $f$,
where the integral on the left-hand side of (\ref{intRS=intL})  is
the Riemann-Stieltjes integral generated by $F$.
Because of this fact often we will use $\mu$ to denote both a measure and its corresponding distribution function.
%Moreover, for continuous functions $f$ (and so in particular for polynomials) the integral on the left
The integral on the left-hand side of (\ref{intRS=intL}) is indeed
the Lebesgue-Stieltjes integral generated by $F$, and this is simply
the Lebesgue integral with respect to the Lebesgue--Stieltjes measure $\mu_F$ generated by $F$.
%Some more detailed explanations on this construction are given in Appendix \ref{Appx3}.

\begin{remark}\em
Notice that the support of a measure $\mu$ and the spectrum of the corresponding distribution
function, $F_\mu$, coincide, i.e.,
\vspace*{-0.5em}
\begin{snugshade}\begin{equation}\label{suppmu=sigmaFmu}
\text{supp}(\mu)=\sigma(F_\mu):=\big\{x\in\mathbb{R}\;:\;F_\mu(x+\delta)-F_\mu(x-\delta)>0\;,\;
\forall\delta>0\big\}\;.
\end{equation}\end{snugshade}
\end{remark}

\begin{snugshade}
\begin{theorem}[spectral theorem for orthogonal polynomials]\label{SpecThmOP}
Let $\{P_n\}_{n\geq0}$ be a monic OPS characterized by the three-term recurrence relation
\begin{equation}\label{TTRRST}
P_{n+1}(x)=(x-\beta_n)P_n(x)-\gamma_nP_{n-1}(x)\;,\quad n=0,1,2,\cdots
\end{equation}
with initial conditions $P_{-1}(x)=0$ and $P_{0}(x)=1$. Suppose that
\begin{equation}\label{TTRRbgn+ST}
\beta_{n-1}\in\mathbb{R}\;,\quad\gamma_n>0\;,\quad n=1,2,3,\ldots\;.
\end{equation}
Then, there exists a positive Borel measure $\mu$ on $\mathbb{R}$, whose support is an infinite set,
and with finite moments of all orders, such that
\begin{equation}\label{upositive1}
\int_\mathbb{R} P_n(x)P_m(x)\,{\rm d}\mu(x)=\zeta_n\delta_{n,m}\;,\quad n,m=0,1,2,\ldots\;,
\end{equation}
where $\zeta_n:=\gamma_1\gamma_2\cdots\gamma_n$ for each $n\in\mathbb{N}_0$ (being $\zeta_0:=1$).
\end{theorem}
\end{snugshade}

{\it Proof.}
By Favard's Theorem, under the given hypothesis
$\{P_n\}_{n\geq0}$ is a monic OPS with respect
to a positive-definite functional ${\bf u}\in\mathcal{P}'$.
Therefore, by the representation Theorem \ref{RepThmOP}, there exists a
distribution function $\psi:\mathbb{R}\to\mathbb{R}$
(a natural representative of ${\bf u}$), whose spectrum
$\sigma(\psi)$ is an infinite subset of $\mathbb{R}$, fulfilling
\begin{equation}\label{SpThmProof1}
u_k:=\langle {\bf u},x^k\rangle=
\int_{-\infty}^{+\infty} x^k\,{\rm d}\psi(x)\;,\quad k=0,1,2,\ldots\;.
\end{equation}
We may assume that $\psi$ is right-continuous without changing its spectrum.
Let $\mu$ be the corresponding Stieltjes-Lebesgue measure (hence it is a positive Borel measure),
so that $\psi$ is the distribution function $F_\mu$ of the measure $\mu$. Then
$$
\psi(x)=F_\mu(x)=\mu\big((-\infty,x]\big)\,,\quad x\in\mathbb{R}
$$
and we deduce $\mbox{\rm supp}(\mu)=\sigma(\psi)$,
%$$
%\mbox{\rm supp}(\mu):=\big\{x\in\mathbb{R}\,:\,
%\mu\big((x-\epsilon,x+\epsilon)\big)>0\;,\;\forall\epsilon>0\big\}=\sigma(\psi)\;,
%$$
hence the support of $\mu$ is an infinite set.
Moreover, from the connection between the
Riemann-Lebesgue and the Stieltjes-Lebesgue integrals, we have
$$
\int_{-\infty}^{+\infty} p(x)\,{\rm d}\psi(x)
=\int_{\mathbb{R}} p(x)\,{\rm d}F_\mu(x)
=\int_{\mathbb{R}} p(x)\,{\rm d}\mu(x)\;,\quad \forall p\in\mathcal{P}\;,
$$
and so, in particular, by (\ref{SpThmProof1}), the moments of $\mu$ all exist and
$$
\int_\mathbb{R} P_n(x)P_m(x)\,{\rm d}\mu(x)=
\int_{-\infty}^{+\infty} P_n(x)P_m(x)\,{\rm d}\psi(x)=
\langle {\bf u},P_n P_m\rangle=\langle{\bf u},P_{n}^2\rangle\delta_{n,m}
$$
for all $n,m\in\mathbb{N}_0$. Now, taking into account (\ref{usefull-R}), we have
$$
\langle{\bf u},P_{n}^2\rangle=u_0\prod_{j=1}^n\gamma_j=u_0\zeta_n\;,\quad n=0,1,2\ldots\;.
$$
Thus, if $u_0=1$, we obtain (\ref{upositive1}); otherwise, we normalize $\mu$
passing to the measure $\widehat{\mu}=u_0^{-1}\mu$,
and so (\ref{upositive1}) holds with $\widehat{\mu}$ instead of $\mu$.
\qed

\begin{remark}\em
Often, we will refer to a measure $\mu$ under the conditions of the spectral theorem as a ``spectral measure'',
or an ``orthogonality measure'' for the given OPS $\{P_n\}_{n\geq0}$. This measure needs not to be unique whenever the
true interval of orthogonality is an unbounded set, as it was observed by Stieltjes.
In section \ref{specmeasuniq} we will analyze this question.
%conditions for the spectral measure to be unique.
\end{remark}

%\begin{remark}\em
%It is worth mentioning that both the representation theorem and
%the spectral theorem are different ways to present one same result.
%Indeed they are equivalent versions of Favard's Theorem in the positive-definite case.
%\end{remark}

%\section{Non-unicity of the spectral measure (Stieltjes example)}

%\section{Conditions for the unicity of the spectral measure}

\section{On the unicity of the spectral measure}\label{specmeasuniq}

As remarked before, the orthogonality measure for an OPS needs not to be unique.
We present an example due to Stieltjes.
Consider the weight function\vspace*{-0.5em}
\begin{snugshade}\begin{equation}\label{StieltjesEx}
w_c(x;\alpha):=\big(1+\alpha\,\sin(2\pi c\ln x)\big)e^{-c\ln^2x}\;,\quad x\in I:=(0,+\infty)\;,
%\exp\big(-c\ln^2x\big)\;,\quad x\in I:=(0,+\infty)\;,
\end{equation}\end{snugshade}
\noindent
where we fix $\alpha$ and $c$ so that $-1<\alpha<1$ and $c>0$. Clearly,
these choices of $c$ and $\alpha$ ensure that $w_c(\cdot;\alpha)$ becomes nonnegative and integrable on $I$.
This weight function defines a positive-definite functional ${\bf u}\in\mathcal{P}'$ given by
$$
\langle{\bf u},p\rangle:=\int_0^{+\infty}p(x)w_c(x;\alpha)\,{\rm d}x\;,\quad p\in\mathcal{P}\;.
$$
An associated distribution function representing ${\bf u}$ is
$$
\psi_c(x;\alpha):=\int_{-\infty}^xw_c(t;\alpha)\chi_{(0,+\infty)}(t)\,{\rm d}t\;,\quad x\in\mathbb{R}\;.
$$
Computing the moments of ${\bf u}$, we deduce (Exercise \ref{Ex-cp4-1b})
\begin{equation}\label{eqMnSti}
u_n:=\int_0^{+\infty}x^n w_c(x;\alpha)\,{\rm d}x= \sqrt{\mbox{$\frac{\pi}{c}$}}\,e^{(n+1)^2/4c}\,,\quad n=0,1,2,\ldots\;.
\end{equation}
Thus, we see that the moments are independent of the choice of $\alpha$,
hence (for fixed $c>0$) by varying $\alpha\in(-1,1)$ we obtain different orthogonality measures with the same moments.
Therefore there are infinitely many orthogonality measures for an OPS with respect to ${\bf u}$.
\smallskip

Next we prove that the orthogonality measure given by the spectral theorem is unique
if both sequences of the $\beta-$parameters and $\gamma-$parameters are bounded.
We begin by stating two preliminary results.

\begin{snugshade}
\begin{lemma}\label{Mxell}
Let $A=[a_{ij}]_{i,j=1}^N$ be a matrix of order $N$, and let $M$ be a positive constant
chosen so that
$$
|a_{i,j}|\leq M\, ,\quad \forall i,j\in\{1,\ldots,N\}\;.
$$
Suppose that each row and each column of $A$ have at most $\ell$ nonzero entries.
Then each eigenvalue $\lambda$ of $A$ satisfies
$$%\begin{equation}\label{eigenvAellM}
|\lambda|\leq\ell M\; .
$$%\end{equation}
\end{lemma}
\end{snugshade}

{\it Proof.}
Take $x=(x_1,\ldots,x_N)$ to be an eigenvector of $A$ corresponding to the eigenvalue $\lambda$, so that
$$
Ax=\lambda x\;,\quad \|x\|=1\;,
$$
where $\|x\|:=\sqrt{\langle x,x\rangle}$, being $\langle x,y\rangle:=\sum_{j=1}^Nx_j\overline{y_j}$
the usual inner product on $\mathbb{C}^N$. Then
$$
\begin{array}{rcl}
|\lambda|^2 &=& \big|\lambda\,\|x\|^2\big|^2=\big|\lambda\langle x,x\rangle\big|^2=\big|\langle\lambda x,x\rangle\big|^2
=\big|\langle Ax,x\rangle\big|^2 \\ [0.5em]
&\leq& \displaystyle \|Ax\|^2\,\|x\|^2=\|Ax\|^2
=\sum_{j=1}^N|(Ax)_j|^2=\sum_{j=1}^N\Big|\sum_{k=1}^N a_{jk}x_k\Big|^2 \\ [0.5em]
&\leq& \displaystyle\sum_{j=1}^N\Big(\sum_{k=1}^N |a_{jk}|^2 \sum_{k=1}^N|x_k|^2\Big)
=\sum_{j=1}^N\sum_{k=1}^N |a_{jk}|^2 \leq\ell^2 M^2\;,
\end{array}
$$
where in the first two inequalities we have applied the Cauchy-Schwartz inequality,
and the last one holds since, by hypothesis,
each row and each column of $A$ have at most $\ell$ nonzero entries and
all the entries of $A$ are bounded by $M$.
\qed

\begin{snugshade}
\begin{theorem}\label{sup-mu-bnd1}
Under the hypothesis of the spectral Theorem \ref{SpecThmOP}, assume further that both
$\{\beta_n\}_{n\geq0}$ and $\{\gamma_n\}_{n\geq1}$ are bounded sequences.
Then the support of the orthogonality measure $\mu$ is a bounded set.
\end{theorem}
\end{snugshade}

{\it Proof.}
By hypothesis,
\begin{equation}\label{bgnA}
\exists C>0\;:\;\forall n\in\mathbb{N}\;,\quad |\beta_n|\leq C\;,\quad |\gamma_n|\leq C\;.
\end{equation}
On the other hand, we know that the zeros $x_{n,j}$ of $P_n$ ($1\leq j\leq n$) are the eigenvalues
of the tridiagonal matrix $J_n$ given by (\ref{Jn}).
Since, by (\ref{bgnA}), the entries of $J_n$ are bounded by $M:=\max\{1,C\}$, and
in each row and each column of $J_n$ there are at most $\ell=3$ nonzero entries,
then Lemma \ref{Mxell} ensures that
\begin{equation}\label{bgnB}
\big|x_{n,j}\big|\leq 3M\;,\quad \forall n\in\mathbb{N}\;,\; \forall j\in\{1,2,\ldots,n\}\;.
\end{equation}
Thus, the spectrum of each distribution function $\psi_n:\mathbb{R}\to\mathbb{R}$
introduced in (\ref{dist-psiN}) is contained in the interval $[-3M,3M]$,
hence the spectrum $\sigma(\psi)$ of any distribution function $\psi$ obtained as
a limit of a subsequence of $\{\psi_n\}_{n\geq1}$ is also contained in $[-3M,3M]$.
Therefore, the orthogonality measure $\mu$ given by the spectral Theorem \ref{SpecThmOP}
satisfies
$$
\mbox{\rm supp}(\mu)=\sigma(\psi)\subseteq[-3M,3M]\;,
$$
so that $\mbox{\rm supp}(\mu)$ is a bounded set.
\qed

\begin{remark}\em
Under the conditions of the spectral Theorem \ref{SpecThmOP}, it follows from the proof of
Theorem \ref{sup-mu-bnd1} that {\it if $\{\beta_n\}_{n\geq0}$ and $\{\gamma_n\}_{n\geq1}$ are bounded sequences
then the true interval of orthogonality of the corresponding monic OPS is bounded}.
\end{remark}

\begin{snugshade}
\begin{theorem}\label{mu-uniq-bnd}
Under the hypothesis of the spectral Theorem \ref{SpecThmOP}, assume further that both
$\{\beta_n\}_{n\geq0}$ and $\{\gamma_n\}_{n\geq1}$ are bounded sequences.
Then the orthogonality measure $\mu$ is unique.
\end{theorem}
\end{snugshade}

{\it Proof.}
By Theorem \ref{sup-mu-bnd1}, there exists at least one orthogonality measure $\mu$ with compact support.
Let $\nu$ be any other orthogonality measure (hence it has the same moments as $\mu$).
Fix $a>0$. Then, for each $n\in\mathbb{N}_0$, we have
\begin{equation}\label{int-mu-nu1}
\int_{|x|\geq a}{\rm d}\nu(x)\leq\int_{|x|\geq a}\Big|\frac{x}{a}\Big|^{2n}\,{\rm d}\nu(x)\leq
a^{-2n}\int_{\mathbb{R}}x^{2n}\,{\rm d}\nu(x)=a^{-2n}\int_{\mathbb{R}}x^{2n}\,{\rm d}\mu(x)\;.
\end{equation}
We have seen in the proof of Theorem \ref{sup-mu-bnd1} that, for each $n\in\mathbb{N}$,
the zeros $x_{n,j}$ of $P_n$ are uniformly bounded, hence
\begin{equation}\label{int-mu-nu2}
\exists r>0\;:\;\forall n\in\mathbb{N}\;,\;\forall j\in\{1,2,\ldots,n+1\}\;,\quad
\big|x_{n+1,j}\big|\leq r\;.
\end{equation}
By the Gauss-Jacobi-Christofell quadrature formula (Theorem \ref{GaussJC}),
with $n$ replaced by $n+1$ and $p(x)=x^{2n}$, we have
$$
\int_{\mathbb{R}}x^{2n}\,{\rm d}\mu(x)=\sum_{j=1}^{n+1}A_{n+1,j}x_{n+1,j}^{2n}\,.
$$
Substituting this in the right-hand side of (\ref{int-mu-nu1}) and taking into
account (\ref{int-mu-nu2}), as well as (\ref{GJC1}), we deduce
$$
\int_{|x|\geq a}{\rm d}\nu(x)%\leq a^{-2n}\sum_{j=1}^{n+1}A_{n+1,j}x_{n+1,j}^{2n}
\leq\Big(\frac{r}{a}\Big)^{2n}\sum_{j=1}^{n+1}A_{n+1,j}=u_0\Big(\frac{r}{a}\Big)^{2n}\;.
$$
Choosing $a>r$ and taking the limit as $n\to+\infty$, we obtain
$\;\int_{|x|\geq a}{\rm d}\nu(x)=0\,$.
Thus $\nu\big(\{x\in\mathbb{R}:|x|\geq a\}\big)=0$, hence
$\mbox{\small\rm supp}\,(\nu)\subseteq[-a,a]$ for each $a>r$, and so
$$
\mbox{\small\rm supp}\,(\nu)\subseteq[-r,r]\;.
$$
This proves that any orthogonality measure has a compact support contained in the interval $[-r,r]$.
%We now prove that $\nu=\mu$ (meaning that $\mu$ is unique).
%For that we will prove that the distribution functions $F_\mu$ and $F_\nu$ coincide on $\mathbb{R}$.
To prove the uniqueness of the orthogonality measure,
take arbitrarily $z\in\mathbb{C}$ and $t\in\mathbb{R}$ such that $|z|\geq 2r$ and $|t|\leq r$.
Then $\left|\frac{t}{z}\right|\leq\frac12$, hence
$$
s_n(t):=\sum_{k=0}^n\frac{t^k}{z^{k+1}}=\frac{1}{z}\sum_{k=0}^n\Big(\frac{t}{z}\Big)^k=
\frac{1-\Big(\frac{t}{z}\Big)^{n+1}}{z-t}\rightarrow\frac{1}{z-t}\;,\quad\mbox{\rm as $n\to+\infty$}\;\,.
$$
Therefore, for any orthogonality measure $\mu$,
\begin{equation}\label{int-mu-nu4}
\int_{\mathbb{R}}\frac{{\rm d}\mu(t)}{z-t}=\int_{\mathbb{R}}\lim_{n\to\infty}s_n(t)\,{\rm d}\mu(t)
=\lim_{n\to\infty}\int_{\mathbb{R}}s_n(t)\,{\rm d}\mu(t)
=\lim_{n\to\infty}\sum_{k=0}^n\frac{u_k}{z^{k+1}}\,,
\end{equation}
where the interchange between the limit and the integral follows by Lebesgue's dominated convergence theorem,
taking into account that
$$
|s_n(t)|=\left|\frac{1-\Big(\frac{t}{z}\Big)^{n+1}}{z-t}\right|\leq\frac{2}{|z-t|}
=:g(t)\in L_\mu^1([-r,r])\;.
$$
Notice that the limit in the right-hand side of (\ref{int-mu-nu4}) depends only of the
moments $u_k$, $k\in\mathbb{N}_0$, hence it has the same value considering
any measure $\mu$ with the same moments and with compact support contained in $[-r,r]$
(i.e., considering any orthogonality measure).
Therefore, the function
$$
F(z):=\int_{\mathbb{R}}\frac{{\rm d}\mu(x)}{z-x}
\; ,\quad z\in\mathbb{C}\setminus[-r,r]
%\; ,\quad z\in\mathbb{C}\setminus\mbox{\rm supp}(\mu)
$$
is uniquely determined by $\mu$ for $z$ outside the circle $|z|=2r$
(meaning that, for $z$ outside this circle, $F(z)$ has the same value for any orthogonality measure $\mu$).
Since $F$ is analytic on $\mathbb{C}\setminus[-r,r]$, then by the identity theorem for analytic functions,\footnote{\,The identity theorem for analytic functions asserts that {\it given functions $f$ and $g$ analytic on a connected open set $D\subseteq\mathbb{C}$, if $f=g$ on some open and non-empty subset of $D$ then $f=g$ on $D$}.} we may conclude that $F(z)$ is uniquely determined by $\mu$ for $z\in\mathbb{C}\setminus[-r,r]$.
Thus the uniqueness of the measure follows from the Perron--Stieltjes inversion
formula:% (to be proved later):
%(which we will prove in the next text---cf. Theorem \ref{Perron-Stieltjes}):
$$
\psi_\mu(t)-\psi_\mu(s)=\lim_{\varepsilon\rightarrow0^+}\,
\frac{1}{\pi}\int_s^t \frac{F(x-i\varepsilon)-F(x+i\varepsilon)}{2 i}\,{\rm d}x\; ,
%\frac{1}{\pi}\int_s^t \Im \left(F(x-i\varepsilon)\right)\,{\rm d}x\; .
$$
where $\psi_\mu$ is an appropriate normalization of the distribution function $F_\mu$, given by
$$
\psi_\mu(-\infty):=0\;,\quad
\psi_\mu(t):=\frac{F_\mu(t+0)+F_\mu(t-0)}{2}\;,\quad t\in\mathbb{R}\;.
$$
Notice that $\psi_\mu$ and $F_\mu$ may be different only at (countably many) points of discontinuity,
hence the integrals of continuous functions with respect to $\psi_\mu$ and $F_\mu$ take the same value.
The Perron--Stieltjes inversion
formula will be proved later.
\qed
\bigskip

We conclude this section stating without proof two results that ensure
uniqueness of the orthogonality measure---see Theorems II-5.1 and II-5.2 in
Freud's book \ref{Freud1971-C3a}.% (pp.\,78--81).

\begin{snugshade}
\begin{theorem}[Riesz uniqueness criterium]\label{mu-uniq-RieszCriterium}
The orthogonality measure $\mu$ is unique whenever its sequence of moments $u_n:=\int_\mathbb{R}t^n\,{\rm d}\mu(t)$
$(n=0,1,2,\ldots)$ satisfies %the condition
\begin{equation}\label{mu-uniq1}
\liminf_{n\to+\infty}\frac{\sqrt[2n]{u_{2n}}}{2n}<\infty\;.
\end{equation}
\end{theorem}
\end{snugshade}

\begin{snugshade}
\begin{corollary}\label{mu-uniq-CorRieszCriterium}
 The orthogonality measure $\mu$ is unique if the condition
\begin{equation}\label{mu-uniq2}
\int_\mathbb{R}e^{\theta|x|}\,{\rm d}\mu(x)<\infty
\end{equation}
holds for some $\theta>0$.
\end{corollary}
\end{snugshade}

%\begin{snugshade}
%\begin{theorem}\label{mu-uniq-RieszCriterium}
%Let $\psi:\mathbb{R}\to\mathbb{R}$ a distribution function, and suppose that its sequence of moments $u_k:=\int_{-\infty}^{+\infty}t^n\,{\rm d}\psi(t)$ ($n=0,1,2,\ldots$) satisfies the condition
%$$
%\liminf_{n\to+\infty}\frac{\sqrt[2n]{u_{2n}}}{2n}<\infty\;.
%$$
%Then the orthogonality measure $\mu$ with distribution function
%$\psi$---i.e., with moments $u_n$ ($n=0,1,2,\ldots$)---is unique.
%\end{theorem}
%\end{snugshade}
%
%\begin{snugshade}
%\begin{corollary}\label{mu-uniq-CorRieszCriterium}
%Let $\psi:\mathbb{R}\to\mathbb{R}$ be a distribution function, and suppose that the condition
%$$
%\int_{-\infty}^{+\infty}e^{\theta|x|}\,{\rm d}\psi(x)<\infty
%$$
%holds. Then the orthogonality measure $\mu$ with distribution function $\psi$ is unique.
%\end{corollary}
%\end{snugshade}

\section*{Exercises}
%\bigskip

{\small
%\noindent
\begin{enumerate}[label=\emph{\bf \arabic*.},leftmargin=*]

\item\label{Ex-cp4-1a}
Show that the conclusion of Helly's convergence Theorem \ref{Helly2}
may not holds whenever $[\xi,\eta]$ is not a bounded interval.
\medskip

\noindent
({\sl Hint.} Consider $\{\phi_n\}_{n\geq0}$ defined by
$\phi_n(x):=0$ if $x<n$, and $\phi_n(x):=1$ if $x\geq n$.)
\medskip

\item\label{Ex-cp4-1b}
Prove (\ref{eqMnSti}).
\medskip

\item\label{Ex-cp4-Charlier}
(Charlier polynomials.)
Define the monic Charlier OPS $\{C_n^{(a)}(x)\}_{n\geq0}$ by the generating function
$$
e^{-aw}(1+w)^x=\sum_{n=0}^\infty C_n^{(a)}(x)\,\frac{w^n}{n!}\;,
$$
being $a\in\mathbb{R}\setminus\{0\}$.
Prove the following assertions:
\begin{enumerate}
\item For each $n\in\mathbb{N}_0$, $C_n^{(a)}(x)$ has the explicit representation
$$
C_n^{(a)}(x)=\sum_{k=0}^{n} {{x\choose k}} {{n\choose k}} k! (-a)^{n-k}\;,
$$
being ${{z\choose 0}}:=1$ and ${{z\choose k}}:=z(z-1)\cdots(z-k+1)/k!$ for $k=1,2,\ldots,n$ ($z\in\mathbb{C}$).
\item $\{C_n^{(a)}(x)\}_{n\geq0}$ is an OPS with respect to the functional ${\bf u}\in\mathcal{P}'$ given by
$$
\big\langle{\bf u},p\big\rangle=\int_{0}^{+\infty}p(x)\,{\rm d}\psi^{(a)}(x)\;,\quad p\in\mathcal{P}\;,
$$
where $\psi:\mathbb{R}\to\mathbb{R}$ is a step function whose jumps are given by
 %$\psi^{(a)}(x+0)-\psi^{(a)}(x-0):=e^{-a}a^x/x!$ at
$e^{-a}a^x/x!$ at the points $x=0,1,2,\ldots$.
The positive-definite case occurs for $a>0$, and in this case $\psi^{(a)}(x)$
is the Poisson distribution function of probability theory.
\item The TTRR for $\{C_n^{(a)}(x)\}_{n\geq0}$ is
$$
C_{n+1}^{(a)}(x)=(x-n-a)C_{n}^{(a)}(x)-anC_{n-1}^{(a)}(x)\;,\quad n=0,1,2,\cdots\;,
$$
being $C_{-1}^{(a)}(x):=0$ and $C_{0}^{(a)}(x)=0$.
\end{enumerate}
\medskip

\item\label{Ex-cp4-Meixner}
(Meixner polynomials.)
Let $\{m_n(x;\beta,c)\}_{n\geq0}$ be the Meixner OPS of the first kind, defined via the generating function
$$
\big(1-c^{-1}w\big)^x(1-w)^{-x-\beta}=\sum_{n=0}^\infty m_n(x;\beta,c)\,\frac{w^n}{n!}\;,
$$
being $c\in\mathbb{R}\setminus\{0,1\}$ and $\beta\in\mathbb{R}\setminus\{0,-1,-2,-3,\cdots\}$.
Prove the following assertions:
\begin{enumerate}
\item For each $n\in\mathbb{N}_0$, $m_n(x;\beta,c)$ has the explicit representation
$$
m_n(x;\beta,c)=(-1)^nn! \sum_{k=0}^{n} {{x\choose k}} {{-x-\beta\choose n-k}} c^{-k}\;.
%\quad n\in\mathbb{N}_0\;.
$$
\item If $0<c<1$ and $\beta>0$, $\{m_n(x;\beta,c)\}_{n\geq0}$ is an OPS with respect to the positive-definite functional ${\bf u}\in\mathcal{P}'$ given by
$$
\big\langle{\bf u},p\big\rangle=\int_{-\infty}^{+\infty}p(x)\,{\rm d}\psi(x)\;,\quad p\in\mathcal{P}\;,
$$
where the distribution function $\psi:\mathbb{R}\to\mathbb{R}$ is a step function
supported on $\mathbb{N}_0$ (i.e., $\sigma(\psi)=\mathbb{N}_0$) whose jump at
the point $x=k$ is given by $c^k(\beta)_k/k!$, $k\in\mathbb{N}_0$.
\end{enumerate}
\smallskip

({\sl Hint.} Proceed as for the Charlier polynomials, using the relation
$$
\sum_{n=0}^\infty\frac{(\beta)_nz^n}{n!}=(1-z)^{-\beta}\,,\quad |z|<1\;.\mbox{\rm )}
$$
\smallskip

\item\label{Ex-cp4-2a}
Let $\{P_n\}_{n\geq0}$ be a monic OPS with respect to a positive Borel measure $\mu$.
%with an infinite support and whose moments of all orders exist.
Prove that each $P_n(x)$ admits the representation
$$
P_n(x)=\frac{1}{n!H_{n-1}} \int_\mathbb{R}\cdots\int_\mathbb{R}\; %\int_{\mathbb{R}^n}
\prod_{i=1}^n(x-x_i)\prod_{1\leq j<k\leq n}^n(x_j-x_k)^2\,{\rm d}\mu(x_1)\cdots{\rm d}\mu(x_n)\;,
$$
where $H_{n-1}$ is the Hankel determinant of order $n$.
\smallskip

{\sl Remark.} This formula goes back at least as far as Heine (from Heine-Borel), 1878.
Nowadays it has important applications in \emph{Random Matrix Theory}.

\end{enumerate}
\medskip
}
%\medskip

\section*{Final remarks}

A concise proof of Lemma \ref{LammaCountableContinuous}
appears in the book \ref{Leoni2009-C3} by Giovanni Leoni.
Lemma \ref{LammaHelly} and Helly's Theorems (under the formulation presented here)
may be found in Chihara's book.
For the main facts concerning the Riemann-Stieltjes integral needed to understand this text we refer the student to the classical book \ref{KolmogorovFomini1968-A2} by Kolmogorov and Fomini.
The sections about the representation theorem (Theorem \ref{RepThmOP}) and the spectral theorem (Theorem \ref{SpecThmOP})
follow closely the presentations appearing in the books by Chihara and Ismail (respectively).
The section where we consider the problem of the uniqueness of the orthogonality measure is also based on Ismail's book.
Theorem \ref{mu-uniq-RieszCriterium} and Corollary \ref{mu-uniq-CorRieszCriterium}
appear in Freud's book \ref{Freud1971-C3a} (cf. Theorems II-5.1 and II-5.2 therein).
Exercises {\bf 1}, {\bf 3} and {\bf 4} were taken from Chihara's book.
Stieltjes example in exercise {\bf 2} may be found in Ismail's book.
Exercise {\bf 5} appear, e.g., in Szeg\"o's book \ref{Szegoo1975-C3a}, being the formulation presented here
based on Percy Deift's book \ref{Deift1997-C3a}.
As a final remark we mention that Deift's book presents another formulation of the spectral theorem, exploring the connections between the theory of OP and the theory of Jacobi operators.
\medskip

%\section*{Notas finais}
%\section*{Coment\'arios finais}
%
%Uma prova concisa e elegante do Lema \ref{LammaCountableContinuous}
%pode ser vista no livro \ref{Leoni2009-C3} de Giovanni Leoni.
%O Lema \ref{LammaHelly} e os Teoremas de Helly, na formula\c c\~ao aqui apresentada,
%encontram-se no livro de Chihara.
%Os factos sobre o integral de Riemann-Stieltjes necess\'arios \`a compreens\~ao deste texto
%podem ser vistos, e.g., no livro cl\'assico de Kolmogorov e Fomini \ref{KolmogorovFomini1968-A2}.
%As sec\c c\~oes sobre os teoremas de representa\c c\~ao (Teorema \ref{RepThmOP}) e espectral (Teorema \ref{SpecThmOP})
%seguem as apresenta\c c\~oes que aparecem nos livros de Chihara e Ismail (repectivamente).
%A sec\c c\~ao acerca da unicidade da medida espectral \'e tamb\'em baseada no livro de Ismail.
%O Teorema \ref{mu-uniq-RieszCriterium} e o Corol\'ario \ref{mu-uniq-CorRieszCriterium}
%aparecem no livro de Freud \ref{Freud1971-C3a} (trata-se dos Teoremas II-5.1 e II-5.2 a\'\i {\,} contidos).
%Os exerc\'\i cios {\bf 1}, {\bf 3} e {\bf 4} s\~ao retirados do livro de Chihara.
%O exemplo de Stieltjes a que se refere o exerc\'\i cio {\bf 2} pode encontrar-se no livro de Ismail.
%O exerc\'\i cio {\bf 5} aparece, e.g., no livro \ref{Szegoo1975-C3a} de Szeg\"o, e a formula\c c\~ao
%apresentada \'e baseada no livro \ref{Deift1997-C3a} de Percy Deift.
%Refira-se, para concluir, que o livro de Deift apresenta ainda outra formula\c c\~ao
%do teorema espectral, explorando a liga\c c\~ao entre a teoria dos PO e a teoria dos operadores de Jacobi.
%\smallskip

\section*{Bibliography}
\smallskip

{\small
\begin{enumerate}[label=\emph{\rm [\arabic*]},leftmargin=*]
\item\label{Chihara1978-C3a} T. S. Chihara, {\sl An introduction to orthogonal polynomials}, Gordon and Breach (1978).
\item\label{Deift1997-C3a} P. Deift, {\sl Orthogonal polynomials and random matrices: A Riemann-Hilbert Approach}, AMS Courant Lecture Notes {\bf 3} (2000).
\item\label{Freud1971-C3a} G. Freud, {\sl Orthogonal polynomials}, Pergamon Press, Oxford (1971).
\item\label{Gautschi2004-C3a} W. Gautschi, {\sl Orthogonal polynomials. Computation and approximation}, Oxford University Press, Oxford (2004).
\item\label{Ismail2004-C3a} M. E. H. Ismail, {\sl Classical and quantum orthogonal polynomials in one variable},
          Cambridge University Press (2005) [paperback edition: 2009].
\item\label{KolmogorovFomini1968-A2} A. Kolmogorov and S. Fomini, {\sl Introductory Real Analysis}, Dover Publications, Inc. (1975).
\item\label{Leoni2009-C3} G. Leoni, {\sl A first course in Sobolev spaces}, Graduate Studies in Math., AMS, {\bf 105}
(2009).
\item\label{Szegoo1975-C3a} G. Szeg\"o, {\sl Orthogonal polynomials}, AMS Colloq. Publ. {\bf 230} (1975), 4th ed.
\end{enumerate}
}

\chapter{Markov's Theorem}\label{MarkovChap}
%\chapter{Associated polynomials and Markov's Theorem}\label{MarkovChap}

\pagestyle{myheadings}\markright{Markov's Theorem}
\pagestyle{myheadings}\markleft{J. Petronilho}

%In this text we consider the question raised in section \ref{SPinvform} of text \ref{}.

According to the spectral Theorem \ref{SpecThmOP}, given a sequence of monic
polynomials $\{P_n\}_{n\geq0}$ fulfilling the three-term recurrence relation
\begin{snugshade}
\begin{equation}\label{TTRRSTa}
xP_n(x)=P_{n+1}(x)+\beta_nP_n(x)+\gamma_nP_{n-1}(x)\;,\quad n=0,1,2,\cdots
\end{equation}
\end{snugshade}\noindent
with initial conditions $P_{-1}(x)=0$ and $P_{0}(x)=1$, and subject to the conditions
\begin{snugshade}
\begin{equation}\label{TTRRbgn+STa}
\beta_{n-1}\in\mathbb{R}\;,\quad\gamma_n>0\;,\quad n=1,2,3,\ldots\;,
\end{equation}
\end{snugshade}\noindent
then there exists a positive Borel measure $\mu$ which is an orthogonality measure with respect to
which $\{P_n\}_{n\geq0}$ is an OPS. Moreover, when the true interval of orthogonality is a compact set
(or when both $\{\beta_n\}_{n\geq0}$ and $\{\gamma_n\}_{n\geq1}$ are bounded sequences)
this measure is unique. Thus, it is natural to ask:
\smallskip

\begin{snugshade}
{\sc Problem.} {\it How to find the orthogonality measure from the three-term
recurrence relation fulfilled by the polynomials?}
\end{snugshade}
\smallskip

This question fits into the study of the so-called {\it inverse problems}
in the Theory of Orthogonal Polynomials.
In this text we will describe a program that leads to the orthogonality measure
starting from the three-term recurrence relation.
The main tools for the success of this program are Markov's Theorem and
the Perron--Stieltjes inversion formula. The latter
allow us to find the measure from the knowledge of its Stieltjes transform,
which in turn is determined by the former.
%In order to state these results we need to introduce

\section{The Perron--Stieltjes inversion formula}\label{SPinvform}

\begin{snugshade}
\begin{definition}\label{Stieltjes-transf}
Let $\mu$ be a finite positive Borel measure, with $\mbox{\rm supp}(\mu)\subseteq\mathbb{R}$.
The {\sl Stieltjes transform} associated with $\mu$ is the complex function $F\equiv F(\cdot;\mu)$
given by
\begin{equation}\label{STmu}
F(z):=\int_{\mathbb{R}}\frac{{\rm d}\mu(x)}{z-x}
\; ,\quad z\in\mathbb{C}\setminus\mbox{\rm supp}(\mu)\; .
\end{equation}
\end{definition}
\end{snugshade}

Recall that $\mbox{\rm supp}(\mu)$ is a closed set. Using this fact we may prove that
$F$ is an analytic function on $\mathbb{C}\setminus\mbox{\rm supp}(\mu)$.
The Perron--Stieltjes inversion formula allow us to recover the measure
from the knowledge of its corresponding Stieltjes transform.

\begin{snugshade}
\begin{theorem}[Perron--Stieltjes inversion formula]\label{Perron-Stieltjes}
Let $\mu$ be a finite positive Borel measure. % with compact support and finite moments of all orders.
Then for every $a,b\in\mathbb{R}$, with $a<b$, the equality
%$$
%F(z):=\int_{\mathbb{R}}\frac{{\rm d}\mu(x)}{z-x}
%\; ,\quad z\in\mathbb{C}\setminus\mbox{\rm supp}(\mu)
%$$
%$F$ be its Stieltjes transform. Then
\begin{equation}\label{SIF}
\lim_{\varepsilon\rightarrow 0^+}\frac{1}{\pi}\int_a^b
\Im \left(F(x-i\varepsilon)\right)dx=\mu((a,b))+
\mbox{$\frac{1}{2}\,$}\mu(\{a\})+\mbox{$\frac{1}{2}\,$}\mu(\{b\})%\; ,
\end{equation}
holds, where $F$ is the Stieltjes transform associated with $\mu$.
\end{theorem}
\end{snugshade}

{\it Proof.}
Observe that $\overline{F(z)}=F(\overline{z})$, hence
$$
\Im (F(z))=\frac{F(z)-F(\overline{z})}{2i}=
-\frac{1}{2i}\int_{\mathbb{R}}\frac{z-\overline{z}}{{|x-z|}^2}\,{\rm d}\mu(x)=
-\int_{\mathbb{R}}\frac{\Im (z)}{{|x-z|}^2}\,{\rm d}\mu(x)\, ,
$$
%$$
%2i\Im (F(z))=F(z)-\overline{F(z)}=F(z)-F(\overline{z})=
%-\int_{\mathbb{R}}\frac{z-\overline{z}}{{|x-z|}^2}d\mu(x)=
%-2i\int_{\mathbb{R}}\frac{\Im (z)}{{|x-z|}^2}d\mu(x)\; ,
%$$
so that, setting $z=x-i\epsilon$, with $x\in\mathbb{R}$ and $\varepsilon>0$, we may write
$$\begin{array}{c} \Im (F(x-i\varepsilon))=
\displaystyle\int_{\mathbb{R}}\frac{\varepsilon}{{|s-x+i\varepsilon|}^2}\,{\rm d}\mu(s)
=\displaystyle\int_{\mathbb{R}}\frac{\varepsilon}{{(s-x)}^2+\varepsilon^2}\,{\rm d}\mu(s)\;.
\end{array}
$$
Integrating on $(a,b)$ with respect to $x$, and then interchanging the order of integration in the last integral
(this is allowed taking into account that the integrand function is positive), we obtain
\begin{equation}\label{4.2}
\int_a^b\Im (F(x-i\varepsilon))\,{\rm d}x
=\int_{\mathbb{R}}\theta_{\varepsilon}(s)\,{\rm d}\mu(s)\; ,
\end{equation}
where
$$
\theta_{\varepsilon}(s):=\int_a^b\frac{\varepsilon}{{(s-x)}^2+\varepsilon^2}\,{\rm d}x
=\arctan\left(\frac{b-s}{\varepsilon}\right)-\arctan\left(\frac{a-s}{\varepsilon}\right)\; .
$$
%Making the change of variables $x\curvearrowright y$ defined by $x=s+\varepsilon y$, we deduce
%$$
%\theta_{\varepsilon}(s)=\arctan\left(\frac{b-s}{\varepsilon}\right)
%-\arctan\left(\frac{a-s}{\varepsilon}\right)\; ,
%$$
Notice that, for each $s\in\mathbb{R}$,
$$
\lim_{\varepsilon\rightarrow 0}
\theta_{\varepsilon}(s)= \left\{
\begin{array}{rcl}
\pi & \mbox{\rm if} & a<s<b \\ [0.25em]
\frac{\pi}{2} & \mbox{\rm if} & s=a \;\;\mbox{\rm or}\;\; s=b \\ [0.25em]%\rule{0pt}{1.0em}
0 & \mbox{\rm if} & s<a  \;\;\mbox{\rm or}\;\; s>b
\end{array}
\right\}=\pi\chi_{(a,b)}(s)+\mbox{$\frac\pi2$}\,\chi_{\{a\}}(s)+\mbox{$\frac\pi2$}\,\chi_{\{b\}}(s)\,.
%=:\theta(s)
$$
Moreover, $|\theta_{\varepsilon}(s)| \leq \pi$ for each $s\in\mathbb{R}$,
and the (constant) function $s\in\mathbb{R}\mapsto\pi$ is integrable with respect to
$\mu$ (since $\mu$ is a finite measure, so that $\mu(\mathbb{R})<\infty$).
Thus, by Lebesgue's dominated convergence theorem,
$$
\lim_{\varepsilon\rightarrow
0}\int_{\mathbb{R}}\theta_{\varepsilon}(s)\,{\rm d}\mu(s)
=\int_{\mathbb{R}}\displaystyle
\lim_{\varepsilon\rightarrow 0}\theta_{\varepsilon}(s)\,{\rm d}\mu(s)
=\pi\mu((a,b))+\mbox{$\frac\pi2$}\,\mu(\{a\})+\mbox{$\frac\pi2$}\,\mu(\{b\})\; ,
$$
and (\ref{SIF}) follows from (\ref{4.2}) taking the limit as $\varepsilon\rightarrow0^+$.
\qed
\medskip

\begin{remark}\em
The Perron--Stieltjes inversion formula (\ref{SIF}) may be stated in terms of the distribution function $F_\mu$
associated to the measure $\mu$, after an appropriate normalization of $F_\mu$.
Indeed, being $\psi:\mathbb{R}\to\mathbb{R}$ defined by %taking an appropriate normalization of $F_\mu$, so that
\begin{snugshade}
$$
%\psi(-\infty):=0\; ,\quad
\psi(x):=\frac{F_\mu(x+0)+F_\mu(x-0)}{2}\; ,$$
\end{snugshade}\noindent
and recalling that $F_\mu(x):=\mu\big((-\infty,x]\big)$ for each $x\in\mathbb{R}$,
(\ref{SIF}) may be rewritten as
\begin{snugshade}
$$
\psi(b)-\psi(a)=\lim_{\varepsilon\rightarrow0^+}\,
\frac{1}{\pi}\int_a^b \Im \left(F(x-i\varepsilon)\right)\,{\rm d}x\; .
$$
\end{snugshade}
\end{remark}

In fact, we deduce
$$
\begin{array}{rl}
\psi(b)-\psi(a)& = \frac12\,\big[F_\mu(b+0)+F_\mu(b-0)-F_\mu(a+0)-F_\mu(a-0)\big] \\ [.7em]
& = \frac12\,\big[\big(F_\mu(b)-F_\mu(a)\big)+\big(F_\mu(b-0)-F_\mu(a-0)\big)\big] \\ [.7em]
& = \frac12\,\big[\mu((a,b])+\mu([a,b))\big] = \mu((a,b))+\frac12\mu(\{a\})+\frac12\mu(\{b\})\;,
\end{array}
$$
where in the second equality we took into account that $F_\mu$ is right-continuous.

\section{Associated polynomials}

Let ${\bf u}\in\mathcal{P}'$, regular (not necessarily positive-definite),
and let $\{P_n\}_{n\geq0}$ be the corresponding monic OPS.
%be the linear functional with respect to which $\{P_n\}_{n\geq0}$ is an OPS.
According with Favard's Theorem, $\{P_n\}_{n\geq0}$ is
characterized by the three-term recurrence relation
(\ref{TTRRSTa}), with $\beta_{n-1}\in\mathbb{C}$ and $\gamma_n\in\mathbb{C}\setminus\{0\}$
for each $n\in\mathbb{N}$. Making a shift on this recurrence relation, we may define a new monic OPS,
$\{P_n^{(k)}\}_{n\geq0}$, being $k\in\mathbb{N}_0$ (fixed), called the (monic)
{\sl associated polynomials of order $k$} corresponding to $\{P_n\}_{n\geq0}$, by
\begin{snugshade}
\begin{equation}\label{Ak}
P_{n+1}^{(k)}(x)=\big(x-\beta_{n+k}\big)P_{n}^{(k)}(x)-
\gamma_{n+k}P_{n-1}^{(k)}(x)\;,\quad
n=0,1,2,\cdots\;
\end{equation}
\end{snugshade}\noindent
with initial conditions $P_{-1}^{(k)}(x)=0$ and $P_0^{(k)}(x)=1$.
Favard's Theorem  ensures that, indeed,
$\{P_n^{(k)}\}_{n\geq0}$ is a monic OPS.

\begin{remark}\em
If $k=0$, then $P_n^{(0)}\equiv P_n$.
When $k=1$, often $\{P_n^{(1)}\}_{n\geq0}$ is called the sequence of (monic)
{\sl associated polynomials of the first kind}, or {\sl numerator polynomials}.
\end{remark}

%Moreover, %from the three-term recurrence relation one easily deduce that
According with (\ref{PndetJn}), for each $k\in\mathbb{N}_0$ and $n\in\mathbb{N}$,
$P_n^{(k)}$ has the following representation as a determinant of order $n$ of a tridiagonal matrix:
\begin{snugshade}
\begin{equation}\label{WassociadosMonic1}
P_n^{(k)}(x)=\left|
\begin{array}{cccccc}
x-\beta_{k} & 1  &  0& \dots & 0 & 0  \\
\gamma_{k+1} & x-\beta_{k+1}  &  1 &\dots & 0 & 0 \\
0 & \gamma_{k+2} & x-\beta_{k+2}&   \dots & 0 & 0 \\
\vdots & \vdots & \vdots   & \ddots & \vdots & \vdots \\
0 &0 & 0 & \cdots& x-\beta_{n+k-2} & 1 \\
0 &0 & 0 &  \cdots & \gamma_{n+k-1} & x-\beta_{n+k-1}
\end{array}\right|\;.
\end{equation}
\end{snugshade}\noindent
%Let ${\bf u}\in\mathcal{P}'$ be the linear functional with respect to which
%$\{P_n\}_{n\geq0}$ is an OPS. Then, another
When $k\geq1$, another representation for $P_n^{(k)}$ is
\begin{snugshade}
\begin{equation}\label{WMonic1a}
P_n^{(k)}(x)=\frac{1}{\langle{\bf u},P_{k-1}^2\rangle}\,
\Big\langle P_{k-1}(\xi){\bf u}(\xi),\frac{P_{n+k}(x)-P_{n+k}(\xi)}{x-\xi}\Big\rangle
\quad(k\in\mathbb{N}\;;\;n\in\mathbb{N}_0)\;,
\end{equation}
\end{snugshade}\noindent
where ${\bf u}(\xi)$ means that ${\bf u}$ acts on polynomials regarded as functions of the variable $\xi$.
The representation (\ref{WMonic1a}) may be easily proved
by checking that the right-hand side defines a polynomial on the variable $x$
which fulfills the three-term recurrence relation (\ref{Ak}), and
for $n=-1$ and $n=0$ the right-hand side of (\ref{WMonic1a}) equals $0$ and $1$, respectively.
(Exercise \ref{Ex-cp5-A3})
A more concise form of writing (\ref{WMonic1a}) is
\begin{snugshade}
\begin{equation}\label{WMonic1b}
P_n^{(k)}(x)=\big\langle(\xi-x)^{-1}{\bf a}_{k-1}(\xi),P_{n+k}(\xi)\big\rangle
\quad(k\in\mathbb{N}\;;\;n\in\mathbb{N}_0)\;,
\end{equation}
\end{snugshade}\noindent
where $\{{\bf a}_n\}_{n\geq0}$ is the dual basis associated
with $\{P_n\}_{n\geq0}$
% ---the notation ${\bf a}_{k-1}(\xi)$ meaning that the functional ${\bf a}_{k-1}$
%acts on functions of the variable $\xi$---,
and the division of a functional
by a polynomial is given by Definition \ref{def-left-mult}.
Finally, considering the operator $\theta_0$ introduced in (\ref{thetac}), defined
for each $p\in\mathcal{P}$ by
$\theta_0p(x):=(p(x)-p(0))/x$
%$$
%\theta_0p(x):=\frac{p(x)-p(0)}{x}\quad (p\in\mathcal{P})
%$$
if $x\neq0$, and $\theta_0p(0):=p'(0)$, and taking into account
the definition of right multiplication of a functional by a polynomial
(Definition \ref{def-right-mult}), %as well as Theorem \ref{dualR},
%then (\ref{WMonic1a}) may be written in a more concise form as
then from (\ref{WMonic1b}) we arrive at a rather elegant representation
for the $n$th degree monic associated polynomial of order $k$:
%Henceforth,
\begin{snugshade}
\begin{equation}\label{WMonic1bb}
P_n^{(k)}={\bf a}_{k-1}\theta_0P_{n+k}
%P_n^{(k)}(x)=\big({\bf a}_{k-1}\theta_0P_{n+k}\big)(x)
\quad(k\in\mathbb{N}\;;\;n\in\mathbb{N}_0)\;.
\end{equation}
\end{snugshade}\noindent
%$(k\in\mathbb{N}\;;\;n\in\mathbb{N}_0)$,
%where $\{{\bf a}_n\}_{n\geq0}$
%is the dual basis associated with $\{P_n\}_{n\geq0}$.

A very useful relation linking the associated polynomials of orders $k$ and $k+1$ is
\begin{snugshade}
\begin{equation}\label{WassociadosMonic}
P_{n}^{(k+1)}(x)P_{n}^{(k)}(x)-P_{n+1}^{(k)}(x)P_{n-1}^{(k+1)}(x)=
\prod_{j=1}^{n}\gamma_{j+k}\;,\quad k,n\in\mathbb{N}_0\;.
\end{equation}
\end{snugshade}\noindent
Indeed, by (\ref{Ak}), we have, for all $k,n\in\mathbb{N}_0$,
$$
\begin{array}{rcl}
P_{n}^{(k+1)}(x)&=&\big(x-\beta_{n+k}\big)P_{n-1}^{(k+1)}(x)-\gamma_{n+k}P_{n-2}^{(k+1)}(x)\;,\quad \\ [0.5em]
P_{n+1}^{(k)}(x)&=&\big(x-\beta_{n+k}\big)P_{n}^{(k)}(x)-\gamma_{n+k}P_{n-1}^{(k)}(x)\;.
\end{array}
$$
Multiplying the first equality by $P_{n}^{(k)}(x)$ and the second one by $-P_{n-1}^{(k+1)}(x)$,
and then adding the resulting equalities, we deduce
$$
P_{n}^{(k+1)}(x)P_{n}^{(k)}(x)-P_{n+1}^{(k)}(x)P_{n-1}^{(k+1)}(x)=\gamma_{n+k}
\big(P_{n-1}^{(k+1)}(x)P_{n-1}^{(k)}(x)-P_{n}^{(k)}(x)P_{n-2}^{(k+1)}(x)\big)\,,
$$
hence (\ref{WassociadosMonic}) follows after repeatedly application of this relation.

We also state the following formulas, close to the Christoffel--Darboux identities:
\begin{snugshade}
\begin{equation}\label{PnxPnyj}
\frac{P_n(x)-P_n(y)}{x-y}=\sum_{k=1}^{n}P_{k-1}(x)
P_{n-k}^{(k)}(y)\;,\quad n\in\mathbb{N}\;.
\end{equation}\vspace*{-0.5em}
\begin{equation}\label{PnxPny2}
P_n'(x)=\sum_{k=1}^{n}P_{k-1}(x)
P_{n-k}^{(k)}(x)\;,\quad n\in\mathbb{N}\;.
\end{equation}
\end{snugshade}\noindent
Clearly, (\ref{PnxPny2}) follows from (\ref{PnxPnyj}) by taking the limit $y\to x$.
To prove (\ref{PnxPnyj}), notice that $(P_n(x)-P_n(y))/(x-y)$ is a polynomial of
degree $n-1$ in the variable $x$ (whose coefficients depend on $y$), so we can write
\begin{equation}\label{PnxPnyj1}
\frac{P_n(x)-P_n(y)}{x-y}=\sum_{k=0}^{n-1}c_{n,k}(y)P_{k}(x)
\end{equation}
for each $n\geq1$. Then for $0\leq k\leq n-1$, we compute the Fourier coefficients $c_{n,k}(y)$:
$$
c_{n,k}(y)=\frac{\Big\langle{\bf u}(x),\frac{P_{n}(x)-P_{n}(y)}{x-y}P_{k}(x)\Big\rangle}{\langle{\bf u},P_{k}^2\rangle}=
P_{n-1-k}^{(k+1)}(y)\;,
$$
where the last equality holds due to (\ref{WMonic1a}).
Substituting the last expression for $c_{n,k}(y)$ into (\ref{PnxPnyj1}) we obtain (\ref{PnxPnyj}).

%Often, specially in the framework of the positive-definite case, it may be useful
%to consider, instead of the monic polynomials, the corresponding orthonormal ones.

\section{Markov's Theorem}

We now return to the positive-definite case, with the purpose to state the celebrated Markov's theorem.
We begin by proving some preliminary results.

%We now return to the positive-definite case, so that we start with a sequence of polynomials
%$\{P_n\}_{n\geq0}$ which is a monic OPS with respect to a positive-definite
%linear functional ${\bf u}\in\mathcal{P}'$, being characterized by the three-term recurrence relation
%(\ref{TTRRSTa}) subject to conditions (\ref{TTRRbgn+STa}).

\begin{snugshade}
\begin{lemma}\label{b.2}
Let ${\bf u}\in\mathcal{P}'$ be positive-definite and $\{P_n\}_{n\geq0}$
the monic OPS with respect to ${\bf u}$. Then
\begin{equation}\label{LemMarkov1}
\frac{u_0P_{n-1}^{(1)}(x)}{P_n(x)}=\sum_{j=1}^n\frac{A_{n,j}}{x-{x_{n,j}}}=
\int_{-\infty}^{+\infty}\frac{{\rm d}\psi_n(t)}{x-t}\; ,\quad
x\in\mathbb{C}\setminus\Lambda_n\quad (n\in\mathbb{N})\; ,
\end{equation}
where $A_{n,1},\ldots,A_{n,n}$ are the coefficients appearing in the quadrature formula $(\ref{GJC1})$,
$\Lambda_n:=\{x_{n,1},\ldots,x_{n,n} \}$, being $x_{n,1}<\cdots<x_{n,n}$ the zeros of $P_n$,
and $\psi_n$ is the distribution function introduced in $(\ref{dist-psiN})$.
\end{lemma}
\end{snugshade}\noindent

{\it Proof.}
Notice first that, by (\ref{WassociadosMonic}) with $k=0$, for each $n\in\mathbb{N}$ the equality
$$
P_{n}^{(1)}(x)P_{n}(x)-P_{n+1}(x)P_{n-1}^{(1)}(x)=\gamma_1\gamma_2\cdots\gamma_{n}
$$
holds. Hence, since $\gamma_1\gamma_2\cdots\gamma_{n}\neq0$,
we see that $P_n$ and $P_{n-1}^{(1)}$ have no common zeros.
Moreover, we know that the zeros of $P_n$ are real and simple. Thus,
for each $n\in\mathbb{N}$, the decomposition of the rational function $P_{n-1}^{(1)}(x)/P_n(x)$
into partial fractions yields
\begin{equation} \label{3.2}
\frac{P_{n-1}^{(1)}(x)}{P_n(x)}=\sum_{j=1}^n
\frac{\lambda_{n,j}}{x-x_{n,j}}\; ,\quad
\lambda_{n,j}:=\frac{P_{n-1}^{(1)}(x_{n,j})}{P_n^\prime(x_{n,j})}\; .
\end{equation}
On the other hand, from (\ref{GJC-Anj}) in the proof of Theorem \ref{GaussJC}, we know that
$$
A_{n,j}=\langle {\bf u},\ell_{j,n}\rangle
=\Big\langle {\bf u},\frac{P_n(x)}{(x-x_{n,j})P_n^\prime(x_{n,j})}\Big\rangle
=\frac{u_0}{P_n^\prime(x_{n,j})}\frac{1}{u_0}
\Big\langle {\bf u},\frac{P_n(x_{n,j})-P_n(x)}{x_{n,j}-x}\Big\rangle\; .
$$
Moreover, by (\ref{WMonic1a}) with $k=1$, and changing $n$ into $n-1$, we see that the relation
$$
P_{n-1}^{(1)}(x)=\frac{1}{u_0}\,\Big\langle {\bf u}_\xi,\frac{P_{n}(x)-P_{n}(\xi)}{x-\xi}\Big\rangle
$$
holds for each $n\in\mathbb{N}$.
Therefore, we conclude that $A_{n,j}=u_0P_{n-1}^{(1)}(x_{n,j})/P_n^\prime(x_{n,j})$, hence
$\lambda_{n,j}=A_{n,j}/u_0$.
Thus the first equality in (\ref{LemMarkov1}) is proved.
The second one is an immediate consequence of the properties
of the Riemann-Stieltjes integral, taking into account that $\psi_n$ is a step function
with spectrum $\sigma(\psi_n)=\Lambda_n$, and with jump equal to $A_{n,j}$ at the point $x_{n,j}$.
\qed

%\begin{remark}
%Os n\'{u}meros $\lambda_{n,j}$ $(1 \leq j \leq n )$ s\~{a}o
%conhecidos por {\sf n\'{u}meros de Christoffel}.
%Decorre da f\'ormula de quadratura de Gauss e de (\ref{Christof})
%que os n\'umeros de Christoffel s\~ao todos positivos. Al\'em disso,
%atendendo a (\ref{ass}) e (\ref{1.2}), podemos observar que
%estes n\'umeros admitem a representa\c c\~ao
%$$
%\lambda_{n,j}=\frac{-1}{a_{n+1}p_n^\prime(x_{n,j})p_{n+1}(x_{n,j})}\; ,\quad
%1\leq j\leq n\, ,\quad n=1,2,3,\cdots\; .
%$$
%\end{remark}

\begin{snugshade}
\begin{lemma}\label{5.1}
Let $\{P_n\}_{n\geq0}$ be a monic OPS with respect to a
positive-definite functional ${\bf u}\in\mathcal{P}'$, and let
$\psi$ be a natural representative for ${\bf u}$,
in the sense of the representation Theorem \ref{RepThmOP}.
Let $s$ be a spectral point of $\psi$. %, i.e., $s\in\sigma(\psi)$.
Then, every neighborhood of $s$ contains
at least one zero of $P_n(x)$ for infinitely many values of $n\in\mathbb{N}$.
In symbols:
\begin{equation}\label{MemSpectral-s}
%\forall s\in\sigma(\psi)\;,\;\;
\forall V\in\mathcal{V}(s)\;,\;\;
\forall N\in\mathbb{N}\;,\;\;
\exists n\in\mathbb{N}\;:\;\;
n>N\;\wedge\; V\cap\Lambda_n\neq\emptyset\;.
\end{equation}
%\begin{equation}\label{MemSpectral-s}
%%\forall s\in\sigma(\psi)\;,\;\;
%\forall V\in\mathcal{V}(s)\;,\;\;
%\forall N\in\mathbb{N}\;,\;\;
%\bigcup_{n>N}V\cap\Lambda_n\neq\emptyset\;.
%\end{equation}
\end{lemma}
\end{snugshade}\noindent

{\it Proof.}
By hypothesis, $s\in\sigma(\psi)$.
If (\ref{MemSpectral-s}) is not true, there exist a neighborhood $V$ of $s$ and $N\in\mathbb{N}$
such that $V$ does not contains zeros of $P_n$ for every $n\geq N$.
%Then, by definition of $\psi_n$, we have $\,\psi_n(x)-\psi_n(z)=0\,$ for all $x,z\in V$ and
%$n\geq N$, hence, passing to the subsequence $\{\psi_{n_j}\}_{j\geq1}$
%which generates the natural representative $\psi$, so that
%$$
%\psi(x):=\lim_{j\to\infty}\psi_{n_j}(x)\;,\quad x\in\mathbb{R}\,,
%$$
%we deduce $\psi(x)-\psi(z)=0\,$ for any $x,z\in V$.
Then, by definition of $\psi_n$, we have $\,\psi_n(x)=\psi_n(z)\,$ for all $x,z\in V$ and
$n\geq N$. Since $\psi$ is a natural representative of ${\bf u}$,
there is a subsequence $\{\psi_{n_j}\}_{j\geq1}$ of $\{\psi_{n}\}_{n\geq1}$ such that
$$
\psi(x):=\lim_{j\to\infty}\psi_{n_j}(x)\;,\quad x\in\mathbb{R}\,.
$$
Therefore, $\psi(x)=\lim_{j\to\infty}\psi_{n_j}(x)=\lim_{j\to\infty}\psi_{n_j}(z)=\psi(z)$ for every $x,z\in V$.
Thus $\psi$ is constant on the neighborhood $V$ of $s$, hence $s\not\in\sigma(\psi)$, contrary to the hypothesis.
\qed
\medskip

It follows from Lemma \ref{5.1} that if $s$ is a spectral point
of a natural representative $\psi$ of ${\bf u}$, then either $s$ is a zero of
$P_n(x)$ for infinitely many $n$, or else $s$
is a limit of a sequence of numbers belonging to the set
\begin{snugshade}
$$Z_1:=\{x_{n,j} \,:\,1 \leq j \leq n \;,\;   n\in\mathbb{N}\}\; .$$
\end{snugshade}\noindent
Therefore, setting
\begin{snugshade}
$$
\begin{array}{c}
X_1:=Z^\prime_1\equiv \{\,\mbox{\rm accumulation points of $Z_1$}\, \}\; , \\ [0.25em]
%\rule{0pt}{1.5em}
X_2:= \{x\in Z_1 \,:\,P_n(x)=0\;\mbox{\rm for infinitely many $n$}\, \}\; ,
\end{array}
$$
\end{snugshade}\noindent
then $\sigma(\psi)\subseteq X_1 \cup X_2\,$.
Moreover, the following holds:
\begin{snugshade}
\begin{equation}\label{SupCoSup}
\sigma(\psi)\subseteq X_1 \cup X_2 \subseteq[\xi,\eta]
=\mbox{\rm co}(\sigma(\psi))\; ,
\end{equation}
\end{snugshade}\noindent
where $\mbox{\rm co}(\sigma(\psi))$  is the convex hull of the set
$\sigma(\psi)$, i.e., it is the smallest closed interval which contains $\sigma(\psi)$,
and, as usual, $[\xi,\eta]$ is the true interval of orthogonality of $\{P_n\}_{n\geq0}$.
Indeed, the second inclusion in (\ref{SupCoSup}) is an immediate consequence
of the definitions of the involved sets, and the last equality holds since {\it the
interval $\mbox{\rm co}(\sigma(\psi))$ is a supporting set for ${\bf u}$}.
To prove this, set $\mbox{\rm co}(\sigma(\psi))=[a,b]$,
and let $p(x)$ be a real polynomial which does not vanish identically on $[a,b]$
and it is non-negative there. Then, we have
$$\langle{\bf u},p\rangle=\int_{-\infty}^{+\infty}p(x)\,{\rm d}\psi(x)
=\int_a^b p(x)\,{\rm d}\psi(x)\,.$$
Since $\sigma(\psi)$ is an infinite set,
there is $x_0\in\sigma(\psi)$ such that $p(x_0)\neq0$, and so,
since $p(x)\geq0$ for each $x\in[a,b]$, by continuity we have $p(x)>0$ for each $x$ on a
neighborhood $(x_0-\delta,x_0+\delta)$ of the point $x_0$
---choosing $\delta>0$ so that $(x_0-\delta,x_0+\delta)\subseteq[a,b]$---, hence,
using the Mean Value Theorem for the Riemann-Stieltjes integral, we may write
$$\langle{\bf u},p\rangle=\int_a^b p(x)\,{\rm d}\psi(x)
\geq\int_{x_0-\delta}^{x_0+\delta} p(x)\,{\rm d}\psi(x)
= p(\xi_0)\big(\psi(x_0+\delta)-\psi(x_0-\delta)\big)\;,$$
for some $\xi_0\in(x_0-\delta,x_0+\delta)$. Henceforth,
since $p(\xi_0)>0$ and $\psi(x_0+\delta)-\psi(x_0-\delta)>0$
(since $x_0\in\sigma(\psi)$), we conclude that $\langle{\bf u},p\rangle>0$, so that
$[a,b]=\mbox{\rm co}(\sigma(\psi))$ is indeed an interval which is a supporting set for ${\bf u}$.
Therefore, by Theorem \ref{OPS-zeros-Thm2}, the (closed) interval
$\mbox{\rm co}(\sigma(\psi))$ contains the zeros of each $P_n$, $n\in\mathbb{N}$,
hence $[\xi,\eta]\subseteq\mbox{\rm co}(\sigma(\psi))$. Now, since
$\sigma(\psi)\subseteq[\xi,\eta]\subseteq\mbox{\rm co}(\sigma(\psi))$,
then $[\xi,\eta]=\mbox{\rm co}(\sigma(\psi))$.
Thus (\ref{SupCoSup}) is proved.
%$$\sigma(\psi_n)\subseteq[\xi,\eta]
%$$\sigma(\psi_n)\subseteq[\xi,\eta]\subseteq\mbox{\rm co}\big(\cup_{j=1}^\infty\Xi_j\big)
%=\mbox{\rm co}(Z_1)\;,\quad n\in\mathbb{N}\;.$$
%and the last one holds since the
%interval $\mbox{\rm co}(\sigma(\psi))$ is a supporting set for ${\bf u}$,
%hence it contains the zeros of each $P_n$, $n\in\mathbb{N}$.
%Indeed, more generally, the following holds:
%\begin{snugshade}
%\begin{equation}\label{SupCoSup}
%\mbox{\rm supp}(\psi)\subset X_1 \cup X_2 \subset[\xi_1,\eta_1]
%\subset \mbox{\rm co}(\mbox{\rm supp}(\psi))\; ,
%\end{equation}
%\end{snugshade}\noindent
%where $\mbox{\rm co}(\mbox{\rm supp}(\psi))$  is the convex hull of the set
%$\mbox{\rm supp}(\psi)$, i.e., it is the smallest closed interval which contains
%$\mbox{\rm supp}(\psi)$.
%The second inclusion in (\ref{SupCoSup}) is an immediate consequence
%of the definitions of the involved sets, and the last one follows since the
%interval $\mbox{\rm co}(\mbox{\rm supp}(\psi))$ is a supporting set for ${\bf u}$,
%hence it contains the zeros of each $P_n$, $n\in\mathbb{N}$.
%and the last one follows by Corollary \ref{CorZn}.

We next show that %observe an fact concerning the sets $X_1$ and $X_2$ is that
{\it $X_1\cup X_2$ is a closed set in $\mathbb{C}$}.
To prove this fact, we will show that the limit of any convergent sequence
of elements in $X_1 \cup X_2 $ belongs to this set. Indeed,
take arbitrarily a sequence $\{x_n\}_{n\geq1}$ such that $x_n\in X_1\cup X_2$ for each $n\in\mathbb{N}$
and $x_n\to x$ in $\mathbb{C}$. We need to prove that $x\in X_1\cup X_2 $.
Since $x_{n}\in X_1\cup X_2$ for each $n\in\mathbb{N}$, then two situations may occur:
$x_{n}\in X_1$ for infinitely many $n$, or $x_{n}\in X_2$ for infinitely many $n$ (or both).
In the first situation (passing, if necessary, to a subsequence),
since $x_n\to x$ and $X_1:=Z_1'$ is a closed set, we have $x\in X_1$;
in the second situation, we have (passing again, if necessary, to a subsequence,
and taking into account that $X_2\subset Z_1$)
$x=\lim_{n\to\infty}x_{n}\in X_2'\subseteq Z_1'=X_1$.
Therefore, in any situation, $x\in X_1\subseteq X_1\cup X_2$,
which proves that $X_1\cup X_2$ is closed.

We also need the following result from Complex Analysis, stated here without proof
(see e.g. Reinhold Remmert's book \ref{Remmert1997-C5a}, pp.\;150--151).
%see also Titchmarsh's book \ref{Titchmarsh1964-C4a}).
%(cf.\,\cite[pg.\,121]{chiara})

%\begin{snugshade}
%\begin{lemma}[Stieltjes-Vitali theorem]\label{Stieltjes-Vitali}
%Let $\{f_n\}_{n\geq1}$ be a sequence of analytic functions
%on a nonempty connected set contained in an open set $\Omega\subseteq\mathbb{C}$.
%Suppose that $\{f_n\}_{n\geq1}$ is uniformly bounded on each compact subset of $\Omega$ and it
%converges on a subset $E\subseteq \Omega$ which contains an accumulation point in $\Omega$.
%Then $\{f_n\}_{n\geq1}$ converges uniformly in any compact subset of $\Omega$.
%\end{lemma}
%\end{snugshade}

%\begin{snugshade}
%\begin{lemma}[Stieltjes-Vitali theorem]\label{Stieltjes-Vitali}
%Let $\{f_n\}_{n\geq1}$ be a sequence of functions, each analytic in a (open) region $D\subseteq\mathbb{C}$.
%Suppose that $\{f_n\}_{n\geq1}$ is uniformly bounded on $D$ and it
%converges on a subset $E\subseteq D$ having an accumulation point inside $D$.
%Then $\{f_n\}_{n\geq1}$ converges uniformly in any region bounded by a contour interior to $D$
%(the limit being, therefore, an analytic function).
%\end{lemma}
%\end{snugshade}

\begin{snugshade}
\begin{lemma}[Vitali's convergence theorem]\label{Stieltjes-Vitali}
Let $G$ be a domain in $\mathbb{C}$ (i.e., $G$ is a nonempty open connected subset of $\mathbb{C}$),
%(de acordo com a def de dominio em Reinhold Remmert, footnote da p.2),
%\footnote{\,As usual, a domain in $\mathbb{C}$ is an open and connected subset of $\mathbb{C}$.},
and let $\{f_n\}_{n\geq1}$ be a sequence of analytic functions in $G$ that is locally bounded in $G$
(equivalently, it is bounded on every compact set in $G$).
Suppose that the set
\begin{equation}\label{setA-Vitali}
A:=\big\{z\in G\,:\,\lim_{n\to\infty}f_n(z)\;\mbox{exists in $\mathbb{C}$}\big\}
\end{equation}
has at least one accumulation point in $G$.
Then the sequence $\{f_n\}_{n\geq1}$ converges uniformly on compact subsets of $G$.
\end{lemma}
\end{snugshade}

%\begin{snugshade}
%\begin{theorem}[Stieltjes-Vitali]
%Let $\{f_n\}_{n\geq1}$ be a sequence of analytic functions
%on a nonempty connected set contained in an open set $\Omega\subseteq\mathbb{C}$.
%If $\{f_n\}_{n\geq1}$ is uniformly bounded on each compact subset of $\Omega$ and it
%converges on a subset $E\subseteq \Omega$ which contains an accumulation point in $\Omega$,
%then $\{f_n\}_{n\geq1}$ converges uniformly in any compact subset of $\Omega$.
%\end{theorem}
%\end{snugshade}

Finally we are ready to state Markov's theorem.

\begin{snugshade}
\begin{theorem}[Markov]\label{Markov}
Let $\{P_n\}_{n\geq0}$ be a monic OPS with respect to a
positive-definite ${\bf u}\in\mathcal{P}'$, and let
$\psi$ be a natural representative of ${\bf u}$,
in the sense of the representation Theorem \ref{RepThmOP}.
Assume further that $\sigma(\psi)$ is a bounded set.
Then
\begin{equation}\label{3.6}
\lim_{n \rightarrow +\infty} \frac{u_0P_{n-1}^{(1)}(z)}{P_n(z)}=
\int_{-\infty}^{+\infty}\frac{{\rm d}\psi(x)}{z-x}\; , \quad
z\in\mathbb{C}\backslash(X_1 \cup X_2) \; ,
\end{equation}
the convergence being uniform on compact subsets of $\mathbb{C}\backslash(X_1 \cup X_2)$.
\end{theorem}
\end{snugshade}

{\it Proof.}
%A demonstra\c c\~ao que se segue \'e baseada nas refer\^encias
%\cite{Colombia}, \cite{chiara} e \cite{van}.
Since $\sigma(\psi)$ is bounded, then it follows immediately from (\ref{SupCoSup})
that the true interval of orthogonality $[\xi,\eta]$ of the sequence $\{P_n\}_{n\geq0}$ is bounded.
According with Lemma \ref{b.2}, we may write
\begin{equation}\label{LemMom}
\frac{u_0P_{n-1}^{(1)}(z)}{P_n(z)}=
\int_{\xi}^{\eta}\frac{{\rm d}\psi_n(x)}{z-x}\; ,\quad
z\in\mathbb{C}\backslash[\xi,\eta]
\end{equation}
for each $n\in\mathbb{N}$.
On the other hand, the representation Theorem \ref{RepThmOP} ensures the existence
of a subsequence $\{\psi_{n_j}\}_{j\geq0}$ which converges
on $[\xi,\eta]$ to the given natural representative $\psi$.
It follows from Helly's convergence Theorem \ref{Helly2} that
\begin{equation}\label{3.6ASub}
\lim_{j\rightarrow+\infty}\frac{u_0P_{n_j-1}^{(1)}(z)}{P_{n_j}(z)}=
\int_{\xi}^{\eta}\frac{{\rm d}\psi(x)}{z-x}\; ,\quad
z\in\mathbb{C} \backslash [\xi, \eta]\; .
\end{equation}
Set $M:=\max\{ |\xi|,|\eta|\}$. We will prove that
\begin{equation}\label{3.6A}
\lim_{n\rightarrow+\infty}\frac{u_0P_{n-1}^{(1)}(z)}{P_n(z)}=
\int_{-\infty}^{+\infty}\frac{{\rm d}\psi(x)}{z-x}\; ,\quad |z|>M\; ,
\end{equation}
the convergence being uniform on each set $\{z\in\mathbb{C}:|z|\geq M'\}$ such that $M'>M$.
We start by noticing that, by (\ref{WassociadosMonic}), with $k=0$,
\begin{equation}\label{Laur}
\frac{u_0P_{n}^{(1)}(z)}{P_{n+1}(z)}-\frac{u_0P_{n-1}^{(1)}(z)}{P_n(z)}=
\frac{C_n}{P_{n+1}(z)P_n(z)}\; ,\quad |z|>M
\end{equation}
%\begin{equation}\label{WassociadosMonic}
%P_{n}^{(k+1)}(x)P_{n}^{(k)}(x)-P_{n+1}^{(k)}(x)P_{n-1}^{(k+1)}(x)=
%\prod_{j=1}^{n}\gamma_{j+k}\;,\quad k,n\in\mathbb{N}_0\;.
%\end{equation}
for each $n=0,1,2,\cdots$, where $C_n:=u_0\prod_{j=1}^{n}\gamma_{j}$.
Since $P_{n+1}(z)P_n(z)$ is a polynomial of degree $2n+1$
with real and simple zeros, %it follows from the previous equality,
then by developing the right-hand side of (\ref{Laur}) in a Laurent series on the (open) annulus $|z|>M$
(taking into account that $1/(z-x)=\sum_{j\geq0}x^j/z^{j+1}$ for $|z|>|x|$),
we see that the Laurent series development of the left-hand side of (\ref{Laur}) takes de form
\begin{equation}\label{Laurent1}
\frac{u_0P_{n}^{(1)}(z)}{P_{n+1}(z)}-\frac{u_0P_{n-1}^{(1)}(z)}{P_n(z)}=
\frac{c_{2n+1}}{z^{2n+1}}+\frac{c_{2n+2}}{z^{2n+2}}+\cdots\;
,\quad |z|>M
\end{equation}
for each $n\in\mathbb{N}_0$. By repeatedly application of (\ref{Laurent1}) we deduce
\begin{equation}\label{Laurent1n1}
\frac{u_0P_{m-1}^{(1)}(z)}{P_{m}(z)}-\frac{u_0P_{n-1}^{(1)}(z)}{P_n(z)}=
\sum_{j=2n}^\infty\frac{c_j^{(m,n)}}{z^{j+1}}\;,\quad |z|>M\;, \quad m\geq n \quad (m,n\in\mathbb{N}_0)\,,
\end{equation}
where the $c_j^{(m,n)}$'s are complex numbers (indeed, we will see that they are real numbers).
Next, for each $n\in\mathbb{N}_0$ consider the Laurent series development
\begin{equation}\label{Laurent2}
\frac{u_0P_{n-1}^{(1)}(z)}{P_{n}(z)}
=\sum_{j=0}^\infty\frac{c_j^{(n)}}{z^{j+1}}\;,\quad |z|>M\;.
\end{equation}
Then, on the first hand, comparing (\ref{Laurent1n1}) and (\ref{Laurent2}), we deduce
%\begin{equation}\label{Laurent3}
%c_j^{(m)}=c_j^{(n)}\; ,\quad j=0,1,\ldots,2n\; ,\quad m\geq n\quad (m,n\in\mathbb{N})\; .
%\end{equation}
\begin{equation}\label{Laurent3a}
c_j^{(m,n)}=c_j^{(m)}-c_j^{(n)}\; ,\quad j\geq 2n\; ,\quad m\geq n\quad (m,n\in\mathbb{N}_0)\; .
\end{equation}
On the other hand, since $1/(z-x)=\sum_{j\geq0}x^j/z^{j+1}$ for
$|z|>|x|$, then (\ref{LemMom}) yields
$$
\frac{u_0P_{n-1}^{(1)}(z)}{P_{n}(z)}=\sum_{j=0}^\infty\frac{1}{z^{j+1}}
\int_{\xi}^{\eta}x^j\,{\rm d}\psi_n(x)\; ,\quad |z|>M\quad (n\in\mathbb{N}_0)\; ,
$$
hence, comparing with (\ref{Laurent2}),
and taking into account the uniqueness of the coefficients of a Laurent development, we obtain
$$
c_j^{(n)} %=\int_{-\infty}^{+\infty} x^j\,{\rm d}\psi_n(x)
=\int_{\xi}^{\eta} x^j\,{\rm d}\psi_n(x)\quad
(j,n\in\mathbb{N}_0)\; .
$$
Therefore, for $m\geq n$ and $j\geq 2n$, we deduce
$$
\big|c_j^{(m,n)}\big|=\big|c_j^{(m)}-c_j^{(n)}\big|
=\Big|\int_{\xi}^{\eta} x^j\,{\rm d}\psi_m(x)-\int_{\xi}^{\eta} x^j\,{\rm d}\psi_n(x)\Big|\leq 2u_0M^j\;,
$$
where the last inequality holds since $[\xi,\eta]\subseteq[-M,M]$ and
$\int_{\xi}^{\eta} \,{\rm d}\psi_n(x)=u_0$.
Thus, from (\ref{Laurent1n1}) we obtain
\begin{equation}\label{Laurent1n1aa}
\left|\frac{u_0P_{m-1}^{(1)}(z)}{P_{m}(z)}-\frac{u_0P_{n-1}^{(1)}(z)}{P_n(z)}\right|\leq
2u_0\sum_{j=2n}^\infty\left(\frac{M}{|z|}\right)^{j+1}\;,\quad |z|>M\;, \quad m\geq n\,.
\end{equation}
Since the series $\sum_{j=0}^\infty\big(\frac{M}{|z|}\big)^{j+1}$ is convergent whenever $|z|>M$,
it follows from (\ref{Laurent1n1aa}) that $\{u_0P_{n-1}^{(1)}(z)/P_n(z)\}_{n\geq0}$
is a Cauchy sequence for each $z$ fulfilling $|z|>M$.
Thus, since by (\ref{3.6ASub}) this sequence has a convergent subsequence,
it follows that the sequence converges (to the same limit as its subsequence).
Hence (\ref{3.6A}) is proved.
Note that the convergence in (\ref{3.6A}) is uniform on each set
$\{z\in\mathbb{C}:|z|\geq M'\}$ with $M'>M$.
In fact, from (\ref{Laurent1n1aa}) we obtain
\begin{equation}\label{LaurentAAA}
\left|\frac{u_0P_{m-1}^{(1)}(z)}{P_{m}(z)}-\frac{u_0P_{n-1}^{(1)}(z)}{P_n(z)}\right|\leq
2u_0\sum_{j=2n}^\infty\left(\frac{M}{M'}\right)^{j+1}\;,\quad |z|\geq M'>M\;, \quad m\geq n\,.
\end{equation}
Therefore, since, clearly,
\begin{equation}\label{LaurentBBB}
\forall\epsilon>0\;,\;\exists n_0\in\mathbb{N}\;:\;\forall n\in\mathbb{N}\;,\; n\geq n_0
\quad\Rightarrow\quad 2u_0\sum_{j=2n}^\infty\left(\frac{M}{M'}\right)^{j+1}<\epsilon\;,
\end{equation}
%Therefore, since for any given $\epsilon>0$ one can choose $n_0\in\mathbb{N}$ such that %$$2u_0\sum_{j=2n}^\infty\left(\frac{M}{M'}\right)^{j+1}<\epsilon$$ for every $m\geq n\geq n_0$,
then $\{u_0P_{n-1}^{(1)}/P_n\}_{n\geq0}$
is a (uniformly) Cauchy sequence on the set $\{z\in\mathbb{C}:|z|\geq M'\}$,
hence it converges uniformly therein,
%ver Elon Lages Lima, p.292, Teorema 1
and so we conclude that, indeed, the convergence in (\ref{3.6A}) holds uniformly on this set.
[This fact can be proved directly as follows:
Fix $\epsilon>0$. By (\ref{LaurentBBB}) and (\ref{LaurentAAA}),
there exists $n_0\in\mathbb{N}$ such that
$$
\forall m,n\in\mathbb{N}\;,\;\forall z\in A_{M'}\;,\quad m\geq n\geq n_0\;\Rightarrow\;
\left|\frac{u_0P_{m-1}^{(1)}(z)}{P_{m}(z)}-\frac{u_0P_{n-1}^{(1)}(z)}{P_n(z)}\right|<\epsilon\;,
$$
where $A_{M'}:=\{z\in\mathbb{C}:|z|\geq M'\}$.
Keeping  $z\in A_{M'}$ and $n$ fixed, and letting $m\to\infty$,
and taking into account that we already proved (\ref{3.6A}) pointwise, it follows that
$$
\forall m\in\mathbb{N}\;,\;\forall z\in A_{M'}\;,\quad m\geq n_0\;\Rightarrow\;
\left| \int_{-\infty}^{+\infty}\frac{{\rm d}\psi(x)}{x-z}-\frac{u_0P_{n-1}^{(1)}(z)}{P_n(z)}\right|\leq\epsilon\;,
$$
hence the convergence in (\ref{3.6A}) holds uniformly on the set $A_{M'}$, for each $M'>M$.]

To complete the proof, we need to show that
the convergence in (\ref{3.6A}) is indeed uniform on compact subsets of
$\mathbb{C}\backslash(X_1 \cup X_2)$. For each $N\in\mathbb{N}$, define
$$Z_N:=\{\,x_{n,j}\,:\,1 \leq j \leq n \; ,\; n\geq N\, \}\;.$$
Let $K$ be a compact subset of $\mathbb{C}\backslash(X_1 \cup X_2)$.
Then $K$ contains at most finitely many zeros of the polynomials in the sequence
$\{P_n\}_{n\geq0}$ (otherwise, $K$ could contain
an infinite subset of points in $Z_1$, hence---since it is compact---it would
contain a point in $Z_1'=X_1$, and so $K\cap X_1\neq\emptyset$,
which contradicts the definition of $K$),
and none of these zeros belong to $X_2$.
Therefore, there exists $N\in\mathbb{N}$ such that
$K\cap(X_1 \cup Z_N)=\emptyset$.
%Notice that $X_1\cup Z_N$ is a closed set
%(this fact can be proved in the same way
%as we did above to prove that $X_1\cup X_2$ is closed).
%Let us prove that $X_1 \cup Z_N $ is closed.
%We will show that the limit of any convergent sequence
%of elements in $X_1 \cup Z_N $ belongs to this set. Indeed,
%take arbitrarily a sequence $\{x_n\}_{n\geq1}$ such that $x_n\in X_1\cup Z_N$ for each $n\in\mathbb{N}$
%and $x_n\to x$ in $\mathbb{C}$. We need to prove that $x\in X_1\cup Z_N $.
%Since $x_{n}\in X_1\cup Z_N$ for each $n\in\mathbb{N}$, then two situations may occur:
%$x_{n}\in X_1$ for infinitely many $n$, or $x_{n}\in Z_N$ for infinitely many $n$ (or both).
%In the first situation (passing, if necessary, to a subsequence),
%since $x_n\to x$ and $X_1:=Z_1'$ is a closed set, we have $x\in X_1$;
%in the second situation, we have (passing again, if necessary, to a subsequence)
%$x=\lim_{n\to\infty}x_{n}\in Z_N'\subseteq Z_1'=X_1$.
%Therefore, in any situation, $x\in X_1\subseteq X_1\cup Z_N$,
%which proves that $X_1\cup Z_N$ is closed.
Next, set
\begin{equation}\label{dist-X1ZN}
\delta:=\mbox{\rm dist}\,\big(K,X_1\cup Z_N\big)=\inf \big\{|z-x| \,:\,  z\in K\,,\, x\in X_1\cup Z_N\big\}\;.
\end{equation}
Since $K$ is compact and $X_1 \cup Z_N $ is closed
(this fact can be proved in the same way
as we did above to prove that $X_1\cup X_2$ is closed),
and $K\cap(X_1 \cup Z_N)=\emptyset$, then $\delta>0$.
Therefore, using (\ref{LemMarkov1}), we deduce, for each $n\geq N$,
$$
\left|\frac{u_0P_{n-1}^{(1)}(z)}{P_n(z)} \right| \leq \displaystyle
\sum_{j=1}^n\frac{A_{n,j}}{|z-x_{n,j}|}\leq
\frac{1}{\delta}\displaystyle\sum_{j=1}^n
A_{n,j}=\frac{u_0}{\delta}\; ,\quad z\in K\; .
$$
(The last inequality follows from (\ref{dist-X1ZN}), taking into account that $x_{n,j}\in Z_N$ if $n\geq N$.)
Therefore, setting $f_n(z):=P_{n+N-1}^{(1)}(z)/P_{n+N}(z)$, we see that the sequence of functions
$\{f_n\}_{n\geq0}$ is bounded in $K$.
Thus, the claimed result follows from (\ref{3.6A}) and Vitali's convergence theorem
(Lemma \ref{Stieltjes-Vitali}), taking therein
$G:=\mathbb{C}\backslash (X_1\cup X_2)$ and noting that the set $A$ defined by (\ref{setA-Vitali}) contains
$\{z\in\mathbb{C}:|z|>M\}$.
%Note also that $G$ is open in $\mathbb{C}$, i.e., $X_1\cup X_2$ is closed.
%In fact, we have already seen that $X_1$ is closed, hence,
%to conclude that $\Omega$ is open we only need to show that $X_2$ is also closed.
%Take $x\in \overline{X_2}$. Then ...
%This is the same as showing that that
%This fact can be proved in the same way
%as we did above to prove that $X_1\cup Z_N$ is closed, noticing that
%$X_2\subset Z_1$, so that if $\{x_n\}_{n\geq1}$ is a sequence in $X_2$ which converges to some element $x$,
%then $x=\lim_{n\to\infty}x_{n}\in X_2'\subseteq Z_1'=X_1$.
\qed

\begin{remark}\em
The statement and proof given here for Markov's theorem is based on references
\ref{CharrisSS1991-C4a}, \ref{CharrisSS2003-C4a}, \ref{Chihara1978-C4a}, and \ref{VanAssche1991-C4a}.
Markov proved Theorem \ref{Markov} for absolutely continuous measures $\mu$
supported on a bounded interval: $\mbox{\rm supp}(\mu)=[a,b]$.
Under such conditions, by (\ref{SupCoSup}), $X_1\cup X_2=[a,b]$.
The result remains true for unbounded intervals,
provided the underlying moment problem is determined (see \ref{CBerg1994-C5}).
Different proofs of Markov's theorem, based on the notion of weak convergence of measures,
appear in \ref{NikishinSorokin1991i} and \ref{CBerg1994-C5}.
\end{remark}

%\begin{remark}\em
%It follows immediately from Theorem \ref{Markov} that the Stieltjes transform
%$F(z)$ is an analytic function on $\mathbb{C}\setminus(X_1\cup X_2)$.
%Moreover, since $\sigma(\psi)$ is a closed set on $\mathbb{R}$ (and so on $\mathbb{C}$),
%it is easy to see that $F$ is continuous on $\mathbb{R}\setminus\sigma(\psi)$.
%Therefore, since $\mathbb{C}\setminus(X_1\cup X_2)$ is an open set in $\mathbb{C}$,
%by analytic continuation we deduce that $F$ in analytic on $\mathbb{C}\setminus\sigma(\psi)$,
%whenever this is a connected set.
%\end{remark}

\begin{remark}\em
The set $X_1\cup X_2$ in Theorem \ref{Markov} cannot be replaced by $\sigma(\psi)$.
For instance, consider the sequence of monic polynomials $\{P_n\}_{n\geq0}$ defined by
$$
P_{2n+1}(x):=2^nxU_n\Big(\mbox{$\frac{x^2-5}{4}$}\big)\,,\quad
P_{2n}(x):=2^n\Big\{U_n\Big(\mbox{$\frac{x^2-5}{4}$}\Big)+2U_{n-1}\Big(\mbox{$\frac{x^2-5}{4}$}\Big)\Big\}\;.
$$
It can be shown (Exercise \ref{Ex-cp5-A6}) that $\{P_n\}_{n\geq0}$ is a monic OPS
with respect to the measure
$$
{\rm d}\mu(x):=\frac{\chi_E(x)}{|x|}\,\sqrt{1-\Big(\mbox{$\frac{x^2-5}{4}$}\Big)^2}\,{\rm d}x\,,
$$
where $E:=[-3,-1]\cup[1,3]$. Clearly,
%$$0\in X_2\subset X_1\cup X_2\;,\quad 0\not\in E= \sigma(\psi)\,,$$
$0\in X_2\subseteq X_1\cup X_2\,$ and $0\not\in E= \sigma(\psi)\,$.
Nevertheless, the ratio
%but the ratio
$P_{n-1}^{(1)}(z)/P_n(z)$ is not well defined at $z=0$ if $n$ is odd, hence
the sequence %$\lim_{n\to\infty}\frac{P_{n-1}^{(1)}(z)}{P_{n}(z)}$ doesn't exist at $z=0$.
$\big\{P_{n-1}^{(1)}(z)/P_{n}(z)\big\}_{n\geq0}$ has no limit as $n\to+\infty$ at $z=0$.
%$$
%\lim_{n\to\infty}\Big|\frac{P_{2n}^{(1)}(0)}{P_{2n+1}(0)}\Big|=+\infty\;.
%$$
\end{remark}

\section*{Exercises}
%\bigskip

{\small
%\noindent
\begin{enumerate}[label=\emph{\bf \arabic*.},leftmargin=*]

%\item\label{Ex-c4-1}
%Let ${\bf u}\in\mathcal{P}'$, regular, and $\{P_n\}_{n\geq0}$ the corresponding monic OPS.
%Let ${\bf u}^{(1)}$ be the (regular) functional corresponding to the monic associated OPS $\{P_n^{(1)}\}_{n\geq0}$.
%Let $S_{\bf u}$ and $S_{{\bf u}^{(1)}}$ be the formal Stieltjes series for ${\bf u}$ and ${\bf u}^{(1)}$, respectively.
%Show that the relation between these formal Stieltjes series is
%$$
%S_{{\bf u}^{(1)}}(z)=\beta_0-z-\frac{u_0}{S_{\bf u}(z)}\;.
%$$
%Find the relation between $S_{{\bf u}^{(k)}}$ and $S_{\bf u}(z)$, being $S_{{\bf u}^{(k)}}$ the formal Stieltjes series for the functional ${\bf u}^{(k)}$ corresponding to the monic associated polynomials of order $k$, $\{P_n^{(k)}\}_{n\geq0}$.
%\medskip

\item\label{Ex-c5-A1}
Prove that the Stieltjes transform $F$ introduced in Definition \ref{Stieltjes-transf}
is an analytic function on $\mathbb{C}\setminus\mbox{\rm supp}(\mu)$.
% ver esta afirmação em Nikishin & Sorokin, p. 60, na linha a seguir à fórmula (6.1)
\medskip

\item\label{Ex-c5-A2aa}
%Show by two different processes (first: use the definition; second: use Markov's theorem)
Let ${\rm d}\mu(x):=\frac{\chi_{(-1,1)}(x)}{\pi\sqrt{1-x^2}}\,{\rm d}x$
(so that $\mu$ is the orthogonality measure
for the Chebyshev polynomials of the first kind, $\{T_n\}_{n\geq0}$).
Show that the Stieltjes transform of $\mu$ is
$$
F(z)=\frac{1}{\sqrt{z^2-1}}\;,\quad z\in\mathbb{C}\setminus[-1,1]\;,
$$
where the branch of the complex square root is chosen so that $\sqrt{z^2-1}$ is an analytic function on $\mathbb{C}\setminus[-1,1]$ and $\sqrt{z^2-1}>0$ if $z>1$.
\medskip

%\item\label{Ex-c5-A2aa}
%%Show by two different processes (first: use the definition; second: use Markov's theorem)
%Show that the Stieltjes transform of the orthogonality measure of the Chebyshev polynomials of the first kind, $\{T_n\}_{n\geq0}$, is given by
%$$
%F(z)=\frac{1}{\sqrt{z^2-1}}\;,\quad z\in\mathbb{C}\setminus[-1,1]\;,
%$$
%where the branch of the complex square root is chosen so that $\sqrt{z^2-1}$ is an analytic function on $\mathbb{C}\setminus[-1,1]$ and $\sqrt{z^2-1}>0$ if $z>1$.
%\medskip

\item\label{Ex-c5-A2}
%Show by two different processes (first: use the definition; second: use Markov's theorem)
Show that the Stieltjes transform of the orthogonality measure
${\rm d}\mu(x):=\frac{\chi_{(-1,1)}(x)}{\pi}\sqrt{1-x^2}\,{\rm d}x$
(for the Chebyshev polynomials of the second kind, $\{U_n\}_{n\geq0}$) is
$$
F(z)=2\big(z-\sqrt{z^2-1}\,\big)\;,\quad z\in\mathbb{C}\setminus[-1,1]\;,
$$
where the branch of the complex square root is chosen as in Exercise \ref{Ex-c5-A2aa}
\medskip

%\item\label{Ex-c5-A2}
%%Show by two different processes (first: use the definition; second: use Markov's theorem)
%Show that the Stieltjes transform of the orthogonality measure of the Chebyshev polynomials of the second kind, $\{U_n\}_{n\geq0}$, is given by
%$$
%F(z)=2\big(z-\sqrt{z^2-1}\,\big)\;,\quad z\in\mathbb{C}\setminus[-1,1]\;,
%%F(z)=\frac{1}{\sqrt{z^2-1}}\;,\quad z\in\mathbb{C}\setminus[-1,1]\;,
%$$
%where the branch of the complex square root is chosen so that $\sqrt{z^2-1}$ is an analytic function on $\mathbb{C}\setminus[-1,1]$ and $\sqrt{z^2-1}>0$ if $z>1$.
%\medskip

\item\label{Ex-cp5-A3}
Prove relations (\ref{WMonic1a}) and (\ref{WMonic1bb}).
\medskip

\item\label{Ex-cp5-A4} % Gautschi2004,p.16
Let ${\rm d}\lambda(x):=\chi_{(-1,1)}{\rm d}x$ be the orthogonality measure for the Legendre polynomials (so that it is the Lebesgue measure on $[-1,1]$). Show that the associated Legendre polynomials of the first kind are orthogonal
with respect to the measure
$$
{\rm d}\lambda^{(1)}(x):=
\frac{2\chi_{(-1,1)}(x)\;}{\;\pi^2+\displaystyle\ln^2\frac{1+x}{1-x}\;}\,{\rm d}x\; .
%,\; x\in\mathbb{R}\; .
$$
%\medskip

\noindent
({\sl Hint.} Denote by $F$ and $F^{(1)}$ the Stieltjes transforms of the orthogonality measures
for the Legendre polynomials and their associated polynomials of the first kind, respectively.
We may start by showing that
$$F(z)=\mbox{\rm Log}\left(\frac{z+1}{z-1}\right)\;,\quad z\in\mathbb{C}\setminus[-1,1]\;.$$
Here we took the principal branch of the logarithm,
so that $F$ is an analytic function on $\mathbb{C}\setminus[-1,1]$.
%Hence, setting $z=x+iy$, with $x,y\in\mathbb{R}$, we deduce
%$$
%\qquad F(x+iy)=\left\{
%\begin{array}{l}
%\displaystyle\ln\left|\frac{1+x}{1-x}\right| \,,\;\mbox{\rm if $\;y=0\;\wedge\; x\not\in[-1,1]$} \\ [1em]
%\displaystyle\ln\sqrt{\frac{(1+x)^2+y^2}{(1-x)^2+y^2}}
%-i\,\left[\arctan\left(\frac{1-x}{y}\right)+\arctan\left(\frac{1+x}{y}\right)\right]
%\,,\,\mbox{\rm if $y\neq0$\,.}
%\end{array}
%\right.
%$$
Hence, setting $z=x-i\epsilon$, with $x\in\mathbb{R}$ and $\epsilon>0$, we deduce
$$\qquad
F(x-i\epsilon)=\ln\sqrt{\frac{(1+x)^2+\epsilon^2}{(1-x)^2+\epsilon^2}}
+2i\,\arctan\left(\sqrt{\left(\frac{x^2-1+\epsilon^2}{2\epsilon}\right)^2+1}
-\frac{x^2-1+\epsilon^2}{2\epsilon}\right)\,.
$$
%$$
%F(x-i\epsilon)=\ln\sqrt{\frac{(1+x)^2+\epsilon^2}{(1-x)^2+\epsilon^2}}
%-i\,\left[\arctan\left(\frac{1-x}{\epsilon}\right)+\arctan\left(\frac{1+x}{\epsilon}\right)\right]\,.
%$$
Next, using the relation
$\,F^{(1)}(z)=z-\beta_0-u_0/F(z)\,$, $z\in\mathbb{C}\setminus[-1,1]$
(as usual, $u_0$ is the moment of order zero for the measure ${\rm d}\lambda$,
and $\beta_0$ is the first $\beta-$parameter appearing in the TTRR for the
monic Legendre polynomials---so that, indeed, we compute $u_0=2$ and $\beta_0=0$),
the orthogonality measure ${\rm d}\lambda^{(1)}$ can be easily computed
using the Perron-Stieltjes inversion formula, noticing that,
writing ${\rm d}\lambda^{(1)}(x)=w^{(1)}(x){\rm d}x$, then
$$
w^{(1)}(x)=\frac{1}{\pi}\,\lim_{\epsilon\to0^+}\Im\big(F^{(1)}(x-i\epsilon)\big)\;,\quad -1<x<1\;.\mbox{\rm )}
$$
\medskip

%\item\label{Ex-cp5-A5}
%Extend Markov's Theorem for functions of the second kind.
%\medskip

\item\label{Ex-cp5-A6}
Let $\{U_n\}_{n\geq0}$ be the Chebyshev OPS of the second kind, which is orthogonal with respect to the positive-definite functional ${\bf u}\in\mathcal{P}'$ defined by
$$
\langle{\bf u},p\rangle:=\frac{2}{\pi}\,\int_{-1}^1p(x)\sqrt{1-x^2\,}\,{\rm d}x\; ,\quad p\in\mathcal{P}\;.
$$
Let $\{P_n\}_{n\geq0}$ be a sequence of polynomials defined by
$$
P_{2n+1}(x):=2^nxU_n\Big(\mbox{$\frac{x^2-5}{4}$}\big)\,,\quad
P_{2n}(x):=2^nU_n\Big(\mbox{$\frac{x^2-5}{4}$}\Big)+2^{n-1}U_{n-1}\Big(\mbox{$\frac{x^2-5}{4}$}\Big)\;,\quad n\in\mathbb{N}_0\;.
$$
\begin{enumerate}
\item
Prove that $\{P_n\}_{n\geq0}$ is a monic OPS with respect to a positive-definite functional, by showing that it fulfills the TTRR
$$
P_{-1}(x)=0\,,\quad P_0(x)=1 \;,\quad P_{n+1}(x)=xP_n(x)-\gamma_nP_{n-1}(x)\;,\quad n\in\mathbb{N}_0\;,
$$
where $\gamma_{2n}=1$ and $\gamma_{2n+1}=4$ for all $n\in\mathbb{N}_0$.
%, and with initial conditions $P_{-1}(x)=0$ and $P_1(x)=1$.
\item
%Show that the (positive-definite) functional ${\bf v}\in\mathcal{P}'$ with respect to which $\{P_n\}_{n\geq0}$ is orthogonal is defined by
%Prove that the orthogonality measure with respect to which  $\{P_n\}_{n\geq0}$ is orthogonal is given by
%$$
%{\rm d}\mu(x):=\frac{\chi_E(x)}{|x|}\,\sqrt{1-\Big(\mbox{$\frac{x^2-5}{4}$}\Big)^2}\,{\rm d}x\,,
%$$
%where $E:=[-3,-1]\cup[1,3]$. Is this measure unique? Why?
Prove that the spectral measure for $\{P_n\}_{n\geq0}$ (appearing in the spectral theorem for orthogonal polynomials) has distribution function $\psi$ given by
$$
\psi(x):=\int_{-\infty}^xw(t)\,{\rm d}t\,,\quad x\in\mathbb{R}\,,
$$
where
$$
\omega(t):=\left\{
\begin{array}{ccl}
\displaystyle\frac{1}{|t|}\,\sqrt{1-\Big(\mbox{$\frac{t^2-5}{4}$}\Big)^2} & \mbox{\rm if} & t\in E\,, \\ [1em]
0 & \mbox{\rm if} & t\not\in E\;,
\end{array}
\right.
$$
being $E:=[-3,-1]\cup[1,3]$. Is the spectral measure unique? Why?
%$$
%{\rm d}\mu(x):=\frac{\chi_E(x)}{|x|}\,\sqrt{1-\Big(\mbox{$\frac{x^2-5}{4}$}\Big)^2}\,{\rm d}x\,,
%$$
%where $E:=[-3,-1]\cup[1,3]$. Is this measure unique? Why?
\item
Use the monic OPS $\{P_n\}_{n\geq0}$ to show that the set $X_1\cup X_2$ in the statement of Markov's theorem cannot be replaced by $\sigma(\psi)$.
\end{enumerate}

\item\label{Ex-cp5-A6aa}
Let $\{P_n\}_{n\geq0}$ be a monic OPS with respect to a positive Borel measure $\mu$.
Denote the zeros of $P_n$ by $x_{n,1},x_{n,2},\ldots,x_{n,n}$, in increasing order, and
let $\{P_n^{(1)}\}_{n\geq0}$ be the sequence of numerator polynomials,
which is a monic OPS with respect to a positive Borel measure $\mu^{(1)}$,
and so $P_n^{(1)}$ has $n$ real and simple zeros for each $n\in\mathbb{N}$.
Denoting these zeros by $x_{n,1}^{(1)},x_{n,2}^{(1)},\ldots,x_{n,n}^{(1)}$, in increasing order,
prove the following interlacing property:
$$
x_{n+1,j}<x_{n,j}^{(1)}<x_{n+1,j+1}\;,\quad j=1,2,\ldots,n\;.
$$
Conclude that $\mbox{\rm co}\big(\mbox{\rm supp}(\mu^{(1)})\big)\subseteq\mbox{\rm co}\big(\mbox{\rm supp}(\mu)\big)$.
%Conclude that $\mbox{\rm co}\big(\psi^{(1)})\big)\subseteq\mbox{\rm co}(\psi)$,
%where $\psi^{(1)}$ and $\psi$ are the distribution functions of $\mu$ and $\mu^{(1)}$, respectively.
\smallskip

\item\label{Ex-cp5-6}
Suppose that ${\bf u}\in\mathcal{P}'$ is regular, normalized so that $u_0:=\langle{\bf u},1\rangle=1$,
and let $\{P_n\}_{n\geq0}$ be the monic OPS with respect to ${\bf u}$.
Let $\lambda\in\mathbb{C}\setminus\{0\}$ and $c\in\mathbb{C}$, and set
$${\bf u}^{\lambda,c}:=\bm{\delta}_c+\lambda(x-c)^{-1}{\bf u}\;.$$
\begin{enumerate}
\item
Prove that ${\bf u}^{\lambda,c}$ is regular if and only if $P_n(c)+\lambda P_{n-1}^{(1)}(c)\neq0$ for all $n\in\mathbb{N}$.
Moreover, under these conditions, setting %$a_0:=0$ and
$$
a_0:=0\;,\quad
a_n\equiv a_n^{\lambda,c}:=
-\frac{P_n(c)+\lambda P_{n-1}^{(1)}(c)}{P_{n-1}(c)+\lambda P_{n-2}^{(1)}(c)}\quad\mbox{\rm if}
\;\; n\geq1\,,$$
the monic OPS $\{P_n^{\lambda,c}\}_{n\geq0}$ with respect to ${\bf u}^{\lambda,c}$ is given by
$$
%P_0^{\lambda,c}(x)=1\;,\quad
%P_n^{\lambda,c}(x)=P_n(x)-\frac{P_n(c)+\lambda P_{n-1}^{(1)}(c)}{P_{n-1}(c)+\lambda P_{n-2}^{(1)}(c)}\,P_{n-1}(x)\;,
P_n^{\lambda,c}(x)=P_n(x)+a_n\,P_{n-1}(x)\;,
\quad n\in\mathbb{N}_0\;,
$$
and $\{P_n^{\lambda,c}\}_{n\geq0}$ fulfills the TTRR $$
P^{\lambda,c}_{n+1}(x)=\big(x-\beta_n^{\lambda,c}\big)P_n^{\lambda,c}(x)-\gamma_n^{\lambda,c}P^{\lambda,c}_{n-1}(x)\;,
\quad n\geq0\;,
$$
where $\beta_n^{\lambda,c}:=\beta_n+a_n-a_{n+1}$ ($n\geq0$),
$\gamma_1^{\lambda,c}:=\lambda a_1$, and $\gamma_n^{\lambda,c}:=\gamma_{n-1}a_n/a_{n-1}$ ($n\geq2$), being
%$a_0=0$ and $a_n:=-[P_n(c)+\lambda P_{n-1}^{(1)}(c)]/[P_{n-1}(c)+\lambda P_{n-2}^{(1)}(c)]$ if $n\geq1$,
$\{\beta_n\}_{n\geq0}$ and $\{\gamma_n\}_{n\geq1}$ the sequences of parameters appearing in the TTRR for $\{P_n\}_{n\geq0}$, so that
$P_{n+1}(x)=(x-\beta_n)P_n(x)-\gamma_nP_{n-1}(x)$, $n\geq0$, with $\beta_n\in\mathbb{C}$ and $\gamma_n\in\mathbb{C}\setminus\{0\}$ for each $n$.
\item
Suppose that the $\beta-$parameters vanish in the TTRR for $\{P_n\}_{n\geq0}$ (i.e., $\beta_n=0$ for each $n\geq0$).
Show that ${\bf u}^{\lambda,0}:=\bm{\delta}+\lambda x^{-1}{\bf u}$ is regular and the
corresponding parameters $a_n\equiv a_n^{\lambda,0}$ defined in (a) are given by
$$
\qquad\qquad
a_{2n}=-\frac{1}{\lambda}\frac{P_{2n}(0)}{P_{2n-2}^{(1)}(0)}
=\frac{1}{\lambda}\prod_{j=0}^{n-1}\frac{\gamma_{2j+1}}{\gamma_{2j}}\;,\quad
a_{2n-1}=-\lambda\frac{P_{2n-2}^{(1)}(0)}{P_{2n-2}(0)}
=-\lambda\prod_{j=1}^{n-1}\frac{\gamma_{2j}}{\gamma_{2j-1}} %\;,\quad n\geq1\;.
\vspace*{-0.25em}
$$
for each $n\geq1$ (with the conventions $\gamma_0:=1$ and empty product equals $1$).
\end{enumerate}
%\medskip
\bigskip

\item\label{Ex-cp5-7}
(Orthogonal polynomials on the semi-circle)
Let ${\bf v}\in\mathcal{P}'$ be defined by
$$
\langle{\bf v},p\rangle:=\frac{1}{\pi}\int_{0}^\pi p\big(e^{i\theta}\big)\,{\rm d}\theta\;,\quad p\in\mathcal{P}\;.
$$
\begin{enumerate}
\item
Show that
$$
{\bf v}=\bm{\delta}-\mbox{$\frac{2}{\pi i}$}\,x^{-1}{\bf u}\; ,
$$
where ${\bf u}\in\mathcal{P}'$ is the (positive-definite) Legendre functional
%with respect to which the Legendre polynomials are an OPS,
normalized so that
$$
\langle{\bf u},p\rangle:=\mbox{$\frac{1}{2}$}\int_{-1}^1 p(x)\,{\rm d}x\;,\quad p\in\mathcal{P}\;.
$$
%\smallskip
\noindent
({\sl Hint.} Note that $\int_\Gamma \frac{p(z)}{z}\,{\rm d}z=0$ for each $p\in\mathcal{P}$, where $\Gamma$ is the closed path on $\mathbb{C}$ defined by $\Gamma:=\Gamma_1+\ell_{\epsilon}^-+\Gamma_\epsilon+\ell_{\epsilon}^+$, $\Gamma_1$ and $\Gamma_\epsilon$ being semicircles on the upper semi-plane of radius $1$ and $\epsilon$, respectively, with $0<\epsilon<1$, $\Gamma_1$ starting at the point  $z=1$ and ending at $z=-1$, and $\Gamma_\epsilon$ starting at $z=-\epsilon$ and ending at $z=\epsilon$, and $\ell_{\epsilon}^{-}$ and $\ell_{\epsilon}^{+}$ are segments on the real line, joining the points $z=-1$ to $z=-\epsilon$, and $z=\epsilon$ to $z=1$, respectively. Consider the integrals along each of the paths
$\Gamma_1$, $\ell_{\epsilon}^-$, $\Gamma_\epsilon$, and $\ell_{\epsilon}^+$, and then take the limit as $\epsilon\to0^+$.)
\smallskip

%\noindent
%({\sl Hint.} Note that $\int_\Gamma \frac{p(z)}{z}\,{\rm d}z=0$ for each $p\in\mathcal{P}$, where $\Gamma$ is the closed path on $\mathbb{C}$ defined by $\Gamma:=\Gamma_1+\ell_{\epsilon}^-+\Gamma_\epsilon+\ell_{\epsilon}^+$, $\Gamma_1$ and $\Gamma_\epsilon$ being semicircles on the upper semi-plane of radius $1$ and $\epsilon$, respectively, with $0<\epsilon<1$, $\Gamma_1$ starting in the point  $(1,0)$ and ending in $(-1,0)$, and $\Gamma_\epsilon$ starting in $(-\epsilon,0)$ and ending in $(\epsilon,0)$, and $\ell_{\epsilon}^{-}$ and $\ell_{\epsilon}^{+}$ are segments on the real line, joining the points $(-1,0)$ to $(-\epsilon,0)$, and $(\epsilon,0)$ to $(1,0)$, respectively.)
%\smallskip

\item
Prove that ${\bf v}$ is regular and the monic OPS $\{Q_n\}_{n\geq0}$ with respect to ${\bf v}$ is given by
$$
Q_n(x)=P_n(x)
-\frac{2i}{2n-1}\left(\frac{\Gamma\Big(\frac{n+1}{2}\Big)}{\Gamma\Big(\frac{n}{2}\Big)}\right)^2\,
P_{n-1}(x)\;,\quad n\geq1\;,
$$
where $\{P_n\}_{n\geq0}$ is the (Legendre) monic OPS with respect to ${\bf u}$.

({\sl Hint.} Use exercise \ref{Ex-cp5-6})
%where
%$$
%a_n:=\frac{-2i}{2n+1}\left(\frac{\Gamma\Big(1+\frac{n}{2}\Big)}{\Gamma\Big(\frac{1+n}{2}\Big)}\right)^2\,,\quad n\geq1\;.
%$$
%\smallskip

%\noindent
%({\sl Hint.} Use exercise \ref{Ex1-hwk6})
\end{enumerate}
\end{enumerate}
%\bigskip
}
\medskip

\section*{Final remarks}

The proof of the Perron--Stieltjes inversion formula (Theorem \ref{Perron-Stieltjes})
is taken from the article \ref{Koelink2002-C4a} by Erik Koelink.
The proofs of Markov's Theorem and the lemmas before it are based on references
\ref{CharrisSS1991-C4a}, \ref{CharrisSS2003-C4a}, \ref{Chihara1978-C4a}, and \ref{VanAssche1991-C4a}.
Alternative proofs of Markov's Theorem appear in Berg's article \ref{CBerg1994-C5}
(based on the concept of weak convergence of measures),
and in the book \ref{NikishinSorokin1991i} by Nikishin and Sorokin.

Exercise {\bf 1} is a statement that appears in the book of Nikishin and Sorokin (p. 60), from where we have also taken exercises {\bf 2} and {\bf 3}. The result of exercise {\bf 4} appears e.g. in the article \ref{Maroni1991-C5a} by Maroni. The result of exercise {\bf 5} may be found in Gautschi's book \ref{Gautschi2004-C5aa}. Exercise {\bf 6} deals with a very special case of an OPS obtained from another one via a polynomial mapping. This kind of polynomial transformation between two families of OP has received considerable attention in the last decades (see e.g. \ref{GeronimoVanAssche1988-C5} and \ref{MarcioPetronilho2010-C5}). Exercise {\bf 7} appears e.g. in Chihara's book \ref{Chihara1978-C4a}. The OP on the semi-circle studied in exercise {\bf 9} were introduced by Gautschi and Milovanovi\'c in \ref{GautschiMilovanovic1986-C5}, where they have studied the main properties of such polynomials by a technique totally different from the one presented in this exercise. The approach considered in exercise {\bf 9} to the OP on the semi-circle, based on exercise {\bf 8}, is due to Maroni (cf. e.g. \ref{Maroni1991-C5a} and references therein).
\medskip

\section*{Bibliography}
\medskip

{\small
\begin{enumerate}[label=\emph{\rm [\arabic*]},leftmargin=*]
\item\label{CBerg1994-C5} C. Berg, {\it Markov's theorem revisited}, {J. Approx. Theory} {\bf 78} (1994), 260--275.
\item\label{CharrisSS1991-C4a} J. Charris, G. Salas, and V. Silva, {\it Polinomios ortogonales relacionados com problemas espectrales}, {Revista Colombiana de Matematicas} {\bf 27} (1991), 35--80.
\item\label{CharrisSS2003-C4a} J. Charris, B. Aldana, and G. Preciado-L\'opez, {\it Recurrence relations, continued fractions, and determining the spectral properties of orthogonal systems of polynomials} (In Spanish), {Rev. Acad. Colombiana Cienc. Exact. Fis. Natur} {\bf 27} (2003), n$^{\rm o}$104, 381--421.
\item\label{Chihara1978-C4a} T. S. Chihara, {\sl An introduction to orthogonal polynomials}, Gordon and Breach (1978).
%\item\label{Deift1997-C4a} P. Deift, {\sl Orthogonal Polynomials and Random Matrices: A Riemann-Hilbert Approach}, AMS Courant Lecture Notes {\bf 3} (2000).
%\item\label{Freud1971-C4a} G. Freud, {\sl Orthogonal polynomials}, Pergamon Press, Oxford (1971).
\item\label{MarcioPetronilho2010-C5} M. N. de Jesus and J. Petronilho, {\it On orthogonal polynomials obtained via polynomial mappings}, {J. Approx. Theory} {\bf 162} (2010) 2243--2277.
\item\label{Gautschi2004-C5aa} W. Gautschi, {\sl Orthogonal polynomials. Computation and Approximation}, Oxford University Press, Oxford (2004).
\item\label{GautschiMilovanovic1986-C5} W. Gautschi and G. V. Milovanovi\'c, {\it Polynomials orthogonal on the semicircle}, {J. Approx. Theory} {\bf 46} (1986), 230--250.
%\item\label{Ismail2004-C4a} M. E. H. Ismail, {\sl Classical and Quantum Orthogonal Polynomials in One Variable},
%          Cambridge University Press (2005) [paperback edition: 2009].
\item\label{GeronimoVanAssche1988-C5} J. Geronimo and W. Van Assche, {\it Orthogonal polynomials on several intervals via a polynomial mapping}, {Trans. Amer. Math. Soc.} {\bf 308} (1988) 559--581.
\item\label{Koelink2002-C4a} E. Koelink, {\sl Spectral theory and special functions}, Laredo Lectures on Orthogonal Polynomials and Special Functions (R. \'Alvarez-Nodarse, F. Marcell\'an, and W. Van Assche Eds.), Nova Science Publishers, 45-84 (2004).
\item\label{Maroni1991-C5a} P. Maroni, {\it Une th\'eorie alg\'ebrique des polyn\^omes orthogonaux. Applications aux polyn\^omes orthogonaux semiclassiques}, In C. Brezinski et al. Eds., Orthogonal Polynomials and Their Applications, Proc. Erice 1990, IMACS, Ann. Comp. App. Math. {\bf 9} (1991) 95--130.
\item\label{Remmert1997-C5a} R. Remmert, {\sl Classical Topics in Complex Function Theory}, Graduate Texts in Mathematics {\bf 172}, Springer (1997).
%\item\label{Szegoo1975-C4a} G. Szeg\"o, {\sl Orthogonal Polynomials}, AMS Colloq. Publ. {\bf 230} (1975), 4th ed.
%\item\label{Titchmarsh1964-C4a} E. C. Titchmarsh, {\sl The Theory of Functions}, Oxford University Press, Oxford (1964), corrected 2th ed.
\item\label{NikishinSorokin1991i} E. M. Nikishin and V. N. Sorokin, {\sl Rational approximations and orthogonality}, Translations of Mathematical Monographs {\bf 92} (1975) AMS.
\item\label{VanAssche1991-C4a} W. Van Assche, {\it Orthogonal polynomials, associated polynomials, and functions of the second kind}, {J. Comput. Appl. Math.} {\bf 37} (1991), 237--249.
\end{enumerate}
}

\chapter{Orthogonal polynomials and dual basis}

%\chapter{Orthogonal polynomials and dual basis.}
%More operations in $\mathcal{P}$ and $\mathcal{P}'$. }

\pagestyle{myheadings}\markright{Orthogonal polynomials and dual basis}
\pagestyle{myheadings}\markleft{J. Petronilho}

Every OPS is a simple set of polynomials, hence it has an
associated dual basis in $\mathcal{P}^\prime$.
In this text we present several properties of the dual basis associated with an OPS.
We also introduce some additional operations on the spaces $\mathcal{P}$ and  $\mathcal{P}'$
(the translation and homothetic operators on these spaces)
which appear as useful tools on the study of several classes of OP,
including the so called classical and semiclassical families.
%From an algebraic viewpoint, often the use of dual basis

\section{Orthogonal polynomials and dual basis}

We begin by establishing some connections between a regular functional
and the dual basis associated with the corresponding monic OPS.

\begin{snugshade}
\begin{theorem}\label{dualR}
Let ${\bf u}\in\mathcal{P}'$ be regular, $\{P_n\}_{n\geq0}$ the
corresponding monic OPS, and $\{{\bf a}_n\}_{n\geq0}$ the associated dual basis.
%Let the TTRR fulfilled by $\{P_n\}_{n\geq0}$ be
%\begin{equation}\label{ttrr-dual-Pn}
%xP_n(x)=P_{n+1}(x)+\beta_nP_n(x)+\gamma_nP_{n-1}(x)\;,\quad n=0,1,2\ldots\;,
%\end{equation}
%with $P_{-1}(x)=0$ and $P_0(x)=1$, being $\beta_n\in\mathbb{C}$ and
%$\gamma_{n}\in\mathbb{C}\setminus\{0\}$.
Then:
\begin{itemize}
\item[{\rm (i)}] For each $n\in\mathbb{N}_0$, ${\bf a}_n$ is explicitly given by
$$
{\bf a}_n=\frac{P_n}{\langle {\bf u},P_n^2\rangle}\, {\bf u}\;.
$$
As a consequence, $\{P_n\}_{n\geq0}$ is a monic OPS with respect to ${\bf a}_0$, being
$$
{\bf u}=u_0\,{\bf a}_0\;.
$$
\item[{\rm (ii)}] Let ${\bf v}\in\mathcal{P}'$ and $N\in \mathbb{N}_0$ such that
%such that $\langle {\bf v},P_n\rangle=0$ for all $n\geq N+1$, then
$$
\langle {\bf v},P_n\rangle=0 \;\;\mbox{\rm if}\;\; n\geq N+1\;.
$$
Then,\vspace*{-0.5em}
%$$
%{\bf v}=\sum_{j=0}^{N}\langle {\bf v},P_j\rangle\, {\bf a}_j\; .
%$$
%\item[(ii)] If ${\bf v}\in\mathcal{P}^*$ and there exists $N\in \mathbb{N}_0$
%such that $\langle {\bf v},P_n\rangle=0$ for all $n\geq N+1$, then
%$$
%{\bf v}=\sum_{j=0}^{N}\langle {\bf v},P_j\rangle\, {\bf a}_j\; .
%$$
%there exists $\phi\in\mathcal{P}$, with
%$\deg\phi\leq N$, such that
$$
{\bf v}=\sum_{j=0}^{N}\langle {\bf v},P_j\rangle\, {\bf a}_j=\phi\,{\bf u}\; ,\quad
%$$
%where $\phi$ is polynomial given explicitly by
%$$
\phi(x):=\sum_{j=0}^{N}\frac{\langle {\bf v},P_j\rangle}{\langle{\bf u},P_j^2\rangle}\, P_j(x)\;.
$$
Further, $\deg\phi\leq N$, and $\deg\phi=N$ if and only if $\langle {\bf v},P_N\rangle\neq0$.\medskip
\item[{\rm (iii)}]
%$\{P_n\}_{n\geq0}$ is a monic OPS with respect to the functional ${\bf a}_0$. Indeed:
%$$
%{\bf u}=u_0\,{\bf a}_0\;.
%$$
%\item[{\rm (iv)}]
Let the TTRR fulfilled by $\{P_n\}_{n\geq0}$ be
$$%\begin{equation}\label{ttrr-dual-Pn}
xP_n(x)=P_{n+1}(x)+\beta_nP_n(x)+\gamma_nP_{n-1}(x)\;,\quad n\in\mathbb{N}_0\;,
$$%\end{equation}
with $P_{-1}(x)=0$, $P_0(x)=1$, $\beta_n\in\mathbb{C}$, and
$\gamma_{n}\in\mathbb{C}\setminus\{0\}$.
Then $\{{\bf a}_n\}_{n\geq0}$ fulfills
$$
x\,{\bf a}_{n}={\bf a}_{n-1}+\beta_n\,{\bf a}_n+\gamma_{n+1}\,{\bf a}_{n+1} \; , \quad n\in\mathbb{N}_0\;,
$$
with initial conditions ${\bf a}_{-1}={\bf 0}$ and ${\bf a}_{0}=u_0^{-1}\,{\bf u}$.
\end{itemize}
\end{theorem}
\end{snugshade}

{\it Proof.}
By Theorem \ref{expDB1}, for each $n\in\mathbb{N}_0$ we may write
$$
P_n{\bf u}=\sum_{j\geq0}\langle P_n{\bf u},P_j\rangle {\bf a}_j=
\sum_{j\geq0}\langle {\bf u},P_nP_j\rangle {\bf a}_j= \langle {\bf u},P_n^2\rangle {\bf a}_n\,,
$$
hence (i) is proved. Statement (ii) follows immediately from (i) using again Theorem \ref{expDB1}.
Finally, for all $n\in\mathbb{N}$ and $j\in\mathbb{N}_0$, we have
$$
\begin{array}{rcl}
x\,{\bf a}_n & = & \displaystyle
\frac{xP_n}{\langle {\bf u},P_n^2\rangle}\,{\bf u}\, = \,
\frac{P_{n+1}+\beta_nP_n+\gamma_{n}P_{n-1}}{\langle {\bf u},P_n^2\rangle}\,{\bf u} \\ [1em]
 & = & \displaystyle
\frac{\langle {\bf u},P_{n+1}^2\rangle}{\langle {\bf u},P_n^2\rangle}\frac{P_{n+1}}{\langle {\bf u},P_{n+1}^2\rangle}{\bf u}
+\beta_n\frac{P_n}{\langle {\bf u},P_n^2\rangle}\,{\bf u}
+\gamma_{n}\frac{\langle {\bf u},P_{n-1}^2\rangle}{\langle {\bf u},P_n^2\rangle}\frac{P_{n-1}}{\langle {\bf u},P_{n-1}^2\rangle}{\bf u} \\ [1em]
& = & \gamma_{n+1}\,{\bf a}_{n+1}+\beta_n\,{\bf a}_n+{\bf a}_{n-1}\;,
\end{array}
$$
where we have used the relation
$\gamma_i=\langle {\bf u},P_i^2\rangle/\langle {\bf u},P_{i-1}^2\rangle$, $i\in\mathbb{N}$
(cf. Corollary \ref{cor-TTRR}).
\qed

%{\it Proof.}
%For every $n,j\in\mathbb{N}_0$, we have
%$$
%\Big\langle\,\frac{P_n}{\langle {\bf u},P_n^2\rangle}\, {\bf u},P_j\,\Big\rangle
%=\frac{1}{\langle {\bf u},P_n^2\rangle}\,\langle\, {\bf u},P_n
%P_j\,\rangle=\delta_{n,j}=\langle {\bf a}_n,P_j\rangle\; ,
%$$
%so the actions of the linear functionals $\; \frac{P_n}{\langle
%{\bf u},P_n^2\rangle}\, {\bf u}\;$ and $\;{\bf a}_n\;$ over each
%element of the basis $\{P_j\}_{j\geq0}$ of $\mathcal{P}$ coincide, hence
%these functional coincide over all the space
%$\mathcal{P}$. This proves (i). Statement (ii) is an
%immediate consequence of Theorem \ref{expDB1} and (i).
%Finally, for all $n\in\mathbb{N}$ and $j\in\mathbb{N}_0$, we have
%$$
%\begin{array}{rcl}
%\langle x\,{\bf a}_n,P_j\rangle & = & \langle {\bf a}_n,xP_j\rangle \, = \,
%\langle {\bf a}_n,P_{j+1}+\beta_j P_j+\gamma_j P_{j-1}\rangle \\
%\rule{0pt}{1.2em}
%&=& \langle \, {\bf a}_{n-1}+\beta_n\,{\bf a}_n+\gamma_{n+1}\,{\bf a}_{n+1} , P_j\, \rangle \; ,
%\end{array}
%$$
%which gives (iii).
%\qed

%As a consequence, we can state the following useful result.

\begin{snugshade}
\begin{corollary}\label{OPSwrtv}
Let $\{P_n\}_{n\geq0}$ be a monic OPS (with respect to some functional in
$\mathcal{P}'$) and let ${\bf v}\in\mathcal{P}'$. Then $\{P_n\}_{n\geq0}$
is a monic OPS with respect to ${\bf v}$ if and only if
\begin{equation}
\langle{\bf v},1\rangle\neq0 \;, \qquad \langle{\bf v},P_n\rangle=0 \;, \;\; n\geq1 \, .
\label{uSPOMregular}
\end{equation}
\end{corollary}
\end{snugshade}

{\it Proof.}
Clearly, if $\{P_n\}_{n\geq0}$ is a monic OPS with respect to ${\bf v}$
then (\ref{uSPOMregular}) holds. Conversely, if (\ref{uSPOMregular}) holds,
then by (ii) in Theorem \ref{dualR},
$${\bf v}=\langle{\bf v},1\rangle\, {\bf a}_0=\frac{\langle{\bf v},1\rangle}{\langle{\bf u},1\rangle}\, {\bf u}\,,$$
where ${\bf a}_0$ is the first element of the dual basis associated with $\{P_n\}_{n\geq0}$, and
${\bf u}$ is the regular functional with respect to which $\{P_n\}_{n\geq0}$ is an OPS.
Since, by hypothesis, $\langle{\bf v},1\rangle\neq0$, it follows
that $\{P_n\}_{n\geq0}$ is a monic OPS with respect to ${\bf v}$.
\qed
\medskip

\section{The translation and homothetic operators}

\begin{snugshade}
\begin{definition}[translation operators]\label{taub}
Let $b\in\mathbb{C}$.
%Let ${\bf u}\in\mathcal{P}'$, $a\in\mathbb{C}\setminus\{0\}$, and $b\in\mathbb{C}$.
\begin{enumerate}
\item[{\rm (i)}]
The {\sl translator operator} on $\mathcal{P}$ is $\tau_b:\mathcal{P}\to\mathcal{P}$ $(p\mapsto\tau_bp)$ defined by
\begin{equation}\label{def-taub}
\tau_bp(x):=p(x-b)\;,\quad p\in\mathcal{P}\;;
\end{equation}
\item[{\rm (ii)}]
The {\sl translator operator} on $\mathcal{P}'$ is $\bm{\tau}_b:=\tau_{-b}^{\,\prime}$, i.e., $\bm{\tau}_b:\mathcal{P}'\to\mathcal{P}'$
%$({\bf u}\mapsto\bm{\tau}_b{\bf u})$,
is the dual operator of $\tau_{-b}$, so that
\begin{equation}\label{def-taub-dual}
\langle\bm{\tau}_b{\bf u},p\rangle:=\langle{\bf u},\tau_{-b}p\rangle=\langle{\bf u},p(x+b)\rangle\;,
\quad {\bf u}\in\mathcal{P}'\;,\quad p\in\mathcal{P}\;.
\end{equation}
\end{enumerate}
\end{definition}
\end{snugshade}

Notice that the moments of the functional $\bm{\tau}_b{\bf u}$ are
\begin{snugshade}
\begin{equation}\label{moments-taubu}
\big(\bm{\tau}_b{\bf u}\big)_n=\sum_{j=0}^n\binom{n}{j}b^{n-j}u_j
=\sum_{i+j=n}\binom{n}{i}b^{i}u_j\;,\quad n\in\mathbb{N}_0\;.
\end{equation}
\end{snugshade}

Indeed, for each $n\in\mathbb{N}_0$,
$$
\big(\bm{\tau}_b{\bf u}\big)_n:=\langle\bm{\tau}_b{\bf u},x^n\rangle
=\langle{\bf u},\tau_{-b}x^n\rangle=\langle{\bf u},(x+b)^n\rangle
=\sum_{j=0}^n\binom{n}{j}b^{n-j}\langle{\bf u},x^j\rangle \,.
$$

\begin{snugshade}
\begin{definition}[homothetic operators]\label{ha}
Let $a\in\mathbb{C}\setminus\{0\}$.
%Let ${\bf u}\in\mathcal{P}'$, $a\in\mathbb{C}\setminus\{0\}$, and $b\in\mathbb{C}$.
\begin{enumerate}
\item[{\rm (i)}]
The {\sl homothetic operator} on $\mathcal{P}$ is $h_a:\mathcal{P}\to\mathcal{P}$ $(p\mapsto h_ap)$ defined by
\begin{equation}\label{def-ha}
h_ap(x):=p(ax)\;,\quad p\in\mathcal{P}\;.
\end{equation}
\item[{\rm (ii)}]
The {\sl homothetic operator} on $\mathcal{P}'$ is $\bm{h}_a:=h_a^{\,\prime}$,
i.e., $\bm{h}_a:\mathcal{P}'\to\mathcal{P}'$ is the dual operator of $h_a$, so that
\begin{equation}\label{def-ha-dual}
\langle\bm{h}_a{\bf u},p\rangle:=\langle{\bf u},h_ap\rangle=\langle{\bf u},p(ax)\rangle\;,
\quad {\bf u}\in\mathcal{P}'\;,\quad p\in\mathcal{P}\;.
\end{equation}
\end{enumerate}
\end{definition}
\end{snugshade}

The moments of the functional $\bm{h}_a{\bf u}$ are
\begin{snugshade}
\begin{equation}\label{moments-hau}
\big(\bm{h}_a{\bf u}\big)_n=a^nu_n\;,\quad n\in\mathbb{N}_0\;.
\end{equation}
\end{snugshade}

In the next proposition we list some useful properties involving the translation and homothetic operators.

\begin{snugshade}
\begin{proposition}\label{prop-taubha}
Let $a\in\mathbb{C}\setminus\{0\}$, $b\in\mathbb{C}$,
${\bf u}\in\mathcal{P}'$, and $p\in\mathcal{P}$.
Then:\vspace*{-0.5em}
\begin{multicols}{2}
\begin{enumerate}
\item[{\rm 1.}] $\tau_0p=h_1p=p$ \vspace*{0.25em}
\item[{\rm 2.}] $\big(\tau_b\circ \tau_{-b}\big)p=\big(\tau_{-b}\circ \tau_b\big)p=p$ \vspace*{0.25em}
\item[{\rm 3.}] $\big(h_a\circ h_{a^{-1}}\big)p=\big(h_{a^{-1}}\circ h_a\big)p=p$ \vspace*{0.25em}
\item[{\rm 4.}] $\big(h_a\circ\tau_b\big)p=\big(\tau_{b/a}\circ h_a\big)p$ \vspace*{0.25em}
\item[{\rm 5.}] $\big(\tau_b\circ h_a\big)p=\big(h_a\circ\tau_{ab}\big)p$ \vspace*{0.25em}
\item[{\rm 6.}\hspace*{3em}]\hspace*{-3em} $\bm{\tau}_0{\bf u}=\bm{h}_1{\bf u}={\bf u}$ \vspace*{0.25em}
\item[{\rm 7.}\hspace*{3em}] \hspace*{-3em} $\big(\bm{\tau}_b\circ \bm{\tau}_{-b}\big){\bf u}=\big(\bm{\tau}_{-b}\circ\bm{\tau}_b\big){\bf u}={\bf u}$ \vspace*{0.25em}
\item[{\rm 8.}\hspace*{3em}] \hspace*{-3em} $\big(\bm{h}_a\circ \bm{h}_{a^{-1}}\big){\bf u}=\big(\bm{h}_{a^{-1}}\circ \bm{h}_a\big){\bf u}={\bf u}$ \vspace*{0.25em}
\item[{\rm 9.}\hspace*{3em}]\hspace*{-3em} $\big(\bm{h}_a\circ\bm{\tau}_{b}\big){\bf u}=\big(\bm{\tau}_{ab}\circ \bm{h}_a\big){\bf u}$ \vspace*{0.25em}
\item[{\rm 10.}\hspace*{3em}] \hspace*{-3em} $\big(\bm{\tau}_b\circ \bm{h}_a\big){\bf u}=\big(\bm{h}_a\circ\bm{\tau}_{b/a}\big){\bf u}$
\end{enumerate}
\end{multicols}
\end{proposition}
\end{snugshade}

{\it Proof.} Properties 1, 2, 3, 6, 7, and 8 follow easily by straightforward computations.
The proof of property 4 can be done as follows: %We have
$$
\begin{array}{rcl}
\big(h_a\circ\tau_b\big)p(x)&=&h_a\big[(\tau_bp)(x)\big]%=h_a\big[p(x-b)\big]
=(\tau_bp)(ax)=p(ax-b)= p\left(a\big(x-\frac{b}{a}\big)\right) \\ [0.5em]
&=&\tau_{b/a}\big[p(ax)\big]=\tau_{b/a}\big[h_ap(x)\big]
=\big(\tau_{b/a}\circ h_a\big)p(x)\;.
\end{array}
$$
%hence property 4 is proved.
Replacing $b$ by $ab$ in property 4 we obtain property 5.
To prove property 9, notice that
$$
\begin{array}{rcl}
\langle\big(\bm{h}_a\circ\bm{\tau}_{b}\big){\bf u},p\rangle&=&
\langle\bm{\tau}_{b}{\bf u},h_ap\rangle = \langle{\bf u},\big(\tau_{-b}\circ h_a\big)p\rangle
% \stackrel{{\rm (2)}}{=}
=\langle{\bf u},\big(h_a\circ\tau_{-ab} \big)p\rangle \\ [0.5em]
&=&\langle\bm{h}_a{\bf u},\tau_{-ab} p\rangle=\langle\big(\bm{\tau}_{ab}\circ \bm{h}_a\big){\bf u},p\rangle\;,
\end{array}
$$
where in the third equality we have used property 5.
Finally, replacing $b$ by $b/a$ in property 9 we obtain property 10.
\qed
\medskip

Properties 2 and 3 show that the operators $\tau_b$ and $h_a$ are invertible in $\mathcal{P}$, being
\begin{snugshade}\vspace*{-0.5em}
\begin{equation}
\tau_b^{-1}p=\tau_{-b}p\;,\quad h_a^{-1}p=h_{a^{-1}}p\;,\quad
a\in\mathbb{C}\setminus\{0\}\;,\;b\in\mathbb{C}\;,\; p\in\mathcal{P}\;.
\end{equation}
\end{snugshade}
\noindent
Similarly, properties 6 and 7 show that $\bm{\tau}_b$ and $\bm{h}_a$ are invertible in $\mathcal{P}'$, being
\begin{snugshade}\vspace*{-0.5em}
\begin{equation}
\bm{\tau}_b^{-1}{\bf u}=\bm{\tau}_{-b}{\bf u}\;,\quad \bm{h}_a^{-1}{\bf u}=\bm{h}_{a^{-1}}{\bf u}\;,\quad
a\in\mathbb{C}\setminus\{0\}\;,\;b\in\mathbb{C}\;,\; {\bf u}\in\mathcal{P}'\;.
\end{equation}
\end{snugshade}
\noindent
As a consequence, we also deduce
\begin{snugshade}\vspace*{-0.5em}
\begin{equation}
\begin{array}{c}
\big(h_a\circ\tau_b\big)^{-1}=\tau_{-b}\circ h_{a^{-1}}\;,\qquad
\big(\tau_b\circ h_a\big)^{-1}=h_{a^{-1}}\circ\tau_{-b}\quad\mbox{\rm on $\mathcal{P}$}\;; \\ [0.5em]
\big(\bm{h}_a\circ\bm{\tau}_b\big)^{-1}=\bm{\tau}_{-b}\circ \bm{h}_{a^{-1}}\;,\qquad
\big(\bm{\tau}_b\circ \bm{h}_a\big)^{-1}=\bm{h}_{a^{-1}}\circ\bm{\tau}_{-b}\quad\mbox{\rm on $\mathcal{P}'$}\;.
\end{array}
\end{equation}
\end{snugshade}

Finally, we point out the following property that one should keep in mind
(it follows immediately from the proof of property 4 above replacing $b$ by $-b$):
\begin{snugshade}\vspace*{-0.5em}
\begin{equation}\label{transfAfim}
\big(h_a\circ\tau_{-b}\big)p(x)=p(ax+b)\;,\quad
a\in\mathbb{C}\setminus\{0\}\;,\;b\in\mathbb{C}\;,\; p\in\mathcal{P}\;.
\end{equation}
\end{snugshade}

The next proposition is of fundamental importance for a
rigorous treatment of the classification problem of semiclassical OP.

\begin{snugshade}
\begin{theorem}\label{Thm-u-equiv-v}
Define a binary relation on $\mathcal{P}'$ as follows: for every ${\bf u},{\bf v}\in\mathcal{P}'$,
\begin{equation}\label{u-equiv-v}
{\bf u}\sim{\bf v}\qquad\mbox{\rm iff}\qquad
\exists a\in\mathbb{C}\setminus\{0\}\;,\;\; \exists b\in\mathbb{C}\;:\;\;
{\bf v}=\big(\bm{h}_{a^{-1}}\circ\bm{\tau}_{-b}\big){\bf u}\;.
\end{equation}
Then, $\;\sim$ is an equivalent relation on $\mathcal{P}'$.
\end{theorem}
\end{snugshade}

{\it Proof.} Let ${\bf u},{\bf v},{\bf w}\in\mathcal{P}'$.
Since $${\bf u}=\big(\bm{h}_{1}\circ\bm{\tau}_{0}\big){\bf u}\;,$$
then ${\bf u}\sim{\bf u}$, so that the binary relation $\;\sim$ is reflexive.
To prove that it is symmetric, assume that ${\bf u}\sim{\bf v}$.
Then (\ref{u-equiv-v}) holds. Therefore, we may write
$$
{\bf u}=\big(\bm{h}_{a^{-1}}\circ\bm{\tau}_{-b}\big)^{-1}{\bf v}
=\big(\bm{\tau}_{b}\circ\bm{h}_{a}\big){\bf v}
=\big(\bm{h}_{a}\circ\bm{\tau}_{b/a}\big){\bf v}
=\big(\bm{h}_{c^{-1}}\circ\bm{\tau}_{-d}\big){\bf v}\;,
$$
where $c:=a^{-1}\in\mathbb{C}\setminus\{0\}$ and $d:=-b/a\in\mathbb{C}$.
(Notice also that the third equality follows from property 10 in Proposition \ref{prop-taubha}.)
Thus, ${\bf v}\sim{\bf u}$.
Finally, to prove that  $\;\sim$ is transitive, suppose that
${\bf u}\sim{\bf v}$ and ${\bf v}\sim{\bf w}$.
Then, there exists $a,c\in\mathbb{C}\setminus\{0\}$ and $b,d\in\mathbb{C}$ such that
$$
{\bf v}=\big(\bm{h}_{a^{-1}}\circ\bm{\tau}_{-b}\big){\bf u}\;,\quad
{\bf w}=\big(\bm{h}_{c^{-1}}\circ\bm{\tau}_{-d}\big){\bf v}\;.
$$
As a consequence, we may write
$$
\begin{array}{rcl}
{\bf w}&=&\big(\bm{h}_{c^{-1}}\circ\bm{\tau}_{-d}\big)\big(\bm{h}_{a^{-1}}\circ\bm{\tau}_{-b}\big){\bf u}
=\big(\bm{h}_{c^{-1}}\circ(\bm{\tau}_{-d}\circ\bm{h}_{a^{-1}})\circ\bm{\tau}_{-b}\big){\bf u} \\ [0.5em]
&=& \big(\bm{h}_{c^{-1}}\circ(\bm{h}_{a^{-1}}\circ\bm{\tau}_{-ad})\circ\bm{\tau}_{-b}\big){\bf u}
=\big((\bm{h}_{c^{-1}}\circ\bm{h}_{a^{-1}})\circ(\bm{\tau}_{-ad}\circ\bm{\tau}_{-b})\big){\bf u} \\ [0.5em]
&=& \big(\bm{h}_{a^{-1}c^{-1}}\circ\bm{\tau}_{-ad-b}\big){\bf u}
=\big(\bm{h}_{\alpha^{-1}}\circ\bm{\tau}_{-\beta}\big){\bf u}\;,
\end{array}
$$
where $\alpha:=ac\in\mathbb{C}\setminus\{0\}$ and $\beta:=b+ad\in\mathbb{C}$.
Therefore,  ${\bf u}\sim{\bf w}$.%, which completes de proof.
\qed

\begin{remark}\em
The relation between ${\bf u}$ and ${\bf v}$ in $(\ref{u-equiv-v})$ may be expressed as
\begin{snugshade}
\begin{equation}\label{u-equiv-v2}
%\langle{\bf v},x^n\rangle=\big\langle{\bf u},(ax+b)^n\big\rangle\;,\quad n=0,1,2,\ldots\;.
\langle{\bf v},x^n\rangle=\big\langle{\bf u},\Big(\mbox{$\frac{x-b}{a}$}\Big)^n\big\rangle\;,\quad n=0,1,2,\ldots\;.
\end{equation}
\end{snugshade}
\end{remark}

\begin{snugshade}
\begin{theorem}\label{uv-equiv}
Let $\{P_n\}_{n\geq0}$ be a simple set in $\mathcal{P}$ and $\{{\bf a}_n\}_{n\geq0}$ its associated dual basis.
Let $a\in\mathbb{C}\setminus\{0\}$ and $b\in\mathbb{C}$. Define
\begin{equation}\label{QnAfimPn}
Q_n:=a^{-n}\big(h_a\circ\tau_{-b}\big)P_n\;,\quad n=0,1,2,\ldots
\end{equation}
Then $\{Q_n\}_{n\geq0}$ is a simple set in $\mathcal{P}$, and
its dual basis, $\{{\bf b}_n\}_{n\geq0}$, is given by
\begin{equation}\label{dual-bas-Qn}
{\bf b}_n=a^n\big(\bm{h}_{a^{-1}}\circ\bm{\tau}_{-b}\big){\bf a}_n\;,\quad n=0,1,2\ldots\;.
\end{equation}
\end{theorem}
\end{snugshade}

{\it Proof.}
It is clear that $\{Q_n\}_{n\geq0}$ is a simple set in $\mathcal{P}$.
Moreover, for every $n,k\in\mathbb{N}_0$,
$$
\begin{array}{rcl}
\langle{\bf b}_n,Q_k\rangle&=&
a^{n-k}\langle\big(\bm{h}_{a^{-1}}\circ\bm{\tau}_{-b}\big){\bf a}_n,
\big(h_a\circ\tau_{-b}\big)P_k\rangle
=a^{n-k}\langle{\bf a}_n,\big(\tau_{b}\circ h_{a^{-1}}\big)
\big(h_a\circ\tau_{-b}\big)P_k\rangle \\ [0.25em]
&=& a^{n-k}\langle{\bf a}_n,P_k\rangle
=a^{n-k}\delta_{n,k}
=\delta_{n,k}\,,
\end{array}
$$
hence $\{{\bf b}_n\}_{n\geq0}$ is the dual basis associated with $\{Q_n\}_{n\geq0}$.
\qed

\begin{remark}\em
By (\ref{transfAfim}) we see that the polynomial $Q_n$ in (\ref{QnAfimPn}) is indeed
\begin{snugshade}\vspace*{-0.25em}
\begin{equation}
Q_n(x):=a^{-n}P_n(ax+b)\;,\quad n=0,1,2,\ldots\;,
\end{equation}
\end{snugshade}
\noindent
so that $Q_n$ is obtained from $P_n$ by an affine change of the variable,
being $Q_n$ normalized so that it becomes a monic polynomial whenever $P_n$ is monic.
\end{remark}

\begin{snugshade}
\begin{theorem}\label{OPSequiv}
Under the hypothesis of Theorem \ref{uv-equiv}, assume further that
$\{P_n\}_{n\geq0}$ is a monic OPS with respect to the functional ${\bf u}\in\mathcal{P}'$,
and let
\begin{equation}\label{ttrr-afim-Pn}
xP_n(x)=P_{n+1}(x)+\beta_nP_n(x)+\gamma_nP_{n-1}(x)\;,\quad n=0,1,2\ldots\;,
\end{equation}
be the TTRR fulfilled by $\{P_n\}_{n\geq0}$, with initial conditions
$P_{-1}(x)=0$ and $P_0(x)=1$, being $\beta_n\in\mathbb{C}$ and
$\gamma_{n}\in\mathbb{C}\setminus\{0\}$. % for each $n\in\mathbb{N}_0$.
Then, $\{Q_n\}_{n\geq0}$ is a monic OPS with respect to %the functional
\begin{equation}\label{Qn-afim-OPS}
{\bf v}:=\big(\bm{h}_{a^{-1}}\circ\bm{\tau}_{-b}\big){\bf u}\;,
\end{equation}
and the TTRR fulfilled by $\{Q_n\}_{n\geq0}$ is
\begin{equation}\label{ttrr-afim-Qn}
xQ_n(x)=Q_{n+1}(x)+\widehat{\beta}_nQ_n(x)+\widehat{\gamma}_nQ_{n-1}(x)\;,\quad n=0,1,2\ldots\;,
\end{equation}
with initial conditions $Q_{-1}(x)=0$ and $Q_0(x)=1$, where
\begin{equation}\label{beta-gamma-ttrr-afim-Qn}
\widehat{\beta}_n:=\frac{\beta_n-b}{a}\;,\quad \widehat{\gamma}_n:=\frac{\gamma_n}{a^2}\;.
\end{equation}
\end{theorem}
\end{snugshade}

{\it Proof.}
Changing $x$ into $ax+b$ in (\ref{ttrr-afim-Pn}) and then multiplying both sides
of the resulting equality by $a^{-n-1}$,
we obtain (\ref{ttrr-afim-Qn}). Since $\{Q_n\}_{n\geq0}$ satisfies (\ref{ttrr-afim-Qn})
and $\widehat{\gamma}_n\neq0$ for each $n\geq1$, then
it is a monic OPS (by Favard's Theorem). By Theorem \ref{uv-equiv},
the dual basis associated with $\{Q_n\}_{n\geq0}$ is given by (\ref{dual-bas-Qn}).
Moreover, by Theorem \ref{dualR}--(i),
$\{Q_n\}_{n\geq0}$ is a monic OPS with respect to ${\bf b}_0$.
Therefore, since
$$
{\bf b}_0=\big(\bm{h}_{a^{-1}}\circ\bm{\tau}_{-b}\big){\bf a}_0=
\big(\bm{h}_{a^{-1}}\circ\bm{\tau}_{-b}\big)u_0^{-1}{\bf u}=u_0^{-1}{\bf v}\;,
$$
so that ${\bf v}=u_0{\bf b}_0$ (being $u_0\neq0$),
we conclude that $\{Q_n\}_{n\geq0}$ is a monic OPS with respect to ${\bf v}$.
\qed

\section*{Exercises}
%\bigskip

{\small
%\noindent
\begin{enumerate}[label=\emph{\bf \arabic*.},leftmargin=*]
%\item
%Given a monic polynomial $P_n$ of degree $n$,
%denote by $P_n^{[k]}$ the monic polynomial of degree $n$ defined by
%$$
%P_n^{[k]}(x):=\frac{{\rm d}^k}{{\rm d}x^k}\,\frac{P_{n+k}(x)}{(n+1)_k}\quad(k,n\in\mathbb{N}_0)\;,
%$$
%where, for a given $\alpha\in\mathbb{C}$, $(\alpha)_n$ is the {\it Pochhammer symbol}, defined as
%$$
%(\alpha)_0:=1\;,\qquad (\alpha)_n:=\alpha(\alpha+1)\cdots(\alpha+n-1)\;,\quad n\in\mathbb{N}\;.
%$$
%Clearly, if $\{P_n\}_{n\geq0}$ is a simple set in $\mathcal{P}$
%(which needs not to be an OPS), then so is $\{P_n^{[k]}\}_{n\geq0}$.
%Let $\{{\bf a}_n\}_{n\geq0}$ and $\big\{{\bf a}_n^{[k]}\,\big\}_{n\geq0}$
%be the dual basis in $\mathcal{P}'$ associated with $\{P_n\}_{n\geq0}$ and $\{P_n^{[k]}\}_{n\geq0}$, respectively.
%Prove that
%$$
%D^{k}\big(\,{\bf a}_n^{[k]}\,\big)=(-1)^k(n+1)_k\,{\bf a}_{n+k} \;,\quad k,n\in\mathbb{N}_0 \, .
%$$
\item
Let ${\bf u}\in\mathcal{P}'$ be regular, $\{P_n\}_{n\geq0}$ the corresponding monic OPS,
and $\{{\bf a}_n\}_{n\geq0}$ its dual basis in $\mathcal{P}'$.
Consider the monic OPS $\{P_n^{(k)}\}_{n\geq0}$ (of the associated polynomials of order $k$, $k\in\mathbb{N}$)
and let $\{{\bf a}_n^{(k)}\}_{n\geq0}$ be its dual basis. Show that
%$$
%{\bf a}_n^{(1)}=u_0\big(x{\bf a}_{n+1}\big){\bf u}^{-1}\;,\quad n\in\mathbb{N}_0\;.
%%(k,n\in\mathbb{N}_0)\;,
%$$
$$x^{-1}\big({\bf a}_n^{(k)} {\bf a}_{k-1}\big) ={\bf a}_{n+k}\;,\quad n\in\mathbb{N}_0\;.$$
%$${\bf a}_n^{(k)} {\bf a}_{k-1}={\bf a}_{n+k}\;,\quad n\in\mathbb{N}_0\;.$$
In particular,
${\bf a}_n^{(1)}=\big(x{\bf a}_{n+1}\big){\bf a}_0^{-1}$ for each $n\in\mathbb{N}_0$,
and so $\{P_n^{(1)}\}_{n\geq0}$ is an OPS with respect to the functional ${\bf u}^{(1)}$ given by
$${\bf u}^{(1)}=c\,(xP_1{\bf u}){\bf u}^{-1}=-cu_0x^2{\bf u}^{-1}\;,\quad
c:=u_0^{(1)}/\gamma_1\in\mathbb{C}\setminus\{0\}\,.\;$$
(This relation suggests taking the normalization $u_0^{(1)}:=\gamma_1$,
which is henceforth considered the standard normalization of the functional ${\bf u}^{(1)}$.)
\end{enumerate}
%\bigskip
}
\medskip

\section*{Final remarks}

This short text is based on the works \ref{Maroni1991-C6} and \ref{Maroni1994-C6} by Pascal Maroni,
although some of the results therein may be found also in Chihara's book.
The equivalence relation (\ref{u-equiv-v}) appears in Maroni's work \ref{Maroni1994-C6} (see p.\;19 therein).
Some detailed computations concerning the results presented in this text may be found in the master thesis \ref{AFLoureiro2000-C6}
(under the co-supervision of Pascal Maroni and Z\'elia da Rocha).

Finally we mention that the results contained in this text are of an elementary nature and they could be left as exercises included in other texts. Our option to include them as an autonomous text is due to the advantage that results from its systematized presentation for reading some subsequent texts.
\medskip

%\section*{Notas finais}
%\section*{Coment\'arios finais}
%
%
%Este breve texto foi elaborado com base nos trabalhos
%\ref{Maroni1991-C6} e \ref{Maroni1994-C6} de Pascal Maroni,
%embora v\'arios dos resultados se encontrem tamb\'em no livro de Chihara.
%A rela\c c\~ao de equival\^encia (\ref{u-equiv-v}), com a formaliza\c c\~ao apresentada,
%\'e usualmente atribu\'\i da a Maroni. V\'arios detalhes apresentados encontram-se
%tamb\'em na tese de mestrado \ref{AFLoureiro2000-C6} de Ana Loureiro
%(realizada sob a orienta\c c\~ao conjunta de Pascal Maroni e de Z\'elia da Rocha).
%
%Refira-se, por \'ultimo, que os resultados contidos neste texto s\~ao elementares,
%pelo que poderiam, na sua maioria, ser deixados como exerc\'\i cios ou
%dilu\'\i dos noutros textos.
%A op\c c\~ao pela sua inclus\~ao num texto aut\'onomo deve-se \`a vantagem
%que resulta da sua apresenta\c c\~ao sistematizada
%para a leitura de alguns textos subsequentes.
%\medskip

\section*{Bibliography}
%\medskip

{\small
\begin{enumerate}[label=\emph{\rm [\arabic*]},leftmargin=*]
\item\label{Chihara1978-C6} T. S. Chihara, {\sl An introduction to orthogonal polynomials}, Gordon and Breach (1978).
\item\label{AFLoureiro2000-C6} A. Loureiro, {\sl Uma nova caracteriza\c c\~ao dos polin\'omios ortogonais cl\'assicos}, MSc thesis, Univ. Porto (2003).
\item\label{Maroni1991-C6} P. Maroni, {\it Une th\'eorie alg\'ebrique des polyn\^omes orthogonaux. Applications aux polyn\^omes
orthogonaux semiclassiques}, In C. Brezinski et al. Eds., Orthogonal Polynomials and Their Applications, Proc. Erice 1990, IMACS, Ann. Comp. App. Math. {\bf 9} (1991) 95--130.
\item\label{Maroni1994-C6} P. Maroni, {\sl Fonctions eul\'eriennes. Polyn\^omes orthogonaux classiques},
T\'echniques de l'Ing\'enieur, trait\'e G\'en\'eralit\'es (Sciences Fondamentales), A {\bf 154} (1994) 1--30.
\end{enumerate}
}

\chapter{Pearson's distributional differential equation}

%\chapter{Orthogonal polynomials and dual basis.}
%More operations in $\mathcal{P}$ and $\mathcal{P}'$. }

\pagestyle{myheadings}\markright{Pearson's distributional differential equation}
\pagestyle{myheadings}\markleft{J. Petronilho}

In this text we start our study of the so called \emph{classical orthogonal polynomials},
which includes four families of OP: Hermite, Laguerre, Jacobi
(including as special cases the Legendre and Chebyshev polynomials), and Bessel polynomials.
Those constitute perhaps the most important class of OP.
The regular functional with respect to which each one of these families is an OPS
satisfies an homogeneous linear distributional differential equation of the first order,
called \emph{Pearson's equation}---see equation (\ref{Pearson-DEq}) in bellow.
Our purpose, here, is the analysis of the solutions ${\bf u}\in\mathcal{P}'$ of this equation.

\section{Pearson's distributional differential equation}

The Pearson's distributional differential equation has the form
%Our purpose in this section is the study
%of the solutions ${\bf u}\in\mathcal{P}'$ of the distributional differential equation
\begin{snugshade}\vspace*{-0.5em}
\begin{equation}\label{Pearson-DEq}
D(\phi{\bf u})=\psi{\bf u}\;,
\end{equation}
\end{snugshade}
\noindent
where $\phi\in\mathcal{P}_2$ and $\psi\in\mathcal{P}_1$,
and ${\bf u}\in\mathcal{P}'$ is the unknown.
Notice that we do not require \emph{a priori} ${\bf u}$ to be a regular functional.
%Since equation (\ref{Pearson-DEq}) is trivially satisfied if
%$\phi\equiv\psi\equiv0$, we will exclude this trivial situation.
%Thus, we may write
We may write
\begin{snugshade}\vspace*{-0.5em}
\begin{equation}\label{Pearson-PhiPsi}
\phi(x)= ax^2+bx+c\;,\quad \psi(x)= px+q\;,
\end{equation}
\end{snugshade}
\noindent
being $a,b,c,p,q\in\mathbb{C}$.
%\begin{snugshade}\vspace*{-0.5em}
%\begin{equation}\label{abcpq1}
%|a|+|b|+|c|+|p|+|q|\neq0\;.
%\end{equation}
%\end{snugshade}
%\noindent
We also define, for each integer or rational number $n$,
\begin{snugshade}\vspace*{-0.5em}
\begin{equation}\label{psi-n}
\psi_n:=\psi+n\phi^{\prime}\;, \quad
d_n:=\psi_{n/2}^{\prime}=na+p\;, \quad e_n:=\psi_n(0)=nb+q\;.
\end{equation}
\end{snugshade}
\noindent
%Notice that $\psi_n\in\mathcal{P}_1$, and $\psi_n$ is well defined even if $n$ is not an integer number.
Notice that $\psi_n(x)=d_{2n}x+e_n\in\mathcal{P}_1$.
Finally, for each ${\bf u}\in\mathcal{P}'$ and each $n\in\mathbb{N}_0$, we set
\begin{snugshade}\vspace*{-0.5em}
\begin{equation}\label{Phi-n-u}
{\bf u}^{[n]}:=\phi^n{\bf u}\;.
\end{equation}
\end{snugshade}

We begin with the following elementary result.

\begin{snugshade}
\begin{lemma}\label{Pearson-lemma1} %{lemarec}
Let ${\bf u}\in\mathcal{P}'$.
Then ${\bf u}$ satisfies the Pearson distributional differential equation $(\ref{Pearson-DEq})$
if and only if the corresponding sequence of moments, $u_n:=\langle{\bf u},x^n\rangle$,
satisfies the second order linear difference equation %three-term recurrence relation
\begin{equation}\label{le1a}
d_nu_{n+1}+e_nu_n+n\phi(0)u_{n-1}=0\;, \quad n\in\mathbb{N}_0\; .
\end{equation}
Moreover, if ${\bf u}$ satisfies $(\ref{Pearson-DEq})$, then ${\bf u}^{[n]}$ satisfies
\begin{equation}\label{le1b}
D\big(\phi{\bf u}^{[n]}\big)=\psi_n{\bf u}^{[n]}\;, \quad n\in\mathbb{N}_0\;.
\end{equation}
\end{lemma}
\end{snugshade}

{\it Proof.}
We have
$$
\begin{array}{rcl}
D(\phi{\bf u})=\psi{\bf u} &\Leftrightarrow&
\langle D(\phi{\bf u}),x^n\rangle=\langle\psi{\bf u},x^n\rangle\;,\quad\forall n\in\mathbb{N}_0 \\ [0.25em]
 &\Leftrightarrow& -n\langle {\bf u},\phi x^{n-1}\rangle=\langle{\bf u},\psi x^n\rangle\;,\quad\forall n\in\mathbb{N}_0\;, \\ [0.25em]
 &\Leftrightarrow& (na+p)u_{n+1}+(nb+q)u_n+ncu_{n-1}=0\;,\quad\forall n\in\mathbb{N}_0\;,
\end{array}
$$
hence the first assertion of the theorem is proved.
To prove (\ref{le1b}) we use mathematical induction.
Since ${\bf u}^{[0]}:={\bf u}$ and $\psi_0=\psi$, then (\ref{le1b}) holds for $n=0$.
Assume that (\ref{le1b}) holds for a certain $n\in\mathbb{N}_0$.
Then $D\big({\bf u}^{[n+1]}\big)=D\big(\phi{\bf u}^{[n]}\big)=\psi_n{\bf u}^{[n]}$, hence
$$
\begin{array}{rcl}
D\big(\phi{\bf u}^{[n+1]}\big)&=&\phi'{\bf u}^{[n+1]}+\phi D\big({\bf u}^{[n+1]}\big)
=\phi'{\bf u}^{[n+1]}+\phi\psi_n{\bf u}^{[n]}
=\big(\phi'+\psi_n\big){\bf u}^{[n+1]}\\ [0.25em]
&=&\psi_{n+1}{\bf u}^{[n+1]}\,.
\end{array}
$$
Thus (\ref{le1b}) is proved.
\qed
\medskip

\begin{snugshade}
\begin{theorem}\label{Pearson-Thm2} %\label{lem2}
Let ${\bf u}\in\mathcal{P}'$, and suppose that
${\bf u}$ satisfies the Pearson distributional differential equation $(\ref{Pearson-DEq})$.
Suppose further that
\begin{equation}\label{dnd0}
d_n\neq0\; , \quad \forall n\in\mathbb{N}_0\; .
\end{equation}
Then, there exists a simple set of polynomials $\{R_n\}_{n\geq0}$ such that
\begin{equation} \label{rodrigues}
R_n{\bf u}=D^n(\phi^n{\bf u})\; , \quad n\in\mathbb{N}_0 .
\end{equation}
Moreover, $\{R_n\}_{n\geq0}$ may be chosen so that it satisfies the three-term recurrence relation
\begin{equation}\label{ttrr1}
R_{n+1}(x)=(\tilde{\alpha}_n x-\tilde{\beta}_n)R_n(x)-\tilde{\gamma}_nR_{n-1}(x)\; , \quad n\in\mathbb{N}_0\;,
\end{equation}
with initial conditions $R_{-1}(x)=0$ and $R_0(x)=1$, being
\begin{equation}\label{coef}
\begin{array}{c}
\displaystyle{ \tilde{\alpha}_n:=\frac{d_{2n-1}d_{2n}}{d_{n-1}}}\; , \quad
\tilde{\beta}_n:=-\tilde{\alpha}_n \frac{d_{-2}q+2bnd_{n-1}}{d_{2n-2}d_{2n}}\; , \quad  n\in\mathbb{N}_0 \; , \\ [1.25em]
\displaystyle{ \tilde{\gamma}_n:=-\tilde{\alpha}_n \frac{n d_{2n-2}}{d_{2n-1}} \, \phi
\left( -\frac{e_{n-1}}{d_{2n-2}}\right) }\; , \quad n\in\mathbb{N} \; .
\end{array}
\end{equation}
\end{theorem}
\end{snugshade}

{\it Proof.}
The proof will be made by mathematical induction on $n$.
Defining
$$R_0:=1 \; , \quad R_1:=\psi \; ,$$
it is clear that (\ref{rodrigues}) holds for $n=0$ and $n=1$.
Suppose now that (\ref{rodrigues}) holds for the indices $n$ and $n-1$, that is,
there exist polynomials $R_n$ and $R_{n-1}$, with degrees $n$ and $n-1$, respectively, such that
\begin{equation}\label{rod2}
R_n{\bf u}=D^n\big(\phi^n{\bf u}\big)\;,\quad
R_{n-1}{\bf u}=D^{n-1}\big(\phi^{n-1}{\bf u}\big)\; .
\end{equation}
%and
%\begin{equation}
%R_{n-1}{\bf u}=D^{n-1}\big(\phi^{n-1}{\bf u}\big)\; .
%\label{rod3}
%\end{equation}
We must show that there is a polynomial $R_{n+1}$ of degree $n+1$ such that
\begin{equation}
R_{n+1}{\bf u}=D^{n+1}\big(\phi^{n+1}{\bf u}\big) \quad .
\label{rod4}
\end{equation}
Using Lemma \ref{Pearson-lemma1}, we have
$$%\begin{equation}
\begin{array}{rcl}
D^{n+1}\big(\phi^{n+1}{\bf u}\big)&=&D^n\big[D\big(\phi{\bf u}^{[n]}\big)\big]
= D^n\big(\psi_n{\bf u}^{[n]}\big)= D^n\big(\psi_n\phi^n{\bf u}\big) \\ [0.5em]
&=& \binom{n}{0}\psi_n D^n(\phi^n{\bf u})+\binom{n}{1}\psi_n^{\prime}D^{n-1}(\phi^n{\bf u})
=\psi_nR_n{\bf u}+n\psi_n^{\prime}D^{n-1}(\phi^n{\bf u})\; ,
\end{array}
$$
where in the fourth equality we applied the distributional Leibniz rule %for the derivative of a product,
and in the last one we used the first relation in the induction hypothesis (\ref{rod2}).
Therefore,
\begin{equation}
D^{n-1}(\phi^{n}{\bf u})=\frac{1}{nd_{2n}}\big(D^{n+1}(\phi^{n+1}{\bf u})-\psi_nR_n{\bf u} \big)\; .
\label{rod6}
\end{equation}
Notice that, according to the hypothesis (\ref{dnd0}), $\psi_{m/2}^{\prime}=d_{m}\neq0$ for each $m\in\mathbb{N}_0$.
We point out that we have deduced (\ref{rod6}) using the first relation in (\ref{rod2}).
Therefore, making the change of indices $n\rightarrow n-1$ in the above reasoning and using the
second relation in (\ref{rod2}), we obtain
\begin{equation}
D^{n-2}(\phi^{n-1}{\bf u}) = \frac{1}{(n-1)d_{2n-2}}(R_n-\psi_{n-1}R_{n-1}){\bf u}\; .
\label{rod7}
\end{equation}
On the other hand, using again Lemma \ref{Pearson-lemma1}, we have
$$
\begin{array}{rcl}
D^{n+1}(\phi^{n+1}{\bf u}) &=& D^{n-1}\big[D(\psi_n\phi^n{\bf u})\big]
%=D^{n-1}\big[\psi_n^{\prime}\phi^n{\bf u}+\psi_nD\big(\phi{\bf u}^{[n-1]}\big)\big] \\ [0.5em]
= D^{n-1}\big[(\psi_n^{\prime}\phi+\psi_n\psi_{n-1})\phi^{n-1}{\bf u}\big] \\ [0.5em]
&=& (\psi_n^{\prime}\phi+\psi_n\psi_{n-1})D^{n-1}(\phi^{n-1}{\bf u})+ \binom{n-1}{1}(\psi_n^{\prime}\phi+\psi_n\psi_{n-1})^{\prime}D^{n-2}(\phi^{n-1}{\bf u}) \\ [0.5em]
&&\quad +\binom{n-1}{2}(\psi_n^{\prime}\phi+\psi_n\psi_{n-1})^{\prime\prime}D^{n-3}(\phi^{n-1}{\bf u})
\end{array}
$$
where in the last equality we have applied again Leibniz rule.
Consequently, using (\ref{rod7}) and the second relation in (\ref{rod2}),
and taking into account the identities
$$
(\psi_n^{\prime}\phi+\psi_n\psi_{n-1})^{\prime}=2d_{2n-1}\psi_n \; , \quad
(\psi_n^{\prime}\phi+\psi_n\psi_{n-1})^{\prime\prime}=2d_{2n}d_{2n-1} \; ,
%(\psi_n^{\prime}\phi+\psi_n\psi_{n-1})^{\prime}=2\psi_{n-\frac{1}{2}}^{\prime}\psi_n \; , \quad
%(\psi_n^{\prime}\phi+\psi_n\psi_{n-1})^{\prime\prime}=2\psi_n^{\prime}\psi_{n-\frac{1}{2}}^{\prime} \; ,
$$
we deduce
%\begin{equation}
%\begin{array}{l}
%\frac{(n-1)(n-2)}{2}D^{n-3}\big(\phi^{n-1}{\bf u}\big)
%\displaystyle=\frac{1}{2\psi_{n-\frac{1}{2}}^{\prime}\psi_n^{\prime}} \\ [0.5em]
%\quad\displaystyle
%\times\Big( D^{n+1}(\phi^{n+1}{\bf u})- (\psi_n^{\prime}\phi+\psi_n\psi_{n-1})R_{n-1}{\bf u}
%-\frac{2\psi_{n-\frac{1}{2}}^{\prime}}{\psi_{n-1}^{\prime}} \psi_n(R_n-\psi_{n-1}R_{n-1}){\bf u} \, \Big)\, .
%\end{array}
%\label{rod12}
%\end{equation}
\begin{equation}
\begin{array}{r}
\frac{(n-1)(n-2)}{2}D^{n-3}\big(\phi^{n-1}{\bf u}\big)
=\displaystyle\frac{1}{2d_{2n}d_{2n-1}}
\Big\{ D^{n+1}(\phi^{n+1}{\bf u})- (\psi_n^{\prime}\phi+\psi_n\psi_{n-1})R_{n-1}{\bf u} \\ [0.5em]
\displaystyle-\frac{2d_{2n-1}}{d_{2n-2}} \psi_n(R_n-\psi_{n-1}R_{n-1}){\bf u} \, \Big\}\, .
\end{array}
\label{rod12}
\end{equation}
Now, consider the left-hand side of (\ref{rod6}). Using Leibniz rule, we have
%$$
%\begin{array}{rl}
%D^{n-1}(\phi^{n}{\bf u}) &= D^{n-1}\big(\phi(\phi^{n-1}{\bf u})\big) \\ [0.25em]
%&=\phi D^{n-1}(\phi^{n-1}{\bf u}) +(n-1)\phi^{\prime}D^{n-2}(\phi^{n-1}{\bf u})+
%\frac{(n-1)(n-2)}{2}\phi^{\prime\prime}D^{n-3}(\phi^{n-1}{\bf u})\; ;
%\end{array}
%$$
$$
D^{n-1}(\phi^{n}{\bf u})
=\phi D^{n-1}(\phi^{n-1}{\bf u}) +(n-1)\phi^{\prime}D^{n-2}(\phi^{n-1}{\bf u})+
(n-1)(n-2)aD^{n-3}(\phi^{n-1}{\bf u})\,.
$$
Substituting into (\ref{rod6}), and using (\ref{rod7}) and the second relation
in (\ref{rod2}), we deduce
\begin{equation}\label{rod13}
\begin{array}{l}
\displaystyle\frac{1}{nd_{2n}}\big( D^{n+1}(\phi^{n+1}{\bf u})-\psi_nR_n{\bf u} \big) \\ [0.5em]
\qquad =\displaystyle\phi R_{n-1}{\bf u} +\frac{\phi^{\prime}}{d_{2n-2}}(R_n-\psi_{n-1}R_{n-1}){\bf u} +a(n-1)(n-2)D^{n-3}(\phi^{n-1}{\bf u}) \; .
\end{array}
\end{equation}
Finally, substituting (\ref{rod12}) in the right-hand side of (\ref{rod13}),
after some computations (we may use {\scriptsize MATHEMATICA}!) we obtain (\ref{rod4}), provided we define
$$
R_{n+1}(x):=\frac{d_{2n-1}d_{2n}}{d_{n-1}} \left[ \left( x+\frac{d_{-2}q+2bnd_{n-1}}{d_{2n-2}d_{2n}} \right) R_n(x)
+\frac{n d_{2n-2}}{d_{2n-1}} \phi\left( -\frac{e_{n-1}}{d_{2n-2}}\right)R_{n-1}(x)\right] .
$$
Since (by the induction hypothesis) $R_n$ and $R_{n-1}$ have degrees $n$ and $n-1$ (respectively),
it follows that $R_{n+1}$ is a polynomial of degree $n+1$. %which fulfills (\ref{rod4}).
Thus the theorem is proved.
\qed
\medskip

%\section{Classical OPS: definition and regularity conditions}
\section{The regular solutions of Pearson's equation}

In the previous section we analyzed Pearson's distributional differential
equation (\ref{Pearson-DEq}) without requiring the regularity condition
on the functional ${\bf u}\in\mathcal{P}'$.
In this section we determine necessary and sufficient conditions, involving only the (coefficients of the)
polynomials $\phi$ and $\psi$, which ensure the regularity of such a functional.
%The OPS with respect to the regular functionals satisfying Pearson's equation
%are the so called \emph{classical orthogonal polynomials}.
%Those form a class of OP that, up to constant factors and affine changes of the variable,
%can be described by four specific families of orthogonal polynomials:
%\emph{Hermite}, \emph{Laguerre}, \emph{Jacobi} and \emph{Bessel} polynomials.

Notice that if both $\phi$ and $\psi$ vanish identically %coincide with the zero polynomial,
then Pearson's equation (\ref{Pearson-DEq}) reduces to a trivial equation,
so we will exclude this situation from our study.
%Remarkably, if at least one of the polynomials $\phi$ and $\psi$ is nonzero
%much more can be said, provided ${\bf u}$ is regular.

\begin{snugshade}
\begin{lemma}\label{poly-grau21}
Let ${\bf u}\in\mathcal{P}'$. Suppose that ${\bf u}$ is regular and satisfies
Pearson's equation $(\ref{Pearson-DEq})$, being
$\phi\in\mathcal{P}_2$ and $\psi\in\mathcal{P}_1$,
and assume that at least one of the polynomials $\phi$ and $\psi$ is nonzero.
Then neither $\phi$ nor $\psi$ is the zero polynomial, and
\begin{equation}\label{grauPsi1}
\deg\psi=1\;.
\end{equation}
\end{lemma}
\end{snugshade}

{\it Proof.}
Since ${\bf u}$ is regular, there is a monic OPS, $\{P_n\}_{n\geq0}$, with respect to ${\bf u}$.
Since ${\bf u}$ fulfills (\ref{Pearson-DEq}),
if $\psi\equiv0$ then $D\left(\phi\textbf{u}\right)=\textbf{0}$;
hence, if $\phi\not\equiv0$, setting $r:=\deg\phi$ and denoting by $k(\neq0)$ the leading coefficient of $\phi$,
we would have
$$\langle \textbf{u},P_r^2\rangle=k^{-1}\langle \phi\textbf{u},P_r\rangle
=-k^{-1}\big\langle D\left(\phi\textbf{u}\right),\int P_r \big\rangle=-k^{-1}\big\langle\textbf{0},\int P_r\big\rangle=0\;,$$
violating the regularity of ${\bf u}$. We conclude that $\psi\equiv0$ implies $\phi\equiv0$.
Suppose now that $\phi\equiv0$. Then, $\psi{\bf u}=0$. If $\psi\not\equiv0$, setting $t:=\deg\psi$ and being $m(\neq0)$ the leading coefficient of $\psi$, we would have
$\langle \textbf{u},P_t^2\rangle=m^{-1}\langle \psi\textbf{u},P_t\rangle=0$,
violating again the regularity of ${\bf u}$.
We conclude that $\phi\equiv0$ implies $\psi\equiv0$.
Finally, suppose that $\psi\equiv\mbox{\rm constant}=q\neq0$. Then
$\langle \textbf{u},1\rangle=q^{-1}\langle \psi\textbf{u},1\rangle
=q^{-1}\langle D(\phi\textbf{u}),1\rangle=0\,$,
violating once again the regularity of ${\bf u}$.
Thus $\deg\psi=1$.
\qed
\smallskip

Given a monic polynomial $P_n$ of degree $n$ (which needs not to belong to an OPS),
we denote by $P_n^{[k]}$ the monic polynomial of degree $n$ defined by
\begin{snugshade}\vspace*{-0.5em}
\begin{equation}\label{Pnk-deriv}
P_n^{[k]}(x):=\frac{{\rm d}^k}{{\rm d}x^k}\,\frac{P_{n+k}(x)}{(n+1)_k}\quad(k,n\in\mathbb{N}_0)\;,
\end{equation}
\end{snugshade}
\noindent
where, for a given $\alpha\in\mathbb{C}$, $(\alpha)_n$ is the {\it Pochhammer symbol}, defined as
\begin{snugshade}\vspace*{-0.5em}
\begin{equation}\label{Pochhammer}
(\alpha)_0:=1\;,\qquad (\alpha)_n:=\alpha(\alpha+1)\cdots(\alpha+n-1)\;,\quad n\in\mathbb{N}\;.
\end{equation}
\end{snugshade}
\noindent
Clearly, if $\{P_n\}_{n\geq0}$ is a simple set in $\mathcal{P}$, then so is $\{P_n^{[k]}\}_{n\geq0}$.
Under such conditions, there is a beautiful relation between the associated dual basis.

\begin{snugshade}
\begin{lemma}\label{teo-rel-dual-an-ank}
Let $\{P_n\}_{n\geq0}$ be a simple set in $\mathcal{P}$,
and let $\{{\bf a}_n\}_{n\geq0}$ and $\big\{{\bf a}_n^{[k]}\,\big\}_{n\geq0}$
be the dual basis in $\mathcal{P}'$ associated with $\{P_n\}_{n\geq0}$ and $\{P_n^{[k]}\}_{n\geq0}$,
respectively. Then,
\begin{equation}\label{rel-dual-an-ank}
D^{k}\big(\,{\bf a}_n^{[k]}\,\big)=(-1)^k(n+1)_k\,{\bf a}_{n+k} \;,\quad k,n\in\mathbb{N}_0 \, .
\end{equation}
\end{lemma}
\end{snugshade}

{\it Proof.}
Fix $j,k,n\in\mathbb{N}_0$. Then,
$$
\begin{array}{rcl}
\big\langle\,D^{k}\big(\,{\bf a}_n^{[k]}\,\big),P_j\,\big\rangle & = &
(-1)^k\big\langle\, {\bf a}_n^{[k]},\frac{{\rm d}^k}{{\rm d}x^k}\,P_j\,\big\rangle \, = \,
(-1)^k (j-k+1)_k\,\big\langle\, {\bf a}_n^{[k]},P_{j-k}^{[k]}\,\big\rangle \\ [0.5em]
&=& (-1)^k(n+1)_k\,\delta_{n,j-k}\,=\,
\big\langle \, (-1)^k(n+1)_k\, {\bf a}_{n+k}, P_j\, \big\rangle \; .
\end{array}
$$
Therefore, (\ref{rel-dual-an-ank}) is proved.
\qed

\begin{snugshade}
\begin{lemma}\label{reg-lemma3}
Let ${\bf u}\in\mathcal{P}'$, and suppose that
${\bf u}$ satisfies the Pearson distributional differential equation $(\ref{Pearson-DEq})$,
with $\phi$ and $\psi$ given by $(\ref{Pearson-PhiPsi})$,
being at least one of these polynomials nonzero.
Suppose further that ${\bf u}$ is regular. Then
\begin{equation}\label{r1}
d_n:=na+p\neq0\; , \quad \forall n\in\mathbb{N}_0\; .
\end{equation}
Moreover, let $\{P_n\}_{n\geq0}$ be the monic OPS with respect to ${\bf u}$, and
$P_n^{[k]}$ defined by $(\ref{Pnk-deriv})$.
Then, for each $k\in\mathbb{N}_0$, ${\bf u}^{[k]}:=\phi^k{\bf u}$ is regular and
$\{P_n^{[k]}\}_{n\geq0}$ is its monic OPS.
\end{lemma}
\end{snugshade}

{\it Proof.}
By Lemma \ref{poly-grau21}, both $\phi$ and $\psi$ are nonzero, and $\deg\psi=1$, hence $p\neq0$.
%Thus $p\neq0$, and so we see immediately that (\ref{r1}) holds if $a=0$.
%(i.e., whenever $\deg\phi\leq1$).
Consider first the case $k=1$.
Write $Q_n=P_n^{[1]}:=P_{n+1}'/(n+1)$.
We will show that
\begin{equation}\label{Nu}
\langle \phi {\bf u} , Q_n Q_m \rangle =-\frac{d_n}{n+1}
\langle {\bf u} , P_{n+1}^2 \rangle \delta_{nm} \; , \quad n,m\in\mathbb{N}_0 \; .
\end{equation}
Indeed, since $D(\phi{\bf u})=\psi{\bf u}$, we may write
$$
\begin{array}{rcl}
(n+1)\langle\phi{\bf u},Q_mQ_n\rangle &=& \langle\phi{\bf u},Q_mP_{n+1}'\rangle
=\langle\phi{\bf u},\big(Q_mP_{n+1}\big)'-Q_m'P_{n+1}\rangle \\ [0.25em]
&=& -\langle D(\phi{\bf u}),Q_mP_{n+1}\rangle-\langle \phi {\bf u},Q_m'P_{n+1}\rangle\\ [0.25em]
&=& -\langle {\bf u},(\psi Q_m+\phi Q_m')P_{n+1}\rangle\;.
\end{array}
$$
Thus, assuming (without loss of generality) $m\leq n$,
since $\deg(\psi Q_m+\phi Q_m')<n+1$ if $m<n$,
%and $\psi Q_m+\phi Q_m'=(na+p)x^{n+1}+[\mbox{\rm lower degree terms}]$ if $m=n$,
and $\psi Q_m+\phi Q_m'=(na+p)x^{n+1}+\pi_n(x)$ if $m=n$, where $\pi_n\in\mathcal{P}_n$,
we obtain (\ref{Nu}).
Let $s:=\deg\phi\in\{0,1,2\}$. For each $n\in\mathbb{N}_0$, write
$P_{n+s}(x)=\sum_{m=0}^{n+s}a_{nm}Q_m(x)$, being $a_{nm}$ complex parameters.
Multiplying both sides of this equality by $\phi Q_n$ and then applying  ${\bf u}$,
and taking into account (\ref{Nu}), we deduce
$$
\langle {\bf u} , \phi Q_n P_{n+s} \rangle
=-\frac{a_{nn}d_n}{n+1} \langle {\bf u} , P_{n+1}^2 \rangle  \; , \quad n\in\mathbb{N}_0 \; .
$$
The left-hand side of this equality never vanishes, since $\{ P_n\}_{n\geq0}$ is an OPS
with respect to ${\bf u}$ and $\deg\phi =s$ (and $\phi$ is not the zero polynomial).
Thus the right-hand side of the equality cannot vanish,
hence $d_n\neq0$ (and also $a_{nn}\neq0$), which proves (\ref{r1}).

It remains to prove that $\{P_n^{[k]}\}_{n\geq0}$ is a monic OPS with respect to
${\bf u}^{[k]}:=\phi^k{\bf u}$, for each $k\in\mathbb{N}$.
If $k=1$, since, by (\ref{Nu}),
\begin{equation}\label{NuP}
\langle {\bf u}^{[1]} , P_n^{[1]} P_m^{[1]} \rangle =-\frac{na+p}{n+1}
\langle {\bf u} , P_{n+1}^2 \rangle \delta_{nm} \; , \quad n,m\in\mathbb{N}_0 \; ,
\end{equation}
then (\ref{r1}) ensures that $\{P_n^{[1]}\}_{n\geq0}$ is a monic OPS with respect to
${\bf u}^{[1]}:=\phi{\bf u}$.
Now, by Lemma \ref{Pearson-lemma1}, ${\bf u}^{[1]}$ fulfills the Pearson's equation
$D\big(\phi{\bf u}^{[1]}\big)=\psi_1{\bf u}^{[1]}$, hence, since
%${\bf u}^{[2]}:=\phi{\bf u}^{[1]}$, and
$P_n^{[2]}=\big(P_{n+1}^{[1]}\big)'/(n+1)$ and $\psi_1(x)=(2a+p)x+b+q$,
from (\ref{NuP}) with ${\bf u}$, $\psi$, and $\{P_n\}_{n\geq0}$ replaced by ${\bf u}^{[1]}$, $\psi_1$, and
$\{P_n^{[1]}\}_{n\geq0}$ (resp.), we deduce,for every $n,m\in\mathbb{N}_0$,
$$%\begin{equation}\label{NuPaa}
\langle {\bf u}^{[2]} , P_n^{[2]} P_m^{[2]} \rangle
=-\frac{(n+2)a+p}{n+1} \langle {\bf u}^{[1]} , \big(P_{n+1}^{[1]}\big)^2 \rangle \delta_{nm}
=\frac{d_{n+1}d_{n+2}}{(n+1)(n+2)} \langle {\bf u} , P_{n+2}^2 \rangle \delta_{nm} \; ,
%\quad n,m\in\mathbb{N}_0 \; ,
$$%\end{equation}
and so $\{P_n^{[2]}\}_{n\geq0}$ is a monic OPS with respect to ${\bf u}^{[2]}$.
Arguing by induction, we prove %are able to state the orthogonality relations
%the orthogonality relations for the monic OPS $\{P_n^{[k]}\}_{n\geq0}$ are explicitly given by
%the orthogonality relations for the monic OPS $\{P_n^{[k]}\}_{n\geq0}$ read as:
\vspace*{-0.5em}
\begin{snugshade}
\begin{equation}\label{NuPaa}
\langle {\bf u}^{[k]} , P_n^{[k]} P_m^{[k]} \rangle
=(-1)^k\frac{\prod_{j=0}^{k-1}d_{n+k+j-1}}{(n+1)_k}\langle{\bf u},P_{n+k}^2\rangle\delta_{nm}
\quad (k,n,m\in\mathbb{N}_0)\;,
\end{equation}
\end{snugshade}
\noindent
hence $\{P_n^{[k]}\}_{n\geq0}$ is a monic OPS with respect to
${\bf u}^{[k]}$, for each $k\in\mathbb{N}_0$.
%the relation
%\begin{equation}\label{NuPaa}
%\langle {\bf u}^{[k]} , P_n^{[k]} P_m^{[k]} \rangle
%=(-1)^k\frac{\prod_{j=k-1}^{2k-2}d_{n+j}}{(n+1)_k}\langle{\bf u},P_{n+k}^2\rangle\delta_{nm}%\;,
%\quad (n,m\in\mathbb{N}_0)
%\end{equation}
%for each $k\in\mathbb{N}_0$. Thus, $\{P_n^{[k]}\}_{n\geq0}$ is a monic OPS with respect to ${\bf u}^{[k]}$.
\qed
\smallskip

We may now establish necessary and sufficient conditions ensuring
the regularity of a given functional ${\bf u}\in\mathcal{P}'$ satisfying (\ref{Pearson-PhiPsi}).
%Pearson's distributional differential equation (\ref{Pearson-PhiPsi}).

\begin{snugshade}
\begin{theorem}\label{vCOP-Thm1}
Let ${\bf u}\in\mathcal{P}^\prime\setminus\{{\bf 0}\}$,
and suppose that ${\bf u}$ satisfies the Pearson equation
\begin{equation}\label{EqDistC1}
D(\phi{\bf u})=\psi{\bf u}\;,
\end{equation}
%where $\phi$ and $\psi$ are polynomials such that $\deg\,\phi\leq2$ and $\deg\,\psi=1$.
where $\phi$ and $\psi$ are nonzero polynomials such that
$\phi\in\mathcal{P}_2$ and $\psi\in\mathcal{P}_1$. Set
$$
\phi(x):=ax^2+bx+c\,,\quad\psi(x):=px+q\,,\quad
d_n:=na+p\,,\quad e_n:=nb+q \quad(n\in\mathbb{N}_0)\,.
$$
Then, ${\bf u}$ is regular if and only if %the conditions
\begin{equation}\label{EqDistC2}
d_n\neq0\,,\quad \phi\Big(-\frac{e_n}{d_{2n}}\Big)\neq0\,,\quad
\forall n\in\mathbb{N}_0\;. %n=0,1,2,\ldots
\end{equation}
Moreover, under these conditions, the monic OPS $\{P_n\}_{n\geq0}$
with respect to ${\bf u}$ is given by the three-term recurrence relation
\begin{equation}\label{ttrrC1}
P_{n+1}(x)=(x-\beta_n)P_n(x)-\gamma_nP_{n-1}(x)\;,\quad  n\in\mathbb{N}_0 %n=0,1,2,\cdots
\end{equation}
with initial conditions $P_{-1}(x)=0$ and $P_0(x)=1$, being
\begin{equation}\label{EqDistC3}
\beta_n=\frac{ne_n}{d_{2n}}-\frac{(n+1)e_{n+1}}{d_{2n+2}}\, ,\quad
\gamma_{n+1}=-\frac{(n+1)d_{n-1}}{d_{2n-1}d_{2n+1}}\phi\Big(-\frac{e_n}{d_{2n}}\Big)
\;,\quad  n\in\mathbb{N}_0\;.
\end{equation}
In addition, for each $n\in\mathbb{N}_0$, $P_n$ satisfies the distributional Rodrigues formula
\begin{equation}\label{EqDistRod}
P_n{\bf u}=k_n\,D^n\big(\phi^n{\bf u}\big)\;,\quad
k_n:=\prod_{i=0}^{n-1}d_{n+i-1}^{-1}\,.
\end{equation}

\end{theorem}
\end{snugshade}

{\it Proof.}
Suppose that ${\bf u}$ is regular.
Let $\{P_n\}_{n\geq0}$ be the monic OPS with respect to ${\bf u}$.
By Lemma \ref{reg-lemma3}, $d_n\neq0$ for each $n\in\mathbb{N}_0$. Moreover,
%, hence the first condition appearing in (\ref{EqDistC2}) holds.
%by Lemma \ref{reg-lemma3},
${\bf u}^{[k]}:=\phi^k{\bf u}$ is regular and
$\{P_n^{[k]}\}_{n\geq0}$ is its monic OPS, for each $k\in\mathbb{N}_0$.
%Let $\{ \beta_{n}^{[k]} \}_{n\geq0}$ and $\{ \gamma_{n}^{[k]} \}_{n\geq1}$
%be the sequences of parameters appearing in the three-term recurrence relation
%fulfilled by $\{P_{n}^{[k]}\}_{n\geq0}$.
This monic OPS satisfies a TTRR
\begin{equation}\label{ttrrC1}
P_{n+1}^{[k]}(x)=(x-\beta_n^{[k]})P_n^{[k]}(x)-\gamma_n^{[k]}P_{n-1}^{[k]}(x)\;,\quad  n\in\mathbb{N}_0
\end{equation}
with initial conditions $P_{-1}^{[k]}(x)=0$ and $P_0^{[k]}(x)=1$, being
$\beta_{n}^{[k]}\in\mathbb{C}$ and $\gamma_{n}^{[k]}\in\mathbb{C}\setminus\{0\}$ for each $n$.
Let us compute $\gamma_{1}^{[n]}$ for each fixed $n\in\mathbb{N}_0$.
We first show that the coefficient $\gamma_1\equiv\gamma_{1}^{[0]}$,
appearing in the TTRR for $\{P_n\}_{n\geq0}$,
%the $\gamma$-parameters of a MOPS corresponding to a functional satisfying
%an equation like (\ref{q1})
is given by
\begin{equation}
\gamma_1=-\frac{1}{p+a}\phi\left(-\frac{q}{p}\right) \; .
\label{g1}
\end{equation}
To prove this relation, take $n=0$ and $n=1$ in the recurrence relation (\ref{le1a})
%remark that equation (\ref{q1}) implies that
for the sequence of moments associated to ${\bf u}$.
This gives
\begin{equation}
u_1=-\frac{q}{p}u_0 \; , \quad u_2=-\frac{1}{p+a}\left[-(b+q)\frac{q}{p}+c \right] u_0 \; .
\label{u1}
\end{equation}
%Since  $P_0(x)=1$ and $P_1(x)=x-\beta_0$ , it follows that $\beta_0=u_1/u_0$ and
%\begin{equation}
%\gamma_1=\frac{\langle {\bf u} , P_1^2 \rangle}{\langle {\bf u} , P_0^2 \rangle}=
%\frac{u_2-2\beta_0u_1+\beta_0^2u_0}{u_0} \; .
%\label{u2}
%\end{equation}
%From (\ref{u1}) and (\ref{u2}), it is easy to deduce (\ref{g1}).
%Now, since equation (\ref{le1b}) is of the same type as equation
%(\ref{EqDistC1}), with the same polynomial $\phi$ and $\psi$ replaced by $\psi_n$,
%we see that the expression of $\gamma_{1}^{[n]}$ may be obtained
%replacing the coefficients $p$ and $q$ of $\psi$ in (\ref{g1})
%by the corresponding coefficients of $\psi_n$. Hence,
On the other hand, by Corollary \ref{cor-TTRR},
%the $\gamma-$parameters may be computed using Hankel determinants, so that
\begin{equation}
\gamma_1=\frac{H_{-1}H_1}{H_0^2}=\frac{u_2u_0-u_1^2}{u_0^2} \; .
\label{u2}
\end{equation}
Substituting $u_1$ and $u_2$ given by (\ref{u1}) into (\ref{u2}) yields (\ref{g1}).
Now, since equation (\ref{le1b}) is of the same type as equation
(\ref{EqDistC1}), with the same polynomial $\phi$ and being $\psi$ replaced by $\psi_n$,
we see that the expression of $\gamma_{1}^{[n]}$ may be obtained
replacing the coefficients $p$ and $q$ of $\psi$ in (\ref{g1})
by the corresponding coefficients of $\psi_n$. Hence,
\begin{equation}
\gamma_{1}^{[n]}=-\frac{1}{d_{2n}+a}\phi\left(-\frac{e_n}{d_{2n}}\right)
=-\frac{1}{d_{2n+1}}\phi\left(-\frac{e_n}{d_{2n}}\right) \;.
\label{g1n}
\end{equation}
Since ${\bf u}^{[n]}$ is regular, then $\gamma_{1}^{[n]}\neq0$.
Thus, the second condition in (\ref{EqDistC2}) holds.
%\smallskip

Conversely, suppose that conditions (\ref{EqDistC2}) hold.
%To prove that ${\bf u}$ is regular we must show that there exists a sequence
%of polynomials $\{P_n)\}_{n\geq0}$ which is an OPS with respect to ${\bf u}$.
%We first remark that, a
According with Theorem \ref{Pearson-Thm2},
there is a simple set of polynomials $\{R_n\}_{n\geq0}$ such that
(\ref{rodrigues}) holds and satisfiyng the TTRR (\ref{ttrr1})--(\ref{coef}).
The hypothesis (\ref{EqDistC2}) ensure that $\tilde{\alpha}_n\neq0$
and $\tilde{\gamma}_n\neq0$ for each $n$.
Thus, by Favard's theorem, $\{R_n\}_{n\geq0}$ is an OPS.
We claim that $\{R_n\}_{n\geq0}$ is an OPS with respect to ${\bf u}$.
By Corollary \ref{OPSwrtv}, we only need to show that
\begin{equation}\label{qq}
u_0\neq0 \; ,\quad \langle{\bf u},R_n\rangle=0 \; , \quad n\geq1\;.
\end{equation}
%(because we already know that $\{R_n\}_{n\geq0}$ is an OPS).
In fact, if $u_0=0$, since (by Lemma \ref{Pearson-lemma1})
Pearson's equation (\ref{EqDistC1}) is equivalent to
the recurrence relation (\ref{le1a}) fulfilled by the moments $u_n$,
and since for $n=0$ (\ref{le1a}) yields $pu_1+qu_0=0$,
we would get $pu_1=0$ and so $u_1=0$ (because $p=d_0\neq0$);
therefore, $u_0=u_1=0$, and it follows recurrently from (\ref{le1a})
that $u_n=0$ for each $n\in\mathbb{N}_0$, hence ${\bf u}=\textbf{0}$,
in contradiction with the hypothesis. Thus, $u_0\neq0$.
On the other hand, by (\ref{rodrigues}), for each $n\geq1$ we have
$\langle{\bf u},R_n\rangle=\langle R_n{\bf u},1\rangle=(-1)^n\langle\phi^n{\bf u},0\rangle=0$.
Thus (\ref{qq}) is proved. Therefore $\{R_n\}_{n\geq0}$ is a monic OPS with respect to ${\bf u}$,
hence ${\bf u}$ is regular.

It remains to prove (\ref{EqDistC3})--(\ref{EqDistRod}).
Since $\{P_n\}_{n\geq0}$ and $\{R_n\}_{n\geq0}$ are both OPS with respect to ${\bf u}$,
then there exist a sequence $\{k_n\}_{n\geq0}$, with $k_n\in\mathbb{C}\setminus\{0\}$, such that
\begin{equation}\label{PnknRn}
P_n(x)=k_nR_n(x)\; , \quad n\in\mathbb{N}_0 \; .
\end{equation}
Multiplying both sides of (\ref{ttrr1}) by $k_n$ we obtain
$$
xP_n(x)=\frac{k_n}{\tilde{\alpha}_n k_{n+1}}P_{n+1}(x)+\frac{\tilde{\beta}_n}{\tilde{\alpha}_n}P_n(x)
+\frac{\tilde{\gamma}_n k_n}{\tilde{\alpha}_n k_{n-1}}P_{n-1}(x) \; , \quad n\geq1\; ,
$$
and  $P_1=x-\tilde{\beta}_0$ where $\tilde{\beta}_0=-q/p$.
Since each $P_n$ is a monic polynomial, we must have $k_n/(\tilde{\alpha}_n k_{n+1})=1$.
Therefore, since $k_0=1$, it follows that
$$
k_n=\prod_{i=0}^{n-1}\tilde{\alpha}_i^{-1}=\prod_{i=0}^{n-1}d_{n+i-1}^{-1} \; , \quad n\geq0\; .
$$
Thus (\ref{EqDistRod}) follows from (\ref{PnknRn}) and (\ref{rodrigues}).
Finally, by (\ref{coef}), the coefficients of the TTRR for $\{P_n\}_{n\geq0}$ are given by
$$
\begin{array}{c}
\displaystyle{ \beta_n= \frac{\tilde{\beta}_n}{\tilde{\alpha}_n}
=-\frac{d_{-2}q+2bnd_{n-1}}{d_{2n}d_{2n-2}} \; , \quad n\geq0 }\;; \\ [1em]
\displaystyle{ \gamma_n=\frac{\tilde{\gamma}_n k_n}{\tilde{\alpha}_n k_{n-1}}
=-\frac{nd_{n-2}}{d_{2n-3}d_{2n-1}} \, \phi \left(-\frac{e_{n-1}}{d_{2n-2}} \right) } \; , \quad n\geq1  \; .
\end{array}
$$
This completes the proof.
\qed

\begin{remark}\em
The regularity conditions (\ref{EqDistC2}) may be expressed as
\vspace*{-0.5em}
\begin{snugshade}
\begin{equation}\label{EqDistC2equiv}
d_n\cdot\big( ae_n^2-be_nd_{2n}+cd_{2n}^2\big)\neq0\,,\quad \forall n\in\mathbb{N}_0\;.
%d_n\neq0\;,\quad ae_n^2-be_nd_{2n}+cd_{2n}^2\neq0\,,\quad \forall n\in\mathbb{N}_0\;.
\end{equation}
\end{snugshade}
\end{remark}

\section*{Exercises}
%\bigskip

{\small
%\noindent
\begin{enumerate}[label=\emph{\bf \arabic*.},leftmargin=*]
\item
Let ${\bf u}\equiv{\bf u}(r_1,r_2)\in\mathcal{P}'$ (with $r_1,r_2\in\mathbb{C}$) be a solution of the Pearson's equation
$$
D\big((x-r_1)(x-r_2){\bf u}\big)=\Big(x-\frac{r_1+r_2}{2}\Big){\bf u}\;.
%\Big(x-\mbox{$\frac{r_1+r_2}{2}$}\Big){\bf u}\;.
$$
\begin{enumerate}
\item
Prove that ${\bf u}$ is regular if and only if $r_1\neq r_2$.
\item
Assuming the condition $r_1\neq r_2$, show that the monic OPS $\{P_n\}_{n\geq0}$
with respect to ${\bf u}$ is given by
$$
P_n(x):=\left(\frac{r_1-r_2}{4}\right)^n\,U_n\left(\frac{2x-r_1-r_2}{r_1-r_2}\right)\;, \quad n\in\mathbb{N}_0\;,
$$
where $\{U_n\}_{n\geq0}$ is the OPS of the Chebyshev polynomials of the second kind.
%\item
%Compute the moments of ${\bf u}$.
\end{enumerate}
%\item
%Dar uma equacao de diferecas para os momentos e resolver?
%\item
%Let ${\bf u}\equiv{\bf u}(\theta)\in\mathcal{P}'$ be a solution of the Pearson equation
%$$
%D\big((x^2-5x+6){\bf u}\big)=(2x+\theta){\bf u}\;.
%$$
%Determine all values of the parameter $\theta$ such that ${\bf u}$ is regular.
%For such values of $\theta$ show that the monic OPS with respect to ${\bf u}$ is ...
\end{enumerate}
}
\medskip

\section*{Final remarks}

This text is based on reference \ref{PacoPetronilho1994-C7}
and the works \ref{Maroni1991-C7}, \ref{Maroni1993-C7}, and \ref{Maroni1994-C7} by Maroni.
As far as we know, the regularity condition (\ref{r1}) was firstly stated  (in a different way) by the Russian mathematician
Ya. L. Geronimus in \ref{Geronimus1940-C7} (cf. Theorem II therein).
%O livro \ref{NikiforovUvarov1988-C7} de Nikiforov e Uvarov
%\'e uma fonte de motiva\c c\~ao para o tratamento dos PO cl\'assicos (e das suas vers\~oes discretas).
\medskip

%\section*{Notas finais}
%\section*{Coment\'arios finais}
%
%Este texto foi elaborado com base no artigo \ref{PacoPetronilho1994-C7},
%com incurs\~oes pelos trabalhos \ref{Maroni1991-C7}, \ref{Maroni1993-C7} e \ref{Maroni1994-C7} de Maroni.
%A condi\c c\~ao de regularidade (\ref{r1}) foi provada
%(por uma via distinta da que aparece nas notas de curso) pelo matem\'atico russo
%Ya. L. Geronimus em \ref{Geronimus1940-C7} (Teorema II).
%%O livro \ref{NikiforovUvarov1988-C7} de Nikiforov e Uvarov
%%\'e uma fonte de motiva\c c\~ao para o tratamento dos PO cl\'assicos (e das suas vers\~oes discretas).
%\medskip

\section*{Bibliography}
%\medskip

{\small
\begin{enumerate}[label=\emph{\rm [\arabic*]},leftmargin=*]
%\item\label{PacoAmilcarPetronilho1994-C5} F. Marcell\'an, A. Branquinho, and J. Petronilho,
%{\it Classical Orthogonal Polynomials: A Functional Approach}, Acta Applicand\ae Mathematic\ae {\bf 34} (1994) 283--303.
\item\label{Geronimus1940-C7} Ya. L. Geronimus,
{\it On polynomials orthogonal with respect to numerical sequences and on Hahn's theorem},
Izv. Akad. Nauk. \textbf{4} (1940) 215--228. (In Russian.)
\item\label{PacoPetronilho1994-C7} F. Marcell\'an and J. Petronilho,
{\it On the solution of some distributional differential equations: existence and characterizations of the classical moment functionals}, Integral Transforms and Special Functions {\bf 2} (1994) 185--218.
\item\label{Maroni1991-C7} P. Maroni, {\it Une th\'eorie alg\'ebrique des polyn\^omes orthogonaux. Applications aux polyn\^omes
orthogonaux semiclassiques}, In C. Brezinski et al. Eds., Orthogonal Polynomials and Their Applications, Proc. Erice 1990, IMACS, Ann. Comp. App. Math. {\bf 9} (1991) 95--130.
\item\label{Maroni1993-C7} P. Maroni, {\it  Variations Around Classical Orthogonal Polynomials. Connected Problems}, J. Comput. Appl. Math. {\bf 48} (1993) 133--155.
\item\label{Maroni1994-C7} P. Maroni, {\sl Fonctions eul\'eriennes. Polyn\^omes orthogonaux classiques},
T\'echniques de l'Ing\'enieur, trait\'e G\'en\'eralit\'es (Sciences Fondamentales), A {\bf 154} (1994) 1--30.
%\item\label{NikiforovUvarov1988-C7} A. F. Nikiforov and V. B. Uvarov,
%{\sl Special Functions of Mathematical Physics}. Birkhauser Verlag, Basel (1988).
\end{enumerate}
}

\chapter{Classical orthogonal polynomials}
%\chapter{Pearson distributional differential equation}

%\chapter{Orthogonal polynomials and dual basis.}
%More operations in $\mathcal{P}$ and $\mathcal{P}'$. }

\pagestyle{myheadings}\markright{Classical orthogonal polynomials}
\pagestyle{myheadings}\markleft{J. Petronilho}

The \emph{classical functionals} are the regular solutions (in $\mathcal{P}'$) of Pearson's equation.
The corresponding OPS are called \emph{classical orthogonal polynomials}.
In this text we present the most significant results concerning this important class of OP.

%\section{Classical OP: definition, existence, and characterizations}
\section{Definition and characterizations}

\begin{snugshade}
\begin{definition}\label{def-u-classical}
Let ${\bf u}\in\mathcal{P}'$.
${\bf u}$ is called a {\sl classical} functional if the following two conditions hold:
\begin{enumerate}
\item[{\rm (i)}]
${\bf u}$ is regular;
\item[{\rm (ii)}]
${\bf u}$ satisfies a Pearson distributional differential equation
\begin{equation}\label{EDClassic1}
D(\phi{\bf u})=\psi{\bf u}\;,
\end{equation}
where $\phi$ and $\psi$ are polynomials fulfilling
\begin{equation}\label{grauPhiPsi}
\deg\phi\leq2\;,\quad\deg\psi=1\;.
\end{equation}
\end{enumerate}
%If $\{P_n\}_{n\geq0}$ is an OPS with respect to a classical functional
%then $\{P_n\}_{n\geq0}$ is called a {\sl classical OPS}.
An OPS $\{P_n\}_{n\geq0}$ with respect to a classical functional
is called a {\sl classical OPS}.
\end{definition}
\end{snugshade}

%\begin{snugshade}
%\begin{definition}\label{def-u-classical}
%Let ${\bf u}\in\mathcal{P}'$.
%${\bf u}$ is a {\sl classical} functional if the following two conditions hold:
%\begin{enumerate}
%\item[{\rm (i)}]
%${\bf u}$ is regular;
%\item[{\rm (ii)}]
%%there exist two nonzero polynomials $\phi$ and $\psi$, being $\phi\in\mathcal{P}_2$ and $\psi\in\mathcal{P}_1$,
%there exist two polynomials $\phi$ and $\psi$, being
%\begin{equation}\label{grauPhiPsi}
%\deg\phi\leq2\;,\quad\deg\psi=1\;,
%\end{equation}
%such that ${\bf u}$ satisfies the Pearson distributional differential equation
%\begin{equation}\label{EDClassic1}
%D(\phi{\bf u})=\psi{\bf u}\;.
%\end{equation}
%\end{enumerate}
%If $\{P_n\}_{n\geq0}$ is an OPS with respect to a classical functional
%then $\{P_n\}_{n\geq0}$ is called a {\sl classical OPS}.
%\end{definition}
%\end{snugshade}

\begin{remark}\em
According with Lemma \ref{poly-grau21}, in the above definition
conditions (\ref{grauPhiPsi}) may be replaced by the weaker conditions
\begin{snugshade}\vspace*{-0.5em}
\begin{equation}\label{phiP2psiP1}
\phi\in\mathcal{P}_2\;,\quad\psi\in\mathcal{P}_1\;,\quad\{\phi,\psi\}\neq\mathcal{P}_{-1}:=\{0\}\;.
%\quad\{\phi,\psi\}\neq\{0\}\;.
%\quad |\phi|+|\psi|\not\equiv0\;.
\end{equation}
\end{snugshade}
%provided $\phi$ and $\psi$ are not simultaneously the zero polynomial.
\end{remark}

%Now, concerning Definition \ref{def-u-classical}, a question of major interest is the following:
%given ${\bf u}\in\mathcal{P}'\setminus\{\textbf{0}\}$ fulfilling (ii),
%%the distributional equation (\ref{EDClassic1}),
%%where $\phi\in\mathcal{P}_2$ and $\psi\in\mathcal{P}_1$, being $\phi$ and $\psi$ not simultaneously zero,
%find necessary and sufficient conditions involving only the (coefficients of the) polynomials $\phi$ and $\psi$
%such that (i) holds. %, i.e, ${\bf u}$ is regular.
%In order to answer to this question we begin by stating two preliminary lemmas.

Theorem \ref{vCOP-Thm1} gives necessary and sufficient conditions for the existence of solutions of Pearson's equation,
characterizing also such functionals (and, in particular, solving the question of the existence of classical functionals).
Thus, we may state:
{\it a functional ${\bf u}\in\mathcal{P}'\setminus\{{\bf 0}\}$ is classical if and only if
there exist $\phi\in\mathcal{P}_2$ and $\psi\in\mathcal{P}_1$ such that the following conditions hold:
\begin{snugshade}
\begin{equation}\label{P-regular1}
\begin{array}{rl}
{\rm (i)} & D(\phi{\bf u})=\psi{\bf u}\,; \\ %[0.25em]
{\rm (ii)} & na+p\neq0\;,\quad\displaystyle\phi\left(-\frac{nb+q}{2na+p}\right)\neq0\;,
\quad \forall n\in\mathbb{N}_0\;,
\end{array}
\end{equation}
\end{snugshade}
\noindent
where we have set $\phi(x)=ax^2+bx+c$ and $\psi(x)=px+q$.}

In the next proposition we state several characterizations of the classical OPS.
For convenience, we introduce the concept of admissible pair of polynomials.

\begin{snugshade}
\begin{definition}
%Let $\phi$ and $\psi$ be two polynomials such that
$(\phi,\psi)$ is called an {\sl admissible pair} if
$$
\phi\in\mathcal{P}_2\; , \quad \psi\in\mathcal{P}_1 \; ,\quad
d_n:=\psi^{\prime}+\mbox{$\frac{n}{2}$}\,\phi^{\prime\prime}\neq0 \; , \;\;\forall n\in\mathbb{N}_0.
$$
\end{definition}
\end{snugshade}

Introducing this concept makes sense, since according with %Theorem \ref{vCOP-Thm1}
conditions (ii) in (\ref{P-regular1}), only
admissible pairs may appear in the framework of the theory of classical OP.

\begin{snugshade}
\begin{theorem}[characterizations of the classical OPS]\label{ThmClassicalOPS}
Let ${\bf u}\in\mathcal{P}^\prime$, regular, and let $\{P_n\}_{n\geq0}$ be its monic OPS.
Then the following properties are equivalent:
\begin{enumerate}
\item[{\rm C1.}]
${\bf u}$ is classical, i.e., there are nonzero polynomials $\phi\in\mathcal{P}_2$ and $\psi\in\mathcal{P}_1$
such that ${\bf u}$ satisfies the distributional Pearson's differential equation
$$D(\phi {\bf u})=\psi {\bf u}\; ;$$
\item[{\rm C1$'$.}]
there is an admissible pair $(\phi,\psi)$ such that ${\bf u}$ satisfies Pearson's equation
$$D(\phi {\bf u})=\psi {\bf u}\; ;$$
\item[{\rm C2.}] {\rm (Al-Salam $\&$ Chihara)} %, \cite{16})}
there exist a polynomial $\phi\in\mathcal{P}_2$ and, for each $n\in\mathbb{N}_0$,
complex parameters $a_n$, $b_n$ and $c_n$, with $c_n\neq0$ if $n\geq1$, such that
$$
\phi(x)P_n^{\prime}(x)=a_nP_{n+1}(x)+b_nP_n(x)+c_nP_{n-1}(x)\;,\quad n\geq0 \; ;
$$
\item[{\rm C3.}] {\rm (Hahn)} %, \cite{24})} %\noindent\;{\rm C3$'$.}
$\Big\{ P_n^{[k]}:=\frac{{\rm d}^k}{{\rm d}x^k}\frac{P_{n+k}}{(n+1)_k}\Big\}_{n\geq0}$ is a monic OPS for some $k\in\mathbb{N}\,$;
%$\Big\{ Q_n:=\frac{P_{n+1}^{\prime}}{n+1}\Big\}_{n\geq0}$ is a monic OPS;
%(being the corresponding regular functional ${\bf v}:=\phi{\bf u}$);
\item[{\rm C3$'$.}] %{\rm (Hahn)} %, \cite{24})} %\noindent\;{\rm C3$'$.}
$\big\{ P_n^{[k]}\big\}_{n\geq0}$ is a monic OPS for each $k\in\mathbb{N}\,$;
\item[{\rm C4.}]
there exist $k\in\mathbb{N}$ and complex parameters $r_{n}^{[k]}$ and $s_{n}^{[k]}$ such that
$$
P_n^{[k-1]}(x)=P_n^{[k]}(x)+r_{n}^{[k]}P_{n-1}^{[k]}(x)+s_{n}^{[k]}P_{n-2}^{[k]}(x) \; , \quad n\geq2\; ;
\eqno(\star)
$$
\item[{\rm C4$'$.}]
for each $k\in\mathbb{N}$, there exist parameters $r_{n}^{[k]}$ and $s_{n}^{[k]}$ such that $(\star)$ holds;
%$$
%P_n^{[k-1]}(x)=P_n^{[k]}(x)+r_{n,k}P_{n-1}^{[k]}(x)+s_{n,k}P_{n-2}^{[k]}(x) \; , \quad n\geq2\; ;
%$$
\item[{\rm C5.}] {\rm (Bochner)} %, \cite{1})}
there exist polynomials $\phi$ and $\psi$ and, for each $n\geq0$, a complex parameter $\lambda_n$,
with $\lambda_n\neq0$ if $n\geq1$, such that $y=P_n(x)$ is a solution of the second order
ordinary differential equation
$$
\phi(x) y^{\prime\prime} + \psi(x) y^{\prime} +\lambda_n y =0 \;,\quad n\geq0\;;
$$
%{\rm (} $\lambda_n$ is given by $\lambda_n=-n[\psi^{\prime}+(n-1)\frac{\phi^{\prime \prime}}{2}] \ , \ n\geq0$ {\rm )} ;
\item[{\rm C6.}] {\rm (Maroni)} %, \cite{14})}
there is an admissible pair $(\phi,\psi)$ so that the formal Stieltjes series associated with ${\bf u}$,
$S_{\bf u}(z):=-\sum_{n=0}^{\infty}u_{n}/z^{n+1}$,
satisfies (formally) %satisfies the formal first order linear ordinary differential equation
$$
\phi(z)S_{\bf u}^{\prime}(z)=[\psi(z)-\phi^{\prime}(z)]S_{\bf u}(z)+(\psi^{\prime}
-\mbox{$\frac12$}\,\phi^{\prime \prime}) u_0 \; ;
$$
\item[{\rm C7.}] {\rm (McCarthy)} %, \cite{22})}
there exists an admissible pair $(\phi,\psi)$ and, for each $n\geq1$, complex parameters $h_n$ and $t_n$ such that
$$
\phi(P_nP_{n-1})^{\prime}(x)=h_nP_n^2(x)-(\psi-\phi^{\prime})P_nP_{n-1}(x)+t_nP_{n-1}^2(x) \; ;
$$
\item[{\rm C8.}] {\rm (distributional Rodrigues formula)} %, \cite{3})}
there exist a polynomial $\phi\in\mathcal{P}_2$ and nonzero complex parameters $k_n$
%, with $k_n\neq0$ for each $n\geq0$,
such that
$$
P_n(x) {\bf u} = k_n D^n\big(\phi^n(x){\bf u}\big) \; , \quad n\geq0 \;.
$$
\end{enumerate}
Moreover, the polynomials $\phi$ and $\psi$ %appearing above %in the different characterizations
may be taken the same in all properties above where they appear.
In addition, let the TTRR %three-term recurrence relation
fulfilled by the monic OPS $\{P_n\}_{n\geq0}$ be
$$%\begin{equation}\label{paramTTRRC}
P_{n+1}(x)=(x-\beta_n)P_n(x)-\gamma_{n}P_{n-1}(x)\;,\quad n\geq0
$$%\end{equation}
($P_{-1}(x)=0$; $P_0(x)=1$).
Write $\phi(x)=ax^2+bx+c$, $\psi(x)=px+q$, $d_n:=na+p$, and $e_n:=nb+q$. Then
$$%\begin{equation}\label{paramTTRRCbeta}
%\beta_n=\frac{ne_n}{d_{2n}}-\frac{(n+1)e_{n+1}}{d_{2n+2}}
\beta_n=-\frac{d_{-2}q+2bnd_{n-1}}{d_{2n}d_{2n-2}}\;,\quad
%\beta_n=-\frac{(-2a+p)q+2bn[(n-1)a+p]}{(2na+p)[(2n-2)a+p]}\;,%\;n\geq0\;;
%$$%\end{equation}
%$$%\begin{equation}\label{paramTTRRCgamma}
\gamma_{n}=-\frac{nd_{n-2}}{d_{2n-3}d_{2n-1}}\phi\Big(-\frac{e_{n-1}}{d_{2n-2}}\Big)\;,
%\gamma_{n+1}=\frac{-(n+1)[(n-1)a+p][a(nb+q)^2-b(nb+q)(2na+p)+c(2na+p)^2]}{[(2n-1)a+p](2na+p)^2[2(n+1)a+p]}%\;,\;n\geq0\;.
$$%\end{equation}
%for each $n\geq0$.
and the parameters appearing in the above characterizations
%Furthermore, the parameters appearing in the above characterizations
%given in Theorem \ref{ThmClassicalOPS}
may be computed explicitly: % as follows:
$$
\begin{array}{c}
\displaystyle
a_n=na\;,\;\; b_n=-\mbox{$\frac12$}\psi(\beta_n)\;,\;\; c_n=-d_{n-1}\gamma_n\;,\;\;
r_n^{[1]}=\mbox{$\frac12$}\frac{\psi(\beta_n)}{d_{n-1}}\;,\;\;
s_n^{[1]}=-\frac{(n-1)a}{d_{n-2}}\gamma_n\;, \\ [1.0em]
\displaystyle
\lambda_n=-nd_{n-1}\;,\;\;
h_n=d_{2n-3}\;,\;\;
t_n=-d_{2n-1}\gamma_n\;,\;\;
k_n=\mbox{$\prod_{i=0}^{n-1}d_{n+i-1}^{-1}$}\;.
\end{array}
$$

%$$%\begin{equation}\label{paramC2}
%a_n=na\;,\quad b_n=-\mbox{$\frac12$}\psi(\beta_n)\;,\quad c_n=-d_{n-1}\gamma_n\;,
%$$%\end{equation}
%$$%\begin{equation}\label{paramC4}
%r_n=\mbox{$\frac12$}\frac{\psi(\beta_n)}{d_{n-1}}\;,\quad
%s_n=-\frac{(n-1)a}{d_{n-2}}\gamma_n\;,
%$$%\end{equation}
%$$%\begin{equation}\label{paramC58}
%\lambda_n=-nd_{n-1}\;,\quad
%k_n=\prod_{i=0}^{n-1}d_{n+i-1}^{-1}\;,
%$$%\end{equation}
%$$%\begin{equation}\label{paramC7}
%h_n=\;,\quad
%t_n=\;.
%$$%\end{equation}
\end{theorem}
\end{snugshade}

{\it Proof.} By Lemma \ref{reg-lemma3} and Theorem \ref{vCOP-Thm1},
C1$\,\Leftrightarrow\,$C1$'$, C1$\,\Rightarrow\,$C3$'$, and C1$'$$\,\Leftrightarrow\,$C8.
%C1$'$ and C8 are equivalent to C1.
Clearly, C3$'$$\,\Rightarrow\,$C3 and C4$'$$\,\Rightarrow\,$C4.
We show that C3$'$$\,\Rightarrow\,$C4$'$ using the same arguments
of the proof of C3$\,\Rightarrow\,$C4 given in bellow.
%We omit the proof of C1$'$$\,\Leftrightarrow\,$C6.
The proof of C1$'$$\,\Leftrightarrow\,$C6 is left to the reader (Exercise \ref{Ex-cp8-0}).
Thus, we only need to show that:
\smallskip
\begin{center}
C1$'$$\,\Rightarrow\,$C2$\,\Rightarrow\,$C3$\,\Rightarrow\,$C4$\,\Rightarrow\,$C1$\,,\;\;$
C1$\,\Leftrightarrow\,$C5$\,,\;\;$ C2$\,\Leftrightarrow\,$C7.
\end{center}
\smallskip
%c1 -> c2 -> c3 -> c4 -> c5 -> c1 <-> cRodrigues

(C1$'$$\,\Rightarrow\,$C2).
Assume that C1$'$ holds. Fix $n\in\mathbb{N}_0$. Since $\deg(\phi P_n')\leq n+1$, then
\begin{equation}\label{ThmChar1}
\phi P'_n=\sum_{j=0}^{n+1}a_{n,j}P_j\,,\quad
a_{n,j}:=\frac{\langle{\bf u},\phi P_n'P_j\rangle}{\langle{\bf u},P_j^2\rangle}\;.
\end{equation}
For each integer number $j$, with $0\leq j\leq n+1$, we deduce
\begin{equation}\label{ThmChar2}
\begin{array}{rl}
\langle{\bf u},\phi P_n'P_j\rangle &=\,\langle\phi{\bf u},(P_nP_j)'-P_nP_j'\rangle=
-\langle D(\phi{\bf u}),P_nP_j\rangle-\langle \phi{\bf u}, P_nP_j'\rangle \\ [0.5em]
&=\, -\langle {\bf u},\psi P_jP_n\rangle-\langle {\bf u},\phi P_j'P_n\rangle\;.
\end{array}
\end{equation}
If $0\leq j\leq n-2$ we obtain $\langle{\bf u},\phi P_n'P_j\rangle=0$, and so $a_{n,j}=0$.
Thus, (\ref{ThmChar1}) reduces to
$$
\phi P'_n=a_nP_{n+1}+b_nP_n+c_nP_{n-1}\;,\quad n\geq0\;,
$$
where, writing $\phi(x)=ax^2+bx+c$ and $\psi(x)=px+q$,
$a_n=na$ (by comparison of coefficients), $b_n=a_{n,n}$, and $c_n:=a_{n,n-1}$.
%\begin{equation}
% a{P_j^2}\,,\quad
%a_{n,j}:=\frac{\langle{\bf u},\phi P_n'P_j\rangle}{\langle{\bf u},P_j^2\rangle}\;.
%\end{equation}
Setting $j=n-1$ in (\ref{ThmChar2}), we deduce
$$
\langle{\bf u},\phi P_n'P_{n-1}\rangle=
-\langle {\bf u},(\psi P_{n-1}+\phi P_{n-1}')P_n\rangle=-d_{n-1}\langle{\bf u}, P_n^2\rangle\;,
$$
hence
$$
c_n:=a_{n,n-1}=\frac{\langle{\bf u},\phi P_n'P_{n-1}\rangle}{\langle{\bf u},P_{n-1}^2\rangle}
=\frac{\langle{\bf u},\phi P_n'P_{n-1}\rangle}{\langle{\bf u},P_{n}^2\rangle}
\frac{\langle{\bf u},P_n^2\rangle}{\langle{\bf u},P_{n-1}^2\rangle}=-d_{n-1}\gamma_n\;,\quad n\geq1\;.
$$
Since, by hypothesis, $(\phi,\psi)$ is an admissible pair,
then we may conclude that $c_n\neq0$ for each $n\geq1$. Thus C1$'$$\,\Rightarrow\,$C2.
Notice that taking $j=n$ in (\ref{ThmChar2}) yields
$$
\langle{\bf u},\phi P_n'P_n\rangle=-\mbox{$\frac12\,$}\langle {\bf u},\psi P_n^2\rangle
=-\mbox{$\frac12\,$}\big(p\langle {\bf u},x P_n^2\rangle+q\langle {\bf u},P_n^2\rangle\big)\;,
$$
hence we deduce the expression for $b_n$ given in the statement of the theorem:
$$
b_n:=a_{n,n}=\frac{\langle{\bf u},\phi P_n'P_{n}\rangle}{\langle{\bf u},P_{n}^2\rangle}
=-\mbox{$\frac12\,$}\Big(p\frac{\langle {\bf u},x P_n^2\rangle}{\langle {\bf u},P_n^2\rangle}+q\Big)
=-\mbox{$\frac12$}\psi(\beta_n)\;.
$$
%\smallskip

(C2$\,\Rightarrow\,$C3).
Suppose that C2 holds.
We will show that $\{ P_n^{[1]}:=P_{n+1}^{\prime}/(n+1)\}_{n\geq0}$
is a monic OPS with respect to ${\bf v}:=\phi{\bf u}$.
Indeed, for each $n\in\mathbb{N}_0$ and $0\leq m\leq n$, %we deduce
$$
\begin{array}{rcl}
(n+1)\langle{\bf v},x^m P_n^{[1]} \rangle &=& \langle\phi{\bf u},x^mP_{n+1}'\rangle
=\langle{\bf u},\big(\phi P_{n+1}'\big)x^m\rangle \\ [0.25em]
&=& \langle {\bf u},(a_{n+1}P_{n+2}+b_{n+1}P_{n+1}+c_{n+1}P_{n})x^m\rangle
=c_{n+1}\langle {\bf u},P_{n}^2\rangle\delta_{m,n}\;.\\ %[0.25em]
%&=& c_{n+1}\langle {\bf u},P_{n}^2\rangle\delta_{m,n}\;.
\end{array}
$$
Therefore, since (by hypothesis) $c_{n+1}\neq0$ for each $n\geq0$,
we conclude that $\{P_n^{[1]}\}_{n\geq0}$
is a monic OPS (with respect to ${\bf v}:=\phi{\bf u}$).
%Thus (C2$\,\Rightarrow\,$C3).
\smallskip

(C3$\,\Rightarrow\,$C4).
By hypothesis, $\{P_n^{[k]}:=\frac{{\rm d}^k}{{\rm d}x^k}\big(\frac{P_{n+k}}{(n+1)_k}\big)\}_{n\geq0}$
is a monic OPS for some (fixed) $k\in\mathbb{N}$. Then there exists $\beta_n^{[k]}\in\mathbb{C}$ and $\gamma_n^{[k]}\in\mathbb{C}\setminus\{0\}$ such that
\begin{equation}\label{PnkTTRR}
xP_n^{[k]}=P_{n+1}^{[k]}+\beta_n^{[k]}P_n^{[k]}+\gamma_n^{[k]}P_{n-1}^{[k]} \;,\quad n\in\mathbb{N}_0 \,.
\end{equation}
Similarly, there exists $\beta_n\in\mathbb{C}$ and $\gamma_n\in\mathbb{C}\setminus\{0\}$ such that
\begin{equation}\label{PnTTRR}
xP_n=P_{n+1}+\beta_nP_n+\gamma_{n}P_{n-1}\;,\quad n\in\mathbb{N}_0\,.
\end{equation}
%where $P_{-1}(x)=0$; $P_0(x)=1$).
Changing $n$ into $n+k$ in (\ref{PnTTRR}), then taking the derivative of order $k$ in both sides of the resulting equation
and using Leibnitz rule on the left-hand side, we find %after dividing by $(n+1)_k$, we find
$$
xP_n^{[k]}+\frac{k}{n+1}P_{n+1}^{[k-1]}=
\frac{n+k+1}{n+1}P_{n+1}^{[k]}+\beta_{n+k}P_n^{[k]}+\frac{n\gamma_{n+k}}{n+k}P_{n-1}^{[k]} \;,\quad n\in\mathbb{N}_0 \,.
$$
In this equation, replacing $xP_n^{[k]}$ by the right-hand side of (\ref{PnkTTRR}), and then changing $n$ into $n-1$,
we obtain $(\star)$, with
$$
r_{n}^{[k]}=\frac{n\,\big(\beta_{n+k-1}-\beta_{n-1}^{[k]}\big)}{k}\;,\quad
s_{n}^{[k]}=\frac{n\,\big((n-1)\gamma_{n+k-1}-(n+k-1)\gamma_{n-1}^{[k]}\big)}{k(n+k-1)}\;.
$$
\smallskip

(C4$\,\Rightarrow\,$C1).
By hypothesis $(\star)$ holds.
Let $\{{\bf a}_n\}_{n\geq0}$ and $\{{\bf a}_n^{[k]}\}_{n\geq0}$ be the dual basis
for $\{P_n\}_{n\geq0}$ and $\{P_n^{[k]}\}_{n\geq0}$, respectively. By Theorem \ref{expDB1},
$\,{\bf a}_n^{[k]}=\sum_{j\geq0}\langle{\bf a}_n^{[k]},P_j^{[k-1]}\rangle {\bf a}_j^{[k-1]}$ for each $n\in\mathbb{N}_0$.
Using $(\star)$, we compute
$$
\;\langle{\bf a}_n^{[k]},P_j^{[k-1]}\rangle
=\langle{\bf a}_n^{[k]},P_j^{[k]}\rangle
+r_j^{[k]}\langle{\bf a}_n^{[k]},P_{j-1}^{[k]}\rangle
+s_j^{[k]}\langle{\bf a}_n^{[k]},P_{j-2}^{[k]}\rangle
=\left\{
\begin{array}{cl}1\;, &\mbox{\rm if}\; j=n \\ [0.5em]
r_{n+1}^{[k]}\;, &\mbox{\rm if}\; j=n+1 \\ [0.5em]
s_{n+2}^{[k]}\;, &\mbox{\rm if}\; j=n+2 \\ [0.5em]
0\;, &\mbox{\rm otherwise}\,, %& j\neq n,n+1,n+2 \\ [0.5em] % j\leq n-1\vee j\geq n+3 \\
\end{array}
\right.
$$
hence
$$
{\bf a}_n^{[k]}= {\bf a}_n^{[k-1]}+r_{n+1}^{[k]}{\bf a}_{n+1}^{[k-1]}+s_{n+2}^{[k]}{\bf a}_{n+2}^{[k-1]}\;,\quad n\in\mathbb{N}_0\;.
$$
Taking the (distributional) derivative of order $k$ in both sides of this equation, and using the relations
$D^j\big({\bf a}_n^{[j]}\big)=(-1)^j(n+1)_j\,{\bf a}_{n+j}$ (see Lemma \ref{teo-rel-dual-an-ank}), %($j,n\in\mathbb{N}_0$),
%$$
%{\bf a}_n=\frac{P_n}{\langle{\bf u},P_n^2\rangle}\,{\bf u}\,,\quad
%D^\ell\big({\bf a}_n^{[\ell]}\big)=(-1)^\ell(n+1)_\ell\,{\bf a}_{n+\ell}\,,\quad\ell,n\in\mathbb{N}_0\;,
%$$
we obtain
$$
D\left(\frac{1}{n+k}\,{\bf a}_{n+k-1}+\frac{r_{n+1}^{[k]}}{n+1}\,{\bf a}_{n+k}
+\frac{(n+k+1)s_{n+2}^{[k]}}{(n+1)(n+2)}\,{\bf a}_{n+k+1}\right)=-{\bf a}_{n+k}\;,
\quad n\in\mathbb{N}_0\;.
$$
Therefore, since, by Theorem \ref{dualR}, ${\bf a}_j=\frac{P_j}{\langle{\bf u},P_j^2\rangle}\,{\bf u}$
for $j\in\mathbb{N}_0$ and, by Corollary \ref{cor-TTRR},
$\gamma_j=\frac{{\langle{\bf u},P_j^2\rangle}}{{\langle{\bf u},P_{j-1}^2\rangle}}$ for $j\in\mathbb{N}$,
being $\gamma_j$ the $\gamma-$parameter appearing in (\ref{PnTTRR}), we deduce
\begin{equation}\label{DPhink1}
D\big(\Phi_{n+k+1}\,{\bf u}\big)=-P_{n+k}\,{\bf u}\;,\quad n\in\mathbb{N}_0\;,
\end{equation}
where $\Phi_{n+k+1}$ is a polynomial of degree at most $n+k+1$, given by
$$%\begin{equation}\label{Phink1}
\Phi_{n+k+1}(x):=\frac{\gamma_{n+k}}{n+k}\,P_{n+k-1}(x)+
\frac{r_{n+1}^{[k]}}{n+1}\,P_{n+k}(x)+\frac{(n+k+1)s_{n+2}^{[k]}}{(n+1)(n+2)\gamma_{n+k+1}}\,P_{n+k+1}(x)\,.
%\quad n\in\mathbb{N}_0\;.
$$%\end{equation}
Since $\Phi_{n+k+1}$ is a (finite) linear combination of polynomials of the simple set $\{P_j\}_{j\geq0}$
and $\gamma_{n+k}\neq0$, then $\Phi_{n+k+1}$ does not vanishes identically, so
$\Phi_{n+k+1}\in\mathcal{P}_{n+k+1}\setminus\{0\}$.
Setting $n=0$ and $n=1$ in (\ref{DPhink1}) we obtain the two equations
\begin{equation}\label{DPhink1a}
D\big(\Phi_{k+1}\,{\bf u}\big)=-P_{k}\,{\bf u}\;,\quad
D\big(\Phi_{k+2}\,{\bf u}\big)=-P_{k+1}\,{\bf u}\;.
\end{equation}
If $k=1$ it follows immediately from the first of these equations that C1 holds. %${\bf u}$ is classical.
Henceforth, assume that $k\geq2$.
%In this case we will prove that C1 holds following essentially the arguments
%presented in the proof of Lemma 3.2 in \ref{LoureiroMaroniZelia2006-C8}.
Setting $n=0$ and $n=1$ in the definition of $\Phi_{n+k+1}$
and using the TTRR (\ref{PnTTRR}), we easily deduce
\begin{equation}\label{DPhink1b}
\left\{
\begin{array}{rcl}
\Phi_{k+1}(x) &=& E_0(x;k)P_{k+1}(x)+F_1(x;k)P_{k}(x)\;, \\ [0.5em]
\Phi_{k+2}(x) &=& G_1(x;k)P_{k+1}(x)+H_0(x;k)P_{k}(x)\;,
\end{array}
\right.
\end{equation}
where $E_0(\cdot;k),H_0(\cdot;k)\in\mathcal{P}_0$ and $F_1(\cdot;k),G_1(\cdot;k)\in\mathcal{P}_1$,
explicitly given by
\begin{equation}\label{DPhink1ba}
\begin{array}{l}
E_0(x;k) := \displaystyle\frac{(k+1)s_{2}^{[k]}}{2\gamma_{k+1}}-\frac{1}{k}\;,\quad
F_1(x;k) := \frac{x-\beta_k}{k}+r_{1}^{[k]} \;, \\ [0.75em]%\in\mathcal{P}_1\;, \\ [0.5em]
G_1(x;k) := \displaystyle\frac{(k+2)s_{3}^{[k]}(x-\beta_{k+1})}{6\gamma_{k+2}}+\frac{r_{2}^{[k]}}{2}\;, \quad
H_0(x;k) := \frac{\gamma_{k+1}}{k+1}-\frac{(k+2)s_{3}^{[k]}\gamma_{k+1}}{6\gamma_{k+2}}\;.
\end{array}
\end{equation}
%%%%%%%%%%%%%%%%%%%%%%%%%%%%%%%%%%%%%%%%%%%%%%%%%%%%%%%%%%%%%%%%%%%%%%%%%%%%%%%%%%
%Consider the system of two equations (\ref{DPhink1b}),
%where we regard $P_k$ and $P_{k+1}$ as unknowns.
%The determinant of such a system is
%$$%\begin{equation}\label{DPhink1Det1}
%\Delta_2(x)\equiv\Delta_2(x;k) :=
%\left|
%\begin{array}{cc}
%E_0(x;k) & F_1(x;k) \\ [0.25em]
%G_1(x;k) & H_0(x;k)
%\end{array}
%\right|\;.
%$$%\end{equation}
%Clearly, $\Delta_2\in\mathcal{P}_2$.
%Solving the system for $P_k$ and $P_{k+1}$ we obtain the two equations
%%%%%%%%%%%%%%%%%%%%%%%%%%%%%%%%%%%%%%%%%%%%%%%%%%%%%%%%%%%%%%%%%%%%%%%%%%%%%%%%%%%%%%%%%%
Let $\Delta_2(x)\equiv\Delta_2(x;k):=E_0(x;k)H_0(x;k)- F_1(x;k)G_1(x;k)$,
the determinant of the system (\ref{DPhink1b}).
Using (\ref{DPhink1a})--(\ref{DPhink1ba}), and taking into account
that ${\bf u}$ is regular, we prove that
$\Delta_2\in\mathcal{P}_2\setminus\{0\}$ (Exercise \ref{Ex-cp8-0}).
Solving (\ref{DPhink1b}) for $P_k$ and $P_{k+1}$ we obtain
%From (\ref{DPhink1b}), we obtain
\begin{eqnarray}
\label{DPhink1Det2}
\Delta_2(x)P_{k+1}(x)=H_0(x;k)\Phi_{k+1}(x)-F_1(x;k)\Phi_{k+2}(x)\,, \\
\label{DPhink1Det3}
\Delta_2(x)P_{k}(x)=E_0(x;k)\Phi_{k+2}(x)-G_1(x;k)\Phi_{k+1}(x)\;.
\end{eqnarray}
Since $P_k$ and $P_{k+1}$ cannot share zeros,
it follows from (\ref{DPhink1Det2})--(\ref{DPhink1Det3})
that any common zero of $\Phi_{k+1}$ and $\Phi_{k+2}$
(if there is some) must be a zero of $\Delta_2$.
Let $\Phi$ be the greatest common divisor of $\Phi_{k+1}$ and $\Phi_{k+2}$, i.e.,
$$
\Phi(x):=\mbox{g.c.d.}\,\{\Phi_{k+1}(x),\Phi_{k+2}(x)\}\;.
$$
Any zero of $\Phi$ is also a zero of both $\Phi_{k+1}$ and $\Phi_{k+2}$,
and so it is a zero of $\Delta_2$. Therefore, $\Phi\in\mathcal{P}_2\setminus\{0\}$.
%$$\Phi\leq\deg\Delta_2\leq2\;.$$
(Notice that indeed $\Phi\not\equiv0$, since $\Phi_{k+1}\not\equiv0$ and $\Phi_{k+2}\not\equiv0$.)
Moreover, there exist polynomials $\Phi_{1,k}$ and $\Phi_{2,k}$, with no common zeros, such that
\begin{eqnarray}
\label{DPhink1A}\Phi_{k+1}=\Phi\,\Phi_{1,k}\;,\quad\Phi_{k+2}=\Phi\,\Phi_{2,k}\;,
\qquad\qquad\qquad \\ [0.25em]
\label{DPhink1Aa}
\Phi_{1,k}\in\mathcal{P}_{k+1-\ell}\setminus\{0\}\;,\quad
\Phi_{2,k}\in\mathcal{P}_{k+2-\ell}\setminus\{0\}\;,\quad
\ell:=\deg\Phi\leq2\;.
\end{eqnarray}
From (\ref{DPhink1a}) and (\ref{DPhink1A}) we deduce
\begin{equation}\label{DPhink1B}
\Phi_{1,k}D(\Phi{\bf u})=-(P_k+\Phi_{1,k}'\Phi){\bf u}\;,\quad
\Phi_{2,k}D(\Phi{\bf u})=-(P_{k+1}+\Phi_{2,k}'\Phi){\bf u}\;.
\end{equation}
Combining these two equations yields
%\begin{equation}\label{DPhink1C}
$\big(\Phi_{1,k}(P_{k+1}+\Phi_{2,k}'\Phi)-\Phi_{2,k}(P_k+\Phi_{1,k}'\Phi)\big){\bf u}={\bf 0}$,
%\end{equation}
and so, since ${\bf u}$ is regular,
%\begin{equation}\label{DPhink1D}
$\Phi_{1,k}(P_{k+1}+\Phi_{2,k}'\Phi)=\Phi_{2,k}(P_k+\Phi_{1,k}'\Phi)$.
%\end{equation}
Therefore, taking into account that $\Phi_{1,k}$ and $\Phi_{2,k}$ have no common zeros
and (\ref{DPhink1Aa}) holds, we may ensure that there exists a polynomial $\Psi\in\mathcal{P}_1$ such that
\begin{equation}\label{DPhink1E}
P_k+\Phi_{1,k}'\Phi=-\Psi\Phi_{1,k}\;,\quad
P_{k+1}+\Phi_{2,k}'\Phi=-\Psi\Phi_{2,k}\;.
\end{equation}
Combining equations (\ref{DPhink1B}) and (\ref{DPhink1E}) we deduce
$$
\Phi_{1,k}\big(D(\Phi{\bf u})-\Psi{\bf u}\big)=
\Phi_{2,k}\big(D(\Phi{\bf u})-\Psi{\bf u}\big)={\bf 0}\;.
$$
From these equations, and using once again the fact that
$\Phi_{1,k}$ and $\Phi_{2,k}$ have no common zeros,
we conclude, by Proposition \ref{prop-uvzero}, that $D(\Phi{\bf u})=\Psi{\bf u}$. %Exercise \ref{Ex-cp8-00}).
Thus C4$\,\Rightarrow\,$C1.
%{\red por o exercicio; referencias; passar $\Phi$ e $\Psi$ a $\phi$ e $\psi$; calcular $r_n$ e $s_n$ apenas para $k=1$ e deixar o %caso geral para os exercicios}
%Setting $n=0$ and taking into account Theorem \ref{dualR}, we obtain
%$$
%{\bf a}_0^{[1]}=\phi{\bf u}\;,\quad
%\phi(x):=\frac{1}{u_0}+\frac{r_1}{\langle{\bf u},P_1^2\rangle}P_1(x)+\frac{s_2}{\langle{\bf u},P_2^2\rangle}P_2(x)\;.
%$$
%Therefore, differentiating and taking into account (\ref{rel-dual-an-an1}), we deduce
%$$
%D(\phi{\bf u})=D{\bf a}_0^{[1]}=-{\bf a}_1=-\frac{P_1}{\langle{\bf u},P_1^2\rangle}{\bf u}=\psi{\bf u}\;,\quad
%\psi(x):=-\frac{P_1(x)}{\langle{\bf u},P_1^2\rangle}\;,
%$$
%and so $D(\phi{\bf u})=\psi{\bf u}$, with $\phi\in\mathcal{P}_2$ and $\deg\psi=1$.
%Thus C4$\,\Rightarrow\,$C1.
%For $k=1$,
The formulas for $r_n^{[1]}$ and $s_n^{[1]}$ given in the statement of
the theorem may be derived as follows. We have already proved that
C4$\,\Rightarrow\,$C1$\,\Rightarrow\,$C1$'$$\,\Rightarrow\,$C2$\,\Rightarrow\,$C3$\,\Rightarrow\,$C4,
and we see that the polynomials $\phi$ and $\psi$ appearing in all these characterizations may be taken the same.
As we have seen, %$b_n=-\frac12\psi(\beta_n)$ and $c_n=-d_{n-1}\gamma_n$.
the formulas for $b_n$ and $c_n$ given in the statement of the theorem hold.
We now use these formulas to obtain the expressions for $r_n^{[1]}$ and $s_n^{[1]}$.
Set $Q_n:=P_n^{[1]}:=P_{n+1}'/(n+1)$.
By C4, $P_n=Q_n+r_n^{[1]}Q_{n-1}+s_n^{[1]}Q_{n-2}$ if $n\geq2$.
Hence, since $\{Q_n\}_{n\geq0}$ is a monic OPS with respect to ${\bf v}:=\phi{\bf u}$,
we deduce, for each $n\geq2$,
$$
\; r_n^{[1]}=\frac{\langle{\bf u},\phi  P_n P_n^{\prime}\rangle }{\langle{\bf u},\phi  P_n^{\prime} P_{n-1}\rangle }
=\frac{\langle{\bf u}, P_{n-1}^2\rangle }{\langle{\bf u},\phi  P_n^{\prime} P_{n-1}\rangle }
\frac{\langle{\bf u}, \phi  P_n^{\prime} P_n\rangle }{\langle{\bf u}, P_n^2\rangle }
\frac{\langle{\bf u}, P_n^2\rangle }{\langle{\bf u}, P_{n-1}^2\rangle }=\frac{1}{c_n} b_n \gamma_n
=\mbox{$\frac{1}{2}$}\frac{\psi(\beta_n)}{d_{n-1}}\,,%\; n\geq2 \, ,
$$
where the third equality holds taking into account C2. %, $\phi P_n^{\prime}=a_nP_{n+1}+b_nP_n+c_nP_{n-1}$.
Similarly, for each $n\geq2$,
$$
s_n^{[1]}=\frac{a\langle{\bf u}, P_n^2\rangle }{\frac{1}{n-1}\langle{\bf u},\phi  P_{n-1}^{\prime} P_{n-2}\rangle}
=\frac{(n-1)a \langle{\bf u}, P_n^2\rangle}{c_{n-1}\langle{\bf u}, P_{n-2}^2\rangle}
=\frac{(n-1)a}{c_{n-1}}\gamma_{n-1}\gamma_n=-\frac{(n-1)a}{d_{n-2}}\gamma_n \;. %\;, \quad n\geq2\; .
$$
%\smallskip

(C1$\,\Rightarrow\,$C5).
By hypothesis, $D(\phi{\bf u})=\psi{\bf u}$, where
$\phi\in\mathcal{P}_2$, $\psi\in\mathcal{P}_1$, and $\deg\psi=1$ (cf. Lemma \ref{poly-grau21}).
Fix $n\in\mathbb{N}$, and write
\begin{equation}\label{C1C5C3eq}
\phi P_n^{\prime\prime}+\psi P_n'=\sum_{j=0}^n\lambda_{n,j}P_j\;.  %,\quad n\geq0\;.
\end{equation}
Then, for each $j$ such that $0\leq j\leq n$,
$$
\begin{array}{rcl}
\langle{\bf u}, P_j^2\rangle\lambda_{n,j} &=&
\langle {\bf u},(\phi P_n^{\prime\prime}+\psi P_n')P_j\rangle =
\langle \phi{\bf u},P_n^{\prime\prime}P_j\rangle +\langle\psi{\bf u}, P_n'P_j\rangle \\ [0.5em]
&=& \langle \phi{\bf u},(P_n'P_j)'\rangle-\langle \phi{\bf u},P_n'P_j'\rangle+\langle\psi{\bf u}, P_n'P_j\rangle
=-\langle \phi{\bf u},P_n'P_j'\rangle\,.
\end{array}
$$
Since by hypothesis C1 holds, and we have already proved that
C1$\,\Rightarrow\,$C1$'$$\,\Rightarrow\,$C2$\,\Rightarrow\,$C3, and in the proof of
C2$\,\Rightarrow\,$C3 we have shown that $\{Q_n:=P'_{n+1}/(n+1)\}_{n\geq0}$
is a monic OPS with respect to ${\bf v}:=\phi{\bf u}$,
then $\langle \phi{\bf u},P_n'P_j'\rangle=0$ if $j\neq n$,
hence (\ref{C1C5C3eq}) reduces to
\begin{equation}\label{C1C5C3eqA}
\phi P_n^{\prime\prime}+\psi P_n'+\lambda_nP_n=0\;,\quad n\geq0\;,
\end{equation}
where $\lambda_n:=-\lambda_{n,n}$.
Comparing leading coefficients in (\ref{C1C5C3eqA}), and setting
$\phi(x)=ax^2+bx+c$ and $\psi(x)=px+q$, we obtain
$\lambda_n=-n\big((n-1)a+p\big)=-nd_{n-1}$, hence $\lambda_n\neq0$ if $n\geq1$
(since C1$\,\Rightarrow\,$C1$'$, so $(\phi,\psi)$ is an admissible pair).
Thus C1$\,\Rightarrow\,$C5.
\smallskip

(C5$\,\Rightarrow\,$C1).
By hypothesis, there extist $\phi,\psi\in\mathcal{P}$, and $\lambda_n\in\mathbb{C}$, with $\lambda_n\neq0$ if $n\geq1$,
such that $-\phi P_{n+1}^{\prime\prime}=\psi P_{n+1}'+\lambda_{n+1} P_{n+1}$.
Taking in this equation $n=0$ and $n=1$ we deduce $\psi=-\lambda_1P_1\in\mathcal{P}_1\setminus\mathcal{P}_0$
and $\phi=-(\psi P_2'+\lambda_2P_2)/2\in\mathcal{P}_2$.
%By hypothesis, there exist $\phi\in\mathcal{P}_2$,
%$\psi\in\mathcal{P}_1$, and $\lambda_n\in\mathbb{C}$, with $\lambda_n\neq0$ if $n\geq1$,
%such that $-\phi P_{n+1}^{\prime\prime}=\psi P_{n+1}'+\lambda_{n+1} P_{n+1}$.
%%(Notice that at least one of the polynomials $\phi$ or $\psi$ is nonzero, since $\lambda_n\neq0$.)
We will prove that $D(\phi{\bf u})=\psi{\bf u}$ by showing that the actions of the functionals
$D(\phi{\bf u})$ and $\psi{\bf u}$ coincide on the simple set $\{Q_n\}_{n\geq0}$. %$\{Q_n:=P'_{n+1}/(n+1)\}_{n\geq0}$.
Indeed,
$$
\begin{array}{rl}
\langle D(\phi{\bf u}),Q_n\rangle =& \mbox{$\frac{1}{n+1}$}\langle D(\phi{\bf u}),P'_{n+1}\rangle
= -\mbox{$\frac{1}{n+1}$}\langle{\bf u},\phi P_{n+1}^{\prime\prime}\rangle
=\mbox{$\frac{1}{n+1}$}\langle{\bf u},\psi P'_{n+1}+\lambda_{n+1}P_{n+1}\rangle \\ [0.5em]
&  = \langle{\bf u}, \psi Q_n\rangle+
\mbox{$\frac{\lambda_{n+1}}{n+1}$}\langle{\bf u},P_{n+1}\rangle=\langle\psi{\bf u}, Q_n\rangle\;.
\end{array}
$$
Since at least one of the polynomials $\phi$ and $\psi$ is nonzero (because $\lambda_n\neq0$), C1 holds.
\smallskip

%%%%%%%%%%%%%%%%%% McCarthy %%%%%%%%%%%%%%%%%%%%%%%%%%%%%%%%%%%%%%%%%%%%%%%%%%%%%%%%%%%%%%%%%
(C2$\,\Rightarrow\,$C7).
Since by hypothesis ({\rm C2}) holds, we may write
\begin{eqnarray}
\phi P_n^{\prime} =a_nP_{n+1}+b_nP_n+c_nP_{n-1} \; , \quad \label{c.27} \\ [0.5em]
\phi P_{n-1}^{\prime} =a_{n-1}P_n+b_{n-1}P_{n-1}+c_{n-1}P_{n-2}\;.
\label{c.28}
\end{eqnarray}
Multiplying (\ref{c.27}) by $P_{n-1}$ and (\ref{c.28}) by $P_n$ and
adding the resulting equalities, we find that $\phi(P_n P_{n-1})^{\prime}$
is a linear combination of the polynomials
$P_n^2$, $P_nP_{n-1}$, $P_{n-1}^2$, $P_{n+1}P_{n-1}$ and $P_nP_{n-2}$.
Substituting $P_{n+1}$ and $P_{n-2}$ by the corresponding expressions given by the TTRR, we deduce
\begin{equation}
\phi(P_nP_{n-1})^{\prime}=A_nP_n^2+(B_nx+C_n)P_nP_{n-1}+D_nP_{n-1}^2 \;, \quad n\geq1\;,
\label{c.29}
\end{equation}
where
$$
\begin{array}{l}
A_n:=a_{n-1}-\frac{c_{n-1}}{\gamma_{n-1}} \; , \quad B_n:=a_n+ \frac{c_{n-1}}{\gamma_{n-1}} \; , \\ [0.5em]
C_n:=-a_n \beta_n +b_n+b_{n-1}-\frac{c_{n-1}}{\gamma_{n-1}}\beta_{n-1} \; , \quad
D_n:=c_n-a_n\gamma_n \; .
\end{array}
$$
Write $\phi(x)=ax^2+bx+c$ and $\psi(x)=px+q$.
We have already seen that C2$\,\Leftrightarrow\,$C1$'$, and while proving
C1$'$$\,\Rightarrow\,$C2 we have shown that the coefficients $a_n$, $b_n$, and $c_n$ appearing in
(\ref{c.27}) are given by $a_n=na$, $b_n=-\frac12\psi(\beta_n)$, and $c_n=-d_{n-1}\gamma_n$.
It follows that
%\begin{equation}\label{AnBnCnDn1}
%\begin{array}{rcl}
%A_n=(2n-3)a+p \; , \quad B_n=2a-p \; , \quad D_n=-[(2n-1)a+p]\gamma_n \; , \\ [0.5em]
%C_n = -\frac{1}{2} \big\{ (2an+p)\beta_n +2a\beta_{n-1} - [2a(n-1)+p]\beta_{n-1} \big\} -q = b-q \; ,
%\end{array}
%\end{equation}
\begin{equation}\label{AnBnCnDn1}
\begin{array}{rcl}
A_n=d_{2n-3} \; , \quad B_n=2a-p \; , \quad D_n=-d_{2n-1}\gamma_n \; , \\ [0.5em]
C_n = -\frac{1}{2} \big( d_{2n}\beta_n -d_{2n-4}\beta_{n-1} \big) -q = b-q \; ,
\end{array}
\end{equation}
where the last equality is easily derived using the expressions
for the $\beta-$parameters given in the statement of the theorem.
%\begin{equation}\label{rel-betanbn2}
%d_{2n+2}\beta_{n+1}-d_{2n-2}\beta_{n}=-2b\;,\quad n\geq0\;.
%\end{equation}
%follows from (72).
Therefore, $B_n x+C_n=\phi^{\prime}-\psi$ (independent of $n$).
%hence $B_nx+C_n$ is independent of $n$.
Finally, substituting (\ref{AnBnCnDn1}) into (\ref{c.29}) yields the equation
appearing in C7, being $h_n=A_n=d_{2n-3}$ and $t_n=D_n=-d_{2n-1}\gamma_n$ for each $n\geq1$.
Thus C2$\,\Rightarrow\,$C7.
%$$%\begin{equation}
%h_n=d_{2n-3} \; , \quad t_n=-d_{2n-1}\gamma_n \; ,\quad n\geq1 \; .
%\label{c.30}
%$$%\end{equation}
%Finally, the admissibility condition implies $h_n \neq0$ and $t_n \neq0$ for each $n\geq1$.
\smallskip

(C7$\,\Rightarrow\,$C2).
Fix an integer $n\geq1$. %$n\in\mathbb{N}$.
For this $n$, rewrite the equation in {\rm (C7)} as
$$%\begin{equation}
\big(\phi P_n^{\prime}+\psi P_n-t_nP_{n-1}\big)P_{n-1}=
\big(-\phi P_{n-1}^{\prime}+\phi' P_{n-1}+h_nP_n\big)P_{n}\;.
$$%\label{Eq1}
%\end{equation}
Therefore, since $P_n$ and $P_{n-1}$ have no common zeros,
there is %of a polynomial
$\pi_{1,n}\in\mathcal{P}_1$ such that
\begin{eqnarray}
\phi P_n^{\prime}+\psi P_n-t_nP_{n-1}=\pi_{1,n} P_n \; , \label{Eq2} \\ %[0.5em]
-\phi P_{n-1}^{\prime}+\phi' P_{n-1}+h_nP_n=\pi_{1,n}P_{n-1}\;.\label{Eq3}
\end{eqnarray}
By comparing the leading coefficients on both sides of equation (\ref{Eq2}) we deduce
$\pi_{1,n}(x)=d_nx+z_n$ for some $z_n\in\mathbb{C}$ (and $d_n:=na+p$).
By hypothesis, $(\phi,\psi)$ is an admissible pair, hence $d_n\neq0$ %(for each $n\geq1$),
and so $\deg\pi_{1,n}=1$. %for each $n\geq1$.
Moreover, by the TTRR for $\{P_n\}_{n\geq0}$, $xP_n=P_{n+1}+\beta_nP_n+\gamma_{n}P_{n-1}$.
Therefore, (\ref{Eq2}) may be rewritten as
$$%\begin{equation}
\phi P_n^{\prime} =a_nP_{n+1}+b_nP_n+c_nP_{n-1} \; , \label{Eq4}
$$%\end{equation}
where $a_n:=na$, $b_n:=na\beta_n+z_n-q$, and $c_n:=na\gamma_n+t_n$.
%$$a_n:=na\;,\quad b_n:=na\beta_n+z_n-q\;,\quad c_n:=na\gamma_n+t_n\;.$$
To conclude the proof we need to show that $c_n\neq0$ for all $n\geq1$.
Indeed, changing $n$ into $n+1$ in (\ref{Eq3}) and adding the resulting equation with (\ref{Eq2}), we obtain
$$%\begin{equation}
(\psi+\phi') P_n-t_nP_{n-1}+h_{n+1}P_{n+1} =
\big((d_n+d_{n+1})x+(z_n+z_{n+1})\big)P_n \; . \label{Eq5}
$$%\end{equation}
Since $\psi+\phi'=(2a+p)x+q+b$ and taking into account once again the TTRR for $\{P_n\}_{n\geq0}$,
the last equation may be rewritten as a trivial linear combination
of the three polynomials $P_{n+1}$, $P_{n}$, and $P_{n-1}$.
Thus, we deduce %and noting that $d_n+d_{n+1}=(2n+1)a+2p$,
%\begin{eqnarray}
%h_{n+1}=d_{2n-1} \; , \label{Eq6} \\ %[0.5em]
%z_{n+1}=-z_n-d_{2n-1}\beta_n+q+b\;,\label{Eq7} \\ %[0.5em]
%t_n=-d_{2n-1}\gamma_n\;.\label{Eq8}
%\end{eqnarray}
$$
h_{n+1}=d_{2n-1}\;,\quad
z_{n+1}=-z_n-d_{2n-1}\beta_n+q+b\;,\quad
t_n=-d_{2n-1}\gamma_n\;.
$$
Therefore, $c_n=na\gamma_n+t_n=-d_{n-1}\gamma_n\neq0$ %, and so $c_n\neq0$
(since $n\geq1$).
This completes the proof.
%%%%%%%%%%%%%%%%%%%%%%%%%%%%%%%%%%%%%%%%%%%%%%%%%%%%%%%%%%%%%%%%%%%%%%%%%%%%%%%%%%%%%%%%%%%%%%
\qed

%\begin{remark}\em
%By the proof of Theorem \ref{ThmClassicalOPS}, the
%polynomials $\phi$ and $\psi$ appearing therein
%may be taken the same in all characterizations where they appear.
%\end{remark}

\begin{remark}\em
The $\beta$ and $\gamma-$parameters in Theorem \ref{ThmClassicalOPS}
may be written explicitly in terms (only) of the coefficients of $\phi$ and $\psi$ as follows (for each $n\geq0$):
\begin{snugshade}\vspace*{-0.5em}
$$
\begin{array}{c}
\beta_n=\displaystyle-\frac{(-2a+p)q+2bn[(n-1)a+p]}{(2na+p)[(2n-2)a+p]}\;, \\ [1.25em]
\gamma_{n+1}=\displaystyle\frac{-(n+1)[(n-1)a+p][a(nb+q)^2-b(nb+q)(2na+p)+c(2na+p)^2]}{[(2n-1)a+p](2na+p)^2[2(n+1)a+p]}\;.
\end{array}
$$
\end{snugshade}
\end{remark}

%\begin{remark}\em
%Characterization C4 appears implicitly in the paper \ref{LoureiroMaroniZelia2006-C8}
%by A. Loureiro, P. Maroni, and Z. da Rocha,
%where these authors gave an alternative proof of Hahn's characterization C3.
%One should mention that the above proof of C4$\,\Rightarrow\,$C1 is implicitly given in \ref{}.
%For $k=1$, C4 was proved in \ref{PacoAmilcarPetronilho1994-C8}.
%For $k\geq2$, the above proof of C4$\,\Rightarrow\,$C1 uses arguments
%originally presented in the proof of Lemma 2.1 of \ref{MaroniZelia2001-C8}.
%(see also \ref{LoureiroMaroniZelia2006-C8}.
%presented in the proof of Lemma 3.2 of \ref{LoureiroMaroniZelia2006-C8}.
%by A. Loureiro, P. Maroni, and Z. da Rocha.
%\end{remark}

\begin{remark}\em
It is worth mentioning that the distributional approach considered
here ---developed mainly by Pascal Maroni--- simplifies dramatically the original proofs
of the characterizations of the classical OPS in Theorem \ref{ThmClassicalOPS}. %presented
The student is invited to look at some of the original proofs --- see references
\ref{Al-SalamChihara1972-C8},\ref{Bochner1929-C8},\ref{Hahn1935-C8},\ref{Hahn1937-C8},\ref{Krall1936-C8},\ref{McCarthy1961-C8}.
\end{remark}

\section{Classification and canonical representatives}

In this section we prove a remarkable property:
{\it up to constant factors and affine changes of variables,
there are only four (parametric) families of classical OP,
namely, the Hermite, Laguerre, Jacobi, and Bessel polynomials}.
%In terms of the corresponding functionals, this means that among the set of classical functionals in $\mathcal{P}'$,
%there are only four distinct equivalent classes. We will denote by
The corresponding regular functionals will be denoted by
${\bf u}_H$, ${\bf u}_L^{(\alpha)}$,
${\bf u}_J^{(\alpha,\beta)}$, and ${\bf u}_B^{(\alpha)}$ (resp.) and these will be called
the {\it canonical representatives} (or {\it canonical forms}) of the classical functionals.
%The description of these four canonical representatives is given in Table \ref{Table1}.
Their description is given in Table \ref{Table1}.
Each one of these functionals fulfils Pearson's equation $(\ref{EDClassic1})$,
being the corresponding pair $(\phi,\psi)\equiv(\Phi,\Psi)$ given in the table.
The regularity conditions in the table are determined by conditions (ii) appearing in
$(\ref{P-regular1})$. %, for each concrete pair $(\Phi,\Psi)$.

%the realization the ensure that the corresponding pair is admissible
%is an admissible one---a fact that is ensured by  for each of the four canonical forms---such that the corresponding functional fulfils Pearson's equation

\begin{table}
\centering
%\begin{tabular}{|>{\columncolor[rgb]{1,1,0}}c|c|c|c|}
\begin{tabular}{|>{\columncolor[gray]{0.95}}c|c|c|c|c|}
\hline \rowcolor[gray]{0.95}
\rule{0pt}{1.2em} Class & ${\bf u}$ & $\Phi$ & $\Psi$ & regularity conditions \\
\hline
\rule{0pt}{1.2em} Hermite & ${\bf u}_H$ & $1$& $-2x$ & ------ \\
%\hline
\rule{0pt}{1.2em} Laguerre & ${\bf u}_L^{(\alpha)}$ & $x$ & $-x+\alpha +1$ &
%$\scriptstyle{\alpha  \neq -n \ , \ n\geq1}$ \\
$-\alpha\not\in\mathbb{N}$ \\
%\hline
\rule{0pt}{1.2em} Jacobi & ${\bf u}_J^{(\alpha,\beta)}$ & $1-x^2$ & $-(\alpha + \beta +2)x +\beta-\alpha $ &
%$\scriptstyle{\alpha  \neq -n \ , \ \beta \neq -n \ , \ \alpha+\beta+1 \neq -n \ , \ n\geq1}$  \\
$-\alpha,-\beta,-(\alpha+\beta+1)\not\in\mathbb{N}$ \\
%\hline
\rule{0pt}{1.2em} Bessel & ${\bf u}_B^{(\alpha)}$ & $x^2$ & $(\alpha +2)x +2$  &
%$\scriptstyle{\alpha  \neq -n \ , \ n\geq2 }$ \\
$-(\alpha+1)\not\in\mathbb{N}$ \\
\hline
\end{tabular}
\medskip
\caption{Classification and canonical forms of the classical functionals}\label{Table1}
\end{table}

%Incluir aqui a descricao das formas canonicas (usando a classe de equivalencia) e o quadro dos parametros para essas formas canonicas.

%Dar um exemplo de uma funcional que satisfaz uma ED de Person e fazer a identificacao (por equivalencia) com a correspondete forma canonica. Por tambem um exercicio sobre isto.

%Ultimately, as a consequence of the results to be stated in this section, we will see that,
Ultimately, denoting by $[{\bf u}]$ the equivalent class
determined by a functional ${\bf u}\in\mathcal{P}'$, and setting
$\mathcal{P}'_C:=\{{\bf u}\in\mathcal{P}'\,|\,\mbox{\rm ${\bf u}$ is classical}\}$,
%(the set of all classical functionals), we will prove that
we will show that
\begin{snugshade}
$$
\mathcal{P}'_C/_\sim:=\big\{\,[{\bf u}]\,|\,{\bf u}\in\mathcal{P}'_C\big\}=
\big\{[{\bf u}_H], [{\bf u}_L^{(\alpha)}],
[{\bf u}_J^{(\alpha,\beta)}], [{\bf u}_B^{(\alpha)}]\,\big\}\;,
$$
\end{snugshade}
\noindent
where the parameters $\alpha$ and $\beta$ vary on $\mathbb{C}$ subject to the regularity conditions
in Table \ref{Table1}, and $\sim$ is the equivalence relation in $\mathcal{P}'$
introduced in Theorem \ref{Thm-u-equiv-v}, defined by
%We begin by recalling (cf.) that the binary relation $\sim$
%defined by %for every functionals ${\bf u},{\bf v}\in\mathcal{P}'$ by
\begin{snugshade}\vspace*{-0.5em}
\begin{equation}\label{u-equiv-v1}
{\bf u}\sim{\bf v}\qquad\mbox{\rm iff}\qquad
\exists A\in\mathbb{C}\setminus\{0\}\;,\;\; \exists B\in\mathbb{C}\;:\;\;
{\bf v}=\big(\bm{h}_{A^{-1}}\circ\bm{\tau}_{-B}\big){\bf u}\;.
\end{equation}
\end{snugshade}
%\noindent
%is an equivalence relation in $\mathcal{P}'$.

%The next proposition allow us to ensure that this equivalence relation
%preserves the classical character of a given classical functional.

We start by proving a proposition that allow us to ensure that this equivalence relation
preserves the classical character of a given classical functional.

%\begin{snugshade}
%\begin{lemma}\label{equiv-u-v}
%Let ${\bf u},{\bf v}\in\mathcal{P}'$ and suppose that
%${\bf u}\sim{\bf v}$, i.e.,
%there exist $A\in\mathbb{C}\setminus\{0\}$ and $B\in\mathbb{C}$ such that
%$${\bf v}=\big(\bm{h}_{A^{-1}}\circ\bm{\tau}_{-B}\big){\bf u}\;.$$
%Suppose that there exist two polynomials $\phi$ and $\psi$ such that
%$$D(\phi{\bf u})=\psi{\bf u}\;.$$
%Let $\Phi(x):=K\phi(Ax+B)$ and $\Psi(x):=KA\psi(Ax+B)$, being $K\in\mathbb{C}\setminus\{0\}$.
%Then
%$$D(\Phi{\bf v})=\Psi{\bf v}\;.$$
%\end{lemma}
%\end{snugshade}

\begin{snugshade}
\begin{lemma}\label{equiv-u-v}
Let ${\bf u},{\bf v}\in\mathcal{P}'$ and suppose that
${\bf u}\sim{\bf v}$, i.e., $(\ref{u-equiv-v1})$ holds.
Suppose that there exist two polynomials $\phi$ and $\psi$ such that
$$D(\phi{\bf u})=\psi{\bf u}\;.$$
Let $\Phi(x):=K\phi(Ax+B)$ and $\Psi(x):=KA\psi(Ax+B)$, being $K\in\mathbb{C}\setminus\{0\}$.
Then
$$D(\Phi{\bf v})=\Psi{\bf v}\;.$$
Moreover, if ${\bf u}$ is a classical functional, then so is ${\bf v}$.
\end{lemma}
\end{snugshade}

{\it Proof.}
Since ${\bf u}$ and ${\bf v}$ fulfill (\ref{u-equiv-v1}), then
$$
\langle {\bf v},x^n\rangle=\Big\langle {\bf u},\Big(\mbox{$\frac{x-B}{A}$}\Big)^n\,\Big\rangle\;,
\quad n\in\mathbb{N}_0\;.
$$
Therefore, for each $n\in\mathbb{N}_0$, we have
$$
\begin{array}{rcl}
\langle D(\Phi{\bf v}),x^n\rangle &=& -n\langle {\bf v},\Phi(x)x^{n-1}\rangle
=-n\big\langle {\bf u},\big(\tau_{B}\circ h_{A^{-1}}\big)\big(\Phi(x)x^{n-1}\big)\big\rangle \\ [0.5em]
 &=& -n\Big\langle {\bf u},\Phi\Big(\mbox{$\frac{x-B}{A}$}\Big)\Big(\mbox{$\frac{x-B}{A}$}\Big)^{n-1}\Big\rangle
=-\Big\langle{\bf u},K\phi(x)\cdot A\frac{{\rm d}}{{\rm d}x}\Big\{\Big(\mbox{$\frac{x-B}{A}$}\Big)^{n}\Big\}\Big\rangle\\ [0.75em]
 &=& KA\Big\langle D\big(\phi(x){\bf u}\big),\Big(\mbox{$\frac{x-B}{A}$}\Big)^{n}\Big\rangle
 =KA\Big\langle \psi(x) {\bf u},\Big(\mbox{$\frac{x-B}{A}$}\Big)^{n}\Big\rangle \\ [0.75em]
 &=& \Big\langle {\bf u},\Psi\Big(\mbox{$\frac{x-B}{A}$}\Big)\Big(\mbox{$\frac{x-B}{A}$}\Big)^{n}\Big\rangle
 =\big\langle {\bf u},\big(\tau_{B}\circ h_{A^{-1}}\big)\big(\Psi(x)x^{n}\big)\big\rangle \\ [0.75em]
 &=& \langle {\bf v},\Psi(x) x^n\rangle=\langle \Psi{\bf v},x^n\rangle \;.
\end{array}
$$
Finally, the last sentence stated in the lemma follows by using Theorem \ref{OPSequiv}.
%Finally, if $(\phi,\psi)$ is an admissible pair, then so is $(\Phi,\Psi)$
%%then also $\Phi\in\mathcal{P}_2$ and $\Psi\in\mathcal{P}_1$
%and, setting
%$d_n:=\psi'+n\phi^{\prime\prime}/2$, $e_n:=\psi(0)+n\phi'(0)$,
%$E_n:=\Psi'+n\Phi^{\prime\prime}/2$, and $D_n:=\Psi(0)+n\Phi'(0)$, the relations
%$$
%D_n=KA^2d_n\;,\quad \Phi(-E_n/D_{2n})=K\phi(-e_n/d_{2n})\;,\quad n\in\mathbb{N}_0\;,
%$$
%allow us to conclude that the last sentence in the statement of the lemma holds,
%taking into account the characterization of the classical functionals given by (\ref{P-regular1}).
\qed

\begin{snugshade}
\begin{theorem}[canonical representatives of the classical functionals]\label{canonic-forms-classical}
Let ${\bf u}\in\mathcal{P}^\prime$ be a classical functional,
so that ${\bf u}$ fulfils the distributional Pearson's equation
\begin{equation}\label{Pears1}
D(\phi{\bf u})=\psi{\bf u}\;,
\end{equation}
where $\phi(x)=ax^2+bx+c$ and $\psi(x)=px+q$, subject to the regularity conditions
\begin{equation}\label{Pears1reg}
na+p\neq0\;,\quad\phi\left(-\frac{nb+q}{2na+p}\right)\neq0\;,
\quad\forall n\in\mathbb{N}_0\;.
\end{equation}
%Then, up to an affine change of the variable and up to a constant factor,
%equation $(\ref{Pears1})$ can be reduced to
Then, there exists a regular functional ${\bf v}\in\mathcal{P}'$ such that
\begin{equation}\label{Pears2}
{\bf u}\sim{\bf v}\;,\quad D(\Phi{\bf v})=\Psi{\bf v}\;,
\end{equation}
where, for each classical functional determined by the pair $(\phi,\psi)$,
the corresponding pair $(\Phi,\Psi)$ is given by Table \ref{Table1}.
More precisely, setting
%$$\Delta:=b^2-4ac\;,\quad d:=\phi'(-q/p)\;,$$
%d:=\phi'\big(-q/p\big)\;,$$
%d:=\phi'\left(-\frac{q}{p}\right)\;,$$
$$\Delta:=b^2-4ac\;;\quad d:=\psi\left(-\mbox{$\frac{b}{2a}$}\right)\;\;\mbox{\rm if}\;\; a\neq0 \;,$$
%\begin{equation}\label{Pears3}
%\Delta:=b^2-4ac\;,
%\end{equation}
the following holds:
\begin{itemize}
\item[\,]\hspace*{-2em}{\rm 1.\,(Hermite)} if $a=b=0$, then:
$$%\begin{equation}\label{Pears4H}
{\bf v}=\big(\bm{h}_{\sqrt{-p/(2c)}}\circ\bm{\tau}_{q/p}\big){\bf u}={\bf u}_{{}_H}\;;
$$%\end{equation}
\item[\,]\hspace*{-2em}{\rm 2.\,(Laguerre)} if $a=0$ and $b\neq0$, then:
$$%\begin{equation}\label{Pears4L}
{\bf v}=\big(\bm{h}_{-p/b}\circ\bm{\tau}_{c/b}\big){\bf u}={\bf u}_{{}_L}^{(\alpha)}\;,\quad\alpha:=-1+(qb-pc)/b^2\;;
$$%\end{equation}
\item[\,]\hspace*{-2em}{\rm 3.\,(Bessel)} if $a\neq0$ and $\Delta=0$, then:
%setting $d:=\psi\left(-\mbox{$\frac{b}{2a}$}\right)$, then:
$$%\begin{equation}\label{Pears4B}
%{\bf v}=\Big(h_{4a^2/(2aq-pb)}\circ\tau_{b/(2a)}\Big){\bf u}\;,\quad\alpha:=-2+p/a\;;
{\bf v}=\big(\bm{h}_{2a/d}\circ\bm{\tau}_{b/(2a)}\big){\bf u}={\bf u}_{{}_B}^{(\alpha)}\;,\quad\alpha:=-2+p/a\;;
$$%\end{equation}
\item[\,]\hspace*{-2em}{\rm 4.\,(Jacobi)} if $a\neq0$ and $\Delta\neq0$, then:
%setting $d:=\psi\left(-\mbox{$\frac{b}{2a}$}\right)$, then:
$$%\begin{equation}\label{Pears4J1}
{\bf v}=\big(\bm{h}_{-2a/\sqrt{\Delta}}\circ\bm{\tau}_{b/(2a)}\big){\bf u}={\bf u}_{{}_J}^{(\alpha,\beta)}\;,
%\quad\alpha:=\beta:=-1+p/(2a)+\eta/\sqrt{\delta}\;;
$$%\end{equation}
%being
$$
%\alpha:=-1+\mbox{$\frac{p}{2a}$}\big(1-d/\sqrt{\Delta}\,\big)\;,\quad
\alpha:=-1+p/(2a)-d/\sqrt{\Delta}\;,\quad
\beta:=-1+p/(2a)+d/\sqrt{\Delta}\;.
$$
%$$
%\begin{array}{lcl}
%\alpha=\beta:=-1+p/(2a) & \mbox{\rm if} & \eta=0\; , \\ [0.5em]
%\alpha:=...\;, \quad \beta:=... & \mbox{\rm if} & \eta\neq0\; .
%\end{array}
%$$
\end{itemize}
%Moreover, all these parameters $\alpha$ and $\beta$ (appearing in cases 2, 3, and 4)
%fulfil the regularity conditions given in table \ref{Table1}.
\end{theorem}
\end{snugshade}

{\it Proof.} Taking into account Lemma \ref{equiv-u-v},
the theorem will be proved if we are able to show that,
for each given pair $(\phi,\psi)$, and for each
corresponding pair $(\Phi,\Psi)$ given by Table \ref{Table1}---where the
``corresponding pair'' $(\Phi,\Psi)$ is the one in the table such that $\phi$ and $\Phi$
have the same degree and their zeros the same multiplicity---,
there exist $A,K\in\mathbb{C}\setminus\{0\}$ and $B\in\mathbb{C}$ such that the relations
\begin{equation}\label{PhiPsiK1}
\Phi(x)=K\phi(Ax+B)\;,\quad \Psi(x)=KA\psi(Ax+B)=KA^2px+KA(Bp+q)
\end{equation}
hold, for appropriate choices of the parameters $\alpha$ and $\beta$
appearing in Table \ref{Table1} for the Laguerre, Bessel, and Jacobi cases.
Indeed, considering the four possible cases %(Hermite, Laguerre, Bessel, and Jacobi)
determined by the polynomial $\phi$, we have:
\smallskip

1. Assume $a=b=0$, i.e., $\phi(x)=c$. The regularity conditions (\ref{Pears1reg})
ensure that $p\neq0$ and $c\neq0$. Therefore, since in this case we require
$(\Phi,\Psi)=(1,-2x)$, from (\ref{PhiPsiK1}) we obtain the equations
$$
1=Kc\;,\quad -2=KA^2p\;,\quad 0=Bp+q\;.
$$
A solution of this system of equations is
$$
K=1/c\;,\quad A=\sqrt{-2c/p}\;,\quad B=-q/p\;,
$$
which gives the desired result for the Hermite case, by Lemma \ref{equiv-u-v}.
\smallskip

2. Assume $a=0$ and $b\neq0$, so that $\phi(x)=bx+c$.
Since in this case we require $(\Phi,\Psi)=(x,-x+\alpha+1)$, from (\ref{PhiPsiK1}) we obtain
$$
1=KAb\;,\quad 0=bB+c\;,\quad -1=KA^2p\;,\quad \alpha+1=KA(Bp+q)\;.
$$
Solving this system we find
$$
K=-p/b^2\;,\quad B=-c/b\;,\quad A=-b/p\;,\quad \alpha=-1+(qb-pc)/b^2\;.
$$
Notice that, in this case,
$$
d_n=p\;,\quad \phi\left(-\frac{nb+q}{2na+p}\right)%=-\frac{b^2}{p}\Big(n+\frac{qb-pc}{b^2}\Big)
=-\frac{b^2}{p}\big(n+\alpha+1\big)\;,
$$
hence the regularity conditions (\ref{Pears1reg}) ensure that $p\neq0$
(and so $K$ and $A$ are well defined, being both nonzero complex numbers) and $-\alpha\not\in\mathbb{N}$.
\smallskip

3. Assume $a\neq0$ and $\Delta=0$. Then $\phi(x)=a\Big(x+\frac{b}{2a}\Big)^2$.
In this case we require $(\Phi,\Psi)=\big(x^2,(\alpha+2)x+2\big)$, hence from (\ref{PhiPsiK1}) we obtain
$$
1=KA^2a\;,\quad 0=B+b/(2a)\;,\quad \alpha+2=KA^2p\;,\quad 2=KA(Bp+q)\;.
$$
Therefore, taking into account that
$d:=\psi\left(-\mbox{$\frac{b}{2a}$}\right)=(2aq-pb)/(2a)$, we deduce
$$
K=4a/d^2\;,\quad B=-b/(2a)\;,\quad A=d/(2a)\;,\quad \alpha=-2+p/a\;.
%K=16a^3/(2aq-pb)^2\;,\quad B=-b/(2a)\;,\quad A=(2aq-pb)/(4a^2)\;,\quad \alpha=-2+p/a\;.
$$
In this case we have
$$
d_n=a(n+\alpha+2)\;,\quad
\phi\left(-\frac{nb+q}{2na+p}\right)=\frac{d^2}{a(2n+\alpha+2)^2}\;,
%=\frac{p^2}{4a}\left(\frac{\phi'(-q/p)}{2na+p}\right)^2\;,
%\phi\left(-\frac{nb+q}{2na+p}\right)=\frac{1}{4a}\left(\frac{2aq-pb}{2na+p}\right)^2\;,
%=\frac{p^2}{4a}\left(\frac{\phi'(-q/p)}{2na+p}\right)^2\;,
$$
%$$
%\begin{array}{c}
%d_n=na+p=a(n+\alpha+2)\;,\quad \\ [1em]
%\displaystyle
%\phi\left(-\frac{nb+q}{2na+p}\right)=\frac{1}{4a}\left(\frac{2aq-pb}{2na+p}\right)^2
%%=\frac{p^2}{4a}\left(\frac{\phi'(-q/p)}{2na+p}\right)^2\;,
%\end{array}
%$$
hence conditions (\ref{Pears1reg}) ensure that $-(\alpha+1)\not\in\mathbb{N}$
and $d\neq0$, and so, in particular, $K$ is well defined, being both $K$ and $A$ nonzero complex numbers.
%Finally, notice that $2aq-pb=-p\phi'(-q/p)$ to obtain the desired result in the Bessel case.
\smallskip

4. Finally, assume $a\neq0$ and $\Delta\neq0$.
Writing $\phi(x)=a\Big[\Big(x+\frac{b}{2a}\Big)^2-\frac{\Delta}{4a^2}\Big]$, since
in this case we require $(\Phi,\Psi)=\big(1-x^2,-(\alpha+\beta+2)x+\beta-\alpha\big)$,
from (\ref{PhiPsiK1}) we obtain
$$
\begin{array}{c}
-1=KA^2a\;,\quad 0=B+b/(2a)\;,\quad 1=Ka\Big[\Big(B+\frac{b}{2a}\Big)^2-\frac{\Delta}{4a^2}\Big]\;, \\ [0.5em]
-(\alpha+\beta+2)=KA^2p\;,\quad \beta-\alpha=KA(Bp+q)\;.
\end{array}
$$
A solution of this system of five equations is\footnote{\,We choose $A$ with the minus sign
since whenever $(\phi,\psi)=(\Phi,\Psi)$ that choice implies $A=1$ and $B=0$, hence ${\bf u}={\bf v}={\bf u}_J^{(\alpha,\beta)}$, and so it is a more natural choice.}
%\footnote{\,We choose $A$ with the minus sign since this choice implies ${\bf v}={\bf u}$ whenever $(\phi,\psi)=(\Phi,\Psi)$.}
$$
\begin{array}{c}
K=-4a/\Delta\;,\quad B=-b/(2a)\;,\quad A=-\sqrt{\Delta}/(2a)\;, \\ [0.25em]
\alpha=-1+p/(2a)-d/\sqrt{\Delta}\;,\quad \beta=-1+p/(2a)+d/\sqrt{\Delta}\;.
\end{array}
$$
Adding and subtracting the last equations for $\alpha$ and $\beta$,
we find $\alpha+\beta+2=p/a$ and $\alpha-\beta=-2d/\sqrt{\Delta}$, hence we deduce
$$
d_n=a(n+\alpha+\beta+2)\;,\quad
\phi\left(-\frac{nb+q}{2na+p}\right)=-\frac{\Delta}{a}\frac{(n+\alpha+1)(n+\beta+1)}{(2n+\alpha+\beta+2)^2}\;,
%=\frac{p^2}{4a}\left(\frac{\phi'(-q/p)}{2na+p}\right)^2\;,
%\phi\left(-\frac{nb+q}{2na+p}\right)=\frac{1}{4a}\left(\frac{2aq-pb}{2na+p}\right)^2\;,
%=\frac{p^2}{4a}\left(\frac{\phi'(-q/p)}{2na+p}\right)^2\;,
$$
%$$
%\begin{array}{c}
%d_n=na+p=a(n+\alpha+2)\;,\quad \\ [1em]
%\displaystyle
%\phi\left(-\frac{nb+q}{2na+p}\right)=\frac{1}{4a}\left(\frac{2aq-pb}{2na+p}\right)^2
%%=\frac{p^2}{4a}\left(\frac{\phi'(-q/p)}{2na+p}\right)^2\;,
%\end{array}
%$$
Therefore, conditions (\ref{Pears1reg}) ensure that $-(\alpha+\beta+1)\not\in\mathbb{N}$,
$-\alpha\not\in\mathbb{N}$, and $-\beta\not\in\mathbb{N}$. This completes the proof.
\qed

\begin{remark}\em
It follows from the proof of Theorem \ref{canonic-forms-classical} that
the parameters $\alpha$ and $\beta$ defined in the statement of this theorem (in cases 2, 3, and 4)
fulfil the regularity conditions appearing in Table \ref{Table1}.
\end{remark}

The preceding theorem allows us to classify each classical functional according with
the degree of the polynomial $\phi$ appearing in Pearson's equation $(\ref{EDClassic1})$.

%\begin{snugshade}
%\begin{corollary}
%Let ${\bf u}$ be a classical functional, fulfilling $(\ref{EDClassic1})$--$(\ref{grauPhiPsi})$.
%\begin{enumerate}
%\item[{\rm (i)}] if $\deg\phi=0$ (hence $\phi$ is a nonzero constant), then ${\bf u}\sim$ is called a {\sl Hermite functional};
%\item[{\rm (ii)}] if $\deg\phi=1$, ${\bf u}$ is called a {\sl Laguerre functional};
%\item[{\rm (iii)}] if $\deg\phi=2$ and $\phi$ has simple zeros, ${\bf u}$ is called a {\sl Jacobi functional};
%\item[{\rm (iv)}] if $\deg\phi=2$ and $\phi$ has a double zero, ${\bf u}$ is called a {\sl Bessel functional}.
%\end{enumerate}
%\end{corollary}
%\end{snugshade}

\begin{snugshade}
\begin{corollary}
Let ${\bf u}$ be a classical functional, fulfilling $(\ref{EDClassic1})$--$(\ref{grauPhiPsi})$.
\begin{enumerate}
\item[{\rm (i)}] if $\deg\phi=0$ (hence $\phi$ is a nonzero constant), then ${\bf u}\sim{\bf u}_H\,$;
\item[{\rm (ii)}] if $\deg\phi=1$, then ${\bf u}\sim{\bf u}_L^{(\alpha)}\,$ for some $\alpha$;
\item[{\rm (iii)}] if $\deg\phi=2$ and $\phi$ has simple zeros, then ${\bf u}\sim{\bf u}_J^{(\alpha,\beta)}\,$ for some pair $(\alpha,\beta)$;
\item[{\rm (iv)}] if $\deg\phi=2$ and $\phi$ has a double zero, then ${\bf u}\sim{\bf u}_B^{(\alpha)}\,$ for some $\alpha$.
\end{enumerate}
\end{corollary}
\end{snugshade}
%\noindent
%This classification relies upon the following facts.

%Finally, we point out that, denoting by $[{\bf u}]$ the equivalent class
%determined by a classical functional ${\bf u}\in\mathcal{P}'$, and setting
%$\mathcal{P}'_C:=\{{\bf u}\in\mathcal{P}'\,|\,\mbox{\rm ${\bf u}$ is classical}\}$
%(the set of all classical functionals), then the above results may be expressed as
%\begin{snugshade}
%$$
%\mathcal{P}'_C/_\sim:=\big\{\,[{\bf u}]\,|\,{\bf u}\in\mathcal{P}'_C\big\}=
%\big\{[{\bf u}_H], [{\bf u}_L^{(\alpha)}],
%[{\bf u}_J^{(\alpha,\beta)}], [{\bf u}_B^{(\alpha)}]\,\big\}\;.
%$$
%\end{snugshade}

%%%%%%%%%%%%%%%%%%%%%%%%

\begin{remark}\em
The monic OPS with respect to the canonical representatives ${\bf u}_H$, ${\bf u}_L^{(\alpha)}$,
${\bf u}_J^{(\alpha,\beta)}$, and ${\bf u}_B^{(\alpha)}$ will be denoted by
$\{\widehat{H}_n\}$, $\{\widehat{L}_n^{(\alpha)}\}$, $\{\widehat{P}_n^{(\alpha,\beta)}\}$,
and $\{\widehat{B}_n^{(\alpha)}\}$ (resp.), and they will be called the (monic)
{\sl Hermite}, {\sl Laguerre}, {\sl Jacobi}, and {\sl Bessel} polynomials.
Table \ref{Table3} summarizes the corresponding parameters appearing in all
characterizations presented in Theorem \ref{ThmClassicalOPS}.
\end{remark}

\begin{remark}\em
%for these classical monic OPS with respect to the canonical representatives given in Table \ref{Table1}.
Note that, in view of Theorem \ref{canonic-forms-classical} and Theorem \ref{OPSequiv},
we may now justify a sentence made at the beginning of the section:
{\it up to constant factors and affine changes of variables,
the Hermite, Laguerre, Jacobi, and Bessel polynomials are the only families of classical OP}.
\end{remark}

\begin{remark}\em
Notice also the following special cases of Jacobi polynomials (up to normalization),
that we have introduced in some previous texts:
%\smallskip
\begin{snugshade}
$\alpha=\beta=0$, {\sl Legendre} polynomials;

$\alpha=\beta=-\frac12$, {\sl Chebyshev polynomials of the first kind}: $\{T_n\}_{n\geq0}$;

$\alpha=\beta=\frac12$, {\sl Chebyshev polynomials of the second kind}: $\{U_n\}_{n\geq0}$;

$\alpha=\beta=:\lambda-\frac12$, {\sl Gegenbauer}
(or {\sl ultraspherical}) polynomials: $\{C_n^\lambda\}_{n\geq0}$  ($-2\lambda\not\in\mathbb{N}$).
\end{snugshade}
%\begin{enumerate}
%\item[$\circ$] for $\alpha=\beta=0$ one obtains {\sl Legendre} polynomials;
%\item[$\circ$] for $\alpha=\beta=\frac12$ one obtains {\sl Chebyshev polynomials of the first kind}: $\{T_n\}_{n\geq0}$;
%\item[$\circ$] for $\alpha=\beta=\frac12$ one obtains {\sl Chebyshev polynomials of the second kind}: $\{U_n\}_{n\geq0}$;
%\item[$\circ$] for $\alpha=\beta=:\lambda-\frac12$ one obtains {\sl Gegenbauer} (or {\sl ultraspherical}) polynomials: $\{C_n^\lambda\}_{n\geq0}$, $\lambda>-\frac12$.
%\end{enumerate}
\end{remark}

\begin{remark}\em
%It is worth mentioning that the
The Legendre polynomials were the first discovered
OPS, and they have been introduced by
the French mathematician Adrien Marie Legendre (1752-1833) in a work published in 1785
entitled ``Sur l'attraction des sph\'eroides''.
Chebyshev polynomials were introduced by the Russian mathematician
Pafnuti Lvovich Chebychev (1821-1894), and
Jacobi polynomials by the German mathematician Carl Gustav Jacob Jacobi (1804-1851)
in a work published in 1859 about the so-called {\it hypergeometric functions},
which are solutions of the ordinary differential equation
$$
x(1-x)\,y^{\prime\prime}+\big(\gamma-(\alpha+\beta+1)x\big)\,y'-\alpha\beta\,y=0\;,
$$
proposed by the German mathematician Carl Friedrich Gauss (1777-1855).
%in a famous work on hypergeometric functions.
\end{remark}

\begin{table}
\centering
%\begin{tabular}{|>{\columncolor[rgb]{1,1,0}}c|c|c|c|}
\begin{tabular}{|>{\columncolor[gray]{0.95}}c|c|c|c|c|}
\hline \rowcolor[gray]{0.95}
\rule{0pt}{1.2em} & $\widehat{H}_n$ & $\widehat{L}_n^{(\alpha)}$ & $\widehat{P}_n^{(\alpha,\beta)}$ & $\widehat{B}_n^{(\alpha)}$\\
\hline
\rule{0pt}{1.2em} $\lambda_n$ & ${\scriptstyle{2n}}$ & ${\scriptstyle{n}}$ & ${\scriptstyle{n(n+\alpha+\beta+1)}}$ & ${\scriptstyle{-n(n+\alpha+1)}}$ \\
%\hline
\rule{0pt}{1.5em} $\beta_n$ & ${\scriptstyle{0}}$ & ${\scriptstyle{2n+\alpha+1}}$ &
    $\frac{\beta^2-\alpha^2}{(2n+\alpha+\beta)(2n+2+\alpha+\beta)}$ &
    $\frac{-2\alpha}{(2n+\alpha)(2n+2+\alpha)}$  \\
%\hline
\rule{0pt}{1.75em} $\gamma_n$ & $\frac{n}{2}$ & ${\scriptstyle{n(n+\alpha)}}$ &
    $\frac{4n(n+\alpha)(n+\beta)(n+\alpha+\beta)}{(2n+\alpha+\beta-1)
    (2n+\alpha+\beta)^2(2n+\alpha+\beta+1)}$ &
    $\frac{-4n(n+\alpha)}{(2n+\alpha-1)(2n+\alpha)^2(2n+\alpha+1)}$ \\
%\hline
\rule{0pt}{1em} $a_n$ & ${\scriptstyle{0}}$ & ${\scriptstyle{0}}$ & ${\scriptstyle{-n}}$ & ${\scriptstyle{n}}$ \\
%\hline
\rule{0pt}{1.25em} $b_n$ & ${\scriptstyle{0}}$ & ${\scriptstyle{n}}$ & $\frac{2(\alpha-\beta)n(n+\alpha+\beta+1)}{(2n+\alpha+\beta)(2n+2+\alpha+\beta)}$ & $\frac{-4n(n+\alpha+1)}{(2n+\alpha)(2n+2+\alpha)}$ \\
%\hline
\rule{0pt}{1.75em} $c_n$ & ${\scriptstyle{n}}$ & ${\scriptstyle{n(n+\alpha)}}$ & $\frac{4n(n+\alpha)(n+\beta)(n+\alpha+\beta)(n+\alpha+\beta+1)}{(2n+\alpha+\beta-1)
    (2n+\alpha+\beta)^2(2n+\alpha+\beta+1)}$ & $\frac{4n(n+\alpha)(n+\alpha+1)}{(2n+\alpha-1)(2n+\alpha)^2(2n+\alpha+1)}$ \\
%\hline
\rule{0pt}{1.5em} $r_n^{[1]}$ & ${\scriptstyle{0}}$ & ${\scriptstyle{n}}$ & $\frac{2(\alpha-\beta)n}{(2n+\alpha+\beta)(2n+2+\alpha+\beta)}$ & $\frac{4n}{(2n+\alpha)(2n+2+\alpha)} $ \\
%\hline
\rule{0pt}{1.75em} $s_n^{[1]}$ & ${\scriptstyle{0}}$ & ${\scriptstyle{0}}$ & $\frac{-4(n-1)n(n+\alpha)(n+\beta)}{(2n+\alpha+\beta-1)(2n+\alpha+\beta)^2(2n+\alpha+\beta+1)}$ & $\frac{4(n-1)n}{(2n+\alpha-1)(2n+\alpha)^2(2n+\alpha+1)}$ \\
%\hline
\rule{0pt}{1.2em} $h_n$ & ${\scriptstyle{-2}}$ & ${\scriptstyle{-1}}$ & ${\scriptstyle{-(2n+\alpha+\beta-1)}} $ & ${\scriptstyle{2n+\alpha-1}} $ \\
%\hline
\rule{0pt}{1.5em} $t_n$ & ${\scriptstyle{n}}$ & ${\scriptstyle{n(n+\alpha)}}$ & $\frac{4n(n+\alpha)(n+\beta)(n+\alpha+\beta)}{(2n+\alpha+\beta-1)(2n+\alpha+\beta)^2}$ & $\frac{4n(n+\alpha)}{(2n+\alpha-1)(2n+\alpha)^2}$ \\
%\hline
\rule{0pt}{1.5em} $k_n$ & $\frac{(-1)^n}{2^n}$ & ${\scriptstyle{(-1)^n}}$ & $\frac{(-1)^n}{(n+\alpha+\beta+1)_n}$ & $\frac{1}{(n+\alpha+1)_n}$  \\ [0.5em]
\hline
\end{tabular}
\medskip
\caption{Parameters for the classical monic OPS
appearing in Theorem \ref{ThmClassicalOPS} with respect to the canonical forms %of $\phi$ and $\psi$
given in Table \ref{Table1}.}\label{Table3}
\end{table}

%\begin{table}
%\centering
%%\begin{tabular}{|>{\columncolor[rgb]{1,1,0}}c|c|c|c|}
%\begin{tabular}{|>{\columncolor[gray]{0.95}}c|c|c|c|}
%\hline \rowcolor[gray]{0.95}
%\rule{0pt}{1.2em} $P_n$ & $\Phi$ & $\Psi$ & regularity conditions \\
%\hline
%\rule{0pt}{1.2em} $H_n$ & $1$& $-2x$ &  \\
%%\hline
%\rule{0pt}{1.2em} $L_n^{(\alpha)}$ & $x$ & $-x+\alpha +1$ &
%$\scriptstyle{\alpha  \neq -n \ , \ n\geq1}$ \\
%%$-\alpha\not\in\mathbb{N}$ \\
%%\hline
%\rule{0pt}{1.2em} $P_n^{(\alpha,\beta)}$ & $1-x^2$ & $-(\alpha + \beta +2)x +\beta-\alpha $ &
%$\scriptstyle{\alpha  \neq -n \ , \ \beta \neq -n \ , \ \alpha+\beta+1 \neq -n \ , \ n\geq1}$  \\
%%$-\alpha,-\beta,-(\alpha+\beta+1)\not\in\mathbb{N}$ \\
%%\hline
%\rule{0pt}{1.2em} $B_n^{(\alpha)}$ & $x^2$ & $(\alpha +2)x +2$  &
%$\scriptstyle{\alpha  \neq -n \ , \ n\geq2 }$ \\
%%$-(\alpha+1)\not\in\mathbb{N}$ \\
%\hline
%\end{tabular}
%\medskip
%\caption{Classification of the classical OPS.}\label{Table1}
%\end{table}

%dar as expressoes dos betas e gamas obtidos a partir da transformacao afim com A e B (e K?)

\section{The positive-definite case}

In this section we analyze the classical functionals in the positive-definite case.
To be more precise, we will find the conditions ensuring that the classical functionals
are positive-definite, and then, under such conditions, we will show that
these functionals may be represented uniquely by simple weight functions
(via proper or improper Riemann integrals).
Of course, up to affine changes of the variables,
we only need to analyze the positive-definiteness of the canonical forms described in Table \ref{Table1}.
We begin by stating the following elementary proposition.

\begin{snugshade}
\begin{lemma}\label{lemmaW}
Let $(\xi,\eta)$ be a bounded or unbounded interval of real numbers.
Let $\omega:(\xi,\eta)\to\mathbb{R}$ be a function fulfilling the following four properties:
\begin{itemize}
\item[{\rm (i)}]
$\omega\in\mathcal{C}^1(\xi,\eta)$ and $\omega(x)>0$ for each $x\in(\xi,\eta)$;
\item[{\rm (ii)}]
$\int_\xi^\eta|x|^k\omega(x)\,{\rm d}x<\infty$ for each $k\in\mathbb{N}_0$;
\item[{\rm (iii)}]
there exist real polynomials $\phi$ and $\psi$
%$\phi\in\mathcal{P}_2$ and $\psi\in\mathcal{P}_1\setminus\mathcal{P}_0$
such that $\omega$ fulfils the first order ODE
\begin{equation}\label{ODE-W1}
\big(\phi\omega\big)'=\psi\omega\quad\mbox{on}\;\;(\xi,\eta)\;;
\end{equation}
\item[{\rm (iv)}]
$\displaystyle\lim_{x\to\xi^+}x^k\phi(x)\omega(x)=\lim_{x\to\eta^-}x^k\phi(x)\omega(x)=0$ for each $k\in\mathbb{N}_0$.
%\item[{\rm (iv)}]
%$x^k\phi(x)\omega(x)\big|_{\xi}^{\eta}=0$ for each $k\in\mathbb{N}_0$
%(taking the limit if $\xi=-\infty$ or $\eta=+\infty$).
\end{itemize}
Define a functional ${\bf u}$ on $\mathcal{P}$ by %${\bf u}:\mathcal{P}\to\mathbb{R}$ by
\begin{equation}\label{uW1}
\langle{\bf u},p\rangle:=\int_\xi^\eta p(x)\omega(x)\,{\rm d}x\;,\quad p\in\mathcal{P}\;.
\end{equation}
Then ${\bf u}$ is a positive-definite functional on $[\xi,\eta]$,
and it fulfils the generalized Pearson's distributional differential equation
\begin{equation}\label{PearsonW1}
D(\phi{\bf u})=\psi{\bf u}\;.
\end{equation}
\end{lemma}
\end{snugshade}

{\it Proof.}
Hypothesis (i)--(ii) ensure that ${\bf u}$ is well defined.
%having finite moments of all orders.
Take arbitrarily $p\in\mathcal{P}$ such that $p(x)\geq0$ on $[\xi,\eta]$ and $p(x)\not\equiv0$.
Since $p$ is continuous on $[\xi,\eta]$ and not vanishes identically there, then
there exist $x_0\in(\xi,\eta)$ and $\delta>0$ so that
$(x_0-\delta,x_0+\delta)\subset(\xi,\eta)$ and $p(x)>\epsilon:=p(x_0)/2>0$ for each $x\in(x_0-\delta,x_0+\delta)$.
Hence
$$
\langle{\bf u},p\rangle=\int_\xi^\eta p(x)\omega(x)\,{\rm d}x\geq \epsilon
\int_{x_0-\delta}^{x_0+\delta}\omega(x)\,{\rm d}x>0\;,
$$
where the last equality follows from hypothesis (i).
Thus ${\bf u}$ is positive-definite on $[\xi,\eta]$.
To prove that ${\bf u}$ satisfies (\ref{PearsonW1}), take $p\in\mathcal{P}$. Then
$$
\begin{array}{rcl}
\langle D(\phi{\bf u}),p\rangle &=& \displaystyle
-\langle{\bf u},\phi p'\rangle=-\int_{\xi}^\eta\phi p'\omega\,{\rm d}x=
-\int_{\xi}^\eta\big[(\phi\omega p)'-(\phi\omega)'p\big]\,{\rm d}x \\ [1em]
&=& \displaystyle
\phi(x)\omega(x)p(x)\big|_{\xi}^{\eta}+\int_{\xi}^\eta\psi\omega p\,{\rm d}x
=\langle\psi{\bf u},p\rangle\;,
\end{array}
$$
%$$
%\begin{array}{rcl}
%\langle D(\phi{\bf u}),p\rangle &=& \displaystyle
%-\langle{\bf u},\phi p'\rangle=-\int_{\xi}^\eta\phi(x) p'(x)\omega(x)\,{\rm d}x=
%-\int_{\xi}^\eta\big[(\phi\omega p)'(x)-(\phi\omega)'(x)p(x)\big]\,{\rm d}x \\ [1em]
%&=& \displaystyle
%\phi(x)\omega(x)p(x)\big|_{\xi}^{\eta}+\int_{\xi}^\eta\psi(x)\omega(x) p(x)\,{\rm d}x
%=\langle\psi{\bf u},p\rangle\;,
%\end{array}
%$$
where we have used (iii) in the fourth equality, and (ii) and (iv) in the last one.
\qed

\begin{remark}\em
Under the conditions of Lemma \ref{lemmaW},
we say that $\omega$ is a {\sl weight function} for ${\bf u}$,
and that ${\bf u}$ is represented by the weight function $\omega$;
and we also say that the OPS with respect to ${\bf u}$
is orthogonal with respect to the weight function $\omega$.
\end{remark}

\subsection{Hermite functional}\label{repIntClass}
By Table \ref{Table3}, the coefficients appearing in the TTRR for the monic OPS
%$\{\widehat{H}_n\}_{n\geq0}$
with respect to the (canonical) Hermite functional, ${\bf u}_H$, satisfy
$$
\beta_n=0\in\mathbb{R}\;,\quad \forall n\in\mathbb{N}_0\;;\qquad
\gamma_n=\mbox{$\frac{n}{2}$}>0\;,\quad \forall n\in\mathbb{N}\;.
$$
Therefore, by Favard's Theorem, ${\bf u}_H$ is positive-definite.
Next we show that ${\bf u}_H$ is represented by a weight function,
in the sense of (\ref{uW1}).
First, we guess the polynomials $\phi$ and $\psi$ from Table \ref{Table1}, so that
$\phi(x)\equiv1$ and $\psi(x)=-2x$. This gives the ODE
$$
\omega'=-2x\omega\,.
$$
The general solution of this equation is $Ce^{-x^2}$, where $C$ is an arbitrary real constant.
Thus we choose the weight function
\begin{snugshade}
\begin{equation}\label{weightH}
\omega(x):=e^{-x^2}\;,\quad x\in\mathbb{R}\;.
\end{equation}
\end{snugshade}
\noindent
Notice that it is quite natural to take $(\xi,\eta):=\mathbb{R}$,
since this is the largest interval where $\omega$ becomes positive, as required on hypothesis (i)
appearing in Lemma \ref{lemmaW}.
%(As a matter of fact, this is the largest interval where $\omega$ is positive.)
Of course, by construction, (iii) is also fulfilled.
Moreover, one immediately sees that $\omega$ satisfies the remaining hypothesis (ii) and (iv).
Thus, by Lemma \ref{lemmaW} and Theorem \ref{canonic-forms-classical},
\begin{snugshade}
\begin{equation}\label{uWH}
\langle{\bf u}_H,p\rangle=\int_{-\infty}^{+\infty} p(x) e^{-x^2}\,{\rm d}x\;,\quad p\in\mathcal{P}\;,
\end{equation}
\end{snugshade}
\noindent
meaning that ${\bf u}_H$ is represented by the weight function (\ref{weightH}).
The corresponding positive Borel measure is supported on $\mathbb{R}$,
and the associated distribution function
$\psi_H:\mathbb{R}\to\mathbb{R}$ is given by
\begin{equation}\label{uWHF}
\psi_H(x):=\int_{-\infty}^{x} e^{-t^2}\,{\rm d}t\;,\quad x\in\mathbb{R}\;.
\end{equation}
Notice that ${\bf u}_H$ is uniquely determined by $\psi_H$.
Indeed, this follows Corollary \ref{mu-uniq-CorRieszCriterium},
by choosing there $\theta=2$ and hence noticing that
$$
\int_{-\infty}^{+\infty}e^{\theta|x|}\,{\rm d}\psi_H(x)=
2\int_{0}^{+\infty}e^{2 x-x^2}\,{\rm d}x<\infty\;.
$$
Thus, in accordance with the spectral Theorem \ref{SpecThmOP}, we conclude that
the Hermite polynomials are orthogonal in the positive-definite
sense with respect to a unique positive Borel measure supported on
$\mathbb{R}$, and characterized by the distribution function (\ref{uWHF}).
The reader would recognizes here, up to normalization,
the Gaussian (or normal) probability distribution function.

%, provided we choose
%$$
%\xi=-\infty\;,\quad \eta=+\infty\;,\quad \phi(x)\equiv1\;,\quad \psi(x)=-2x\;.
%$$

\subsection{Laguerre functional}
By Table \ref{Table3}, the coefficients appearing in the TTRR for the monic OPS
%$\{\widehat{H}_n\}_{n\geq0}$
with respect to the Laguerre functional, ${\bf u}_L^{(\alpha)}$, satisfy
$$
\beta_n\in\mathbb{R}\;,\; \forall n\in\mathbb{N}_0\;\;\Leftrightarrow\;\; \alpha\in\mathbb{R}\;;\qquad
\gamma_n>0\;,\; \forall n\in\mathbb{N}\;\;\Leftrightarrow\;\;\alpha>-1\;.
$$
Therefore, ${\bf u}_L^{(\alpha)}$ is positive-definite if and only if $\alpha>-1$.
To show that ${\bf u}_L^{(\alpha)}$ is represented by a weight function (if $\alpha>-1$),
consider the corresponding polynomials $\phi(x)= x$ and $\psi(x)=-x+\alpha+1$,
given by Table \ref{Table1}. This gives the ODE
$$
(x\omega)'=(-x+\alpha+1)\omega\,.
$$
The general solution of this equation is $Cx^{\alpha}e^{-x}$, $C\in\mathbb{R}$.
Thus we choose
\begin{snugshade}
\begin{equation}\label{weightL}
\omega(x):=x^{\alpha}e^{-x}\;,\quad x\in(0,+\infty)\;.
\end{equation}
\end{snugshade}
\noindent
As before, we chose $(\xi,\eta):=(0,+\infty)$
since this is the largest interval where $\omega$ becomes positive.
Thus, hypothesis (i) and (iii) in Lemma \ref{lemmaW} are fulfilled.
Moreover, since $\alpha>-1$, we have, for each $k\in\mathbb{N}_0$,
$$
\int_0^{+\infty}|x|^kw(x)\,{\rm d}x=
\int_0^1x^{\alpha+k}e^{-x}\,{\rm d}x+\int_1^{+\infty}x^{\alpha+k}e^{-x}\,{\rm d}x<\infty\;.
$$
The last two integrals converge. Indeed, on one hand,
$\int_0^1x^{\alpha+k}e^{-x}\,{\rm d}x\leq\int_0^1x^{\alpha+k}\,{\rm d}x<\infty$
(because $\alpha+k>-1$ for each $k\in\mathbb{N}_0$); on the other hand,
$\int_1^{+\infty}x^{\alpha+k}e^{-x}\,{\rm d}x<\infty$,
since $\int_1^{+\infty}\frac{1}{x^s}\,{\rm d}x<\infty$ for an arbitrarily fixed $s>1$,
and %, we can use a comparison test for convergent improper integrals
$$
\frac{x^{\alpha+k}e^{-x}}{\frac{1}{x^s}}=x^{\alpha+k+s}e^{-x}\to0\quad(x\to+\infty)\;.
$$
Thus, $\omega$ satisfies hypothesis (ii). Of course, $\omega$ also satisfies (iv).
Thus, by Lemma \ref{lemmaW} and Theorem \ref{canonic-forms-classical},
\vspace*{-0.5em}
\begin{snugshade}
\begin{equation}\label{uWL}
\langle{\bf u}_L^{(\alpha)},p\rangle=\int_{0}^{+\infty} p(x) x^{\alpha}e^{-x}\,{\rm d}x\;,\quad p\in\mathcal{P}\;,
\end{equation}
\end{snugshade}
\noindent
so ${\bf u}_L^{(\alpha)}$ is represented by the weight function (\ref{weightL}).
The corresponding positive Borel measure is supported on the closed interval $[0,+\infty)$,
and the associated distribution function
$\psi_L^{(\alpha)}:\mathbb{R}\to\mathbb{R}$ is given by
\begin{equation}\label{uWLF}
\psi_L^{(\alpha)}(x):=\int_{-\infty}^{x} t^{\alpha}e^{-t}\chi_{(0,+\infty)}(t)\,{\rm d}t\;,\quad x\in\mathbb{R}\;.
\end{equation}
Notice that ${\bf u}_L^{(\alpha)}$ is uniquely determined by $\psi_L^{(\alpha)}$.
This follows e.g. from Corollary \ref{mu-uniq-CorRieszCriterium},
by choosing there $0<\theta<1$ and hence noticing that
$$
\int_{-\infty}^{+\infty}e^{\theta|x|}\,{\rm d}\psi_L^{(\alpha)}(x)=
\int_0^1x^{\alpha}e^{(\theta-1)x}\,{\rm d}x+\int_1^{+\infty}x^{\alpha}e^{(\theta-1)x}\,{\rm d}x<\infty\;.
$$
Thus, in accordance with the spectral Theorem \ref{SpecThmOP}, we conclude that
if $\alpha>-1$ the Laguerre polynomials are orthogonal in the positive-definite
sense with respect to a unique positive Borel measure supported on
$[0,+\infty)$, and characterized by the distribution function (\ref{uWLF}).
The reader would recognizes here, up to normalization,
the gamma probability distribution function.

\subsection{Jacobi functional}
By Table \ref{Table3}, the coefficients appearing in the TTRR for the monic OPS
%$\{\widehat{H}_n\}_{n\geq0}$
with respect to the Jacobi functional, ${\bf u}_J^{(\alpha,\beta)}$, satisfy
$$
\beta_n\in\mathbb{R}\;,\; \forall n\in\mathbb{N}_0\;\;\Leftrightarrow\;\; \alpha,\beta\in\mathbb{R}\;;\qquad
\gamma_n>0\;,\; \forall n\in\mathbb{N}\;\;\Leftrightarrow\;\;\alpha,\beta>-1\;.
$$
Therefore, ${\bf u}_J^{(\alpha,\beta)}$ is positive-definite if and only if $\alpha>-1$ and $\beta>-1$.
To show that ${\bf u}_J^{(\alpha,\beta)}$ is represented by a weight function (if $\alpha>-1$ and $\beta>-1$),
consider the corresponding polynomials $\phi(x)= 1-x^2$ and $\psi(x)=-(\alpha+\beta+2)x+\beta-\alpha$,
given by Table \ref{Table1}. This gives the ODE
$$
\big((1-x^2)\omega\big)'=\big(-(\alpha+\beta+2)x+\beta-\alpha\big)\omega\,.
$$
The general solution of this equation is $C(1-x)^{\alpha}(1+x)^\beta$, $C\in\mathbb{R}$.
Thus we choose
\vspace*{-0.5em}
\begin{snugshade}
\begin{equation}\label{weightJ}
\omega(x):=(1-x)^{\alpha}(1+x)^\beta\;,\quad x\in(-1,1)\;.
\end{equation}
\end{snugshade}
\noindent
Clearly, hypothesis (i), (iii) and (iv) appearing in Lemma \ref{lemmaW} are fulfilled, where $(\xi,\eta):=(-1,1)$.
Moreover, since $\alpha>-1$ and $\beta>-1$, we have, for each $k\in\mathbb{N}_0$,
%$$
%\begin{array}{rcl}
%\int_{-1}^1|x|^kw(x)\,{\rm d}x &=&
%\int_{-1}^0|x|^k(1-x)^{\alpha}(1+x)^\beta\,{\rm d}x+\int_0^1|x|^k(1-x)^{\alpha}(1+x)^\beta\,{\rm d}x \\ [0.5em]
%&\leq& 2^{\alpha}\int_{-1}^0(1+x)^\beta\,{\rm d}x+2^\beta\int_0^1(1-x)^{\alpha}\,{\rm d}x<\infty\;,
%\end{array}
%$$
$$
\int_{-1}^1|x|^kw(x)\,{\rm d}x
\leq 2^{\alpha}\int_{-1}^0(1+x)^\beta\,{\rm d}x+2^\beta\int_0^1(1-x)^{\alpha}\,{\rm d}x<\infty\;,
$$
and so $\omega$ satisfies hypothesis (ii).
Thus, by Lemma \ref{lemmaW} and Theorem \ref{canonic-forms-classical},
\begin{snugshade}
\begin{equation}\label{uWJ}
\langle{\bf u}_J^{(\alpha,\beta)},p\rangle=\int_{-1}^{1} p(x) (1-x)^{\alpha}(1+x)^\beta\,{\rm d}x\;,\quad p\in\mathcal{P}\;,
\end{equation}
\end{snugshade}
\noindent
hence ${\bf u}_J^{(\alpha,\beta)}$ is represented by the weight function (\ref{weightJ}).
The corresponding positive Borel measure is supported on the closed interval $[-1,1]$,
and the associated distribution function
$\psi_J^{(\alpha,\beta)}:\mathbb{R}\to\mathbb{R}$ is given by
\begin{equation}\label{uWJF}
\psi_J^{(\alpha,\beta)}(x):=\int_{-\infty}^{x} (1-t)^{\alpha}(1+t)^\beta\chi_{(-1,1)}(t)\,{\rm d}t\;,\quad x\in\mathbb{R}\;.
\end{equation}
Notice that, since the sequences $\{\beta_n\}_{n\geq0}$ and $\{\gamma_n\}_{n\geq1}$ are bounded,
then, by Theorem \ref{mu-uniq-bnd}, ${\bf u}_J^{(\alpha,\beta)}$ is uniquely determined by $\psi_J^{(\alpha,\beta)}$.
Therefore, we conclude that if $\alpha>-1$ and $\beta>-1$ then
the Jacobi polynomials are orthogonal in the positive-definite
sense with respect to a unique positive Borel measure supported on
$[-1,1]$, and characterized by the distribution function (\ref{uWJF}).
The reader would recognizes here, up to normalization,
the beta probability distribution function.

\subsection{Bessel functional}
Consider the coefficients given by Table \ref{Table3} for the TTRR of the monic OPS
with respect to the Bessel functional, ${\bf u}_B^{(\alpha)}$.
%We have $\beta_0=-2/(2+\alpha)$, hence
We see that the condition $\alpha\in\mathbb{R}\setminus\{-2,-3,-4,\cdots\}$ is
necessary for $\beta_n\in\mathbb{R}$ for all $n\in\mathbb{N}_0$.
Under this condition, %for each $n\in\mathbb{N}$, we have
%$\gamma_n>0$ for all $n\in\mathbb{N}$ if and only if $\alpha$ fulfils the property
%$\gamma_n>0$ if and only if $\alpha<-(2n+1)$ or $-(2n-1)<\alpha<-n$.
%Hence
we see that $\gamma_n>0$ for all $n\in\mathbb{N}$ if and only if $\alpha$ fulfils the property
$$
\alpha<-(2n+1)\;\vee\;-(2n-1)<\alpha<-n\;,\quad\forall n\in\mathbb{N}\;.
$$
Clearly, there is no $\alpha$ fulfilling this property.
Therefore, ${\bf u}_B^{(\alpha)}$ is not positive-definite whatever the choice of the parameter $\alpha$.

%\begin{remark}\em
%It is worth mentioning that the Bessel polynomials $\{\widehat{B}_n^{(\alpha)}\}_{n\geq0}$
%are orthogonal on the unit circle $\mathbb{S}^1:=\{z\in\mathbb{C}:|z|=1\}$ with respect to the function
%\begin{snugshade}
%$$
%\rho^{(\alpha)}(z):=\frac{1}{2\pi i}\sum_{k=0}^\infty\frac{1}{(\alpha+1)_k} \left(-\frac{2}{z}\right)^k\;,
%\quad z\in\mathbb{C}\setminus\{0\}
%$$
%\end{snugshade}
%%$$
%%w_\alpha(z)=z^{\alpha}e^{-2/z} \quad (\alpha\neq -2,-3,-4, \ldots)\;,
%%$$
%\noindent
%(this series converges absolutely and uniformly on $\mathbb{C}\setminus\{0\}$),
%although $\{\widehat{B}_n^{(\alpha)}\}_{n\geq0}$ is not an OPS in the positive-definite sense.
%Indeed noting that $\rho\equiv\rho^{(\alpha)}(z)$ fulfills
%\begin{snugshade}
%\begin{equation}\label{1stODEBessel}
%(z^2\rho)'=\big((\alpha+2)z+2\big)\rho-\frac{\alpha z}{2\pi i}\;,\quad z\in\mathbb{C}\setminus\{0\}\;,
%\end{equation}
%\end{snugshade}
%\noindent
%one may prove the following orthogonality relations (Exercise \ref{Ex-cp8-5})
%\begin{snugshade}
%\begin{equation}\label{ortBessel}
%\int_{\mathbb{S}^1}\widehat{B}_m^{(\alpha)}(z)\widehat{B}_n^{(\alpha)}(z)\rho^{(\alpha)}(z)\,{\rm d}z
%=\frac{2^{2n+1}(-1)^{n+1}\,n!}{(\alpha+1)_{2n}(n+\alpha+1)_{n+1}}\delta_{m,n} \; (m,n\in\mathbb{N}_0)\,.
%\end{equation}
%\end{snugshade}
%\noindent
%These relations were stated (up to normalization) by Krall and Frink \ref{KrallFrink1955-C8}.
%\end{remark}

\begin{remark}\em
Table \ref{Table4} summarizes the weight functions representing
the classical functionals on the positive-definite case.
\end{remark}

%\begin{table}
%\centering
%%\begin{tabular}{|>{\columncolor[rgb]{1,1,0}}c|c|c|c|}
%\begin{tabular}{|c|c|c|c|}
%\hline \rowcolor[gray]{0.95}
%\rule{0pt}{1.5em} $P_n(x)$ & Interval of orthogonality & $w(x)$ & Restrictions \\
%\hline
%\rule{0pt}{1.25em} $P_n^{(\alpha,\beta)}(x)$ & $[-1,1]$ & $(1-x)^{\alpha}(1+x)^{\beta} $& $\alpha>-1\;,\;\beta>-1$ \\
%%\hline
%\rule{0pt}{1.25em} $L_n^{(\alpha)}(x)$ & $[0,+\infty[$ & $x^{\alpha}e^{-x}$ & $\alpha >-1$  \\
%%\hline
%\rule{0pt}{1.25em}  $H_n(x)$ & $]-\infty,+\infty[$ & $e^{-x^2}$ & ---  \\ [0.25em]
%\hline
%\end{tabular}
%\medskip
%\caption{Classical OP: the positive-definite case (weight functions).}\label{Table4}
%\end{table}

\begin{table}
\centering
%\begin{tabular}{|>{\columncolor[rgb]{1,1,0}}c|c|c|c|}
\begin{tabular}{|c|c|c|c|}
\hline \rowcolor[gray]{0.95}
\rule{0pt}{1.5em} ${\bf u}$ & Interval of orthogonality & $w(x)$ & Restrictions \\
\hline
\rule{0pt}{1.25em} ${\bf u}_J^{(\alpha,\beta)}$ & $[-1,1]$ & $(1-x)^{\alpha}(1+x)^{\beta} $& $\alpha>-1\;,\;\beta>-1$ \\
%\hline
\rule{0pt}{1.25em} ${\bf u}_L^{(\alpha)}$ & $[0,+\infty[$ & $x^{\alpha}e^{-x}$ & $\alpha >-1$  \\
%\hline
\rule{0pt}{1.25em}  ${\bf u}_H$ & $]-\infty,+\infty[$ & $e^{-x^2}$ & ---  \\ [0.25em]
\hline
\end{tabular}
\medskip
\caption{Weight functions ($w$) representing the canonical forms
(presented in Table \ref{Table1}) in the positive-definite case.}\label{Table4}
\end{table}

%%%%%%%%%%%%%%%%%%%%%%%%%%%%%%%%%%%%%%%%%%%%%%%%%%%%%%%%%%%%%%%%%%%%%%%%%%%%%%%%%%%%%%%%%%%%%%%%%%%%%%%%%%%

\section{Orthogonality of the Bessel polynomials on $\mathbb{S}^1$}

We have seen that the Bessel OPS $\{\widehat{B}_n^{(\alpha)}\}_{n\geq0}$
is not an OPS with respect to a positive-definite functional.
Despite this fact, Krall and Frink \ref{KrallFrink1955-C8} proved that $\{\widehat{B}_n^{(\alpha)}\}_{n\geq0}$
fulfills the orthogonality relations (\ref{ortBessel}) in bellow, where the integration is over the unit circle
$\mathbb{S}^1:=\{z\in\mathbb{C}:|z|=1\}$ (or any closed contour around the origin)
%(in fact, instead of $\mathbb{S}^1$, we can consider any closed smooth contour around the origin)
and the ``\,weight'' function is given by
%\begin{snugshade}
%\begin{equation}\label{rhoBessel}
%\rho^{(\alpha)}(z):=\frac{1}{2\pi i}\sum_{k=0}^\infty\frac{1}{(\alpha+1)_k} \left(-\frac{2}{z}\right)^k\;,
%\quad z\in\mathbb{C}\setminus\{0\}\;.
%\end{equation}
%\end{snugshade}
\begin{snugshade}\vspace*{-0.5em}
\begin{equation}\label{rhoBessel}
\rho^{(\alpha)}(z):=\frac{1}{2\pi i}\left[1+\alpha
+\sum_{k=1}^\infty\frac{1}{(\alpha+2)_{k-1}} \left(-\frac{2}{z}\right)^k\right]\;,
\quad z\in\mathbb{C}\setminus\{0\}\;.
\end{equation}
\end{snugshade}
%$$
%w_\alpha(z)=z^{\alpha}e^{-2/z} \quad (\alpha\neq -2,-3,-4, \ldots)\;,
%$$
\noindent
%Notice that, since
%$$
%\left|\frac{\frac{1}{(\alpha+1)_{k+1}}}{\frac{1}{(\alpha+1)_{k}}}\right|
%=\frac{1}{\big|\alpha+k+1\big|}\to0\quad\mbox{\rm (as $k\to+\infty$)}\;,
%$$
The ratio test ensures that the series in (\ref{rhoBessel}) converges absolutely
on $\mathbb{C}\setminus\{0\}$ and uniformly on each compact subset of this set.
This function $\rho\equiv\rho^{(\alpha)}(z)$ fulfills
\begin{snugshade}\vspace*{-0.5em}
\begin{equation}\label{1stODEBessel}
(z^2\rho)'=\big((\alpha+2)z+2\big)\rho-\frac{\alpha(\alpha+1)}{2\pi i}\,z\;,\quad
z\in\mathbb{C}\setminus\{0\}\;.
\end{equation}
\end{snugshade}
\noindent
%Indeed, computing $\rho'(z)$ we have:
%$$\red
%\begin{array}{rl}
%\rho'(z)&\,=\displaystyle
%\frac{1}{4\pi i}\sum_{k=1}^\infty\frac{(\alpha+k)-\alpha}{(\alpha+1)_k} \left(-\frac{2}{z}\right)^{k+1} \\ [1em]
%&\,=\displaystyle\frac{1}{4\pi i}\sum_{k=1}^\infty\frac{1}{(\alpha+1)_{k-1}} \left(-\frac{2}{z}\right)^{k+1}-
%\frac{\alpha}{4\pi i}\sum_{k=1}^\infty\frac{1}{(\alpha+1)_{k}} \left(-\frac{2}{z}\right)^{k+1}\\ [1em]
%&\,=\displaystyle\frac{2}{z^2}\frac{1}{2\pi i}\sum_{k=0}^\infty\frac{1}{(\alpha+1)_{k}} \left(-\frac{2}{z}\right)^{k}+
%\frac{\alpha}{z}\left(\frac{1}{2\pi i}\sum_{k=0}^\infty\frac{1}{(\alpha+1)_{k}}
%\left(-\frac{2}{z}\right)^{k}-\frac{1}{2\pi i}\right)\\ [1em]
%&\,=\displaystyle\frac{2}{z^2}\rho(z)+
%\frac{\alpha}{z}\left(\rho(z)-\frac{1}{2\pi i}\right)\;,
%\end{array}
%$$
%hence, after multiplying both sides by $z^2$, one obtains (\ref{1stODEBessel}).
Indeed, since
$$
\rho'(z)=\frac{1}{4\pi i}\sum_{k=0}^\infty\frac{k+1}{(\alpha+2)_k} \left(-\frac{2}{z}\right)^{k+2}\;,
$$
we deduce
$$
\begin{array}{l}
(z^2\rho)'(z)-\big[(\alpha+2)z+2\big]\rho(z)=z^2\rho'(z)-(\alpha z+2)\rho(z) \\  [0.25em] \quad\displaystyle
=\frac{1}{\pi i}\sum_{k=0}^\infty\frac{k+1}{(\alpha+2)_k}\left(-\frac{2}{z}\right)^{k}
-\frac{\alpha z+2}{2\pi i}\left[\alpha+1-\frac{2}{z}
\sum_{k=0}^\infty\frac{1}{(\alpha+2)_{k}}\left(-\frac{2}{z}\right)^{k}\right]\\ [1.25em]\quad\displaystyle
=-\frac{\alpha(\alpha+1)}{2\pi i}\,z+
\frac{1}{\pi i}\sum_{k=0}^\infty\frac{k+1+\alpha}{(\alpha+2)_{k}} \left(-\frac{2}{z}\right)^{k}
-\frac{1}{\pi i}\left[\alpha+1-\frac{2}{z}
\sum_{k=0}^\infty\frac{1}{(\alpha+2)_{k}}\left(-\frac{2}{z}\right)^{k}\right] \\ [1em]\quad\displaystyle
=-\frac{\alpha(\alpha+1)}{2\pi i}\,z\;.
\end{array}
$$
%hence (\ref{1stODEBessel}) is proved.
\begin{snugshade}
\begin{theorem}[{\rm Krall $\&$ Frink}]\label{KF-Bessel1}
Let $\{\widehat{B}_n^{(\alpha)}\}_{n\geq0}$ be the Bessel monic OPS
(being $\alpha\in\mathbb{C}\setminus\{-2,-3,-4,\ldots\}$), and let
$\rho^{(\alpha)}$ be defined as in $(\ref{rhoBessel})$. Then
%\begin{equation}\label{ortBessel}
%\int_{\mathbb{S}^1}\widehat{B}_m^{(\alpha)}(z)\widehat{B}_n^{(\alpha)}(z)\rho^{(\alpha)}(z)\,{\rm d}z
%=\frac{2^{2n+1}(-1)^{n+1}\,n!}{(\alpha+1)_{2n}(n+\alpha+1)_{n+1}}\delta_{m,n} \; (m,n\in\mathbb{N}_0)\,.
%\end{equation}
\begin{equation}\label{ortBessel}
\int_{\mathbb{S}^1}\widehat{B}_m^{(\alpha)}(z)\widehat{B}_n^{(\alpha)}(z)\rho^{(\alpha)}(z)\,{\rm d}z
=\frac{2^{2n+1}(-1)^{n+1}\,n!}{(\alpha+2)_{2n}(n+\alpha+1)_{n}}\delta_{m,n} \quad (m,n\in\mathbb{N}_0)\,.
\end{equation}
\end{theorem}
\end{snugshade}
%\noindent

{\it Proof.}
From C5 in Theorem \ref{ThmClassicalOPS}
and Tables \ref{Table1} and \ref{Table3},
$y_n:=\widehat{B}_n^{(\alpha)}(z)$ fulfills %the second order ODE
\begin{equation}\label{ODEBn}
z^2y_n^{\prime\prime}+\big((\alpha+2)z+2\big)y_n'=n(n+\alpha+1)y_n\;,
\quad n\in\mathbb{N}\;.
\end{equation}
Multiplying both sides of (\ref{ODEBn}) by $\rho$
and taking into account (\ref{1stODEBessel}), we deduce
$$
\big(z^2\rho y_n'\big)'+\frac{\alpha(\alpha+1)}{2\pi i} zy_n'=n(n+\alpha+1)\rho y_n \;,
\quad n\in\mathbb{N}\;.
$$
Multiplying both sides of this equality by $y_m$
and then integrating around $\mathbb{S}^1$,
$$
\int_{\mathbb{S}^1}\big(z^2\rho y_n'\big)'y_m\,{\rm d}z
+\frac{\alpha(\alpha+1)}{2\pi i}\int_{\mathbb{S}^1}zy_n'y_m\,{\rm d}z
=n(n+\alpha+1)\int_{\mathbb{S}^1}y_ny_m\rho\,{\rm d}z\quad (n,m\in\mathbb{N}_0)\;.
$$
Clearly, by Cauchy's theorem, $\int_{\mathbb{S}^1}zy_n'y_m\,{\rm d}z=0$.
Also, integrating by parts\footnote{ Recall that if $f$ and $g$ are complex functions
holomorphic on a neighborhood of the image of a differentiable and closed path $\gamma$,
then $\int_\gamma f'(z)g(z)\,{\rm d}z=-\int_\gamma f(z)g'(z)\,{\rm d}z$
(integration by parts formula). Indeed,
considering a path parametrization $\gamma:[0,1]\to\mathbb{C}$, we may write
%$$
%\int_\gamma f'(z)g(z)\,{\rm d}z
%=\int_0^1f'\big(\gamma(t)\big)g\big(\gamma(t)\big)\gamma'(t)\,{\rm d}t
%=\int_0^1\big(f\circ\gamma\big)'(t)\cdot \big(g\circ\gamma\big)(t)\,{\rm d}t\;.
%$$
%Hence using the integration by parts formula for the Riemann integral and taking into account that
%$\gamma(0)=\gamma(1)$ (since $\gamma$ is closed), one easily obtains the integration by parts formula.},
$$
\begin{array}{l}
\displaystyle\int_\gamma f'(z)g(z)\,{\rm d}z
=\int_0^1f'\big(\gamma(t)\big)g\big(\gamma(t)\big)\gamma'(t)\,{\rm d}t
=\int_0^1\big(f\circ\gamma\big)'(t)\cdot \big(g\circ\gamma\big)(t)\,{\rm d}t \\ %[0.5em]
\qquad=\displaystyle\;-\int_0^1\big(f\circ\gamma\big)(t)\cdot \big(g\circ\gamma\big)'(t)\,{\rm d}t
=-\int_0^1f\big(\gamma(t)\big)\cdot g'\big(\gamma(t)\big)\gamma'(t)\,{\rm d}t
=-\int_\gamma f(z)g'(z)\,{\rm d}z\;,
\end{array}
$$
where in the third equality we used the integration by parts formula
for the Riemann integral and took into account that the relation
$\gamma(0)=\gamma(1)$ holds (since $\gamma$ is closed).},
we deduce
$\int_{\mathbb{S}^1}\big(z^2\rho y_n'\big)'y_m\,{\rm d}z=
-\int_{\mathbb{S}^1} z^2\rho y_n' y_m'\,{\rm d}z$.
Hence the above equality reduces to
\begin{equation}\label{int1ODEBn}
n(n+\alpha+1)\int_{\mathbb{S}^1}y_ny_m\rho\,{\rm d}z=
-\int_{\mathbb{S}^1} z^2y_n'y_m'\rho \,{\rm d}z
\quad (n,m\in\mathbb{N}_0)\;.
\end{equation}
Interchanging $n$ and $m$ and subtracting the resulting equality to (\ref{int1ODEBn}), yields
\begin{equation}\label{int2ODEBn}
(n-m)(n+m+\alpha+1)\int_{\mathbb{S}^1}y_ny_m\rho\,{\rm d}z=0
\quad (n,m\in\mathbb{N}_0)\;.
\end{equation}
Since $-(\alpha+1)\not\in\mathbb{N}_0$ then $n+m+\alpha+1\neq0$ if $n\neq m$.
Thus (\ref{int2ODEBn}) gives us
\begin{equation}\label{int3ODEBn}
\int_{\mathbb{S}^1}y_ny_m\rho\,{\rm d}z=0\quad\mbox{\rm if}\quad n\neq m
\quad (n,m\in\mathbb{N}_0)\;.
\end{equation}
This proves (\ref{ortBessel}) whenever $n\neq m$.
If $n=m$, from (\ref{int1ODEBn}) we find
\begin{equation}\label{int4ODEBn}
I_n:=\int_{\mathbb{S}^1}y_n^2\rho\,{\rm d}z=
-\frac{1}{n(n+\alpha+1)}\int_{\mathbb{S}^1} z^2y_n'\cdot y_n'\rho \,{\rm d}z\;,
\quad n\in\mathbb{N}\;.
\end{equation}
By C2 in Theorem \ref{ThmClassicalOPS}, we have
$z^2y_n'=a_ny_{n+1}+b_ny_n+c_ny_{n-1}$. Moreover, clearly,
$y_n'=ny_{n-1}+\sum_{j=0}^{n-2}a_{nj}y_j$
for some complex numbers $a_{nj}$. Substituting these expressions into
the integrand on the right-hand side of (\ref{int4ODEBn}) and using (\ref{int3ODEBn}),
we obtain
\begin{equation}\label{int5ODEBn}
I_n:=\int_{\mathbb{S}^1}y_n^2\rho\,{\rm d}z=
-\frac{c_n}{n+\alpha+1}\int_{\mathbb{S}^1} y_{n-1}^2\rho \,{\rm d}z=\gamma_nI_{n-1}\;,
\quad n\in\mathbb{N}\;,
\end{equation}
where the last equality holds since $c_n=-d_{n-1}\gamma_n=-(n+\alpha+1)\gamma_n$
(see Theorem \ref{ThmClassicalOPS} and Table \ref{Table1}).
Iterating (\ref{int5ODEBn}) we deduce $I_n=\gamma_n\gamma_{n-1}\cdots\gamma_1 I_0$,
hence using the expression for $\gamma_n$ given in Table \ref{Table3}, we deduce
\begin{equation}\label{int6ODEBn}
I_n=\frac{(-4)^nn!}{(\alpha+2)_{2n}(n+\alpha+1)_{n}}I_0\;,
\quad n\in\mathbb{N}_0\;.
\end{equation}
It remains to compute $I_0$. Since $\rho$ is given by the Laurent series (\ref{rhoBessel}),
one see by the definition of residue that ${\rm Res}\,(\rho;z=0)=-\frac{1}{\pi i}\,$, hence, by the residue theorem,
\begin{equation}\label{int7ODEBn}
I_0:=\int_{\mathbb{S}^1}\rho\,{\rm d}z=2\pi i\,\mbox{\rm Res}(\rho;z=0)=-2\;.
\end{equation}
%It remains to compute $I_0$. By (\ref{rhoBessel}),
%\begin{equation}\label{int7ODEBn}
%I_0:=\int_{\mathbb{S}^1}\rho\,{\rm d}z=
%\frac{1}{2\pi i}\sum_{k=1}^\infty\frac{(-2)^k}{(\alpha+2)_{k-1}}\int_{\mathbb{S}^1}z^{-k}\,{\rm d}z=-2\;,
%%\frac{-2}{\alpha+1}\;,
%\end{equation}
%where the last equality holds since, making $z=e^{i\theta}$, $0\leq\theta\leq2\pi$, one has
%$$
%\int_{\mathbb{S}^1}z^{-k}\,{\rm d}z=i\int_0^{2\pi}e^{i(1-k)\theta}\,{\rm d}\theta=
%\left\{
%\begin{array}{ccl}
%2\pi i & \mbox{\rm if} & k=1\;, \\ %[0.5em]
%0 & \mbox{\rm if} & k\neq1\;.
%\end{array}
%\right.
%$$
Inserting (\ref{int7ODEBn}) into (\ref{int6ODEBn}) yields (\ref{ortBessel}) for $n=m$.
This completes the proof.
\qed

\begin{remark}\em
The function $\omega_\alpha(z):=z^\alpha e^{-2/z}$ fulfills the Pearson's type equation
\begin{snugshade}\vspace*{-0.5em}
\begin{equation}\label{walphaBessel}
\big(z^2\omega_\alpha(z)\big)'=\big((\alpha+2)z+2\big)\omega_\alpha(z)\;,\quad
z\in\mathbb{C}\setminus\{0\}\;.
\end{equation}
\end{snugshade}
\noindent
This suggests using $\omega_\alpha$ instead of $\rho^{(\alpha)}$ as ``\,weight'' function in the
orthogonality relations (\ref{ortBessel}).
However, $\omega_\alpha$ is a multivalued function if $\alpha$ is not an integer number
and this may be inconvenient for integration around the origin $0$. %$\mathbb{S}^1$.
%This is the reason why in (\ref{ortBessel}) the ``\,weight'' $\rho^{(\alpha)}$
%is considered instead of $\omega_\alpha$.
In general, $\rho^{(\alpha)}$ and $\omega_\alpha$ yield different orthogonality relations
for the Bessel polynomials, unless $\alpha=0$ or $\alpha=-1$
(compare equations (\ref{1stODEBessel}) and (\ref{walphaBessel}), and see Exercise \ref{Ex-cp8-5}).
\end{remark}

%\begin{remark}\em
%Although ${\bf u}_B^{(\alpha)}$ is not a positive-definite functional,
%the Bessel polynomials $\widehat{B}_n^{(\alpha)}$ fulfill the (ordinary) Rodrigues-type formula
%\begin{snugshade}\vspace*{-0.5em}
%\begin{equation}\label{RodrBessel}
%\widehat{B}_n^{(\alpha)}(x)=\frac{1}{(n+\alpha+1)_n}\frac{1}{\omega_\alpha(x)} \frac{{\rm d}^n}{{\rm d}x^n}\big(x^{2n}\omega_\alpha(x)\big)\;,\quad \omega_\alpha(x):=x^{\alpha}e^{-2/x}
%\end{equation}
%\end{snugshade}
%\noindent
%(for $x\in\mathbb{C}\setminus\{0\}$).
%This formula can be used to derive the explicit expression (\ref{ExpB1}).
%\end{remark}

\section{Explicit expressions for the classical OP}

\subsection{The Rodrigues-type formula}
Theorem \ref{vCOP-Thm1} states that the monic OPS $\{P_n\}_{n\geq0}$ with respect to a classical functional ${\bf u}$
(which needs not to be positive-definite) satisfies a distributional Rodrigues formula, involving $P_n$ and ${\bf u}$.
In the next we prove that if, in addition, ${\bf u}$ is (classical and) represented by the weight function $\omega$,
then also a Rodrigues-type formula involving $P_n$ and $\omega$ holds.

\begin{snugshade}
\begin{theorem}[Rodrigues-type formula]\label{ThmRodW}
Assume the hypothesis of Lemma \ref{lemmaW},
so that ${\bf u}$ is positive-definite and represented by the weight function $\omega$,
as in $(\ref{uW1})$.
Assume further that $\omega\in\mathcal{C}^\infty(\xi,\eta)$,
and $\phi$ and $\psi$ are nonzero real polynomials,
$\phi\in\mathcal{P}_2$ and $\psi\in\mathcal{P}_1$.
Let $\{P_n\}_{n\geq0}$ be the monic OPS with respect to ${\bf u}$.
%Under the conditions of Lemma \ref{lemmaW},
%assume further that $\phi$ and $\psi$ are nonzero (real) polynomials such that
%$\phi\in\mathcal{P}_2$ and $\psi\in\mathcal{P}_1$.
%%$(\phi,\psi)$ an admissible pair of (real) polynomials.
%Let $\{P_n\}_{n\geq0}$ be the monic OPS with respect to ${\bf u}$
%(note that, according with $(\ref{uW1})$, ${\bf u}$ is represented by the weight function $\omega$).
Then, for each $n\in\mathbb{N}_0$,
\begin{equation}\label{PearsonRodfunc}
P_n(x)=\frac{k_n}{w(x)} \frac{{\rm d}^n}{{\rm d}x^n}\big(\phi^n(x)\omega(x)\big)\;,\quad \xi<x<\eta\,,
\end{equation}
where $k_n:=\prod_{i=0}^{n-1}d_{n+i-1}^{-1}$,
being $d_k:=\frac{k}{2}\,\phi^{\prime\prime}+\psi'$.
\end{theorem}
\end{snugshade}

{\it Proof.}
By Lemma \ref{lemmaW}, ${\bf u}$ fulfills the distributional Pearson's equation
\begin{equation}\label{PearsonRod1}
D(\phi{\bf u})=\psi{\bf u}\,. %\quad\mbox{\rm in $\mathcal{P}'$}\,,
\end{equation}
Therefore, since ${\bf u}$ is regular (because it is positive-definite),
Theorem \ref{vCOP-Thm1} ensures that $d_k\neq0$ for all $k\in\mathbb{N}_0$,
and the distributional Rodrigues formula holds:
\begin{equation}\label{PearsonRod2}
P_n{\bf u}=k_n D^n(\phi^n{\bf u})\;,\quad n\in\mathbb{N}_0\;. %\quad\mbox{\rm in $\mathcal{P}'$}\,,
\end{equation}
While proving Theorem \ref{Pearson-Thm2}, and up to normalization (being $P_n=k_nR_n$),
we have deduced (\ref{PearsonRod2}) from (\ref{PearsonRod1})
by purely algebraic arguments, the essential tool in the proof being the distributional Leibnitz rule
for the derivative of order $n$ of the functional $\phi{\bf u}$
(the left product of the functional ${\bf u}$ by the polynomial $\phi$).
Therefore, since the weight function $\omega$ fulfills the ODE
(which can be regarded as an analogue version for ordinary functions
of the distributional differential equation (\ref{PearsonRod1}))
\begin{equation}\label{PearsonRod3}
(\phi\omega)'=\psi\omega\quad\mbox{\rm on $(\xi,\eta)$}\,,
\end{equation}
%using Leibnitz rule for the product of ordinary functions,
we see without effort, {\it mutatis mutandis},
that the steps of the proof of Theorem \ref{Pearson-Thm2} may be followed
(replacing therin ${\bf u}$ by $\omega$, and ${\bf u}^{[n]}$ by $\omega^{[n]}:=\phi^n\omega$,
and considering the ordinary derivative instead of the distributional derivative)
allowing us to deduce (\ref{PearsonRodfunc}) from (\ref{PearsonRod3})
and the Leibnitz rule for the product of ordinary functions.
\qed
\medskip

The Rodrigues formula (\ref{PearsonRodfunc}) gives an explicit representation
for $P_n$ as a derivative of order $n$ of a simple real function, divided by the weight function.
This representation is very useful in many areas, e.g., in Number Theory, or in Physics.
In the next we use (\ref{PearsonRodfunc}) to derive explicit expressions for $P_n$
as a linear combination of powers of $x$.
Up to an affine change of variables, we can restrict our study to
the canonical forms described in the previous sections.

%Let ${\bf u}$ be a positive-definite classical functional, so that, up to an affine change of variables,
%it is one of the functionals of Hermite, Laguerre (with $\alpha>-1$), or Jacobi (with $\alpha>-1$ and $\beta>-1$),
%all of them admitting an integral representation of the form (\ref{uW1}), where in each case the
%polynomials $\phi$ and $\psi$ are given in Table \ref{Table1}, and the weight function
%$\omega$ as well as the interval of orthogonality $(\xi,\eta)$ are given in Table \ref{Table4}.

\subsection{Explicit formula for Jacobi polynomials}
Consider $\alpha>-1$ and $\beta>-1$.
Substituting in (\ref{PearsonRodfunc}) the explicit expression of
$\phi$, $\psi$, $k_n$, and $\omega$ appearing in Tables \ref{Table1}, \ref{Table3}, and \ref{Table4},
we may write, for each $n\in\mathbb{N}_0$ and $-1<x<1$,
\begin{snugshade}
$$%\begin{equation}\label{ExpJ0}
\widehat{P}_n^{(\alpha,\beta)}(x)=\frac{(-1)^n}{(n+\alpha+\beta+1)_n}\frac{1}{(1-x)^\alpha(1+x)^\beta}
\frac{{\rm d}^n}{{\rm d}x^n}\big\{(1-x)^{n+\alpha}(1+x)^{n+\beta}\big\}\;.
$$%\end{equation}
\end{snugshade}
\noindent
By Leibniz's rule for the $n$th derivative of a product, and making use of
the generalized binomial coefficient, defined by
%\begin{snugshade}
%\begin{equation}\label{gBinomial}
%{\alpha\choose 0}:=1\;,\quad
%{\alpha\choose k}:=\frac{\alpha(\alpha-1)\cdots(\alpha-k+1)}{k!}\;,
%\quad \alpha\in\mathbb{C}\;,\quad k\in\mathbb{N}\;.
%\end{equation}
%\end{snugshade}
%\noindent
\begin{snugshade}
$$%\begin{equation}\label{gBinomial}
{z\choose 0}:=1\;,\quad
{z\choose k}:=\frac{z(z-1)\cdots(z-k+1)}{k!}\;,
\quad z\in\mathbb{C}\;,\;\; k\in\mathbb{N}\;,
$$%\end{equation}
\end{snugshade}
\noindent
we deduce
\begin{snugshade}
\begin{equation}\label{ExpJ1}
\widehat{P}_n^{(\alpha,\beta)}(x)=\frac{1}{{2n+\alpha+\beta\choose n}}
\sum_{k=0}^{n} {{n+\alpha\choose n-k}} {{n+\beta\choose k}} (x-1)^{k}(x+1)^{n-k}\;,
\quad n\in\mathbb{N}_0\;.
\end{equation}
\end{snugshade}
\noindent
In the literature on OP
%(see e.g., Szego \ref{}, Freud \ref{}, Chihara \ref{}, and Ismail \ref{})
it is usual to consider (non monic) Jacobi polynomials $\{P_n^{(\alpha,\beta)}\}_{n\geq0}$
normalized so that %
\begin{snugshade}
\begin{equation}\label{ExpJ2}
P_n^{(\alpha,\beta)}(1)={n+\alpha\choose n}\;,
\quad n\in\mathbb{N}_0\;,
\end{equation}
\end{snugshade}
\noindent
called the {\it standard normalization} for Jacobi polynomials.
Therefore, since, by (\ref{ExpJ1}),
\begin{snugshade}
\begin{equation}\label{ExpJ2a}
\widehat{P}_n^{(\alpha,\beta)}(1)
=\frac{2^n}{{2n+\alpha+\beta\choose n}}{n+\alpha\choose n}
=\frac{2^n (\alpha+1)_n}{(n+\alpha+\beta+1)_n}\;,
\end{equation}
\end{snugshade}
\noindent
we conclude that the relation between the Jacobi polynomials with
standard normalization %(\ref{ExpJ2a}),
and the monic Jacobi polynomials is
\begin{snugshade}
\begin{equation}\label{ExpJ3}
P_n^{(\alpha,\beta)}(x)=2^{-n}{2n+\alpha+\beta\choose n}\,\widehat{P}_n^{(\alpha,\beta)}(x)\;,
\quad n\in\mathbb{N}_0\;.
\end{equation}
\end{snugshade}
\noindent
This together with (\ref{ExpJ1}) leads to the explicit expression for the Jacobi
polynomials with the standard normalization (\ref{ExpJ3}):
\begin{snugshade}
$$%\begin{equation}\label{ExpJ4}
P_n^{(\alpha,\beta)}(x)=
2^{-n}\sum_{k=0}^{n} {{n+\alpha\choose n-k}} {{n+\beta\choose k}} (x-1)^k(x+1)^{n-k}\;,
\quad n\in\mathbb{N}_0\;.
$$%\end{equation}
\end{snugshade}
\noindent
We also point out the following useful relation:
\begin{snugshade}
\begin{equation}\label{ExpJ5}
P_n^{(\alpha,\beta)}(-x)=(-1)^nP_n^{(\beta,\alpha)}(x)\;,
\quad n\in\mathbb{N}_0\;.
\end{equation}
\end{snugshade}
\noindent
(Clearly this relation holds also for the monic polynomials.)
Also, from (\ref{ExpJ2}), (\ref{ExpJ2a}), and (\ref{ExpJ5}),
\begin{snugshade}
$$%\begin{equation}\label{ExpJ6}
P_n^{(\alpha,\beta)}(-1)=(-1)^n{n+\beta\choose n}\;,\quad
\widehat{P}_n^{(\alpha,\beta)}(-1)=\frac{(-2)^n(\beta+1)_n}{(n+\alpha+\beta+1)_n}\;,
\quad n\in\mathbb{N}_0\;.
$$%\end{equation}
\end{snugshade}
\noindent
Finally, for each $k\in\mathbb{N}_0$, the following formula holds (Exercise \ref{Ex-cp8-2}):
\begin{snugshade}
\begin{equation}\label{ExpJ7}
\frac{{\rm d}^k}{{\rm d}x^k}\big\{\widehat{P}_n^{(\alpha,\beta)}(x)\big\}
=(n-k+1)_k\,\widehat{P}_{n-k}^{(\alpha+k,\beta+k)}(x)\;,\quad n\geq k\;.
\end{equation}
\end{snugshade}

\subsection{Explicit formula for Laguerre polynomials}
Let $\alpha>-1$.
Substituting in (\ref{PearsonRodfunc}) the explicit expression of
$\phi$, $\psi$, $k_n$, and $\omega$ appearing in Tables \ref{Table1}, \ref{Table3}, and \ref{Table4},
we may write, for each $n\in\mathbb{N}_0$ and $x>0$,
\begin{snugshade}
$$%\begin{equation}\label{ExpL0}
\widehat{L}_n^{(\alpha)}(x)=(-1)^nx^{-\alpha}e^x
\frac{{\rm d}^n}{{\rm d}x^n}\big\{x^{n+\alpha}e^{-x}\big\}\;.
$$%\end{equation}
\end{snugshade}
\noindent
By Leibniz's rule for the $n$th derivative of a product, we deduce
\begin{snugshade}
\begin{equation}\label{ExpL1}
\widehat{L}_n^{(\alpha)}(x)=(-1)^n n!
\sum_{k=0}^{n} {{n+\alpha\choose n-k}} \frac{(-x)^{k}}{k!}\;,
\quad n\in\mathbb{N}_0\;.
\end{equation}
\end{snugshade}
\noindent
Considering Laguerre polynomials $\{L_n^{(\alpha)}\}_{n\geq0}$
with {\it standard normalization}, i.e., %
\begin{snugshade}
$$%\begin{equation}\label{ExpL2}
L_n^{(\alpha)}(0)={n+\alpha\choose n}\;,
\quad n\in\mathbb{N}_0\;,
$$%\end{equation}
\end{snugshade}
\noindent
then, since, by (\ref{ExpL1}),
\begin{snugshade}
$$%\begin{equation}\label{ExpL2a}
\widehat{L}_n^{(\alpha)}(0)
=(-1)^n n! {n+\alpha\choose n}
=(-1)^n (\alpha+1)_n\;,
$$%\end{equation}
\end{snugshade}
\noindent
we see that the relation between the Laguerre polynomials with
standard normalization %(\ref{ExpL2a}),
and the monic Laguerre polynomials is
\begin{snugshade}
$$%\begin{equation}\label{ExpL3}
L_n^{(\alpha)}(x)=\frac{(-1)^n}{n!}\,\widehat{L}_n^{(\alpha)}(x)
\quad n\in\mathbb{N}_0\;.
$$%\end{equation}
\end{snugshade}
\noindent
This together with (\ref{ExpL1}) leads to the explicit expression for the Laguerre
polynomials with the standard normalization: %(\ref{ExpL3}):
\begin{snugshade}
$$%\begin{equation}\label{ExpL4}
L_n^{(\alpha)}(x)=
\sum_{k=0}^{n} {{n+\alpha\choose n-k}} \frac{(-x)^k}{k!}\;,
\quad n\in\mathbb{N}_0\;.
$$%\end{equation}
\end{snugshade}
\noindent
We also point out that, for each $k\in\mathbb{N}_0$, the following formula holds (Exercise \ref{Ex-cp8-2}):
\begin{snugshade}
\begin{equation}\label{ExpL7}
\frac{{\rm d}^k}{{\rm d}x^k}\big\{\widehat{L}_n^{(\alpha)}(x)\big\}
=(n-k+1)_k\,\widehat{L}_{n-k}^{(\alpha+k)}(x)\;,\quad n\geq k\;.
\end{equation}
\end{snugshade}

\subsection{Explicit formula for Hermite polynomials}

Substituting in (\ref{PearsonRodfunc}) the explicit expression of
$\phi$, $\psi$, $k_n$, and $\omega$ appearing in Tables \ref{Table1}, \ref{Table3}, and \ref{Table4},
we may write, for each $n\in\mathbb{N}_0$ and $x\in\mathbb{R}$,
\begin{snugshade}\vspace*{-0.25em}
$$%\begin{equation}\label{ExpH0}
\widehat{H}_n(x)=\frac{(-1)^n}{2^n}\,e^{x^2}\,
\frac{{\rm d}^n}{{\rm d}x^n}\big\{e^{-x^2}\big\}\;.
$$%\end{equation}
\end{snugshade}
\noindent
Using this formula we can derive the explicit expression for the Hermite polynomials.
Nevertheless we will obtain such a formula by a different way.
By Lemma \ref{reg-lemma3}, for each $k\in\mathbb{N}_0$,
$\big\{\widehat{H}_n^{[k]}:=\frac{1}{(n+1)_k}\frac{{\rm d}^k\widehat{H}_{n+k}}{{\rm d}x^k}\big\}_{n\geq0}$
is a monic OPS with respect to the functional ${\bf u}^{[k]}:=\phi^k{\bf u}_H={\bf u}_H$
(since $\phi\equiv1$, by Table \ref{Table1}), hence $\widehat{H}_n^{[k]}\equiv H_n$, and so
\begin{snugshade}\vspace*{-0.5em}
\begin{equation}\label{ExpH7}
\frac{{\rm d}^k}{{\rm d}x^k}\big\{\widehat{H}_n(x)\big\}
=(n-k+1)_k\,\widehat{H}_{n-k}(x)\;,\quad n\geq k\;.
\end{equation}
\end{snugshade}
\noindent
Therefore, using McLaurin formula, and since $(n-k+1)_k=k!{n\choose k}$, we may write
\begin{equation}\label{ExpHn1}
\widehat{H}_n(x)
=\sum_{k=0}^n\frac{ \frac{{\rm d}^k\widehat{H}_n}{{\rm d}x^k}(0)}{k!}\,x^k
=\sum_{k=0}^n{n\choose k}\widehat{H}_{n-k}(0)\,x^k
=\sum_{j=0}^n{n\choose j}\widehat{H}_{j}(0)\,x^{n-j}\;.
\end{equation}
To compute $\widehat{H}_{j}(0)$,
we start with the TTRR for $\{\widehat{H}_{n}\}_{n\geq0}$ (see Table \ref{Table3}):
$$
\widehat{H}_{n+1}(x)=x\,\widehat{H}_{n}(x)-\mbox{$\frac{n}{2}$}\,\widehat{H}_{n-1}(x)\;,\quad n\geq0
$$
($\widehat{H}_{-1}(x)=0$, $\widehat{H}_{0}(x)=1$). Thus
$\widehat{H}_{n+1}(0)=-\frac{n}{2}\,\widehat{H}_{n-1}(0)$ for each $n\geq0$,
hence
\begin{snugshade}\vspace*{-0.5em}
\begin{equation}\label{ExpHn0}
\widehat{H}_{2n-1}(0)=0\;,\quad \widehat{H}_{2n}(0)=
(-1)^n\frac{(2n-1)!!}{2^n}\;,\quad n\geq 1\;.
\end{equation}
\end{snugshade}
\noindent
Inserting (\ref{ExpHn0}) into (\ref{ExpHn1}) we easily deduce the desired explicit expression:
\begin{snugshade}\vspace*{-0.5em}
\begin{equation}\label{ExpH1}
\widehat{H}_n(x)=\frac{n!}{2^n}
\sum_{k=0}^{\lfloor n/2\rfloor} \frac{(-1)^k(2x)^{n-2k}}{(n-2k)!k!}\;,
\quad n\in\mathbb{N}_0\;.
\end{equation}
\end{snugshade}
\noindent
The standard normalization for the Hermite polynomials is $\{H_n\}_{n\geq0}$ given by
\begin{snugshade}\vspace*{-0.5em}
\begin{equation}\label{ExpH1standard}
H_n(x)=2^n\widehat{H}_n(x)=n!
\sum_{k=0}^{\lfloor n/2\rfloor} \frac{(-1)^k(2x)^{n-2k}}{(n-2k)!k!}\;,
\quad n\in\mathbb{N}_0\;.
\end{equation}
\end{snugshade}

\subsection{Explicit formula for Bessel polynomials}

The Bessel functional ${\bf u}_B^{(\alpha)}$
fulfils Pearson's equation $D\big(\phi{\bf u}_B^{(\alpha)}\big)=\psi\,{\bf u}_B^{(\alpha)}$,
where $\phi(x)=x^2$ and $\psi(x)=(\alpha+2)x+2$ (see Table \ref{Table1}).
Hence, by Lemma \ref{reg-lemma3}, the sequence
$\big\{\big[\widehat{B}_n^{(\alpha)}\big]^{[k]}
:=\frac{1}{(n+1)_k}\frac{{\rm d}^k\widehat{B}_{n+k}^{(\alpha)}}{{\rm d}x^k}\big\}_{n\geq0}$
is a monic OPS with respect to the functional ${\bf u}^{[k]}:=x^{2k}{\bf u}_B^{(\alpha)}$,
for each $k\in\mathbb{N}_0$.
Moreover, by Lemma \ref{Pearson-lemma1},
${\bf u}^{[k]}$ fulfills Pearson's equation $D\big(x^{2}{\bf u}^{[k]}\big)=\psi_k\,{\bf u}^{[k]}$,
where $\psi_k:=\psi+k\phi'=(\alpha+2k+2)x+2$. Thus, ${\bf u}^{[k]}={\bf u}_B^{(\alpha+2k)}$,
$\big[\widehat{B}_n^{(\alpha)}\big]^{[k]}\equiv\widehat{B}_n^{(\alpha+2k)}$, and so
\begin{snugshade}\vspace*{-0.5em}
\begin{equation}\label{ExpB7}
\frac{{\rm d}^k}{{\rm d}x^k}\big\{\widehat{B}_n^{(\alpha)}(x)\big\}
=(n-k+1)_k\,\widehat{B}_{n-k}^{(\alpha+2k)}(x)\;,\quad n\geq k\;.
\end{equation}
\end{snugshade}
\noindent
Therefore, using McLaurin formula, we may write
\begin{equation}\label{ExpBn1}
\widehat{B}_n^{(\alpha)}(x)
=\sum_{k=0}^n\frac{ \frac{{\rm d}^k\widehat{B}_n^{(\alpha)}}{{\rm d}x^k}(0)}{k!}\,x^k
=\sum_{k=0}^n{n\choose k}\widehat{B}_{n-k}^{(\alpha+2k)}(0)\,x^k\;.
\end{equation}
Thus we need to compute $\widehat{B}_{n}^{(\alpha)}(0)$ for each $n\in\mathbb{N}_0$.
%We start with the ordinary differential equation (\ref{ODEBn}) for $y_n:=\widehat{B}_{n}^{(\alpha)}(x)$.
%given by C5 in Theorem \ref{ThmClassicalOPS} (cf. Tables \ref{Table1}, \ref{Table3}):
%\begin{equation}\label{ODEBn}
%x^2y_n^{\prime\prime}+\big((\alpha+2)x+2\big)y_n'=n(n+\alpha+1)y_n\;,\quad n\geq0\;.
%\end{equation}
Taking $z=x=0$ in the ODE (\ref{ODEBn}) %ordinary differential equation (\ref{ODEBn})
for $y_n:=\widehat{B}_{n}^{(\alpha)}(x)$, and since, by (\ref{ExpB7}) with $k=1$,
$y_n'(x)=n\widehat{B}_{n-1}^{(\alpha+2)}(x)$, we obtain (for each $n\geq1$)
$\widehat{B}_{n}^{(\alpha)}(0)=\frac{2}{n+\alpha+1}\,\widehat{B}_{n-1}^{(\alpha+2)}(0)$,
hence, by iteration of this identity, we deduce
$\widehat{B}_{n}^{(\alpha)}(0)=\frac{2^n}{(n+\alpha+1)_n}\,\widehat{B}_{0}^{(\alpha+2n)}(0)$, i.e.,
\begin{snugshade}\vspace*{-0.5em}
\begin{equation}\label{ExpBn0}
\widehat{B}_{n}^{(\alpha)}(0)=\frac{2^n}{(n+\alpha+1)_n}\,,\quad n\geq0\,.
\end{equation}
\end{snugshade}
\noindent
Finally, inserting (\ref{ExpBn0}) into (\ref{ExpBn1}) we obtain
\begin{snugshade}\vspace*{-0.5em}
\begin{equation}\label{ExpB1}
\widehat{B}_n^{(\alpha)}(x)=\frac{2^n}{(n+\alpha+1)_n}
\sum_{k=0}^{n} {n\choose k}(n+\alpha+1)_k\left(\frac{x}{2}\right)^k\;,
\quad n\in\mathbb{N}_0\;.
\end{equation}
\end{snugshade}
\noindent
A standard normalization for the Bessel polynomials
(cf. Chihara \ref{Chihara1978-C8}, p.\,182--183) is $\{Y_n^{(\alpha)}\}_{n\geq0}$
chosen so that $Y_n^{(\alpha)}(0)=1$, and so
\begin{snugshade}\vspace*{-0.5em}
\begin{equation}\label{ExpB1standard}
Y_n^{(\alpha)}(x)=\sum_{k=0}^{n} {n\choose k}(n+\alpha+1)_k\left(\frac{x}{2}\right)^k\;,
\quad n\in\mathbb{N}_0\;.
\end{equation}
\end{snugshade}

\begin{remark}\em
Although ${\bf u}_B^{(\alpha)}$ is not a positive-definite functional,
the Bessel polynomials $\widehat{B}_n^{(\alpha)}$ fulfill the (ordinary) Rodrigues-type formula
(for $x\in\mathbb{C}\setminus\{0\}$)
\begin{snugshade}\vspace*{-0.5em}
\begin{equation}\label{RodrBessel}
\widehat{B}_n^{(\alpha)}(x)=\frac{1}{(n+\alpha+1)_n}\frac{1}{\omega_\alpha(x)} \frac{{\rm d}^n}{{\rm d}x^n}\big(x^{2n}\omega_\alpha(x)\big)\;,\quad \omega_\alpha(x):=x^{\alpha}e^{-2/x}
\end{equation}
\end{snugshade}
\noindent
(Exercise \ref{Ex-cp8-5}).
This formula can be used to derive the explicit expression (\ref{ExpB1}).
\end{remark}

\begin{remark}\em
The explicit expressions (\ref{ExpJ1}) and (\ref{ExpL1}) for Jacobi and Laguerre polynomials
can be deduced using the technique we have applied to derive (\ref{ExpB1}).
Thus we may remove the restrictions $\alpha>-1$ and $\beta>-1$
considered on the proof of (\ref{ExpJ1}) and (\ref{ExpL1}), and hence these
explicit formulas remain true requiring only that the corresponding functionals are regular
(not necessarily positive-definite).
\end{remark}

%\section{$\nabla_\omega-$classical and $D_{q}-$classical orthogonal polynomials}
%
%
%\section{Remarks and historical notes}

\section*{Exercises}
%\bigskip

{\small
%\noindent
\begin{enumerate}[label=\emph{\bf \arabic*.},leftmargin=*]
\item\label{Ex-cp8-0}
Complete the proof of Theorem \ref{ThmClassicalOPS} by proving that:
\begin{enumerate}
\item the polynomial $\Delta_2$ introduced in the proof of C4$\,\Rightarrow\,$C1 fulfills
$\Delta_2\in\mathcal{P}_2\setminus\{0\}$.
\item C1$'$$\,\Leftrightarrow\,$C6.
\end{enumerate}
\medskip
%\item
%Let ${\bf u},{\bf v}\in\mathcal{P}'$ such that ${\bf u}\sim{\bf v}$.
%Prove that if ${\bf u}$ is classical, then so is ${\bf v}$, and they belong to the same equivalence class, i.e.,
%$[{\bf u}]=[{\bf v}]$.
%%${\bf u}$ and ${\bf v}$ belong to the same equivalence class.
%%Prove that if the functional ${\bf u}$ in Lemma \ref{equiv-u-v} is classical, then so is ${\bf v}$.
%\item
%Mostrar a representacao explicita para os polin de Legendre normalizados para [0,1]
%\item
%In the positive-definite case, find the values of the classical monic OPS of Laguerre, Jacobi and Bessel at the extreme points of their (true) intervals of orthogonality.
%\medskip
%\item
%Appel ortogonais = Hermite (ver texto de Maroni das funcoes eulerianas)
%\medskip
%\item
%Colocar a caracterizacao de Ya L Geronimus sobre os classicos (equiv com os momentos)
%%%%%%%%%%%%%%%%%%%%%%%%%%%%%%%%%%%%%%%
%\item\label{Ex-cp8-0}
%Let $\{P_n\}_{n\geq0}$ be a monic OPS. Show that $\{P_n\}_{n\geq0}$ is classical if and only if
%the following property holds (another characterization of the classical OPS):
%\smallskip
%
%\noindent\;{\rm C3$'$.}
%$\Big\{ P_n^{[k]}:=\frac{{\rm d}^k}{{\rm d}x^k}\frac{P_{n+k}}{(n+1)_k}\Big\}_{n\geq0}$ is a monic OPS for some $k\in\mathbb{N}$.
%\medskip
%%%%%%%%%%%%%%%%%%%%%%%%%%%%%%%%%%%%%%%
\item\label{Ex-cp8-01}
For arbitrary $k\in\mathbb{N}$, find expressions for the parameters $r_n^{[k]}$ and $s_n^{[k]}$ appearing in characterization C4 of Theorem \ref{ThmClassicalOPS}, only in terms of the coefficients of the polynomials $\phi$ and $\psi$
appearing in the Pearson's equation for ${\bf u}$. Compute these expressions for the classical canonical forms of Hermite, Laguerre, Jacobi, and Bessel.
\medskip
\item\label{Ex-cp8-1}
%Calcular os determinantes de Hankel para os classicos e relacionar com as random matrices
Let ${\bf u}$ be a classical functional, so that it is a regular functional on $\mathcal{P}$
which fulfills Pearson's equation
$D(\phi{\bf u})=\psi{\bf u}$, being $\phi\in\mathcal{P}_2$ and $\psi\in\mathcal{P}_1\setminus\mathcal{P}_0$.
\begin{enumerate}
\item
Find a closed formula for the Hankel determinant  $H_n:=\det\big\{[u_{i+j}]_{i,j=0}^{n}\big\}$ (of order $n+1$),
involving only the (coefficients of the) polynomials $\phi$ and $\psi$.
\item
Compute $H_n$ for the classical canonical forms of Hermite, Laguerre, Jacobi, and Bessel.
%Compute $H_n$ for the canonical forms of the classical OPS (Hermite, Laguerre, Jacobi, and Bessel).
Give also expressions for the moments in each case.
\end{enumerate}
\smallskip

%\item
%State (using simple arguments) the following relations among the monic classical OPS
%and their derivatives (being $n\in\mathbb{N}$ and $0\leq k\leq n$):
%$$
%\begin{array}{clcl}
%\quad\mbox{\rm (i)} &
%\displaystyle\frac{{\rm d}^k}{{\rm d}x^k}\,\widehat{H}_n
%=(n-k+1)_k\,\widehat{H}_{n-k}\, &
%\mbox{\rm (iii)} &
%\displaystyle\frac{{\rm d}^k}{{\rm d}x^k}\,\widehat{B}_n^{(\alpha)}
%=(n-k+1)_k\,\widehat{B}_{n-k}^{(\alpha+2k)} \\ [1em]
%%=\,k!{n\choose k}\widehat{P}_{n-k}^{(\alpha+k,\beta+k)}(x)
%\quad\mbox{\rm (ii)} & \displaystyle\frac{{\rm d}^k}{{\rm d}x^k}\,\widehat{L}_n^{(\alpha)}
%=(n-k+1)_k\,\widehat{L}_{n-k}^{(\alpha+k)} &
%\mbox{\rm (iv)} &
%\displaystyle\frac{{\rm d}^k}{{\rm d}x^k}\,\widehat{P}_n^{(\alpha,\beta)}
%=(n-k+1)_k\,\widehat{P}_{n-k}^{(\alpha+k,\beta+k)}\,.
%\end{array}
%$$
%\noindent
%%({\sl Hint.} Use Lemmas 7.5 and 7.1, together with Table 1.)
%({\sl Hint.} Use Lemmas \ref{reg-lemma3} and \ref{Pearson-lemma1} together with Table \ref{Table1}.)
%\smallskip

\item\label{Ex-cp8-2}
Prove relations (\ref{ExpJ7}) and (\ref{ExpL7}).
({\sl Hint.} Proceed as we did for proving (\ref{ExpB7}).)
%, by using Lemmas \ref{reg-lemma3} and \ref{Pearson-lemma1}, together with Table \ref{Table1}.)
\smallskip

\item\label{Ex-cp8-3}
Prove that the (standard) Jacobi polynomials admit the explicit representation
$$
P_n^{(\alpha,\beta)}(x)={2n+\alpha+\beta\choose n}\,
\sum_{k=0}^{n} \frac{{n\choose k}{n+\alpha\choose n-k}}{{2n+\alpha+\beta\choose n-k}} \left(\frac{x-1}{2}\right)^{k}\;,
\quad n\in\mathbb{N}_0\;.
$$
%$$
%\widehat{P}_n^{(\alpha,\beta)}(x)=2^n\,
%\sum_{k=0}^{n} \frac{{n\choose k}{n+\alpha\choose n-k}}{{2n+\alpha+\beta\choose n-k}} \left(\frac{x-1}{2}\right)^{k}\;,
%\quad n\in\mathbb{N}_0\;.
%$$

\item\label{Ex-cp8-4}
\begin{enumerate}
\item
Suppose that ${\bf u}\in\mathcal{P}'$ is regular, and let $\{P_n\}_{n\geq0}$ be its monic OPS.
Let $c\in\mathbb{C}$ and set ${\bf v}:=(x-c){\bf u}$. Prove that ${\bf v}$ is regular if and only if $P_n(c)\neq0$ for all $n\in\mathbb{N}$. Under such conditions, $\{Q_n\}_{n\geq0}$ being the monic OPS with respect to ${\bf v}$, show that
$$
Q_n(x)=\frac{1}{x-c}\,\Big[P_{n+1}(x)-\frac{P_{n+1}(c)}{P_n(c)}P_n(x)\Big]\;,\quad n\in\mathbb{N}_0\;.
$$
\item
Using the results in (a), prove that the following relation among (standard) Jacobi polynomials
holds for each $n\in\mathbb{N}_0$:
$$
\quad\qquad(2n+\alpha+\beta+2)(x+1)P_n^{(\alpha,\beta+1)}(x)
=2(n+\beta+1)P_n^{(\alpha,\beta)}(x)+2(n+1)P_{n+1}^{(\alpha,\beta)}(x)\,.
$$
% $P_n^{(\alpha,\beta)}$being normalized so that its leading coefficient is
%$2^{-n}\Big({{2n+\alpha+\beta}\atop{n}}\Big)$.
%You may assume $\alpha,\beta>-1$. % and $\beta>-1$,
%(It may be useful to recall that $P_n^{(\alpha,\beta)}(-1)=(-1)^n\Big({{n+\beta}\atop{n}}\Big)$.)
\end{enumerate}
%\medskip

\item\label{Ex-cp8-5}
Let ${\bf u}_B^{(\alpha)}$, with $\alpha\in\mathbb{C}\setminus\{-2,-3,-4,\cdots\}$,
be the (canonical) Bessel functional,
and $\{\widehat{B}_n^{(\alpha)}\}_{n\geq0}$ the monic OPS with respect to ${\bf u}_B^{(\alpha)}$.
Show that, although ${\bf u}_B^{(\alpha)}$ is not a positive-definite functional, the following holds:
\begin{enumerate}
\item For each $n\in\mathbb{N}_0$, $\widehat{B}_n^{(\alpha)}$ fulfills the Rodrigues-type formula
(\ref{RodrBessel}).
\item
If $\alpha\in\mathbb{N}_0\cup\{-1\}$, then $\{\widehat{B}_n^{(\alpha)}\}_{n\geq0}$ satisfies the orthogonality relations
$$
\qquad\frac{1}{2\pi i}\,\int_{\mathbb{S}^1}\widehat{B}_m^{(\alpha)}(z)\widehat{B}_n^{(\alpha)}(z)\omega_\alpha(z)\,{\rm d}z
=\frac{(-1)^{n+\alpha+1}\,2^{2n+\alpha+1}\,n!}{(2n+\alpha+1)!\,(n+\alpha+1)_{n}}\delta_{m,n}
$$
for all $m,n\in\mathbb{N}_0$, where $\mathbb{S}^1:=\{z\in\mathbb{C}:|z|=1\}$ (the unit circle).
\end{enumerate}
\medskip

\item\label{Ex-cp8-6}
Let ${\bf u}\in\mathcal{P}'$ be a regular functional fulfilling the generalized Pearson's distributional differential equation
$$2D\big((x^2+2x+1){\bf u}\big)=(-2x^2-x+1){\bf u}$$ and such that $u_1=-u_0/2$ (where, as usual, $u_n:=\langle{\bf u},x^n\rangle$, $n\in\mathbb{N}_0$). Show that ${\bf u}$ is a classical functional,
identifying ${\bf u}$ as well as the corresponding monic OPS.
Conclude that ${\bf u}$ is a positive-definite
functional uniquely represented by a positive Borel measure $\mu$ with
finite moments of all orders and $\mbox{\rm supp}(\mu)=[-1,+\infty)$.
Determine $\mu$ explicitly.
\medskip

\item\label{Ex-cp8-7}
\begin{enumerate}
\item
Prove that the integral representation (\ref{uWJ}) for the Jacobi functional
${\bf u}_J^{(\alpha,\beta)}$ is still valid provided that $\Re\alpha>-1$ and $\Re\beta>-1$.
%\begin{equation}\label{uWJ}
%\langle{\bf u}_J^{(\alpha,\beta)},p\rangle=\int_{-1}^{1}p(x)(1-x)^{\alpha}(1+x)^\beta\,{\rm d}x\;,\quad p\in\mathcal{P}\;,
%\end{equation}
\item
Prove that the integral representation (\ref{uWL}) for the Laguerre functional
${\bf u}_L^{(\alpha)}$ is still valid provided that $\Re\alpha>-1$.
\end{enumerate}
\noindent
({\sl Hint.} Use the identity principle for complex analytic functions
regarding $\alpha$ and $\beta$ as complex variables.)
\medskip

\end{enumerate}
%\medskip
}

\section*{Final remarks}

%, and \ref{MaroniZelia2001-C8}.
As we already mentioned, the (distributional) approach considered here to the classical OP is due to Pascal Maroni. This approach simplifies considerably most of the original proofs of the characterization properties presented in Theorem \ref{ThmClassicalOPS}.
The statement and proof of this theorem is based, essentially, on the articles \ref{Maroni1993-C8}, \ref{PacoAmilcarPetronilho1994-C8}, and \ref{PacoPetronilho1994-C8}.
We did not found characterizations C4 and C4$'$ (see Theorem \ref{ThmClassicalOPS}) in the available literature, for arbitrary $k$.
For $k=1$, C4 was proved in \ref{PacoAmilcarPetronilho1994-C8}.
For $k\geq2$, the proof of C4$\,\Rightarrow\,$C1 uses arguments
originally presented in %the proof of Lemma 2.1 in 
the article \ref{MaroniZelia2001-C8} by Maroni and da Rocha (see also \ref{LoureiroMaroniZelia2006-C8}).
%, except for $k=1$, case in which C4 was stated in \ref{PacoAmilcarPetronilho1994-C8}. For $k\geq2$ the proof of C4 uses ideas %contained in the article \ref{MaroniZelia2001-C8} by Maroni and da Rocha (see also \ref{LoureiroMaroniZelia2006-C8}).
It is a well known fact that any classical functional is equivalent to one of the canonical forms presented in Table \ref{Table1} (see \ref{Maroni1994-C8}, p.\;19). This fact is expressed by Theorem \ref{canonic-forms-classical}, whose explicit statement we also have not found in the literature.
%The positive-definite case considered in Section 3 is very important. Indeed,
Many authors consider that classical OPS only include Hermite, Laguerre, and Jacobi OP, with appropriate constraints on the involved parameters ensuring that their orthogonality occurs in the positive-definite sense.
The content of Section 4, about Bessel polynomials, is taken from the original article by Krall and Frink \ref{KrallFrink1955-C8} (although here we made a minor simplification in the proof of Theorem \ref{ortBessel}).
%It is worth mention that in \ref{RonveauxMawhin2005-C8} Andr\'e Ronveaux and Jean Mawhin discuss several historical aspects %concerning Rodrigues formula.

The relations presented in exercises {\bf 4} and {\bf 5} are very well known and they appear in many texts about OP.
Exercise {\bf 6} may be found in Chihara's book \ref{Chihara1978-C8}. Exercise {\bf 7} is a result presented in the article \ref{KrallFrink1955-C8} by Krall and Frink, considering the normalization adopted in Chihara's book \ref{Chihara1978-C8}.
The result expressed by exercise {\bf 9} appears in the text \ref{Maroni1994-C8} by Maroni (but notice that the hint given here leads to an alternative proof).

%Historically, the Legendre polynomials constitute the first family of OP that appeared in the literature, introduced  by Adrien %Marie Legendre (1752-1833) in a paper published in 1785 entitled ``Sur l'attraction des sph\'eroides''. Nowadays this family of %polynomials is very useful in several domains, including Number Theory, especially due to its Rodrigues formula.
%In \ref{RonveauxMawhin2005-C8} Andr\'e Ronveaux and Jean Mawhin discuss several historical aspects concerning Rodrigues formula. Jacobi OP were introduced by Carl Gustav Jacob Jacobi (1804-1851) in a work published in 1859 about {\it hypergeometric} functions, studied by Carl Friedrich Gauss (1777-1855) in a famous work on these functions. Such functions will be the subject of our next chapter.
\medskip

\section*{Bibliography}
%\medskip

{\small
\begin{enumerate}[label=\emph{\rm [\arabic*]},leftmargin=*]
%\item\label{AlSalam1990-C8} W. Al--Salam, {\it Characterization theorems for orthogonal polynomials},
%In P. Nevai Ed., Orthogonal Polynomials: Theory and Practice, Kluwer Academic Publishers, Dordrecht (1990) 1--24.
\item\label{Al-SalamChihara1972-C8} W. Al-Salam, T. S. Chihara, {\it Another characterization of the classical orthogonal polynomials}, SIAM J. Math. Anal. \textbf{3}(1) (1972) 65--70.
\item\label{Bochner1929-C8} S. Bochner, {\it Uber Sturm-Liouvillesche Polynomsysteme}, Math. Zeit. \textbf{29} (1929) 730--736.
\item\label{Chihara1978-C8} T. S. Chihara, {\sl An introduction to orthogonal polynomials}, Gordon and Breach (1978).
\item\label{Freud1971-C8} G. Freud, {\sl Orthogonal polynomials}, Pergamon Press, Oxford (1971).
\item\label{Geronimus1940-C8} Ya. L. Geronimus,
{\it On polynomials orthogonal with respect to numerical sequences and on Hahn's theorem},
Izv. Akad. Nauk. \textbf{4} (1940) 215--228. (In Russian.)
\item\label{Hahn1935-C8} W. Hahn, {\it Uber die Jacobischen polynome und zwei verwandte polynomklassen},
Math. Zeit. \textbf{39} (1935) 634--638.
\item\label{Hahn1937-C8} W. Hahn, {\it Uber h\"ohere ableitungen von orthogonal polynomen},
Math. Zeit. \textbf{43} (1937) 101.
\item\label{Ismail2004-C8} M. E. H. Ismail, {\sl Classical and Quantum Orthogonal Polynomials in One Variable},
          Cambridge University Press (2005) [paperback edition: 2009].
\item\label{Krall1936-C8} H. L. Krall, {\it On derivatives of orthogonal polynomials},
Bull. Amer. Math. Soc. \textbf{42} (1936) 867--870.
\item\label{Krall1941-C8} H. L. Krall, {\it On higher derivatives of orthogonal polynomials II},
Bull. Amer. Math. Soc. \textbf{47} (1941) 261--264.
\item\label{KrallFrink1955-C8} H. L. Krall and O. Frink,
{\it A new class of orthogonal polynomials: the Bessel polynomials},
Trans. Amer. Math. Soc. \textbf{65} (1949) 100--115.
\item\label{LoureiroMaroniZelia2006-C8} A. Loureiro, P. Maroni, and Z. da Rocha,
{\it The generalized Bochner condition about classical orthogonal polynomials revisited},
J. Math. Anal. Appl. {\bf 322} (2006) 645--667.
%\item\label{Maroni1985-C5} P. Maroni, {\it Sur quelques espaces de distributions qui sont des formes lin\'eaires sur l'espace
%vectoriel des polyn\^omes}, In C. Brezinski et al. Eds., Simposium Laguerre, Bar-le-Duc, Lecture
%Notes in Math. {\bf 1171}, Springer-Verlag (1985) 184--194.
%\item\label{Maroni1987-C5} P. Maroni, {\it Prol\'egom\`enes \`a l'\'etude des polyn\^omes orthogonaux semiclassiques}, Ann. Mat. Pura Appl. {\bf 149} (4) (1987) 165--184.
%\item\label{Maroni1988-C5} P. Maroni, {\it Le calcul des formes lin\'eaires et les polyn\^omes orthogonaux semiclassiques}, In M. Alfaro et al. Eds., Orthogonal Polynomials and Their Applications, Lecture Notes in Math.
%{\bf 1329}, Springer-Verlag (1988) 279--290.
\item\label{PacoAmilcarPetronilho1994-C8} F. Marcell\'an, A. Branquinho, and J. Petronilho,
{\it Classical orthogonal polynomials: a functional approach}, Acta Applicand\ae\, Mathematic\ae\, {\bf 34} (1994) 283--303.
\item\label{PacoPetronilho1994-C8} F. Marcell\'an and J. Petronilho,
{\it On the solution of some distributional differential equations: existence and characterizations of the classical moment functionals}, Integral Transforms and Special Functions {\bf 2} (1994) 185--218.
\item\label{Maroni1991-C8} P. Maroni, {\it Une th\'eorie alg\'ebrique des polyn\^omes orthogonaux. Applications aux polyn\^omes
orthogonaux semiclassiques}, In C. Brezinski et al. Eds., Orthogonal Polynomials and Their Applications, Proc. Erice 1990, IMACS, Ann. Comp. App. Math. {\bf 9} (1991) 95--130.
\item\label{Maroni1993-C8} P. Maroni, {\it  Variations around classical orthogonal polynomials. Connected problems}, J. Comput. Appl. Math. {\bf 48} (1993) 133--155.
\item\label{Maroni1994-C8} P. Maroni, {\sl Fonctions eul\'eriennes. Polyn\^omes orthogonaux classiques},
T\'echniques de l'Ing\'enieur, trait\'e G\'en\'eralit\'es (Sciences Fondamentales), A {\bf 154} (1994) 1--30.
\item\label{MaroniZelia2001-C8} P. Maroni and Z. da Rocha, {\it A new characterization of classical forms},
Comm. Appl. Anal. {\bf 5} (2001) 351--362.
\item\label{McCarthy1961-C8} P. J. McCarthy, {\it Characterizations of classical polynomials}. Port. Math. {\bf 20} (1961) 47--52.
\item\label{NikiforovUvarov1988-C8} A. F. Nikiforov and V. B. Uvarov,
{\sl Special Functions of Mathematical Physics}. Birkhauser Verlag, Basel (1988).
%\item\label{ReedSimon1972-C5} M. Reed and B. Simon, {\sl Methods of Modern Mathematical Physics I: Functional Analysis},
%Academic Press (1972).
%\item\label{RonveauxMawhin2005-C8} A. Ronveaux and J. Mawhin: {\it Rediscovering the contributions of Rodrigues on the
%representation of special functions}, Expo. Math. {\bf 23} (2005) 361--369.
\item\label{Szegoo1975-C8} G. Szeg\"o, {\sl Orthogonal Polynomials}, AMS Colloq. Publ. {\bf 230} (1975), 4th ed.
%\item\label{Treves1967-C5} F. Tr\`eves, {\sl Topological Vector Spaces, Distributions and Kernels}, Academic Press (1967).
\end{enumerate}
}

\chapter{Introduction to hypergeometric series}
%\chapter{The gamma and beta functions}

%\chapter{Orthogonal polynomials and dual basis.}
%More operations in $\mathcal{P}$ and $\mathcal{P}'$. }

\pagestyle{myheadings}\markright{Introduction to hypergeometric functions}
\pagestyle{myheadings}\markleft{J. Petronilho}

%In this text we introduce two basic functions which appear as fundamental tools in
%all the area of Special Functions, being also extremely useful in
%many branches of Mathematics: the gamma and the beta functions.

In this text we give a short introduction to {\it hypergeometric series and functions}.
Our presentation is mainly based in chapters 1 and 2 in the book \ref{AndrewsAskeyRoy1999-C16}
by G. Andrews, R. Askey and R. Roy (which contains much more information concerning
this topic), although in some points of the presentation we also
had supported in the books \ref{Rainville1965-C16} by Rainville,
\ref{Bailey1935-C16} by Bailey, \ref{WhittakerWatson1927-C16}
by Whittaker and Watson, \ref{Lebedev1965-C16} by Lebedev,  as well as in the
Batman Manuscript Project \ref{Erdelyi1955-C16} (directed by A. Erd\'elyi),
and Maroni's monograph \ref{Maroni1994-C16}.
The hypergeometric series (and functions) are fundamental tools in all the area of Special Functions,
being also extremely useful in many branches of Mathematics and its applications.
Before introducing such series, we need to review
two other basic functions, namely the gamma and the beta functions.

\section{The gamma and beta functions}

\begin{snugshade}
\begin{definition}[Gauss]\label{GammaF-def}
The {\sl gamma function} is defined as
\begin{equation}\label{GammaF1}
\Gamma(z):=\lim_{\substack{n\to+\infty \\ (n\in\mathbb{N})}}\frac{n!\, n^{z-1}}{(z)_n}
\;,\quad z\in\mathbb{C}\setminus\{0,-1,-2,-3,\cdots\}\;.
\end{equation}
\end{definition}
\end{snugshade}

%Historically, the gamma function was first defined by Euler, as in (\ref{GammaF4}) in bellow,
%being the notation ``$\,\Gamma(z)\,$'' introduced by Legendre in 1814.

The gamma function is a generalization of the factorial.
Indeed, assuming momentarily that the above limit exists, we may write
%, for each $z\in\mathbb{C}\setminus\{0,-1,-2,-3,\cdots\}$,
$$
\Gamma(z+1):=\lim_{n\to\infty}\frac{n!\, n^{z}}{(z+1)_n}
=\lim_{n\to\infty}\frac{n!\, n^{z-1}}{(z)_n}\frac{zn}{z+n}=z\Gamma(z)\,,
$$
hence the following property holds (difference equation for the gamma function):
\begin{snugshade}\vspace*{-0.5em}
\begin{equation}\label{GammaF2}
\Gamma(z+1)=z\Gamma(z)\;,\quad z\in\mathbb{C}\setminus\{0,-1,-2,-3,\cdots\}\;.
\end{equation}
\end{snugshade}
\noindent
In particular, and since $\Gamma(1)=\lim_{n\to\infty}\,n!\, n^{0}/(1)_n=\lim_{n\to\infty}\,n!/n!=1$,
we deduce
\begin{snugshade}\vspace*{-0.5em}
\begin{equation}\label{GammaF3}
\Gamma(n+1)=n!\;,\quad n\in\mathbb{N}_0\;.
\end{equation}
\end{snugshade}

%\noindent
Notice also that (\ref{GammaF2}) allow us write the following useful identity
\begin{snugshade}\vspace*{-0.5em}
\begin{equation}\label{Poch-Gamma}
(a)_n=\frac{\Gamma(n+a)}{\Gamma(a)}\;,\quad
a\in\mathbb{C}\setminus\{0,-1,-2,-3,\cdots\}\;,\;\;n\in\mathbb{N}_0\;.
\end{equation}
\end{snugshade}

The next theorem shows that indeed the limit defining the gamma function exists.
We need to recall the definition of the {\sl Euler-Mascheroni}
constant\,\footnote{\,This constant is often referred as {\sl Euler's constant}. Its arithmetic nature --- to know whether $\gamma$ is a rational or transcendental number --- is unknown. It is conjectured that ``\,$\gamma\not\in\mathbb{Q}\,$'' (indeed, it is expected that ``\,$\gamma$ is a transcendental number\,'', but a proof (or disproof) has been resisting along the times. This is an old and important conjecture in Number Theory, that fits into the class of problems related with Hilbert's seventh problem appearing in the famous list of open problems presented by David Hilbert on the occasion of the International Congress of Mathematics held in Paris in 1900.}:
\begin{snugshade}\vspace*{-0.5em}
\begin{equation}\label{Euler-Mascheroni}
\gamma:=\lim_{n\to+\infty}\left(1+\frac12+\cdots+\frac{1}{n}-\ln n\right)=0.5772156\ldots\;.
\end{equation}
\end{snugshade}
\noindent
Notice that this limit exists.
In fact, setting $u_n:=\int_0^1\frac{t}{n(n+t)}\,{\rm d}t=\frac{1}{n}-\ln\frac{n+1}{n}$, we have
$0<u_n\leq\int_0^1\frac{1}{n(n+0)}\,{\rm d}t=\frac{1}{n^2}$ for each $n\in\mathbb{N}$, hence
$\sum_{n=1}^\infty u_n$ converges, and so
$$
1+\frac12+\cdots+\frac{1}{n}-\ln n=\sum_{k=1}^nu_k+\ln\frac{n+1}{n}
\xlongrightarrow[n\rightarrow+\infty]{}\sum_{k=1}^\infty u_k=\gamma\;.
$$

\begin{snugshade}
\begin{theorem}\label{GammaF-exist}
The limit $(\ref{GammaF1})$ exists and is never zero.
Moreover, $\Gamma$ is an analytic function in all its domain $\mathbb{C}\setminus\{0,-1,-2,-3,\cdots\}$,
with simple poles at the points $0,-1,-2,-3,\ldots$.
In addition, the identities
\begin{eqnarray}
\label{GammaF4}
\qquad\Gamma(z)=\frac{1}{z}\prod_{n=1}^\infty\left[\left(1+\frac{z}{n}\right)^{-1}\left(1+\frac{1}{n}\right)^z\,\right]
\qquad\qquad\mbox{\rm (Euler)} \\
\label{GammaF4a}
\qquad\frac{1}{\Gamma(z)}=ze^{\gamma z}\prod_{n=1}^\infty\left[\left(1+\frac{z}{n}\right)e^{-z/n}\,\right]
\qquad\qquad\mbox{\rm (Schl\"omilch)}%, 1844)
\end{eqnarray}
hold for each $z\in\mathbb{C}\setminus\{0,-1,-2,-3,\cdots\}$,
where $\gamma$ is the Euler-Mascheroni constant.
%\begin{equation}\label{Euler-Mascheroni}
%\gamma:=\lim_{n\to+\infty}\left(1+\frac12+\cdots+\frac{1}{n}-\ln n\right)=0.5772156\ldots\;.
%\end{equation}
\end{theorem}
\end{snugshade}

{\it Proof.}\footnote{\,We present a proof that does not assume knowledge of the theory of
infinite products, following the exposition at the begin of chapter XII in Whittaker and Watson's book \ref{WhittakerWatson1927-C16}.
(Indeed, assuming some basic facts concerning this theory, a more concise proof could be done.)}
%Recall that $(z)_0:=1$ and $(z)_n:=z(z+1)(z+2)\cdots(z+n-1)$ if $n\in\mathbb{N}$.
Let $N\in\mathbb{N}$ and $z\in\mathbb{C}$ such that $|z|\leq\frac12\,N$.
Recall that, taking the principal value of $\log(1+w)$,
we have $\log(1+w)=\sum_{k=1}^\infty\frac{(-1)^{k-1}}{k}w^k$
if $|w|<1$, hence %, for each $n>N$,
%Taking the principal value of $\log\left(1+\frac{z}{n}\right)$,
%we have $\log\left(1+\frac{z}{n}\right)=\sum_{k=1}^\infty\frac{(-1)^{k-1}}{k}\left(\frac{z}{n}\right)^k$
%if $\left|\frac{z}{n}\right|<1$, hence %, for each $n>N$,
%(so $\left|\frac{z}{n}\right|<\frac12$),
%$$
%\left|\log\left(1+\frac{z}{n}\right)-\frac{z}{n}\right|=
%\left|-\frac12\left(\frac{z}{n}\right)^2+\frac13\left(\frac{z}{n}\right)^3-\cdots\right|
%\leq\frac{|z|^2}{n^2}\left(1+\left|\frac{z}{n}\right|+\left|\frac{z^2}{n^2}\right|+\cdots\right)
%$$
$$
\left|\log\left(1+\frac{z}{n}\right)-\frac{z}{n}\right|\leq
\frac{|z|^2}{n^2}\sum_{k=0}^\infty\left|\frac{z}{n}\right|^k\leq
\frac{N^2}{4n^2}\sum_{k=0}^\infty\left(\frac12\right)^k=\frac{N^2}{2n^2}\quad\mbox{\rm if}\quad n>N\;.
$$
Since the series $\sum_{n=N+1}^\infty\frac{N^2}{2n^2}$ is convergent,
then Weierstrass $M-$test ensures that
$\Sigma_N(z):=\sum_{n=N+1}^\infty\left[\log\left(1+\frac{z}{n}\right)-\frac{z}{n}\right]$
is an absolutely and uniformly convergent series in the region $|z|\leq\frac12\,N$,
and so, since its terms are analytic functions in this region,
then $\Sigma_N$ is an analytic function in the same region. %, for each $N\in\mathbb{N}$.
Consequently, its exponential
$$
\Lambda_N(z):=e^{\Sigma_N(z)}
=\lim_{m\to\infty}\prod_{n=N+1}^m\left[\left(1+\frac{z}{n}\right)e^{-\frac{z}{n}}\,\right]=:
\prod_{n=N+1}^\infty\left[\left(1+\frac{z}{n}\right)e^{-\frac{z}{n}}\,\right]
$$
is an analytic function in the region $|z|\leq\frac12\,N$ which is never zero there
(since it is the exponential of a finite complex value, for each $N$ and $z$),
and so
$$
\Lambda_N(z) \,\prod_{n=1}^N\left[\left(1+\frac{z}{n}\right)e^{-\frac{z}{n}}\,\right]
=\lim_{m\to\infty}\prod_{n=1}^m\left[\left(1+\frac{z}{n}\right)e^{-\frac{z}{n}}\,\right]=:
\prod_{n=1}^\infty\left[\left(1+\frac{z}{n}\right)e^{-\frac{z}{n}}\,\right]=:\Lambda(z)
$$
is an analytic function in the region $|z|\leq\frac12\,N$ which fulfills
$\Lambda(z)\neq0$ for each $z$ in this region
% such that $z\neq-1,-2,\ldots,-N$.
that does not coincide with a nonnegative integer number.
Therefore, since we can take $N$ arbitrarily large, we conclude that $\Lambda$
is analytic in $\mathbb{C}$ (an entire function) and fulfills
$\Lambda(z)\neq0$ for each $z\in\mathbb{C}\setminus\{-1,-2,-3,\ldots\}$.
Clearly, the zeros of $\Lambda(z)$ are precisely the numbers
$-1,-2,-3,\ldots$, which are simple zeros.
%(note that the factor $\prod_{n=1}^N\left[\left(1+\frac{z}{n}\right)e^{-z/n}\,\right]$
%appearing in the definition of $\Lambda$ has its zeros at the points $z=-1,-2,\ldots,-N$).
Now, we may write
$$
\begin{array}{rcl}
e^{\gamma z}\Lambda(z)&=&\displaystyle
\lim_{n\to\infty}e^{\left(1+\frac12+\cdots+\frac{1}{n}-\ln n\right)z}\cdot
\lim_{n\to\infty}\prod_{j=1}^n\left[\left(1+\frac{z}{j}\right)e^{-\frac{z}{j}}\,\right] \\ [1em]
&=&\displaystyle
\lim_{n\to\infty}\left\{e^{\left(1+\frac12+\cdots+\frac{1}{n}-\ln n\right)z}\cdot
\prod_{j=1}^n\left[\left(1+\frac{z}{j}\right)e^{-\frac{z}{j}}\,\right]\right\} \\ [1em]
&=&\displaystyle
\lim_{n\to\infty}\,n^{-z}\cdot\prod_{j=1}^n\left(1+\frac{z}{j}\right)\;,
%\,=\,\lim_{m\to\infty}\,\left[\,\prod_{j=1}^{m-1}\left(1+\frac{1}{j}\right)^{-z}
%\cdot\prod_{j=1}^m\left(1+\frac{z}{j}\right)\,\right]
\end{array}
$$
%$$
%\begin{array}{rcl}
%e^{\gamma z}\Lambda(z)&=&\displaystyle
%\lim_{n\to\infty}e^{\left(1+\frac12+\cdots+\frac{1}{n}-\ln n\right)z}\cdot
%\lim_{n\to\infty}\prod_{j=1}^n\left[\left(1+\frac{z}{j}\right)e^{-\frac{z}{j}}\,\right] \\ [1em]
%&=&\displaystyle
%\lim_{n\to\infty}\left\{e^{\left(1+\frac12+\cdots+\frac{1}{n}-\ln n\right)z}\cdot
%\prod_{j=1}^n\left[\left(1+\frac{z}{j}\right)e^{-\frac{z}{j}}\,\right]\right\}
%=\lim_{n\to\infty}\,n^{-z}\cdot\prod_{j=1}^n\left(1+\frac{z}{j}\right)\;,
%\end{array}
%$$
and since, for each $n\in\mathbb{N}$ and $z\in\mathbb{C}\setminus\{-1,-2,-3,\ldots,-(n-1)\}$,
\begin{equation}\label{GammaF5a2}
\frac{n!\, n^{z-1}}{(z)_n}=\frac{1}{z}\frac{n^z}{\frac{1+z}{1}\frac{2+z}{2}\cdots\frac{n-1+z}{n-1}}
=\frac{1}{z}\left(1+\frac{z}{n}\right)\,\Big[\,n^{-z}\,\prod_{j=1}^n\left(1+\frac{z}{j}\right)\Big]^{-1}\;,
\end{equation}
we conclude that there exists, and it is nonzero, the limit
%$$
%\lim_{n\to+\infty}\frac{n!\, n^{z-1}}{(z)_n}=\frac{1}{z e^{\gamma z}\Lambda(z)}
%\,,\quad z\in\mathbb{C}\setminus\{0,-1,-2,-3,\cdots\}\;.
%$$
\begin{equation}\label{GammaF5a}
\Gamma(z):=\lim_{n\to+\infty}\frac{n!\, n^{z-1}}{(z)_n}=\frac{1}{z e^{\gamma z}\Lambda(z)}\,,\quad z\in\mathbb{C}\setminus\{0,-1,-2,-3,\cdots\}\;.
\end{equation}
Moreover, taking into account the facts proved above about the function $\Lambda$,
it follows immediately that $\Gamma$ is an analytic function on $\mathbb{C}\setminus\{0,-1,-2,-3,\cdots\}$
and it has simple poles at the points $0,-1,-2,-3,\cdots$.
Notice that (\ref{GammaF5a}) also proves (\ref{GammaF4a}).
Finally, for each $z\in\mathbb{C}\setminus\{-1,-2,-3,\ldots,-(n-1)\}$,
(\ref{GammaF5a2}) can be rewritten as
$$
\frac{n!\, n^{z-1}}{(z)_n}=\left(\frac{n}{n+1}\right)^z\left(1+\frac{z}{n}\right)\,
\frac{1}{z}\,\prod_{j=1}^n\left[\left(1+\frac{z}{j}\right)^{-1}\left(1+\frac{1}{j}\right)^z\right]\;,
$$
hence taking the limit as $n\to\infty$ we obtain (\ref{GammaF4}).
\qed
\medskip

\begin{remark}\em
Historically, the gamma function was first defined by Euler, as in (\ref{GammaF4}) above,
being the notation ``$\,\Gamma(z)\,$'' introduced by Legendre in 1814.
\end{remark}
%\medskip

%\section{The beta function}

%Next we introduce the beta integral and series.

\begin{snugshade}
\begin{definition}\label{BetaF-def}
The {\sl beta integral} is
\begin{equation}\label{BetaF1}
B(x,y):=\int_0^1 t^{x-1}(1-t)^{y-1}\,{\rm d}t
\;,\quad \Re x>0\;,\;\; \Re y>0\;.
\end{equation}
The {\sl beta function} is obtained from the beta integral by analytic continuation, and
we still denote it by $B(x,y)$.
\end{definition}
\end{snugshade}

Notice that the integral in (\ref{BetaF1}) is symmetric in $x$ and $y$, i.e.,
\begin{snugshade}\vspace*{-0.5em}
\begin{equation}\label{BetaF2}
B(x,y)=B(y,x)
\;,\quad \Re x>0\;,\;\; \Re y>0\;.
\end{equation}
\end{snugshade}
\noindent
This identity follows immediately from (\ref{BetaF1}) making the change of variables $t=1-s$.

\begin{snugshade}
\begin{theorem}\label{BetaF-Gamma}
The beta function fulfills
\begin{equation}\label{BetaF3}
B(x,y)=\frac{\Gamma(x)\Gamma(y)}{\Gamma(x+y)}
\qquad \big(x,y,x+y\in\mathbb{C}\setminus\{0,-1,-2,-3,\cdots\}\big)\;.
\end{equation}
%for all $x$ and $y$ such that $\;x,y,x+y\in\mathbb{C}\setminus\{0,-1,-2,-3,\cdots\}$.
\end{theorem}
\end{snugshade}

{\it Proof.}
If $\Re x>0$ and $\Re y>0$, the beta function $B(x,y)$ is given by (\ref{BetaF1}).
On the other hand, by Theorem \ref{GammaF-exist},
$\Gamma(x)$, $\Gamma(y)$, and $\Gamma(x+y)$ are well defined and non-zero
for all $x$ and $y$ such that $\;x,y,x+y\in\mathbb{C}\setminus\{0,-1,-2,-3,\cdots\}$.
Thus the right-hand side of (\ref{BetaF3}) is well defined, and we
only need to prove (\ref{BetaF3}) for $x$ and $y$ such that $\Re x>0$ and $\Re y>0$,
since --- taking into account Theorem \ref{GammaF-exist} again --- the right-hand side of (\ref{BetaF3})
provides the analytic continuation of the beta integral.
Assuming $\Re x>0$ and $\Re y>0$, we start by proving that $B$ fulfills the functional equation
\begin{equation}\label{BetaF4}
B(x,y)=\frac{x+y}{y}\,B(x,y+1)\;. %,\quad \Re x>0\;,\;\; \Re y>0\;.
\end{equation}
By (\ref{BetaF1}), we may write %for $\Re x>0$ and $\Re y>0$,
$B(x,y+1)%=\int_0^1 t^{x-1}(1-t)(1-t)^{y-1}\,{\rm d}t
=\int_0^1 t^{x-1}(1-t)^{y-1}\,{\rm d}t-\int_0^1 t^{x}(1-t)^{y-1}\,{\rm d}t$,
i.e.,
\begin{equation}\label{BetaF5}
B(x,y+1)=B(x,y)-B(x+1,y)\;.
\end{equation}
On the other hand, integration by parts yields
$B(x,y+1)=\int_0^1 t^{x-1}(1-t)^{y}\,{\rm d}t
=\left.\frac{t^x}{x}(1-t)^y\right|_{t=0}^1+\frac{y}{x}\int_0^1 t^{x}(1-t)^{y-1}\,{\rm d}t
=\frac{y}{x}B(x+1,y)$,
hence
\begin{equation}\label{BetaF6}
B(x+1,y)=\frac{x}{y}\,B(x,y+1)\;.
\end{equation}
Inserting (\ref{BetaF6}) in the right-hand side of (\ref{BetaF5}) proves (\ref{BetaF4}).
Iterating (\ref{BetaF4}), % yields
\begin{equation}\label{BetaF7}
B(x,y)=\frac{(x+y)_n}{(y)_n}\,B(x,y+n)\;,\quad n\in\mathbb{N}\;.
\end{equation}
One sees (making the change of variables $t=s/n$) that
$$
B(x,y+n)=\int_0^1 t^{x-1}(1-t)^{y+n-1}\,{\rm d}t=\frac{1}{n^x}
\int_0^n s^{x-1}\left(1-\frac{s}{n}\right)^{y+n-1}\,{\rm d}s\;,
$$
and so (\ref{BetaF7}) may be rewritten as
\begin{equation}\label{BetaF8}
B(x,y)=\frac{(x+y)_n}{n!\,n^{x+y-1}}\,\frac{n!\,n^{y-1}}{(y)_n}\,
\int_0^n s^{x-1}\left(1-\frac{s}{n}\right)^{y+n-1}\,{\rm d}s\;,\quad n\in\mathbb{N}\;.
\end{equation}
Now, by definition of the gamma function, we have
\begin{equation}\label{BetaF9}
\lim_{n\to\infty}\,\frac{(x+y)_n}{n!\,n^{x+y-1}}=\frac{1}{\Gamma(x+y)}\;,\qquad
\lim_{n\to\infty}\,\frac{n!\,n^{y-1}}{(y)_n}=\Gamma(y)\;.
\end{equation}
Moreover, using (for instance) Lebesgue's dominated convergence theorem
and the connections between the Lebesgue and the Riemann integrals,
we deduce (Exercise \ref{Ex-cp16-1})
\begin{equation}\label{BetaF10}
\lim_{n\to\infty}\,\int_0^n s^{x-1}\left(1-\frac{s}{n}\right)^{y+n-1}\,{\rm d}s
=\int_0^{+\infty} s^{x-1}e^{-s}\,{\rm d}s\;.
\end{equation}
Therefore, taking the limit as $n\to+\infty$ in (\ref{BetaF8}),
from (\ref{BetaF9}) and (\ref{BetaF10}) we obtain
\begin{equation}\label{BetaF11}
B(x,y)=\frac{\Gamma(y)}{\Gamma(x+y)}\,
\int_0^{+\infty} t^{x-1}e^{-t}\,{\rm d}t\;.
\end{equation}
Taking $y=1$, and since $\Gamma(1)=1$, we deduce
\begin{equation}\label{BetaF12}
\int_0^{+\infty} t^{x-1}e^{-t}\,{\rm d}t=\Gamma(x+1)B(x,1)
=x\,\Gamma(x)\,\int_0^{1} t^{x-1}\,{\rm d}t=\Gamma(x)\;.
\end{equation}
Therefore, inserting (\ref{BetaF12}) into (\ref{BetaF11}) gives (\ref{BetaF3})
for $\Re x>0$ and $\Re y>0$.
Thus by the considerations at the begin of the proof, (\ref{BetaF3}) follows
for all $x$ and $y$ such that $\;x,y,x+y\in\mathbb{C}\setminus\{0,-1,-2,-3,\cdots\}$.
\qed
\medskip

\begin{snugshade}
\begin{corollary}\label{FGamma+}
The gamma function fulfills
\begin{equation}\label{FGamma+1}
\Gamma(x)=\int_0^{+\infty} t^{x-1}e^{-t}\,{\rm d}t\;,
\qquad \Re x>0\;.
\end{equation}
\end{corollary}
\end{snugshade}

{\it Proof.}
Indeed, (\ref{FGamma+1}) is relation (\ref{BetaF12}) stated in the proof of Theorem \ref{BetaF-Gamma}.
\qed

\begin{snugshade}
\begin{corollary}[Euler's reflection formula]\label{EulerReflection}
The gamma function fulfills
\begin{equation}\label{FGammaReflect}
\Gamma(z)\Gamma(1-z)=\frac{\pi}{\sin(\pi z)}\;,
\qquad z\in\mathbb{C}\setminus\mathbb{Z}\;.
\end{equation}
\end{corollary}
\end{snugshade}

{\it Proof.}
The proof is left to the reader (Exercise \ref{Ex-cp16-2}).
\qed

\section{Hypergeometric series}

%\subsection{Definition and analysis of the convergence}

%A series $\sum_{n=0}^\infty c_n$ is called {\it hypergeometric} if
An {\it hypergeometric series} is a series $\sum_{n=0}^\infty c_n$ where
\begin{snugshade}%\vspace*{-0.5em}
$$\frac{c_{n+1}}{c_n}\quad\mbox{\rm is a rational function of $n$.}$$
\end{snugshade}
\noindent
On factorizing the polynomials in $n$, we may write
\begin{equation}\label{HF1}
\frac{c_{n+1}}{c_n}=\frac{(n+a_1)(n+a_2)\cdots(n+a_p)}{(n+b_1)(n+b_2)\cdots(n+b_q)}\frac{x}{n+1}\;,
\end{equation}
where $x$ is a complex number (which appears because the polynomials may be non monic) and
$a_j$ and $b_j$ are complex parameters such that $b_j\in\mathbb{C}\setminus\{0,-1,-2,-3,\ldots\}$.
Therefore, for each $n\in\mathbb{N}$,
$$
c_{n}=c_{n-1}\,\frac{(a_1+n-1)(a_2+n-1)\cdots(a_p+n-1)}{(b_1+n-1)(b_2+n-1)\cdots(b_q+n-1)}\frac{x}{n}\;,
$$
and by iterating this relation we obtain
\begin{equation}\label{HF2}
c_{n}=c_{0}\,\frac{(a_1)_n(a_2)_n\cdots(a_p)_n}{(b_1)_n(b_2)_n\cdots(b_q)_n}\frac{x^n}{n!}\;,\quad n\in\mathbb{N}_0\;.
\end{equation}

Thus (up to a constant factor) an {\sl hypergeometric series} is a series of the form
%form\,\footnote{\,Of course, the factor $c_0$ can be inserted into $x^n$.}
\begin{snugshade}\vspace*{-0.5em}
\begin{equation}\label{HF3}
\pFq{p}{q}{a_1,\cdots,a_p}{b_1,\cdots,b_q}{x}:=
\sum_{n=0}^\infty \frac{(a_1)_n(a_2)_n\cdots(a_p)_n}{(b_1)_n(b_2)_n\cdots(b_q)_n}\frac{x^n}{n!}\;,
\end{equation}
\end{snugshade}
\noindent
being $x\in\mathbb{C}$ and, for all possible $j$,
\begin{snugshade}\vspace*{-0.5em}
\begin{equation}\label{HF4}
a_j\in\mathbb{C}\;,\quad b_j\in\mathbb{C}\setminus\{0,-1,-2,-3,\ldots\}\;.
\end{equation}
\end{snugshade}

\begin{remark}\em
Often, instead of the left-hand side of (\ref{HF3}), the notations
\begin{snugshade}\vspace*{-0.5em}
\begin{equation}\label{HF5}
{}_{p}F_{q}(a_1,\ldots,a_p;\,b_1,\ldots,b_q;\,x)\;,\quad
{}_{p}F_{q}(x)\;,\quad {}_{p}F_{q}
\end{equation}
\end{snugshade}
\noindent
are used, provided concerning the last two ones there is no danger of misunderstanding.
Moreover, it may happens that in the numerator or in the denominator (or in both)
of the fraction defining the general term of an hypergeometric series, no parameters $a_j$ or $b_j$ appear
(this situation takes place when the number of corresponding parameters is $p=0$ or $q=0$, respectively).
In this situation we write ``\,$\mbox{---}$'' instead of the $a_j$ or $b_j$ parameters, to indicate their absence.
For instance,
$$
\pFq{0}{1}{\mbox{---}}{b_1}{x}=\sum_{n=0}^\infty\frac{x^n}{(b_1)_n\,n!}\;,\quad
\pFq{3}{0}{a_1,a_2,a_3}{\mbox{---}}{x}=\sum_{n=0}^\infty\frac{(a_1)_n(a_2)_n(a_3)_n}{n!}\,x^n\;.
$$
\end{remark}

\begin{remark}\em
Notice that if $a_j$ is zero or a negative integer number for some $j\in\{1,2,\ldots,p\}$,
then the series on the right-hand side of (\ref{HF3}) terminates,
i.e., it reduces to a finite sum (and thus it becomes a polynomial in the variable $x$).
In such a case, (\ref{HF3}) is called a {\sl terminating hypergeometric series}.
\end{remark}

Next we analyze the convergence of the hypergeometric series.

\begin{snugshade}
\begin{theorem}\label{pFq-conv}
Let ${}_{p}F_{q}$ be the hypergeometric series defined by $(\ref{HF3})$. Then:
\begin{enumerate}
\item[{\rm (i)}] if $p\leq q$, then ${}_{p}F_{q}(x)$ converges absolutely for each $x\in\mathbb{C}$;
\item[{\rm (ii)}] if $p=q+1$, then ${}_{p}F_{q}(x)$ converges absolutely if $|x|<1$,
and it diverges if $|x|>1$ and the series does not terminates.
%for each $x\in\mathbb{D}:=\{z\in\mathbb{C}:|z|<1\}$;
\item[{\rm (iii)}] if $p>q+1$, then ${}_{p}F_{q}(x)$ diverges for each $x\in\mathbb{C}\setminus\{0\}$,
provided that the series does not terminates.
\end{enumerate}
\end{theorem}
\end{snugshade}

{\it Proof.}
We will apply the ratio test.
We may assume that the series does not terminates (otherwise it converges trivially).
By (\ref{HF1}), we may write
\begin{equation}\label{HF1a}
\left|\frac{c_{n+1}}{c_n}\right|=
\frac{\left|1+\frac{a_1}{n}\right|\cdots\left|1+\frac{a_p}{n}\right|n^{p-(q+1)}}
{\left|1+\frac{b_1}{n}\right|\cdots\left|1+\frac{b_q}{n}\right|\,\left|1+\frac{1}{n}\right|}\,|x|\;,
\quad n\in\mathbb{N}_0\;.
\end{equation}
Therefore, the following holds:
\smallskip

(i) Suppose $p\leq q$.
Then, by (\ref{HF1a}),
$\lim_{n\to+\infty}\left|\frac{c_{n+1}}{c_n}\right|=0<1$,
hence, by the ratio test, the series defining
${}_{p}F_{q}(x)$ converges absolutely for each $x\in\mathbb{C}$.
%\smallskip

(ii) Suppose $p=q+1$.
Then, by (\ref{HF1a}),
$\lim_{n\to+\infty}\left|\frac{c_{n+1}}{c_n}\right|=|x|$,
hence the series ${}_{p}F_{q}(x)$ converges absolutely if $|x|<1$,
and it diverges if $|x|>1$.
%\smallskip

(iii) Suppose $p>q+1$. %(and the series does not terminates).
By (\ref{HF1a}),
$\lim_{n\to+\infty}\left|\frac{c_{n+1}}{c_n}\right|=+\infty$,
hence the series ${}_{p}F_{q}(x)$ diverges for each $x\in\mathbb{C}\setminus\{0\}$.
\qed\smallskip

The case whether $|x|=1$ when $p=q+1$ is of great interest.
The next theorem gives the conditions for convergence in this case.
Its proof requires the following %convergence test for positive series.

\begin{snugshade}
\begin{lemma}[Gauss's test]\label{Gauss-test}
Let $\{a_n\}_{n\geq1}$ be a sequence of positive numbers.
Suppose that there exist $r>1$, $N\in\mathbb{N}$, and a bounded sequence $\{C_n\}_{n\geq 1}$
such that
\begin{equation}\label{GaussT1}
\frac{a_n}{a_{n+1}}=1+\frac{\epsilon}{n}+\frac{C_n}{n^r}\;,\quad n\geq N\;.
\end{equation}
Then the series $\sum_{n=1}^\infty a_n$ is convergent if $\epsilon>1$, and it is divergent if $\epsilon\leq1$.
\end{lemma}
\end{snugshade}

%{\it Proof.}
%We prove the result only for the case $\epsilon\leq1$
%(in fact, we only will need this case for proving
%Theorem \ref{pFq-conv-thm2}).\footnote{\,As it is well known, Gauss's test is a consequence of Raabes's test if $\epsilon<1$,
%and of Bertrand's test if $\epsilon=1$. Here we adapt the proofs of these tests to give
%a ``direct'' proof of Gauss's test, in order to maintain the proof of
%Theorem \ref{pFq-conv-thm2} more self contained.}
%The case $\epsilon>1$ can be proved using similar arguments
%(in this case it may be useful also using the binomial theorem).
{\it Proof.}\footnote{\,Gauss's test is a consequence of Raabes's test if $\epsilon\neq1$,
and of Bertrand's test if $\epsilon=1$. Here we adapt the proofs of these tests to give
a ``direct'' proof of Gauss's test, in order to maintain the proof of
Theorem \ref{pFq-conv-thm2} more self contained.}
Assume first $\epsilon<1$.
%Since $\{C_n\}_{n\geq1}$ is bounded, then there exists $C>0$ such that $|C_n|\leq C$ for each $n\in\mathbb{N}$,
Since $\{C_n\}_{n\geq1}$ is bounded and $r>1$, there exists an integer $N_1\geq N$ such that
$|C_n|/n^{r-1}\leq(1-\epsilon)/2$ if $n\geq N_1$.
Therefore, from (\ref{GaussT1}),
$\frac{a_n}{a_{n+1}}\leq 1+\frac{1+\epsilon}{2}\frac{1}{n}<\frac{n+1}{n}$ for $n\geq N_1$,
hence, $(n+1)a_{n+1}>na_n$ for each $n\geq N_1$.
By repeatedly application of this inequality we deduce
$$
a_n\geq\frac{c_1}{n}\quad\mbox{\rm if}\quad n\geq N_1\,,
$$
where $c_1:=N_1a_{N_1}>0$. Thus, since the series $\sum_{n=1}^\infty\frac{1}{n}$ is divergent,
then so is $\sum_{n=1}^\infty a_n$.
Assume now $\epsilon=1$.
Since $r>1$ then %(by L'H\^opital rule, e.g.)
$\frac{\ln x}{x^{r-1}}\to0$ as $x\to+\infty$,
hence, since $\{C_n\}_{n\geq1}$ is bounded, also
$C_n\ln n/n^{r-1}\to0$ as $n\to+\infty$.
%$\lim_{n\to\infty}\frac{C_n\ln n}{n^{r-1}}=0\;.$
Consequently, there  exists an integer $N_2> N$ such that
$C_n\ln n/n^{r-1}\leq1$ if $n\geq N_2$, and so, from (\ref{GaussT1}),
\begin{equation}\label{GaussT2}
\frac{a_n}{a_{n+1}}\leq1+\frac{1}{n}+\frac{1}{n\ln n}
\leq\frac{(n+1)\ln(n+1)}{n\ln n}\quad\mbox{\rm if}\quad n\geq N_2\,.
\end{equation}
The second inequality in (\ref{GaussT2}) holds since it is equivalent to the inequality $1\leq f(n)$,
being $f(x):=(x+1)\ln\frac{x+1}{x}$; and this last inequality holds since $f$ is (strictly)
decreasing on $(0,+\infty)$, and so $f(n)\geq\lim_{x\to+\infty}f(x)=1$ for each $n>1$.
From (\ref{GaussT2}),
$$
a_n\geq\frac{c_2}{n\ln n}\quad\mbox{\rm if}\quad n\geq N_2\,,
$$
where $c_2:=a_{N_2}N_2\ln N_2>0$, and since the series $\sum_{n=2}^\infty\frac{1}{n\ln n}$ is divergent
(use the integral test:
$\int_2^{+\infty}\frac{{\rm d}x}{x\ln x}=\left.\ln\ln x\right|_2^{+\infty}=+\infty$),
then so is $\sum_{n=1}^\infty a_n$.
Finally, assume $\epsilon>1$.
Arguing as before, there exists $N_1'\geq N$ such that
$|C_n|/n^{r-1}\leq(\epsilon-1)/2$ if $n\geq N_1'$.
Hence, setting $q:=(\epsilon+1)/2$ and taking $s$ such that $1<s<q$, we deduce
\begin{equation}\label{GaussT3}
\frac{a_n}{a_{n+1}}\geq 1+\frac{\epsilon}{n}-\frac{\epsilon-1}{2}\frac{1}{n}=1+\frac{q}{n}\geq\left(1+\frac{1}{n}\right)^s
\quad\mbox{\rm if}\quad n\geq N_3\,,
\end{equation}
being $N_3$ an integer chosen so that $N_3\geq\max\big\{N_1',\frac{2^{s+2}}{q-s}\big\}$.
The last inequality in (\ref{GaussT3}) holds by the binomial theorem\,\footnote{\,The binomial theorem states:
$(1+z)^\alpha=\sum_{k=0}^\infty{{\alpha}\choose{k}}z^k=\sum_{k=0}^\infty\frac{(-\alpha)_k}{k!}(-z)^k$ if $\alpha,z\in\mathbb{C}$, $|z|<1$.},
which allow us writing, for each $n\geq2$,
$$
\left(1+\frac{1}{n}\right)^s
=1+\frac{s}{n}+\sum_{k=2}^\infty{{s}\choose{k}}\left(\frac{1}{n}\right)^k
\leq1+\frac{s}{n}+\frac{4}{n^2}\sum_{k=0}^\infty\frac{(s)_k}{k!}\left(\frac{1}{2}\right)^k
=1+\frac{s}{n}+\frac{2^{s+2}}{n^2}
$$
(where we have used the inequality $\big|{{\alpha}\choose{k}}\big|\leq(|\alpha|)_k/k!$,
valid for all $\alpha\in\mathbb{C}$ and $k\in\mathbb{N}$),
and so the last inequality in (\ref{GaussT3}) follows taking into account that
$1+\frac{q}{n}=1+\frac{s}{n}+\frac{(q-s)n}{n^2}$.
From (\ref{GaussT3}) we obtain $(n+1)^sa_{n+1}\leq n^sa_n$ if $n\geq N_3$, hence
$$
a_n\leq\frac{c_3}{n^s}\quad\mbox{\rm if}\quad n\geq N_3\,,
$$
where $c_3:=a_{N_3}N_3^s>0$, and since $\sum_{n=1}^\infty\frac{1}{n^s}$ is convergent,
then so is $\sum_{n=1}^\infty a_n$.
\qed
\medskip

Before proving the theorem we also point out the following fact\,\footnote{\,Recall that,
given two sequences of real or complex numbers $\{a_n\}_{n\geq0}$ and $\{b_n\}_{n\geq0}$,
the notation ``$\, a_n\sim b_n$ as $n\to+\infty\,$'' means that $a_n/b_n\to1$ as $n\to+\infty$.}:
{\it the coefficient of $x^n$ in the series ${}_{q+1}F_{q}(x)$ is}
\begin{snugshade}\vspace*{0.25em}
\begin{equation}\label{HF6}
\frac{(a_1)_n(a_2)_n\cdots(a_{q+1})_n}{(b_1)_n(b_2)_n\cdots(b_q)_n\, n!}
\sim\frac{\Gamma(b_1)\cdots\Gamma(b_q)}{\Gamma(a_1)\cdots\Gamma(a_{q+1})}\,n^{\sum a_j-\sum b_j-1}\;,\;\;
\mbox{as $\;n\to+\infty$}\;.
\end{equation}
\end{snugshade}
\noindent
(Here we use the abbreviations $\sum a_j:=\sum_{j=1}^{q+1}a_j$ and $\sum b_j:=\sum_{j=1}^{q}b_j$.)
Relation (\ref{HF6}) follows at once from the definition (\ref{GammaF1}) of the gamma function, which gives
$$
(z)_n\sim\frac{n!\,n^{z-1}}{\Gamma(z)}\;,\;\; \mbox{as $\;n\to+\infty$}\;.
$$
%\begin{snugshade}
%\begin{lemma}\label{pFq-conv-lemma}
%The coefficient of $x^n$ in the series ${}_{q+1}F_{q}(x)$ is
%\begin{equation}\label{HF6}
%\frac{(a_1)_n(a_2)_n\cdots(a_{q+1})_n}{(b_1)_n(b_2)_n\cdots(b_q)_n\, n!}
%\sim\frac{\Gamma(b_1)\cdots\Gamma(b_q)}{\Gamma(a_1)\cdots\Gamma(a_{q+1})}\,n^{\sum a_j-\sum b_j-1}\;,\;\;
%\mbox{as $\;n\to+\infty$}\;.
%\end{equation}
%\end{lemma}
%\end{snugshade}

\begin{snugshade}
\begin{theorem}\label{pFq-conv-thm2}
Let $|x|=1$ and $p=q+1$ in the hypergeometric series defined by $(\ref{HF3})$,
and suppose that it is a nonterminating series.
\begin{enumerate}
\item[{\rm (i)}] if $\Re\big(\sum a_j-\sum b_j\big)<0$, then ${}_{q+1}F_{q}(x)$ converges absolutely;
\item[{\rm (ii)}] if $0\leq\Re\big(\sum a_j-\sum b_j\big)<1$ and $x\neq1$, then ${}_{q+1}F_{q}(x)$ converges conditionally;
\item[{\rm (iii)}] if $\Re\big(\sum a_j-\sum b_j\big)\geq1$, then ${}_{q+1}F_{q}(x)$ diverges.
%, provided it is nonterminating.
\end{enumerate}
\end{theorem}
\end{snugshade}

{\it Proof.}
Since $|x|=1$, then $x=e^{i\theta}$ for some $\theta\in\mathbb{R}$.
Define
$$
\alpha_n\equiv\alpha_n(\theta):=x^n=e^{i n\theta}\;, \;\;
\beta_n:=\frac{(a_1)_n\cdots(a_{q+1})_n}{(b_1)_n\cdots(b_q)_n\,n!}\quad (n\in\mathbb{N}_0)\,.
$$
Thus, we may write ${}_{q+1}F_{q}(x)=\sum_{n=0}^\infty f_n(x)$, where $f_n(x):=\alpha_n\beta_n$.
Define also
$$
\gamma:=\frac{\Gamma(b_1)\cdots\Gamma(b_q)}{\Gamma(a_1)\cdots\Gamma(a_{q+1})}\;,\quad
\epsilon:=1-\Re\big(\mbox{$\sum a_j-\sum b_j$}\big)\;.
$$
Notice that, taking into account (\ref{HF6}), we have
\begin{equation}\label{HF7}
\big| f_n(x)\big|=|\beta_n|
%\left|\frac{(a_1)_n(a_2)_n\cdots(a_{q+1})_n}{(b_1)_n(b_2)_n\cdots(b_q)_n\, n!}\right|
\sim\frac{|\gamma|}{n^{\epsilon}}\quad\mbox{\rm as}\quad n\to+\infty\;.
\end{equation}
The three cases (i), (ii), and (iii) in the statement of the theorem correspond, respectively, to
$\epsilon>1$, $0<\epsilon\leq1$, and $\epsilon\leq0$.
So we will analyze the convergence of the series ${}_{q+1}F_{q}(x)$
considering separately these three cases.\footnote{\,Note that (\ref{HF7}) gives us
$\,\left|\frac{f_{n+1}(x)}{f_n(x)}\right|
\sim\left(\frac{n}{n+1}\right)^{\epsilon}\longrightarrow1\,$
as $n\to+\infty$, and so one sees that the ratio test is inconclusive.}
\smallskip

If $\epsilon\leq0$, then, by (\ref{HF7}), if $\lim_{n\to\infty} f_n(x)$ exists, it cannot be zero,
% $|f_n(x)|\to+\infty$ as $n\to+\infty$,
hence the series $\sum_{n=0}^\infty f_n(x)\equiv{}_{q+1}F_{q}(x)$ is divergent. This proves (iii).
\smallskip

%If $\epsilon>1$, then $\delta:=(\epsilon-1)/2>0$, and so, by (\ref{HF7}),
%$$
%\quad\frac{\;\;\big|f_n(x)\big|\;\;}{\displaystyle\frac{1}{\;\;n^{1+\delta}\;\;}\;\;}
%=\frac{\;\;\big|f_n(x)\big|\;\;}{\displaystyle\frac{|\gamma|}{\;\;n^{\epsilon}\;\;}\;\;}\,\cdot
%\frac{|\gamma|}{n^{\delta}}
%\xlongrightarrow[n\rightarrow+\infty]{}1\times0=0\,.
%$$
%Consequently, since the series $\sum_{n=0}^\infty \frac{1}{n^{1+\delta}}$ is convergent, then so is
%$\sum_{n=0}^\infty |f_n(x)|$, hence
%$\sum_{n=0}^\infty f_n(x)\equiv{}_{q+1}F_{q}(x)$ is absolutely convergent.
%This proves (i).
%\smallskip

If $\epsilon>1$, the series $\sum_{n=1}^\infty \frac{1}{n^\epsilon}$ is convergent,
and then, by (\ref{HF7}), so is $\sum_{n=0}^\infty |f_n(x)|$, hence
$\sum_{n=0}^\infty f_n(x)\equiv{}_{q+1}F_{q}(x)$ is absolutely convergent.
This proves (i).
\smallskip

At last, suppose that $0<\epsilon\leq1$ and $x\neq1$. Since $x=e^{i\theta}$, we may take $0<\theta<2\pi$.
Thus, setting  $S_n\equiv S_n(\theta):=\sum_{k=0}^{n-1}\alpha_k\beta_k$,
we may ensure the convergence of the series ${}_{q+1}F_q(x):=\lim_{n\to\infty}S_n$
provided that we are able to show that this last limit exists.
Indeed, by the summation by parts formula,
$$
S_n=\sum_{k=0}^{n-1}\alpha_k\beta_k=A_n \beta_{n-1}-\sum_{k=1}^{n-1}A_k(\beta_k-\beta_{k-1})\;,\quad
A_n:=\sum_{j=0}^{n-1}\alpha_j\;.
$$
The sequence $\{A_n\}_{n\geq0}$ is bounded, since for each $n\geq1$,
$$|A_n|=\left|\sum_{k=0}^{n-1}e^{ik\theta}\right|
=\left|\frac{1-e^{in\theta}}{1-e^{i\theta}}\right|
=\left|\frac{\sin(n\theta/2)}{\sin(\theta/2)}\right|\leq\frac{1}{\sin(\theta/2)}\;.$$
Moreover, taking into account (\ref{HF7}), the sequence $\{\beta_n\}_{n\geq0}$ converges to zero. %$$|\beta_n|\sim\frac{|\gamma|}{n^\epsilon}\to0\quad\mbox{\rm as}\quad n\to\infty\;.$$
It follows that $A_n\beta_{n-1}\to0$ as $n\to\infty$.
Thus, to conclude that $\lim_{n\to\infty}S_n$ exists, we need to show that
$\sum_{n=1}^{+\infty}A_{n+1}(\beta_{n+1}-\beta_{n})$ is a convergent series.
Indeed, we have
\begin{equation}\label{Hypbeta}%$$
\beta_{n+1}-\beta_{n}
=\beta_{n+1}\left(1-\frac{(n+b_1)\cdots(n+b_{q})(n+1)}{(n+a_1)\cdots(n+a_q)(n+a_{q+1})}\right)
=\beta_{n+1}\,\frac{c\, n^q+\pi_{q-1}(n)}{n^{q+1}+\pi_q(n)}\;,
%\quad c:=\sum_{j=1}^{q+1} a_j-\sum_{j=1}^{q} b_j-1\;,
\end{equation}%$$
where $c:=\sum_{j=1}^{q+1} a_j-\sum_{j=1}^{q} b_j-1$, $\pi_{q-1}\in\mathcal{P}_{q-1}$,
and $\pi_{q}\in\mathcal{P}_{q}$.
Notice that $c\neq0$, since $|c|\geq|\Re c|=\epsilon>0$. Therefore, we deduce
%Indeed, we have
%$$
%\beta_{n+1}-\beta_{n}
%=\beta_n\left(\frac{(n+a_1)\cdots(n+a_q)(n+a_{q+1})}{(n+b_1)\cdots(n+b_{q})(n+1)}-1\right)
%=\beta_n\,\frac{c\, n^q+\pi_{q-1}(n)}{n^{q+1}+\pi_q(n)}\;,
%%\quad c:=\sum_{j=1}^{q+1} a_j-\sum_{j=1}^{q} b_j-1\;,
%$$
%where $c:=\sum_{j=1}^{q+1} a_j-\sum_{j=1}^{q} b_j-1$, $\pi_{q-1}\in\mathcal{P}_{q-1}$, and $\pi_{q}\in\mathcal{P}_{q}$.
%Notice that $c\neq0$, since $|c|\geq|\Re c|=\delta>0$. Therefore, we deduce
%Indeed, since
%$$
%\beta_{n+1}-\beta_{n}=\beta_n\,
%\frac{c\, n^q+\pi_{q-1}(n)}{n^{q+1}+\pi_q(n)}\;,\quad
%c:=\sum_{j=1}^{q+1} a_j-\sum_{j=1}^{q} b_j-1\;,
%$$
%where $\pi_{q-1}\in\mathcal{P}_{q-1}$ and $\pi_{q}\in\mathcal{P}_{q}$, we deduce
%$$
%\beta_{k+1}-\beta_{k}=\beta_k\,
%\frac{c\, k^q+\mbox{\rm lower degree terms}}{k^{q+1}+\mbox{\rm lower degree terms}}\;,\quad
%c:=\sum_{j=1}^{q+1} a_j-\sum_{j=1}^{q} b_j-1
%$$,
$$
|A_{n+1}(\beta_{n+1}-\beta_{n})|\leq\frac{|\beta_{n+1}|}{\sin(\theta/2)}\,
\Big|\frac{c\, n^q+\pi_{q-1}(n)}{n^{q+1}+\pi_{q}(n)}\Big|\sim\frac{M}{n^{1+\epsilon}}\;,\quad
M:=\frac{|\gamma\,c|}{\sin(\theta/2)}>0\;,
$$
%(note that $c\neq0$, since $|c|\geq|\Re c|=\delta>0$),
hence the (absolute) convergence of the series $\sum_{n=1}^{+\infty}A_{n+1}(\beta_{n+1}-\beta_{n})$
follows from the convergence of the series $\sum_{n=1}^{+\infty}\frac{1}{n^{1+\epsilon}}$.

To prove that the convergence of the series ${}_{q+1}F_q(x)$ is not absolute,
we need to show that the series $\sum_{n=0}^\infty|\beta_n|$ is divergent.
This can be done using Gauss's test (Lemma \ref{Gauss-test}), according to which (since $\epsilon\leq1$)
we may conclude that this series diverges
if we can show that there exist $N\in\mathbb{N}$
and a bounded sequence $\{C_n\}_{n\geq 1}$, such that
$|\beta_n/\beta_{n+1}|=1+\epsilon/n+C_n/n^2$ for each $n\geq N$.
Indeed, taking into account (\ref{Hypbeta}) and the equality $|1-z|^2=1-2\Re z+|z|^2$,
valid for any complex number $z$, we deduce
$$
\left|\frac{\beta_n}{\beta_{n+1}}\right|^2
=\left|1-\frac{c\, n^q+\pi_{q-1}(n)}{n^{q+1}+\pi_q(n)}\right|^2
=1+\frac{2\epsilon}{n}+\frac{B_n}{n^2}\;,
$$
%$$
%\left|\frac{\beta_n}{\beta_{n+1}}\right|
%=\left|1-\frac{c\, n^q+\pi_{q-1}(n)}{n^{q+1}+\pi_q(n)}\right|
%=\sqrt{1+\frac{2\delta}{n}+\frac{B_n}{n^2}\;}
%%=\sqrt{1+\frac{2\delta}{n}+O\left(\frac{1}{n^2}\right)\;}
%=1+\frac{\delta}{n}+\frac{C_n}{n^2}\;,
%$$
%where $B_n:=|c|^2|1+F_n/n|^2-2\Re(cF_n)$, being
%$F_n:=[n\pi_{q-1}(n)/c-\pi_q(n)]/[n^q+\pi_q(n)/n]$,
%where $B_n:=|c-c_n/n|^2+2\Re c_n$, being
%$c_n:=n[c\pi_q(n)-n\pi_{q-1}(n)]/[n^{q+1}+\pi_q(n)]$.
where
$$
B_n:=\Big|c-\frac{c_n}{n}\Big|^2+2\Re c_n\,,\quad
c_n:=\frac{n[c\,\pi_q(n)-n\pi_{q-1}(n)]}{n^{q+1}+\pi_q(n)}\;.
$$
Note that $\{B_n\}_{n\geq1}$ is a bounded sequence (in fact, it is convergent),
so there exists $B>0$ such that $|B_n|\leq B$ for each $n\in\mathbb{N}$.
Finally, using the binomial theorem,
%\,\footnote{\,The binomial theorem states that
%$(1+z)^\alpha=\sum_{k=0}^\infty{{\alpha}\choose{k}}z^k$ for $\alpha,z\in\mathbb{C}$ and $|z|<1$.},
we obtain
$$
\left|\frac{\beta_n}{\beta_{n+1}}\right|
=\sqrt{1+\frac{2\epsilon}{n}+\frac{B_n}{n^2}\;}
=1+\frac{\epsilon}{n}+\frac{C_n}{n^2}\quad\mbox{\rm if}\quad n\geq N\;,
$$
where $N$ is an integer number choosen large enough such that
$|2\epsilon/n+B_n/n^2|\leq1/2$ for each $n\geq N$, and %$C_n:=B_n/2+(2\epsilon+B_n/n)^2\sum_{k=2}^\infty{{1/2}\choose{k}}(2\epsilon/n+B_n/n^2)^{k-2}$.
$$
C_n:=\frac{B_n}{2}+\left(2\epsilon+\frac{B_n}{n}\right)^2\,\sum_{k=2}^\infty{{1/2}\choose{k}}
\left(\frac{2\epsilon}{n}+\frac{B_n}{n^2}\right)^{k-2}\,,\quad n\geq N\,.
$$
(It doesn't matter how to define $C_1,\ldots,C_{N-1}$.)
Clearly, $\{C_n\}_{n\geq 1}$ is bounded, since %for each $n\geq N$,
%$$
%|C_n|\leq \frac{B}{2}+(2\epsilon+B)^2\,\sum_{k=2}^\infty{{1/2}\choose{k}}L^{k-2}
%\leq\frac{B}{2}+\left(\frac{2\epsilon+B}{L}\right)^2\,\sqrt{1+L\,}
%\quad\mbox{\rm if}\quad n\geq N\,.
%$$
$$
|C_n|
%\leq \frac{B}{2}+\left(2+\frac{B}{N}\right)^2\,
%\sum_{k=2}^\infty\Big|{{1/2}\choose{k}}\Big|\left(\frac{2}{N}+\frac{B}{N^2}\right)^{k-2}
\leq\frac{B}{2}+4(2+B)^2\sum_{k=0}^\infty\frac{\big(\frac12\big)_k}{k!}\,\left(\frac12\right)^k
=\frac{B}{2}+4\sqrt{2}\,(2+B)^2\quad\mbox{\rm if}\quad n\geq N\,.
$$
Thus the proof of (ii) is complete.
\qed
\smallskip

%If $0\leq\Re\big(\sum a_j-\sum b_j\big)<1$ and $x\neq1$,
%we leave the proof of statement (ii) to the reader.
%A sketch of the proof is given in Exercise \ref{Ex-cp16-3}
%\qed
%\smallskip

\begin{remark}\label{pFq-elem1}\em
Many elementary functions have representations as hypergeometric series.
We present some simple examples (Exercise \ref{Ex-cp16-3a}):
\begin{snugshade}
$$\begin{array}{ll}
\mbox{\rm (i)}\quad e^x=\pFq{0}{0}{\mbox{---}}{\mbox{---}}{-x} &
\mbox{\rm (ii)}\quad\log(1-x)=-x\,\pFq{2}{1}{1,1}{2}{x} \\ [0.5em]
\mbox{\rm (iii)}\quad\sin x=x\,\pFq{0}{1}{\mbox{---}}{3/2}{-x^2/4} &
\mbox{\rm (iv)}\quad\cos x=\pFq{0}{1}{\mbox{---}}{1/2}{-x^2/4} \\ [0.75em]
\mbox{\rm (v)}\quad\arcsin x=x\,\pFq{2}{1}{1/2,1/2}{3/2}{x^2} &
\mbox{\rm (vi)}\quad\arctan x=x\,\pFq{2}{1}{1/2,1}{3/2}{-x^2}\,.
\end{array}$$
\end{snugshade}
\end{remark}

Finally, we note that the {\it binomial theorem} can be written in hypergeometric form:
\begin{snugshade}\vspace*{-0.5em}
\begin{equation}\label{binThm}
(1-x)^{-a}=\pFq{1}{0}{a}{\mbox{---}}{x}=\sum_{n=0}^\infty \frac{(a)_n}{n!}\,x^n
\;,\quad |x|<1\;,\quad a\in\mathbb{C}\,.
\end{equation}
\end{snugshade}

\section{The hypergeometric function ${}_{2}F_{1}$} %${}_{2}F_{1}(a,b;\,c;\,x)$}

\subsection{Definition}
The preceding example (ii) involving $\log(1-x)$
shows that although the series converges for $|x|<1$, it has an analytic continuation
as a single-valued function in the complex plane from which a line joining $1$ to $\infty$ is deleted.
We will see that this behavior describes the general situation, i.e.,
a ${}_{2}F_{1}$ series has a continuation to the complex plane with branch points at $1$ and $\infty$.

\begin{snugshade}
\begin{definition}\label{2F1-def}
The {\sl hypergeometric function} is defined by the series
\begin{equation}\label{2F1-1}
\pFq{2}{1}{a,b}{c}{x}:=
\sum_{n=0}^\infty \frac{(a)_n(b)_n}{(c)_n}\frac{x^n}{n!}
\end{equation}
for $|x|<1$, and by analytic continuation elsewhere.
\end{definition}
\end{snugshade}

\begin{remark}\em
Of course, in definition (\ref{2F1-1}) it is implicitly assumed that
%the parameters $a$, $b$, and $c$ fulfill
$a,b\in\mathbb{C}$ and $c\in\mathbb{C}\setminus\{0,-1,-2,\cdots\}$.
\end{remark}

\begin{remark} \em
Usually we reserve the use of the words {\it hypergeometric function}
for ${}_{2}F_{1}$, and {\it hypergeometric series} will be the series
${}_{p}F_{q}$ defined by (\ref{HF3}) --- which includes ${}_{2}F_{1}$,
but will not necessarily mean just ${}_{2}F_{1}$.
\end{remark}
%\smallskip

Notice that Theorems \ref{pFq-conv} and \ref{pFq-conv-thm2}
applied to the specific ${}_{2}F_{1}$ series yield:

\begin{snugshade}
\begin{theorem}\label{2F1-conv-thm2}
Consider the hypergeometric series $(\ref{2F1-1})$,
and suppose that it is a nonterminating series. Then:
\begin{enumerate}
\item[{\rm 1.}]
If $|x|<1$, then the series is absolutely convergent.
\item[{\rm 2.}]
If $|x|>1$, then the series is divergent.
\item[{\rm 3.}]
If $|x|=1$, then the following holds:
\begin{enumerate}
\item[{\rm (i)}] if $\Re(a+b-c)<0$, then the series converges absolutely;
\item[{\rm (ii)}] if $0\leq\Re(a+b-c)<1$ and $x\neq1$, then the series converges conditionally;
\item[{\rm (iii)}] if $\Re(a+b-c)\geq1$, then the series diverges.
\end{enumerate}
\end{enumerate}
\end{theorem}
\end{snugshade}

Consider, for instance, $a=b=c=1$. Then, since $(1)_n=n!$, one has
$$
\pFq{2}{1}{1,1}{1}{x}=\sum_{n=0}^\infty x^n=\frac{1}{1-x}\quad\mbox{\rm if}\;\; |x|<1\;.
$$
In this case, the series is convergent if $|x|<1$, and it is divergent otherwise.
Clearly the function $1/(x-1)$
provides the analytic continuation to $\mathbb{C}\setminus\{1\}$,
and thus the hypergeometric function $\pFq{2}{1}{1,1}{1}{x}$ becomes
defined for each $x\in\mathbb{C}\setminus\{1\}$.
\smallskip

Next we state some important results concerning the hypergeometric function ${}_{2}F_{1}$,
including Euler's integral representation and Gauss theorem, as well as two other results
involving terminating series ${}_{2}F_{1}(x)$ and ${}_{3}F_{2}(x)$ at the point $x=1$, namely the
Chu-Vandermonde and Pfaff-Saalsch\"utz identities.

%We begin by stating
\subsection{Euler's integral representation}
Euler's integral representation may be viewed as the analytic
continuation of (\ref{2F1-1}), %${}_{2}F_{1}(a,b;c;x)$,
provided that the condition $\Re c>\Re b>0$ is satisfied.
This condition involves only the parameters $b$ and $c$,
and not the parameter $a$, which is involved in the function $(1-xt)^{-a}$
that appears in the integrand of the integral representation --- see (\ref{2F1-2}) in bellow.
Regarded as a function of the complex variable $x$ (and being $0\leq t\leq1$, fixed),
this function is in general multivalued
(it is single-valued if $a$ is an integer number --- see (\ref{1-xta}) in bellow).
Taking its principal value, we obtain a single-valued
function which is analytic in the $x-$plane cut along the real axis from $1$ to $\infty$,
i.e., it is an analytic function of the variable $x$ in $\mathbb{C}\setminus[1,+\infty)$.
To see why this holds, we recall that, %$\alpha$ being a fixed complex number,
if $\alpha$ is a (fixed) complex number,
the function defined for $z\in\mathbb{C}\setminus\{0\}$ by
$z^\alpha:=e^{\alpha\,\log_r z}$,
%($\alpha$ being a fixed complex number),
where $\log_r z:=\ln|z|+i\arg_r z$, $\arg_r z\in]r,r+2\pi]$, and $r\in\mathbb{R}$
(fixing the branch of the logarithm),
is an analytic function on $\mathbb{C}\setminus\ell_r\,$, where $\ell_r\,$ is the ray
$\ell_r:=\{\rho e^{ir}\,|\,\rho\geq0\}$.
Using this fact one sees that for its principal value (which is obtained for $r=-\pi$),
the function $(1-xt)^{-a}$, for fixed $t\in[0,1]$, is analytic outside the range of
values $x\in\mathbb{C}$ such that the condition $1-xt\in\ell_{-\pi}:=(-\infty,0]$ holds.
This condition is impossible if $t=0$, hence $(1-xt)^{-a}$
is analytic in $\mathbb{C}$ if $t=0$. If $t\in(0,1]$,
then $1-xt\in(-\infty,0]$ if and only if $x\geq1/t$,
and so $(1-xt)^{-a}$ is analytic in $\mathbb{C}\setminus[1/t,+\infty)$.
%and if this condition holds for all $t\in(0,1]$ then $x\geq1$.
%then $1-xt\in(-\infty,0]$ if and only if
%$1-xt=-\rho$ and $\rho\geq0$, i.e., $x=(1+\rho)/t$ and $\rho\geq0$,
%and this condition holds for all $t\in(0,1]$ provided $x\geq1$.
%Thus, $(1-xt)^{-a}$ with its principal value (being $0\leq t\leq1$) defines a single-valued
%function analytic in the $x-$plane cut along the real axis from $1$ to $+\infty$.
%If $t=0$, obviously $(1-xt)^{-a}\equiv1$ is analytic.
%Notice also that the choice of the branch implies $(1-xt)^{-a}\to1$ as $t\to0^+$.
The choice of the (principal) branch implies
%for each $x\in\mathbb{C}\setminus[1,+\infty)$ and $t\in[0,1]$,
the following explicit expression of $(1-xt)^{-a}$
as single-valued function (of the variable $x$):\footnote{\,Notice also that $1-xt=1$ if $t=0$ and $1-xt\neq0$ for each
$x\in\mathbb{C}\setminus[1,+\infty)$ if $0< t\leq1$, hence $1-xt\neq0$ for each
$x\in\mathbb{C}\setminus[1,+\infty)$ if $t\in[0,1]$, and so $(1-xt)^{-a}$ is well defined
for every $x\in\mathbb{C}\setminus[1,+\infty)$ and $t\in[0,1]$, whatever the choice of $a\in\mathbb{C}$.}
\vspace*{-0.5em}
%$$
%(1-xt)^{-a}:=|1-xt|^{-\Re a} e^{\Im a\cdot\arg(1-xt)-i\big(\Re a\cdot\arg(1-xt)+\Im a\cdot\ln|1-xt|\big)}\;,
%\;\arg(1-xt)\in[-\pi,\pi[\,.
%$$
\begin{snugshade}\vspace*{-0.5em}
\begin{equation}\label{1-xta}
\begin{array}{c}
(1-xt)^{-a}:=|1-xt|^{-\Re a} e^{\Im a\cdot\arg(1-xt)-i\big(\Re a\cdot\arg(1-xt)+\Im a\cdot\ln|1-xt|\big)}\;, \\ [0.5em]
x\in\mathbb{C}\setminus[1,+\infty)\;,\quad t\in[0,1]\;,\quad\arg(1-xt)\in]-\pi,\pi]\,.
\end{array}
\end{equation}
\end{snugshade}
\noindent
(Note that from this we obtain $(1-xt)^{-a}\equiv1$ if $t=0$ and $(1-xt)^{-a}\to1$ as $t\to0^+$.)
In conclusion: %for every $t\in[0,1]$,
$(1-xt)^{-a}$ with its principal value defines a single-valued
function analytic in the $x-$plane cut along the real axis from $1$ to $+\infty$,
whatever the value of $t\in[0,1]$.

%\begin{Snugshade}[236,236,0]
%\begin{theorem}[Euler's integral representation]\label{2F1-Euler-int-rep}
%If $\Re c>\Re b>0$, then
%\begin{equation}\label{2F1-2}
%\pFq{2}{1}{a,b}{c}{x}=\frac{\Gamma(c)}{\Gamma(b)\Gamma(c-b)}
%\int_0^1 t^{b-1}(1-t)^{c-b-1}(1-xt)^{-a}\,{\rm d}t\;,
%\end{equation}
%in the $x-$plane cut along the real axis from $1$ to $\infty$.
%Here it should be understood that $\arg t=\arg(1-t)=0$ and $(1-xt)^{-a}$ as its principal value.
%\end{theorem}
%\end{Snugshade}

\begin{snugshade}
\begin{theorem}[Euler's integral representation]\label{2F1-Euler-int-rep}
If $\Re c>\Re b>0$, then
\begin{equation}\label{2F1-2}
\pFq{2}{1}{a,b}{c}{x}=\frac{\Gamma(c)}{\Gamma(b)\Gamma(c-b)}
\int_0^1 t^{b-1}(1-t)^{c-b-1}(1-xt)^{-a}\,{\rm d}t\;,
\end{equation}
in the $x-$plane cut along the real axis from $1$ to $+\infty$.
Here it should be understood that $\arg t=\arg(1-t)=0$ and $(1-xt)^{-a}$ as its principal value.
\end{theorem}
\end{snugshade}

{\it Proof.}
Fix $x\in\mathbb{C}$ such that $|x|<1$.
According with the binomial theorem (\ref{binThm}), %we may write
\begin{equation}\label{binEulerA}
t^{b-1}(1-t)^{c-b-1}(1-xt)^{-a}%=\sum_{n=0}^\infty \frac{(a)_n}{n!}\,x^nt^{n+b-1}(1-t)^{c-b-1}
=\sum_{n=0}^\infty f_n(t)\;,
\end{equation}
where $f_n\equiv f_n(\cdot;x):(0,1)\to\mathbb{C}$ (regarded as a function of $t$) is defined by
%(it doesn't matter how we define $f_n$ at $t=0,1$)
$$
f_n(t):=\frac{(a)_n\,x^n}{n!}\,t^{n+b-1}(1-t)^{c-b-1}\;. %,\quad 0<t<1\;.
$$
Notice that $f_n\in L^1(0,1)$. Indeed, for each $t\in(0,1)$, we may write
$$
\big|f_n(t)\big|=\frac{|(a)_n|\,|x|^n}{n!}\,|t^n|\,t^{\Re b-1}(1-t)^{\Re(c-b)-1}
\leq\frac{(|a|)_n\,|x|^n}{n!}\,t^{\Re b-1}(1-t)^{\Re(c-b)-1}\;,
$$
the last inequality being justified by the obvious inequality $|(\alpha)_n|\leq(|\alpha|)_n$,
which holds for all $\alpha\in\mathbb{C}$ and $n\in\mathbb{N}_0$.
Since, by assumption, $\Re b>0$ and $\Re(c-b)>0$, then the
function $t\in(0,1)\mapsto t^{\Re b-1}(1-t)^{\Re(c-b)-1}$ is in $L^1(0,1)$.
To see why this holds, notice simply that the integral of such a function is
the beta integral (cf. Definition \ref{BetaF-def})
$$
\int_0^1\,t^{\Re b-1}(1-t)^{\Re(c-b)-1}\,{\rm d}t=B\big(\Re b,\Re(c-b)\big)\;.
$$
Moreover, for each $n\in\mathbb{N}_0$, we may write
$$
\int_0^1|f_n(t)|\,{\rm d}t\leq\frac{(|a|)_n\,|x|^n}{n!}\,B\big(\Re b,\Re(c-b)\big)\;,
$$
and so, summing up for $n=0,1,2,\cdots$, and noticing that,
taking into account (ii) in Theorem \ref{pFq-conv}, the series
$\sum_{n=0}^\infty\frac{(|a|)_n}{n!}\,|x|^n=\pFq{1}{0}{|a|}{\mbox{---}}{|x|}$
is convergent, we obtain
\begin{equation}\label{2F1-exch-int}
\sum_{n=0}^\infty\int_0^1|f_n(t)|\,{\rm d}t\leq
%B\big(\Re b,\Re(c-b)\big)\,\sum_{n=0}^\infty\frac{(|a|)_n}{n!}\,|x|^n=
B\big(\Re b,\Re(c-b)\big)\,\pFq{1}{0}{|a|}{\mbox{---}}{|x|}<\infty\;.
\end{equation}
Now, integrating both sides of (\ref{binEulerA}) with respect to the variable $t$,
(\ref{2F1-exch-int}) allow us to perform the change in the order of
integration and summation.\footnote{\,Recall the following well known (consequence of the Monotone Convergence Theorem)

{\sc Theorem:} {\it Let $(\Omega,\mathcal{A},\mu)$ be a measure space
and $\{f_n\}_{n\geq1}$ a sequence of (complex) functions in $L^1(\Omega,\mu)$
%such that $f_n\in L^1(\Omega,\mu)$ for each $n\in\mathbb{N}$ and
fulfilling $\sum_{n=1}^\infty\int_\Omega |f_k|\,{\rm d}\mu<\infty$.
Then $\sum_{n=1}^\infty |f_n|<\infty$ $\mu$-a.e. in $\Omega$, $\sum_{n=1}^\infty f_n\in L^1(\Omega,\mu)$, and
$$\int_\Omega \Big(\sum_{n=1}^\infty f_n\Big)\,{\rm d}\mu=
\sum_{n=1}^\infty \int_\Omega f_n\,{\rm d}\mu\;.$$}}
This yields %, and thus we obtain
%$$
%\begin{array}{rl}
%\displaystyle\int_0^1 t^{b-1}(1-t)^{c-b-1}(1-xt)^{-a}\,{\rm d}t &
%=\displaystyle\int_0^1 \sum_{n=0}^\infty f_n(t;x)\,{\rm d}t=\sum_{n=0}^\infty \int_0^1 f_n(t;x)\,{\rm d}t \\ [1em]
%&=\displaystyle \sum_{n=0}^\infty \frac{(a)_n}{n!}x^n\,\int_0^1 t^{n+b-1}(1-t)^{c-b-1}\,{\rm d}t \\ [0.5em]
%&=\displaystyle\sum_{n=0}^\infty \frac{(a)_n}{n!}x^n\,B(n+b,c-b)\,,
%\end{array}
%$$
\begin{equation}\label{2F1-intB}
\int_0^1 t^{b-1}(1-t)^{c-b-1}(1-xt)^{-a}\,{\rm d}t
=\sum_{n=0}^\infty \frac{(a)_n}{n!}x^n\,\int_0^1 t^{n+b-1}(1-t)^{c-b-1}\,{\rm d}t\,.
\end{equation}
By (\ref{BetaF1}) and (\ref{BetaF3}), and taking into account (\ref{Poch-Gamma}), we may write
$$
\int_0^1 t^{n+b-1}(1-t)^{c-b-1}\,{\rm d}t=B(n+b,c-b)=\frac{\Gamma(n+b)\Gamma(c-b)}{\Gamma(n+c)}
=\frac{(b)_n\Gamma(b)\Gamma(c-b)}{(c)_n\Gamma(c)}\;.
$$
Inserting this into the right-hand side of (\ref{2F1-intB}) yields (\ref{2F1-2}) for $|x|<1$.
%Finally, since the integral is analytic in the cut plane, (\ref{2F1-2}) holds in this region as well.
%Finally, since for each fixed $t\in[0,1]$ the function $x\mapsto(1-xt)^{-a}$ is analytic
%%and bounded in each compact set of
%in the cut plane $\mathbb{C}\setminus[1,+\infty)$,
%then one easily sees that the integral on the right-hand side of (\ref{2F1-2}) is also an analytic function
%of the variable $x$ in the cut plane (Exercise \ref{Ex-cp16-2a}), hence (\ref{2F1-2}) holds in this region as well.
%\qed
%\smallskip
%
%%%%%%%%%%%%% Proof that the Euler's integral is analytic %%%%%%%%%%%%%%%%%%%%%
To prove that (\ref{2F1-2}) holds in the cut plane $\mathbb{C}\setminus[1,+\infty)$,
we will show that the integral on the right-hand side of (\ref{2F1-2}) is
an analytic function of $x$ in the cut plane.
Indeed, set\footnote{\,We will apply the following general theorem, taking therein
$\Omega=[0,1]$, $G=\mathbb{C}\setminus[1,+\infty)$, and $\mu$ the Lebesgue measure in $\mathbb{R}$
restricted to the interval $[0,1]$.

{\sc Theorem}\,\ref{LutzMattner2001-C16}: {\it Let $(\Omega,\mathcal{A},\mu)$ be a measure space,
let $G\subseteq\mathbb{C}$ be open, and let $f:G\times \Omega\to\mathbb{C}$ be
a function fulfilling the following three properties:
{\rm (i)} $f(z,\cdot)$ is $\mathcal{A}-$measurable for each $z\in G$;
{\rm (ii)} $f(\cdot,t)$ is analytic in $G$ for each $t\in\Omega$; and
{\rm (iii)} $\int_\Omega|f(\cdot,t)|\,{\rm d}\mu(t)$ is locally bounded, that is,
$$
\forall z_0\in G\;,\quad\exists\delta\equiv\delta(z_0)>0\quad:\quad
\sup_{\substack{|z-z_0|\leq\delta \\ (z\in G)}}\,\int_\Omega|f(z,t)|\,{\rm d}\mu(t)<\infty\;.
$$
Then $\int_\Omega f(\cdot,t){\rm d}\mu(t)$ is analytic in $G$ and it may be differentiated under the integral}.}
$$
f(x,t):=t^{b-1}(1-t)^{c-b-1}(1-xt)^{-a}\;,\quad x\in\mathbb{C}\setminus[1,+\infty)\;,\quad t\in(0,1)\;.
$$
(It doesn't matter how we define $f$ for $t=0$ or $t=1$,
provided it remains analytic in the variable $x$.)
We have already seen that for each fixed $t\in[0,1]$ the function $x\mapsto(1-xt)^{-a}$ is analytic
in the cut plane $\mathbb{C}\setminus[1,+\infty)$, and so the same holds for $f$,
regarded as a function of the variable $x$.
On the other hand, it is not difficult to see that
for each fixed $x\in\mathbb{C}\setminus[1,+\infty)$, the function $t\mapsto(1-xt)^{-a}$
is continuous on $[0,1]$, hence it is measurable there, and then so is $f$,
regarded as a function of the variable $t$.
Moreover, from (\ref{1-xta}), it is straightforward to show (Exercise \ref{Ex-cp16-2a})
that for each $z_0\in\mathbb{C}\setminus[1,+\infty)$,
there exists $\delta\equiv\delta(z_0)>0$ such that
\begin{equation}\label{1-xtaCdelta}
\big|(1-xt)^{-a}\big|\leq C(a,z_0,\delta)e^{\pi|\Im a|}\;,\quad\forall x\in \overline{B}(z_0,\delta)\;,\;
\forall t\in[0,1]\;,
\end{equation}
being $\overline{B}(z_0,\delta):=\{x\in\mathbb{C}:|x-z_0|\leq\delta\}\subset\mathbb{C}\setminus[1,+\infty)$ and
$C(a,z_0,\delta)$ is a constant that depends only of $a$, $z_0$, and $\delta$,
%$$
%\begin{array}{l}
%\delta:=\frac12\,|\Im z_0|\chi_{\mathbb{R}\setminus\{0\}}(\Im z_0) \\ [0.5em]
%\qquad\quad +\big[\frac12\,(1-z_0)\chi_{(-\infty,0]}(\Re a)+
%\frac13\,\min\{1,1-z_0\}\chi_{(0,\infty)}(\Re a)\big]\chi_{\{0\}}(\Im z_0)\,,
%\end{array}
%$$
%$$
%\begin{array}{l}
%C(a,z_0,\delta):=
%(1+|z_0|+\delta)^{-\Re a}\,\chi_{(-\infty,0]}(\Re a) \\ [0.2em]
%\qquad\qquad\qquad +\big[1/\delta^{\Re a}\chi_{\{0\}}(\Im z_0)+
%(1+|z_0|/\delta)^{\Re a}\chi_{\mathbb{R}\setminus\{0\}}(\Im z_0)\big]\,\chi_{(0,\infty)}(\Re a)\,,
%\end{array}
%$$
and so, we obtain
$$
\sup_{|x-z_0|\leq\delta}\int_0^1 \big|f(x,t)\big|\,{\rm d}t\leq
C(a,z_0,\delta)\,e^{\pi|\Im a|}\,B\big(\Re b,\Re(c-b)\big)<\infty\;.
$$
Thus, we conclude that the right-hand side of (\ref{2F1-2}) is an analytic function
of the variable $x$ in the cut plane,
%\footnote{\,Recall the following general theorem, where we have to take
%%by applying the following general theorem, taking therein
%$\Omega=[0,1]$, $G=\mathbb{C}\setminus[1,+\infty)$, and $\mu$ the Lebesgue measure in $\mathbb{R}$
%restricted to the interval $[0,1]$.
%
%{\sc Theorem:} {\it Let $(\Omega,\mathcal{A},\mu)$ be a measure space, $G\subseteq\mathbb{C}$ open,
%and $f:G\times \Omega\to\mathbb{C}$ a function fulfilling the following three properties:
%{\rm (i)} $f(z,\cdot)$ is $\mathcal{A}-$measurable for each $z\in G$;
%{\rm (ii)} $f(\cdot,t)$ is analytic in $G$ for each $t\in\Omega$; and
%{\rm (iii)} $\int_\Omega|f(\cdot,t)|\,{\rm d}\mu(t)$ is locally bounded, that is,
%$$
%\forall z_0\in G\;,\quad\exists\delta\equiv\delta(z_0)>0\quad:\quad
%\sup_{\substack{|z-z_0|\leq\delta \\ (z\in G)}}\int_\Omega|f(z,t)|\,{\rm d}\mu(t)<\infty\;.
%$$
%Then $\int_\Omega f(\cdot,t){\rm d}\mu(t)$ is analytic in $G$ and it may be differentiated under the integral}.} %$\mathbb{C}\setminus[1,+\infty)$,
hence, since we have already proved that (\ref{2F1-2}) holds if $|x|<1$,
then if follows by analytic continuation that it holds in the cut $x-$plane
$\mathbb{C}\setminus[1,+\infty)$ as well.
\qed
\medskip
%%%%%%%%%%%%%%%%%%%%%%%%%%%%%%%%%%%%
%$$
%\delta:=\left\{
%\begin{array}{ccl}
%\big|\Im z_0\big|/2 & \mbox{\rm if} & \Im z_0\neq0\;, \\ [0.25em]
%(1-z_0)/2 & \mbox{\rm if} & \Im z_0=0\,\wedge\,\Re a\leq0\;, \\ [0.25em]
%\min\{1,1-z_0\}/3 & \mbox{\rm if} & \Im z_0=0\,\wedge\,\Re a>0\;,
%\end{array}
%\right.
%$$
%and
%$$
%C(a,z_0,\delta):=\left\{
%\begin{array}{ccl}
%\big(1+|z_0|+\delta\big)^{-\Re a} & \mbox{\rm if} & \Re a\leq0\;, \\ [0.25em]
%\delta^{-\Re a} & \mbox{\rm if} & \Re a>0\,\wedge\,\Im z_0=0\;, \\ [0.25em]
%\big(1+|z_0|/\delta\big)^{\Re a} & \mbox{\rm if} & \Re a>0\,\wedge\,\Im z_0\neq0\;.
%\end{array}
%\right.
%$$
%%%%%%%%%%%%%%%%%%% End %%%%%%%%%%%%%%%%%%%%%%%%%%%%%%%%%%%%%%%%%%%%%%%%%%%%%%%%%%%%%%%%%%%%%

As a first application of Euler's integral representation we derive
two transformation formulas of hypergeometric functions.

%\begin{snugshade}
%\begin{corollary}\label{2F1-Pfaff-Euler}
%If $\Re c>\Re b>0$, then the following transformation formulas hold
%in the $x-$plane cut along the real axis from $1$ to $\infty$:
%\begin{equation}\label{2F1-PE1}
%\pFq{2}{1}{a,b}{c}{x}=(1-x)^{-a}\,\pFq{2}{1}{a,c-b}{c}{\frac{x}{x-1}}\;,\;\;\qquad{\rm (Pfaff, 1797)}
%\end{equation}
%\begin{equation}\label{2F1-PE2}
%\pFq{2}{1}{a,b}{c}{x}=(1-x)^{c-a-b}\,\pFq{2}{1}{c-a,c-b}{c}{x}\;.\qquad{\rm (Euler, 1794)}
%\end{equation}
%Here it should be understood that $(1-x)^{-a}$ and $(1-x)^{c-a-b}$ have their principal values.
%\end{corollary}
%\end{snugshade}

\begin{snugshade}
\begin{corollary}\label{2F1-Pfaff-Euler}
If $|x|<1$ and $|x/(x-1)|<1$, then the following transformation formula holds:
%in the $x-$plane cut along the real axis from $1$ to $\infty$:
\begin{equation}\label{2F1-PE1}
\pFq{2}{1}{a,b}{c}{x}=(1-x)^{-a}\,\pFq{2}{1}{a,c-b}{c}{\frac{x}{x-1}}\;;\;\;\qquad{\rm (Pfaff, 1797)}
\end{equation}
and if $|x|<1$, then
\begin{equation}\label{2F1-PE2}
\pFq{2}{1}{a,b}{c}{x}=(1-x)^{c-a-b}\,\pFq{2}{1}{c-a,c-b}{c}{x}\;.\qquad{\rm (Euler, 1794)}
\end{equation}
Here, it should be understood that $(1-x)^{-a}$ and $(1-x)^{c-a-b}$ have their principal values.
Moreover, these formulas are valid for all complex parameters $a$, $b$, and $c$,
provided that $c$ is not zero neither a negative integer number.
\end{corollary}
\end{snugshade}

{\it Proof.}
Assume first $\Re c>\Re b>0$. To prove Pfaff's transformation,
make the substitution $t=1-s$ in Euler's integral (\ref{2F1-2}). Then
$$
\begin{array}{rcl}
\displaystyle\pFq{2}{1}{a,b}{c}{x}&=&
\displaystyle\frac{\Gamma(c)}{\Gamma(b)\Gamma(c-b)}
\int_0^1 (1-s)^{b-1}s^{c-b-1}(1-x+xs)^{-a}\,{\rm d}s \\ [1em]
&=& \displaystyle\frac{(1-x)^{-a}\Gamma(c)}{\Gamma(b)\Gamma(c-b)}
\int_0^1 (1-s)^{b-1}s^{c-b-1}\left(1-\frac{xs}{x-1}\right)^{-a}\,{\rm d}s \\ [1em]
&=& \displaystyle (1-x)^{-a}\,\pFq{2}{1}{a,c-b}{c}{\frac{x}{x-1}}\;.
\end{array}
$$
To prove Euler's transformation, we consider Pfaff's transformation and note that
the hypergeometric series is symmetric in the parameters appearing in the numerator.
Therefore, we may write
$$
\pFq{2}{1}{a,b}{c}{x}%=(1-x)^{-a}\,\pFq{2}{1}{a,c-b}{c}{\frac{x}{x-1}}
=(1-x)^{-a}\,\pFq{2}{1}{c-b,a}{c}{\frac{x}{x-1}}\;.
$$
Applying again Pfaff's transformation (to the last ${}_{2}F_{1}$), we obtain
$$
\begin{array}{rcl}
\displaystyle\pFq{2}{1}{a,b}{c}{x}&=&
\displaystyle (1-x)^{-a}\left(\mbox{$1-\frac{x}{x-1}$}\right)^{-c+b}\,
\pFq{2}{1}{c-b,c-a}{c}{\frac{\frac{x}{x-1}}{\frac{x}{x-1}-1}} \\ [1.25em]
&=&\displaystyle (1-x)^{-a}(1-x)^{c-b}\,
\pFq{2}{1}{c-a,c-b}{c}{x}\;.
\end{array}
$$
So (\ref{2F1-PE1}) and (\ref{2F1-PE2}) hold under the assumption $\Re c>\Re b>0$.
Analytic continuation in the parameters $b$ and $c$ (Exercise \ref{Ex-cp16-5aa})
gives (\ref{2F1-PE1}) and (\ref{2F1-PE2})
for all complex values of $a$, $b$ and $c$, with $c\neq0,-1,-2,-3,\ldots\,$.
\qed

\begin{remark}\em
The hypergeometric ${}_2F_1$ series defined on the right-hand side of (\ref{2F1-PE1})
converges for $|x/(x-1)|<1$. Thus, since this condition is equivalent to $\Re x<\frac12$,
the right-hand side of Pfaff's transformation gives the analytic continuation
to the region $\Re x<\frac12$ (via Euler's integral representation) of the series defined
by $\pFq{2}{1}{a,b}{c}{x}$.
\end{remark}
\medskip

%A a curious special case of Pfaff's transformation is the well know formula
%$$
%\arctan x = \arcsin\frac{x}{\sqrt{1+x^2\,}} \quad (x\in\mathbb{R})\;,
%$$
%which follows immediately from (\ref{2F1-PE1}) taking into account examples
%(v) and (vi) given in Remark \ref{pFq-elem1}.

\subsection{Gauss's summation formula}
Our next result is a celebrated theorem by Gauss. It is convenient to state firstly the following

\begin{snugshade}
\begin{lemma}\label{2F1-Gauss-lemma}
%If $-c\not\in\mathbb{N}_0$ and $\Re (c-a-b)>0$, then
%If $c\in\mathbb{C}\setminus\{0,-1,-2,\cdots\}\;$ and $\Re (c-a-b)>0$, then
If $\Re (c-a-b)>0$, then
\begin{equation}\label{2F1-3}
\pFq{2}{1}{a,b}{c}{1}=\frac{(c-a)(c-b)}{c(c-a-b)}\,\pFq{2}{1}{a,b}{c+1}{1}\;.
\end{equation}
\end{lemma}
\end{snugshade}

{\it Proof.}
Set
$$
A_n:=\frac{(a)_n(b)_n}{n!(c)_n}\;,\quad
B_n:=\frac{(a)_n(b)_n}{n!(c+1)_n}\;,\quad n\in\mathbb{N}_0\;.
$$
After straightforward computations we deduce
$$
c(c-a-b)A_n=(c-a)(c-b)B_n+cnA_n-c(n+1)A_{n+1}\;,\quad n\in\mathbb{N}_0\;.
$$
Therefore, summing up from $n=0$ to $n=N$, we obtain
\begin{equation}\label{ANBNlimit}
c(c-a-b)\sum_{n=0}^NA_n=(c-a)(c-b)\sum_{n=0}^NB_n-c(N+1)A_{N+1}\;,\quad N\in\mathbb{N}\;.
\end{equation}
Now, as $N\to+\infty$,
$$
\sum_{n=0}^NA_n\rightarrow\pFq{2}{1}{a,b}{c}{1}\;,\quad
\sum_{n=0}^NB_n\to\pFq{2}{1}{a,b}{c+1}{1}\;.
%\;,\quad\mbox{\rm as}\; N\to+\infty\;.
$$
Moreover, by (\ref{HF6}), as $N\to+\infty$,
$$
\frac{(a)_{N+1}(b)_{N+1}}{(c)_{N+1}\, (N+1)!}
\sim\frac{\Gamma(c)}{\Gamma(a)\Gamma(b)}\,(N+1)^{a+b-c-1}\;,
$$
and so, since $\Re(a+b-c)<0$, we obtain
$$
(N+1)A_{N+1}=(N+1)\,\frac{(a)_{N+1}(b)_{N+1}}{(c)_{N+1}\, (N+1)!}
\sim\frac{\Gamma(c)}{\Gamma(a)\Gamma(b)}\,(N+1)^{a+b-c}\to0\;.
$$
Therefore, taking $N\to+\infty$ in (\ref{ANBNlimit}) yields (\ref{2F1-3}).
\qed

\begin{snugshade}
\begin{theorem}[Gauss's summation formula, 1812]\label{2F1-Gauss}
If $\Re(c-a-b)>0$, then
\begin{equation}\label{2F1-4}
\pFq{2}{1}{a,b}{c}{1}=\frac{\Gamma(c)\Gamma(c-a-b)}{\Gamma(c-a)\Gamma(c-b)}\;.
\end{equation}
\end{theorem}
\end{snugshade}

{\it Proof.}
Iterating (\ref{2F1-3}) $n$ times yields
\begin{equation}\label{2F1-Gauss1}
\pFq{2}{1}{a,b}{c}{1}=\frac{(c-a)_n(c-b)_n}{(c)_n(c-a-b)_n}\,\pFq{2}{1}{a,b}{c+n}{1}\;,\quad n\in\mathbb{N}\,.
\end{equation}
By  (\ref{HF6}), as $n\to+\infty$,
$$
\frac{(c-a)_n(c-b)_n}{(c)_n(c-a-b)_n}=\frac{(c-a)_n(c-b)_n(1)_n}{(c)_n(c-a-b)_n\,n!}
\sim\frac{\Gamma(c)\Gamma(c-a-b)}{\Gamma(c-a)\Gamma(c-b)\Gamma(1)}\;,
$$
i.e., recalling that $\Gamma(1)=1$,
$$%\begin{equation}\label{2F1-Gauss2}
\lim_{n\to+\infty}\frac{(c-a)_n(c-b)_n}{(c)_n(c-a-b)_n}=\frac{\Gamma(c)\Gamma(c-a-b)}{\Gamma(c-a)\Gamma(c-b)}\;.
$$%\end{equation}
Therefore, (\ref{2F1-4}) will be proved taking the limit in (\ref{2F1-Gauss1}) as $n\to\infty$,
provided we are able to show that
\begin{equation}\label{2F1-Gauss3}
\lim_{n\to+\infty}\,\pFq{2}{1}{a,b}{c+n}{1}=1\;.
\end{equation}
Let $u_k(a,b,c)$ denote the coefficient of $x^k$ in $\pFq{2}{1}{a,b}{c}{x}$, i.e., write
$$
\pFq{2}{1}{a,b}{c}{x}=\sum_{k=0}^\infty u_k(a,b,c)x^k\;,\quad
u_k(a,b,c):=\frac{(a)_k(b)_k}{(c)_k\,k!}\;. %,\quad k\in\mathbb{N}_0\,.
$$
For each $k\in\mathbb{N}_0$ and $n\in\mathbb{N}$ such that $n>|c|$,
we have $|(a)_k|\leq(|a|)_k$, $|(b)_k|\leq(|b|)_k$, and
$|(c+n)_k|\geq(n-|c|)_k$, hence
$$
\big|u_k(a,b,c+n)\big|=\left|\frac{(a)_k(b)_k}{(c+n)_k\,k!}\right|
\leq\frac{(|a|)_k(|b|)_k}{(n-|c|)_k\,k!}=u_k\big(|a|,|b|,n-|c|\big)\;,
$$
and so we may write
%$$
%\begin{array}{rcl}
%\displaystyle\left|\pFq{2}{1}{a,b}{c+n}{1}-1\right| &=&
%\displaystyle\left|\sum_{k=1}^\infty u_k(a,b,c+n)\right|\leq
%\sum_{k=1}^\infty \big|u_k(a,b,c+n)\big| \\ [1.25em]
%&\leq& \displaystyle\sum_{k=1}^\infty u_{k}\big(|a|,|b|,n-|c|\big)=
%\sum_{k=0}^\infty u_{k+1}\big(|a|,|b|,n-|c|\big)\;.
%\end{array}
%$$
$$
\left|\pFq{2}{1}{a,b}{c+n}{1}-1\right| =
\left|\sum_{k=1}^\infty u_k(a,b,c+n)\right|\leq
\sum_{k=1}^\infty u_{k}\big(|a|,|b|,n-|c|\big)\;.
$$
Thus
\begin{equation}\label{2F1-Gauss4}
\left|\pFq{2}{1}{a,b}{c+n}{1}-1\right|\leq\sum_{k=0}^\infty u_{k+1}\big(|a|,|b|,n-|c|\big)\;,
\quad n>|c|\;.
\end{equation}
Next, notice that, for each $k\in\mathbb{N}_0$ and $n>|c|$,
$$
\begin{array}{rl}
u_{k+1}\big(|a|,|b|,n-|c|\big)&=\;\displaystyle\frac{(|a|)_{k+1}(|b|)_{k+1}}{(n-|c|)_{k+1}\,(k+1)!}
=\frac{1}{k+1}\frac{|ab|}{n-|c|}\frac{(|a|+1)_{k}(|b|+1)_{k}}{(n+1-|c|)_{k}\,k!} \\ [1em]
&\leq\;\displaystyle\frac{|ab|}{n-|c|}u_{k}\big(|a|+1,|b|+1,n+1-|c|\big)\,.
\end{array}
$$
Therefore, from (\ref{2F1-Gauss4}) we obtain
%$$
%\left|\pFq{2}{1}{a,b}{c+n}{1}-1\right|\leq\frac{|ab|}{n-|c|}\,
%\sum_{k=0}^\infty u_{k}\big(|a|+1,|b|+1,n+1-|c|\big)\;,\quad n>|c|\,.
%$$
%This inequality may be rewritten as
$$
\left|\pFq{2}{1}{a,b}{c+n}{1}-1\right|\leq\frac{|ab|}{n-|c|}\,
\pFq{2}{1}{|a|+1,|b|+1}{n+1-|c|}{1}\;,\quad n>|c|\,.
$$
According to (i) in Theorem \ref{pFq-conv-thm2},
the series $\pFq{2}{1}{|a|+1,|b|+1}{n+1-|c|}{1}$ converges (absolutely) if $n>|a|+|b|+|c|+1$.
This series is, clearly, a decreasing function of $n$, hence it is bounded
by a positive number independent of $n$, say, $M\equiv M(a,b,c)>0$, and so
$$
\left|\pFq{2}{1}{a,b}{c+n}{1}-1\right|\leq\frac{|ab|M}{n-|c|}\quad\mbox{\rm if}\quad n>|a|+|b|+|c|+1\,.
$$
Therefore, taking the limit as $n\to+\infty$ we obtain (\ref{2F1-Gauss3}).
\qed

\begin{snugshade}
\begin{corollary}[Chu-Vandermonde]\label{2F1-Chu-Vandermonde}
For each $n\in\mathbb{N}_0$, %and $\Re (c-a)>0$,
\begin{equation}\label{2F1-ChuVand}
\pFq{2}{1}{-n,a}{c}{1}=\frac{(c-a)_n}{(c)_n}\;.
\end{equation}
\end{corollary}
\end{snugshade}

{\it Proof.}
The result follows immediately taking $b=-n$ in (\ref{2F1-4}) and using (\ref{Poch-Gamma}).
\qed

\section{The Pfaff-Saalsch\"utz, Dixon's, and Dougall's identities}

The Chu-Vandermonde identity (\ref{2F1-ChuVand}) gives a closed formula
for a terminating ${}_{2}F_{1}$ hypergeometric series.
Similarly, the Pfaff-Saalsch\"utz identity gives a closed formula
for a terminating ${}_{3}F_{2}$ hypergeometric series.
These kind of formulas are very useful on the computation of binomial sums in closed form,
as we will see in the next section.

\begin{snugshade}
\begin{theorem}[Pfaff-Saalsch\"utz]\label{3F2-Pfaff-S}
For each $n\in\mathbb{N}$,
\begin{equation}\label{2F1-7}
\pFq{3}{2}{-n,a,b}{c,1+a+b-c-n}{1}=\frac{(c-a)_n(c-b)_n}{(c)_n(c-a-b)_n}\;.
\end{equation}
\end{theorem}
\end{snugshade}

{\it Proof.}
By Euler's transformation formula (\ref{2F1-PE2}) and the binomial theorem (\ref{binThm}),
$$
\pFq{1}{0}{c-a-b}{\mbox{---}}{x}\cdot\pFq{2}{1}{a,b}{c}{x}=\pFq{2}{1}{c-a,c-b}{c}{x}\;,\quad |x|<1\;.
$$
Rewrite this equation as
$$
\left(\,\sum_{n=0}^\infty \frac{(c-a-b)_n}{n!}\,x^n\right)
\left(\,\sum_{n=0}^\infty \frac{(a)_n(b)_n}{(c)_n\,n!}\,x^n\right)
=\sum_{n=0}^\infty \frac{(c-a)_n(c-b)_n}{(c)_n\,n!}\,x^n\;,\quad |x|<1\;.
$$
Form the Cauchy product of the series on the left-hand side and then equate the coefficients of $x^n$
in both sides of the resulting equality. This yields
\begin{equation}\label{2F1-7a}
\sum_{j=0}^n \frac{(a)_j(b)_j(c-a-b)_{n-j}}{j!\,(c)_j(n-j)_j}=\frac{(c-a)_n(c-b)_n}{(c)_n\,n!}\;,
\quad n\in\mathbb{N}_0\;.
\end{equation}
Now, taking into account the equalities
$$
(\alpha)_{n-j}=\frac{(-1)^j(\alpha)_n}{(1-\alpha-n)_j}\;,
\quad \frac{(-1)^jn!}{(n-j)!}=(-n)_j\;,\quad 0\leq j\leq n\;,
$$
%$$
%(\alpha)_{n-j}=\frac{(-1)^j(\alpha)_n}{(1-\alpha-n)_j}\;,\quad 0\leq j\leq n\;,
%%\quad\alpha\in\mathbb{C}\setminus\mathbb{N}\;,
%$$
%and since $(-1)^jn!/(n-j)!=(-n)_j$,
the sum on the left-hand side of (\ref{2F1-7a}) becomes
$$
%\sum_{j=0}^n \frac{(a)_j(b)_j(c-a-b)_{n-j}}{j!\,(c)_j(n-j)_j}=
\sum_{j=0}^n \frac{(a)_j(b)_j(-n)_j(c-a-b)_n}{j!(c)_j(1+a+b-c-n)_j\,n!}
=\frac{(c-a-b)_n}{n!}\pFq{3}{2}{-n,a,b}{c,1+a+b-c-n}{1}\; .
$$
Thus the theorem is proved.
\qed
\medskip

In the next section we will present examples illustrating how
Chu-Vandermonde and Pfaff-Saalsch\"utz identities can be useful
to obtain closed formulas for sums involving binomial coefficients.
In the applications to such {\it binomial identities}, often the case $p=q+1$ occur,
the success of the procedure depending upon certain relations fulfilled by the parameters
$a_1,\ldots,a_{q+1}$ and $b_1,\ldots,b_{q}$ appearing in the definition of
\begin{equation}\label{pFq-sumsBin}
\pFq{q+1}{q}{a_1,\cdots,a_{q+1}}{b_1,\cdots,b_{q}}{x}\;.
\end{equation}
The series (\ref{pFq-sumsBin}) is called {\sl $k-$balanced} at $x=1$ if
one of the $a_j$'s is a negative integer number, and the following condition holds:
\begin{snugshade}\vspace*{-0.5em}
\begin{equation}\label{pFq-balanced-def}
%k+a_1+\cdots+a_{q+1}=b_1+\cdots+b_q\;.
k+\sum_{j=1}^{q+1}a_j=\sum_{j=1}^{q}b_j\;.
\end{equation}
\end{snugshade}
\noindent
The condition that one of the $a_j$'s is a negative integer number means that the series terminates.
This condition seem artificial, but without it many results do not hold.
An {\sl $1-$balanced} series is also called {\sl Saalsch\"utzian}.
(\ref{pFq-sumsBin}) is called {\sl well-poised} if
\begin{snugshade}\vspace*{-0.5em}
\begin{equation}\label{pFq-well-poised-def}
1+a_1=b_1+a_2=\cdots=b_q+a_{q+1}\;.
\end{equation}
\end{snugshade}
\noindent
We conclude by stating without proof (the proofs can be founded in several of the textbooks
presented in the Bibliography) two theorems involving two identities of these types.
Dixon's identity applies to a well-poised ${}_3F_2$ series,
while Dougall's identity applies to a well-poised $2-$balanced ${}_7F_6$ series.
%The parameters appearing in these identities are supposed to be
%such that all the involved expressions are well defined.
%We omit the proofs (Exercise \ref{Ex-cp16-4a}).
%The proofs are left to the reader (Exercise \ref{Ex-cp16-4a}).

%\begin{snugshade}
%\begin{theorem}[Dixon]\label{3F2-Dixon}
%The identity
%\begin{equation}\label{3F2-Dix1}
%\begin{array}{l}
%\displaystyle\pFq{3}{2}{a,-b,-c}{1+a+b,1+a+c}{1} \\ [1.5em]
%\qquad
%=\displaystyle\frac{\Gamma(1+a/2)\Gamma(1+a+b)\Gamma(1+a+c)\Gamma(1+a/2+b+c)}
%{\Gamma(1+a)\Gamma(1+a/2+b)\Gamma(1+a/2+c)\Gamma(1+a+b+c)}
%\end{array}
%\end{equation}
%holds, where the condition $a+2b+2c+2>0$ is assumed whenever the left-hand side is an infinite series.
%\end{theorem}
%\end{snugshade}

\begin{snugshade}
\begin{theorem}[Dixon]\label{3F2-Dixon}
The identity
\begin{equation}\label{3F2-Dix1}
\pFq{3}{2}{a,-b,-c}{1+a+b,1+a+c}{1}
=\frac{\Gamma\big(1+\frac{a}{2}\big)\Gamma(1+a+b)\Gamma(1+a+c)\Gamma\big(1+\frac{a}{2}+b+c\big)}
{\Gamma(1+a)\Gamma\big(1+\frac{a}{2}+b\big)\Gamma\big(1+\frac{a}{2}+c\big)\Gamma(1+a+b+c)}
\end{equation}
holds, where the condition $\Re(a+2b+2c+2)>0$ is assumed whenever the left-hand side is an infinite series.
\end{theorem}
\end{snugshade}

\begin{snugshade}
\begin{theorem}[Dougall]\label{7F6-Dougall}
For each $n\in\mathbb{N}$, % and $a,b,c,d,e\in\mathbb{C}$,
\begin{equation}\label{7F6-Doug1}
\begin{array}{r}
\displaystyle\pFq{7}{6}{-n,a,1+a/2,-b,-c,-d,-e}{a/2,1+a+b,1+a+c,1+a+d,1+d+e,1+a+n}{1}\qquad\qquad \\ [1.5em]
%\qquad\qquad
=\displaystyle\frac{(1+a)_n(1+a+b+c)_n(1+a+b+d)_n(1+a+c+d)_n}{(1+a+b)_n(1+a+c)_n(1+a+d)_n(1+a+b+c+d)_n}\;,
\end{array}
\end{equation}
provided that $\;1+2a+b+c+d+e+n=0$.
\end{theorem}
\end{snugshade}

%Taking $n\to+\infty$ we obtain another identity due to Dougall:
%
%\begin{snugshade}
%\begin{theorem}[Dougall]\label{5F4-Dougall}
%For each $n\in\mathbb{N}$, % and $a,b,c,d,e\in\mathbb{C}$,
%\begin{equation}\label{5F4-Doug1}
%\begin{array}{r}
%\displaystyle\pFq{5}{4}{a,1+a/2,-b,-c,-d}{a/2,1+a+b,1+a+c,1+a+d}{1}\qquad\qquad \\ [1.5em]
%%\qquad\qquad
%=\displaystyle\frac{\Gamma(1+a+b)\Gamma(1+a+c)\Gamma(1+a+d)\Gamma(1+a+b+c+d)}
%{\Gamma(1+a)\Gamma(1+a+b+c)\Gamma(1+a+b+d)\Gamma(1+a+c+d)}\;,
%\end{array}
%\end{equation}
%provided that $\;\Re(a+b+c+d+1)>0$.
%\end{theorem}
%\end{snugshade}
%
%Dixon's identity is the special case $d=-a/2$.

It is worth mentioning that many other interesting identities are given in the books presented in the Bibliography.

\section{Binomial sums}\label{SumBin}
One area where hypergeometric identities are very useful is in the evaluation
of sums of products of binomial coefficients.
The main idea behind this procedure is writing such a sum as an hypergeometric series.
In this section we present three examples illustrating the power of this technique.
While working on examples of this type, we need to compute quotients involving binomial coefficients, so
often it is useful to make use of the following identities (easy to check),
which hold for $\alpha\in\mathbb{C}$ and $k\in\mathbb{N}_0$: %that one should keep in mind:
\begin{snugshade}
\begin{equation}\label{BinEx0}
\begin{array}{ll}
\displaystyle\mbox{\rm (i)}\quad\frac{{{\alpha+1}\choose{k+1}}}{{{\alpha}\choose{k}}}
=\frac{\alpha+1}{k+1}\;, & \qquad
\displaystyle\mbox{\rm (iii)}\quad\frac{{{\alpha+2}\choose{k+1}}}{{{\alpha}\choose{k}}}
=\frac{(\alpha+2)(\alpha+1)}{(k+1)(\alpha-k+1)}\;, \\ [1.25em]
\displaystyle\mbox{\rm (ii)}\quad\frac{{{\alpha}\choose{k+1}}}{{{\alpha}\choose{k}}}=\frac{\alpha-k}{k+1}\;, & \qquad
\displaystyle\mbox{\rm (iv)}\quad\frac{{{\alpha+1}\choose{k+2}}}{{{\alpha}\choose{k}}}
=\frac{(\alpha+1)(\alpha-k)}{(k+2)(k+1)} \;.
\end{array}
\end{equation}
\end{snugshade}
\medskip

\noindent
It is also useful to keep in mind the relations
\begin{snugshade}
$$
{{\alpha}\choose{n}}:=\frac{\alpha(\alpha-1)\cdots(\alpha-n+1)}{n!}
=\frac{(-1)^n(-\alpha)_n}{n!}=\frac{(\alpha-n+1)_n}{n!}\;,
$$
\end{snugshade}
\noindent
which hold for $\alpha\in\mathbb{C}$ and $n\in\mathbb{N}$.

{\sc Example 1.} As a first example, we show that
\begin{snugshade}
\begin{equation}\label{BinEx1}
\sum_{j=0}^n(-1)^j\frac{{{\alpha}\choose{j}} {{\alpha-1-j}\choose{n-j}}}{j+1}
=\frac{{{\alpha}\choose{n+1}}+(-1)^n}{\alpha+1}\;,\quad
\alpha\in\mathbb{C}\setminus\{-1\}\;,\;\; n\in\mathbb{N}_0\;.
\end{equation}
\end{snugshade}

{\it Proof.}
Denote the sum of the left-hand side by $S$, so that
$$
S:=\sum_{j=0}^n c_j\;,\quad
c_j:=(-1)^j\frac{{{\alpha}\choose{j}} {{\alpha-1-j}\choose{n-j}}}{j+1}\;.
$$
To write this sum as an hypergeometric series, we first compute the ratio $c_{j+1}/c_j$
and then we put it in the form (\ref{HF1}):
$$
\frac{c_{j+1}}{c_j}
=-\frac{(j+1){{\alpha}\choose{j+1}} {{\alpha-2-j}\choose{n-1-j}}}
{(j+2){{\alpha}\choose{j}} {{\alpha-1-j}\choose{n-j}}}
=\frac{(j-n)(j-\alpha)(j+1)}{(j-\alpha+1)(j+2)}\frac{1}{j+1}\;,
$$
the last equality following immediately by (ii) and (i) in (\ref{BinEx0}).
Thus --- cf. (\ref{HF2}) ---,
\begin{equation}\label{BinEx1f1}
S:=\sum_{j=0}^n c_j=c_0\,\sum_{j=0}^n\frac{(-n)_j(-\alpha)_j(1)_j}{(-\alpha+1)_j(2)_jj!}\;,\quad
c_0:={{\alpha-1}\choose{n}}\;,
\end{equation}
and so the given binomial sum can be written in hypergeometric form as
\begin{equation}\label{BinEx1f2}
S={{\alpha-1}\choose{n}}\,
\pFq{3}{2}{-n,-\alpha,1}{-\alpha+1,2}{1}\;.
\end{equation}
At this point, one could try to apply the Pfaff-Saalsch\"utz identity (\ref{2F1-7}).
However, the ${}_3F_2$ series in (\ref{BinEx1f2}) is not of the form
of the ${}_3F_2$ appearing in (\ref{2F1-7}). (Indeed, if
$a=-\alpha$, $b=1$, and $c=-\alpha+1$, then $1+a+b-c-n=1-n\neq2$.)
Thus, the Pfaff-Saalsch\"utz identity does not apply.
Nevertheless, returning to (\ref{BinEx1f1}), and noting that $(1)_j=j!$,
$(2)_j=(1)_{j+1}=(j+1)!$, and $(z)_j=(z-1)_{j+1}/(z-1)$, we may write
\begin{equation}\label{BinEx1f3}
\begin{array}{rl}
S&=\;\displaystyle -{{\alpha-1}\choose{n}}\frac{\alpha}{(n+1)(\alpha+1)}
\sum_{j=0}^n \frac{(-n-1)_{j+1}(-\alpha-1)_{j+1}}{(-\alpha)_{j+1}(j+1)!} \\ [1.25em]
&=\;\displaystyle -{{\alpha-1}\choose{n}}\frac{\alpha}{(n+1)(\alpha+1)}
\left(\,\sum_{j=0}^{n+1} \frac{(-n-1)_{j}(-\alpha-1)_{j}}{(-\alpha)_{j}\,j!}-1\right) \\ [1.25em]
&=\;\displaystyle -{{\alpha-1}\choose{n}}\frac{\alpha}{(n+1)(\alpha+1)}
\left[\,\pFq{2}{1}{-(n+1),-\alpha-1}{-\alpha}{1}-1\right] \;.
\end{array}
\end{equation}
The last ${}_2F_1$ may be computed by the Chu-Vandermonde identity
(\ref{2F1-ChuVand}), and so
$$
\pFq{2}{1}{-(n+1),-\alpha-1}{-\alpha}{1}=\frac{(1)_{n+1}}{(-\alpha)_{n+1}}=\frac{(n+1)!}{(-\alpha)_{n+1}}\;.
$$
Inserting this expression into (\ref{BinEx1f3}) and simplifying the result, we obtain (\ref{BinEx1}).
%Inserting this expression into the right-hand side of (\ref{BinEx1f3})
%and simplifying the resulting equality, we obtain (\ref{BinEx1}).
\qed
\medskip

{\sc Example 2.} As a second example, let us show that
\begin{snugshade}
\begin{equation}\label{BinEx2}
\sum_{j=0}^{n-m}\frac{(-1)^j}{j+1}\,{{n+j}\choose{m+2j}} {{2j}\choose{j}}={{n-1}\choose{m-1}}\;,
\quad m,n\in\mathbb{N}\,.
%\quad 1\leq m\leq n\quad (m,n\in\mathbb{N})\;.
\end{equation}
\end{snugshade}

{\it Proof.}
It is clear that (\ref{BinEx2}) holds if $n<m$, since in such case both sides of (\ref{BinEx2}) are equal to zero
%(with the usual convention that empty sum equals zero---alternatively,
%observe that ${{n+j}\choose{m+2j}}=0$ if $n+j<m+2j$, i.e., if $j>n-m$).
(this holds because ${{k}\choose{\ell}}=0$ if $k,\ell\in\mathbb{N}$ with $k<\ell$).
Henceforth we assume $1\leq m\leq n$. Let
$$
S:=\sum_{j\geq0} c_j\;,\quad
c_j:=\frac{(-1)^j}{j+1}\,{{n+j}\choose{m+2j}}{{2j}\choose{j}}\;.
$$
Since ${{n+j}\choose{m+2j}}=0$ if $n+j<m+2j$, i.e., if $j>n-m$,
then $c_j=0$ if $j>n-m$, and so $S$ is indeed a finite sum.
Using (iii) and (iv) in (\ref{BinEx0}), we compute
%To write this sum as an hypergeometric series, we first compute the ratio $c_{j+1}/c_j$
%and then we put it in the form (\ref{HF1}):
$$
\frac{c_{j+1}}{c_j}
=-\frac{(j+1){{n+j+1}\choose{m+2j+2}}{{2j+2}\choose{j+1}}}
{(j+2){{{n+j}\choose{m+2j}} {{2j}\choose{j}}}}
=\frac{(j+n+1)(j-n+m)\big(j+\frac12\big)}{\big(j+\frac{m}{2}+1\big)\big(j+\frac{m+1}{2}\big)(j+2)}\,.
$$
Thus, since $c_0={{n}\choose{m}}$, we obtain --- cf. (\ref{HF2}) ---,
%$$
\begin{equation}\label{Examp2aa}
S={{n}\choose{m}}\sum_{j\geq0}\frac{(n+1)_j(m-n)_j\big(\frac12\big)_j}
{\big(\frac{m}{2}+1\big)_j\big(\frac{m+1}{2}\big)_j(j+1)!}
={{n}\choose{m}}\sum_{j\geq1}\frac{(n+1)_{j-1}(m-n)_{j-1}\big(\frac12\big)_{j-1}}
{\big(\frac{m}{2}+1\big)_{j-1}\big(\frac{m+1}{2}\big)_{j-1}\,j!}\,.
\end{equation}
%$$
If $m>1$, using  $(z+1)_{j-1}=(z)_j/z$, the last sum can be written as
\begin{equation}\label{BinEx2f1}
S=-\mbox{$\frac12$}\,{{n}\choose{m}}\,\frac{m(m-1)}{n(m-n-1)}\,
\left[\,\pFq{3}{2}{n,m-n-1,-\frac12}{\frac{m}{2},\frac{m-1}{2}}{1}-1\,\right]\,.
\end{equation}
Since $m<n+1$, by the Pfaff-Saalsch\"utz identity (\ref{2F1-7}) the last ${}_3F_2$ series becomes
$$
%\pFq{3}{2}{n,m-n-1,-\frac12}{\frac{m}{2},\frac{m-1}{2}}{1}=
\pFq{3}{2}{-(n+1-m),n,-\frac12}{\frac{m-1}{2},\frac{m}{2}}{1}
=\frac{\big(\frac{m-1}{2}-n\big)_{n+1-m}\big(\frac{m}{2}\big)_{n+1-m}}
{\big(\frac{m-1}{2}\big)_{n+1-m}\big(\frac{m}{2}-n\big)_{n+1-m}}\;.
$$
Inserting this expression into the right-hand side of (\ref{BinEx2f1}) and simplifying
the resulting equality --- in this simplification process the relation
$$(\alpha-n)_k=(-1)^k(-\alpha+n-k+1)_k$$ may be useful ---, we obtain (\ref{BinEx2}) whenever $m>1$.
If $m=1$, noting that $(1)_j=j!$ and $(j+1)!=(2)_j$, the first equality in (\ref{Examp2aa}) gives
$$
%S={{n}\choose{1}}\sum_{j\geq0}\frac{(n+1)_j(1-n)_j\big(\frac12\big)_j}
%{\big(\frac{3}{2}\big)_j(2)_j\,j!}
S={{n}\choose{1}}\,\pFq{3}{2}{-(n-1),n+1,\frac12}{\frac{3}{2},2}{1}=
n\,\frac{\big(\frac{1}{2}-n\big)_{n-1}(1)_{n-1}}
{\big(\frac{3}{2}\big)_{n-1} (-n)_{n-1}}=1\;,
$$
where in the second equality we have used again the Pfaff-Saalsch\"utz identity (\ref{2F1-7}).
This proves (\ref{BinEx2}) for $m=1$.
Notice that since $c_j=0$ if $j>n-m$ then in (\ref{BinEx2})
one may replace $\sum_{j=0}^{n-m}$ by $\sum_{j\geq0}$.
\qed
\bigskip

{\sc Example 3.}
As a last example, we show that if $n,p\in\mathbb{N}$ %and $\ell:=\min\{n,p\}$,
\begin{snugshade}\vspace*{-0.5em}
\begin{equation}\label{BinEx3}
\sum_{k=1}^{\ell}2k\,{{2p}\choose{k+p}}{{2n}\choose{k+n}}=
\frac{4np}{n+p}\,{{2p-1}\choose{p}} {{2n-1}\choose{n}}\;,\quad
\ell:=\min\{n,p\}\;.
\end{equation}
\end{snugshade}

{\it Proof.}
Denote the sum of the left-hand side of (\ref{BinEx3}) by $S$. Then
$$
S=\sum_{k\geq0} c_k\;,\quad
c_k:=2(k+1)\,{{2p}\choose{k+1+p}}{{2n}\choose{k+1+n}}\;.
$$
The last equality holds, indeed, since $c_k=0$ if $k\geq\min\{n,p\}=:\ell$.
To write this sum as a (terminating) hypergeometric series, we compute the ratio $c_{k+1}/c_k\,$:
%(for $k=1,2,\ldots,\ell$):
%and then we put it in the form (\ref{HF1}):
$$
\frac{c_{k+1}}{c_k}
=\frac{(k+2)(k+1-p)(k+1-n)}{(k+2+p)(k+2+n)}\frac{1}{k+1}\;,\quad
k=0,1,\ldots,\ell-1\;.
$$
Thus --- cf. (\ref{HF2}) ---, the given sum can be written in hypergeometric form as
$$
S=c_0\,\sum_{k\geq0}\frac{(2)_k(1-p)_k(1-n)_k}{(2+p)_k(2+n)_k}
=2\,{{2p}\choose{1+p}}{{2n}\choose{1+n}}\,
\pFq{3}{2}{2,1-p,1-n}{2+p,2+n}{1}\,.
$$
This ${\,}_3F_2$ series can be computed using Dixon's identity (\ref{3F2-Dix1}),
taking therein $a=2$, $b=p-1$, and $c=n-1$,
% $d=2+p$, and $e=2+n$. In fact, we see that $a+b+1=d$ and $a+c+1=e$, hence Dixon's identity (\ref{3F2-Dix1}) is applicable,
and so
$$
\begin{array}{rcl}
S&=&\displaystyle 2\,{{2p}\choose{1+p}}{{2n}\choose{1+n}}\,
\frac{\Gamma(2)\Gamma(2+p)\Gamma(2+n)\Gamma(p+n)}{\Gamma(3)\Gamma(1+p)\Gamma(1+n)\Gamma(p+n+1)} \\ [1.25em]
&=&\displaystyle {{2p}\choose{1+p}}{{2n}\choose{1+n}}\,
\frac{(1+p)(1+n)}{p+n}\,,
\end{array}
$$
where the last equality follows from (\ref{GammaF2}).
Therefore, taking into account (i) in (\ref{BinEx0}), we obtain (\ref{BinEx3}).
\qed
\medskip

\begin{remark}\em
The reader is invited to read the very interesting article \ref{Roy1987-C16},
where Ranjan Roy presented (the above and) several other examples,
pointing out the power of this technique to compute intricate binomial sums.
\end{remark}

\begin{remark}\em
A powerful technique to prove identities between hypergeometric functions
was developed by Zeilberger and Wilf, called the {\it creative telescoping} method.
This method is also referred to as the {\it W--Z method}, and it is described e.g.
in the books \ref{AndrewsAskeyRoy1999-C16} and \ref{PetWilfZeilb1996-C16}.
It is worth mentioning that in \ref{AndrewsAskeyRoy1999-C16} (p.\,175), Andrews, Askey, and Roy wrote:
``\,The discoveries of Wilf and Zeilberger truly revolutionized the study
of summations of terminating hypergeometric series.''
As a matter of fact, many further developments of these techniques have appeared since then,
including extensions to the so-called (basic) $q-$hypergeometric series as well as
fully algorithm implementation on the computer.
\end{remark}

%\section{The hypergeometric equation}

\section*{Exercises}
%\bigskip

{\small
%\noindent
\begin{enumerate}[label=\emph{\bf \arabic*.},leftmargin=*]
\item\label{Ex-cp16-1}
Prove the limit relation (\ref{BetaF10}).
%You should use Lebesgue's convergence dominated theorem and the connections between
%the Lebesgue integral and the (proper and improper) Riemann integral.
\smallskip

\noindent
({\sl Hint:}
Use Lebesgue's convergence dominated theorem and the connections between
the Lebesgue integral and the proper and improper Riemann integrals.
It may be useful to notice that for real $t$ and $a>0$,
$$\Big(1-\frac{t}{n}\Big)^{a+n-1}\chi_{[0,n]}(t)\leq e^{1-t}\;,\quad \forall n\in\mathbb{N}\,.)\smallskip$$

\item\label{Ex-cp16-2}
Prove Euler's reflection fomula (\ref{FGammaReflect}).
%$$
%\Gamma(z)\Gamma(1-z)=\frac{\pi}{\sin(\pi z)}\,,\quad z\in\mathbb{C}\setminus\mathbb{Z}\;.
%$$%
\smallskip

\noindent
({\sl Hint.} Set $t=s/(1+s)$ in the definition of the beta integral to obtain
$$
\int_0^\infty\frac{s^{x-1}}{(1+s)^{x+y}}\,{\rm d}s=\frac{\Gamma(x)\Gamma(y)}{\Gamma(x+y)}\;,\quad
\Re x>0\;,\;\;\Re y>0\;;
$$
this gives
$$
\Gamma(x)\Gamma(1-x)=\int_0^\infty\frac{t^{x-1}}{1+t}\,{\rm d}t\;,\quad 0<x<1\;.
$$
The last integral can be computed by the residue theorem using the contour integral
$$
\int_{\Gamma_{\epsilon,R}}\frac{z^{x-1}}{1-z}\,{\rm d}z\;,
$$
where $\Gamma_{\epsilon,R}:=C_{\epsilon,R}\cup\ell_{\epsilon,R}^+\cup C_{\epsilon}\cup\ell_{\epsilon,R}^-$ is a closed path, $C_{\epsilon,R}$ is an incomplete circle around the origin of radius $R$ with starting and ending points at $z=-R\cos\theta\pm i\epsilon$, not containing the point $z=-R$, being $0<\epsilon<1<R$ and $\theta:=\arcsin(\epsilon/R)$, $C_{\epsilon}$ is the semicircle around the origin of radius $\epsilon$ joining the points $z=\pm i\epsilon$ and containing $z=\epsilon$,
%in the semi-plane $\Re z\geq0$,
and $\ell_{\epsilon,R}^{\pm}$ are two segments parallel to the negative real axis, one of them starting at $z=-R\cos\theta+i\epsilon$ and ending at $z=i\epsilon$, and the other one starting at $z=-i\epsilon$ and ending at $z=-R\cos\theta-i\epsilon$.)

%where $C\equiv C_{\epsilon,R}$ is a closed path consisting of two (incomplete) circles about the origin of radii $R$ and $\epsilon$, being $0<\epsilon<R$, which are joined by two straight lines parallel to the negative real axis, on of them from $-R\cos\theta+i\epsilon$ to $i\epsilon$, and the other one from $-i\epsilon$ to $-R\cos\theta-i\epsilon$, being $\theta:=\arcsin(\epsilon/R)\in]0,\pi/2[$.)
%along the negative real axis from $-R$ to $-\epsilon$.)
\medskip

\item\label{Ex-cp16-3aaa}
% Rainville, pp.23-24
Prove Legendre's duplication formula:
$$
\quad\sqrt{\pi}\,\Gamma(2x)=2^{2x-1}\,\Gamma(x)\Gamma\big(x+\mbox{$\frac12$}\big)\;,\quad
x\in\mathbb{C}\setminus\mbox{$\big\{-\frac{k}{2}\,:\,k\in\mathbb{N}_0\big\}$}\;.
%x\in\mathbb{C}\setminus\mbox{$\big\{0,-\frac12,-1,-\frac32,-2,-\frac52,\ldots\big\}$}\;.
$$
%\smallskip

\noindent
({\sl Hint.}
Use $(2x)_{2n}=2^{2n}(x)_n\big(x+\frac12\big)_n$ together with (\ref{Poch-Gamma})
and the definition of $\Gamma$.) % the gamma function.)
%Using (\ref{Poch-Gamma}) and $(2x)_{2n}=2^{2n}(x)_n\big(x+\frac12\big)_n$,
%we may start by showing that, for each $n\in\mathbb{N}_0$, the identity
%$\frac{\Gamma(2x+2n)}{2^{2n}\,\Gamma(x+n)\Gamma\big(x+\frac12+n\big)}
%=\frac{\Gamma(2x)}{\Gamma(x)\Gamma\big(x+\frac12\big)}$
%holds. Then, notice that the right-hand side of this identity is independent of $n$,
%hence after inserting appropriate factors
%in the left-hand side ---to permit us to make use of (\ref{GammaF1})---and
%then taking the limit as $n\to+\infty$, we obtain
%$$
%\lim_{n\to+\infty}\frac{2^{2x}(2n-1)!}{2^{2n}\sqrt{n}\,[(n-1)!]^2}=
%\frac{\Gamma(2x)}{\Gamma(x)\Gamma\big(x+\frac12\big)}\;.
%$$
\medskip

\item\label{Ex-cp16-3a}
Prove the hypergeometric series representations (i)---(vi) given in Remark \ref{pFq-elem1}.
\medskip

\item\label{Ex-cp16-2a}
Show that the estimative (\ref{1-xtaCdelta}) holds.
\smallskip

\noindent
({\sl Hint.} For each $z_0\in\mathbb{C}\setminus[1,+\infty)$,
(\ref{1-xtaCdelta}) holds if we define $\delta$ and $C(a,z_0,\delta)$ as
$$
\begin{array}{l}
\delta:=\frac12\,|\Im z_0|\chi_{\mathbb{R}\setminus\{0\}}(\Im z_0) \\ [0.4em]
\qquad\quad +\big[\frac12\,(1-z_0)\chi_{(-\infty,0]}(\Re a)+
\frac16\,(2-z_0-|z_0|)\chi_{(0,+\infty)}(\Re a)\big]\chi_{\{0\}}(\Im z_0)\,,
%\frac13\,\min\{1,1-z_0\}\chi_{(0,\infty)}(\Re a)\big]\chi_{\{0\}}(\Im z_0)\,,
\end{array}
$$
$$
\begin{array}{l}
C(a,z_0,\delta):=(1+|z_0|+\delta)^{-\Re a}\,\chi_{(-\infty,0]}(\Re a) \\ [0.25em]
\qquad\qquad\qquad +\big[1/\delta^{\Re a}\chi_{\{0\}}(\Im z_0)+
(1+|z_0|/\delta)^{\Re a}\chi_{\mathbb{R}\setminus\{0\}}(\Im z_0)\big]\,\chi_{(0,+\infty)}(\Re a)\,.
\end{array}
$$
In case $\Re a>0$, it may be useful to notice that
$$
\qquad
|1-xt|^2=\left\{
\begin{array}{lcl}
(1-z_0t)[(1-z_0t)+2t(z_0-\Re x)]+|x-z_0|^2t^2 &\mbox{\rm if}& \Im z_0=0 \\ [0.2em]
|x|^2\Big(t-\frac{\Re x}{|x|^2}\Big)^2+\Big(\frac{\Im x}{|x|}\Big)^2 &\mbox{\rm if}& \Im z_0\neq0
\end{array}
\right.
$$
holds for each $x\in\overline{B}(z_0,\delta)$ and $t\in[0,1]$, and so
%(in case $\Re a>0$), for each $x\in\overline{B}(z_0,\delta)$ and $t\in[0,1]$,
$$
\qquad
|1-xt|\geq\left\{
\begin{array}{lcl}
\big\{(1-z_0t)[(1-z_0t)+2t(z_0-\Re x)]\big\}^{1/2}\geq\delta%\{\delta(3\delta-2\delta)\}^{1/2}
 &\mbox{\rm if}& \Im z_0=0 \\ [0.5em]
|\Im x|/|x|\geq\delta/(\delta+|z_0|) &\mbox{\rm if}& \Im z_0\neq0\,.
\end{array}
\right.
$$
Then use (\ref{1-xta}) to obtain the desired estimative.)
\medskip

\item\label{Ex-cp16-5aa}
% see Rainville: Thm. 17, p.48; Thm. 19, p.56
Prove the following statements:
\begin{enumerate}
\item If $x$ is fixed in $\mathbb{C}$ and $|x|<1$,
then $\pFq{2}{1}{a,b}{c}{x}$ is an analytic function of the
variables $a$, $b$, and $c$ for all finite (complex) values of $a$, $b$, and $c$,
except for simple poles at $c\in\{0,-1,-2,\ldots\}$.
\item $\pFq{2}{1}{a,b}{c}{1}$ is an analytic function of
$a$, $b$, and $c$ for all finite values of $a$, $b$, and $c$
such that $\Re(c-a-b)>0$ and $c\in\mathbb{C}\setminus\{0,-1,-2,\ldots\}$.
\end{enumerate}
\medskip

\item\label{Ex-cp16-5ab}
% see Rainville: Thm. 18, p.49
Give an alternative proof to Gauss's summation formula (Theorem \ref{2F1-Gauss}),
by using Euler's integral representation (or the technique of its proof)
to firstly state Gauss's formula for $\Re c>\Re b>0$,
and then removing this constraint by analytic continuation on the parameters,
using statement (b) in exercise \ref{Ex-cp16-5aa}
\medskip

\item\label{Ex-cp16-5abc}
% see Rainville, p.49
%Show that, for each $n\in\mathbb{N}$ and $\Re b>0$,
Use Gauss's summation formula and Legendre's duplication formula to show that
$$
\pFq{2}{1}{-\frac{n}{2},-\frac{n-1}{2}}{\frac{2b+1}{2}}{1}=\frac{2^n\,(b)_n}{(2b)_n}\;,
\quad n\in\mathbb{N}\;,\quad \Re b>0\;.
$$

\item\label{Ex-cp16-4}
%AndrewAskeyRoy,p.108 (questao alternativa: livro A=B,p.44-45)
Let $m$ and $n$ be nonnegative integer numbers. Prove that
$$
\sum_{k=0}^n{{m}\choose{k}}{{m+n-k}\choose{m}}\frac{(-1)^k}{m+n+1-k}=\frac{n!}{(m+1)_{n+1}}\;.
$$
({\sl Hint.} This sum can be written as
${{m+n}\choose{m}}\;\pFq{3}{2}{-n,-m,-m-n-1}{-m-n,-m-n}{1}/(m+n+1)$.)
\medskip

%\item\label{Ex-cp16-4a}
%Prove Dougall's identity (\ref{7F6-Doug1}) and use it to derive Dixon's identity (\ref{3F2-Dix1}).
%\medskip

\item\label{Ex-hwk6-4b}
Use Dixon's identity and Euler's reflection formula to show that
$$
\sum_{j=-\ell}^\ell(-1)^j{{2\ell}\choose{\ell+j}}{{2m}\choose{m+j}}{{2n}\choose{n+j}}=
\frac{{{2\ell}\choose{\ell}}{{2m}\choose{m}}{{2n}\choose{n}}{{\ell+m+n}\choose{m+n}}}
{{{\ell+m}\choose{\ell}}{{\ell+n}\choose{\ell}}}\;,
%\frac{(\ell+m+n)!\,(2\ell)!\,(2m)!\,(2n)!}{(\ell+m)!\,(m+n)!\,(n+\ell)!\,\ell !\,m!\,n!\,}\;,
$$
where $m,n\in\mathbb{N}_0$ and $\ell:=\min\{m,n\}$.
\smallskip

\noindent
({\sl Hint.} Write the sum as
$\,(-1)^\ell\,{{2m}\choose{m-\ell}}{{2n}\choose{n-\ell}}\,
\pFq{3}{2}{-2\ell,-m-\ell,-n-\ell}{m-\ell+1,n-\ell+1}{1}$.
To compute this ${}_3F_2$ apply Dixon's identity to
$\,\pFq{3}{2}{-2\ell-2\epsilon,-m-\ell-\epsilon,-n-\ell-\epsilon}{m-\ell+1-\epsilon,n-\ell+1-\epsilon}{1}\,$
for small $\epsilon>0$, to the result apply Euler's reflection formula,
and then take the limit as $\epsilon\to0^+$.)
\medskip

%%\noindent
%({\sl Hint.} Write the above sum as
%$\;\frac{(-1)^\ell(2m)!\,(2n)!}{(m-\ell)!\,(m+\ell)!\,(n-\ell)!\,(n+\ell)!}\;
%\pFq{3}{2}{-2\ell,-m-\ell,-n-\ell}{m-\ell+1,n-\ell+1}{1}$.
%To compute this  ${}_3F_2$, apply Dixon's identity to
%$\,\pFq{3}{2}{-2\ell-2\epsilon,-m-\ell-\epsilon,-n-\ell-\epsilon}{m-\ell+1-\epsilon,n-\ell+1-\epsilon}{1}\,$
%for small $\epsilon>0$, to the result apply Euler's reflection formula, and then let $\epsilon\to0^+$.)
%\medskip

\item\label{Ex-cp16-5}
%Rainville pp. 53-54
Show that $y:=\pFq{2}{1}{a,b}{c}{x}$  fulfills the {\it hypergeometric differential equation}
$$
x(1-x)\,y^{\prime\prime}+[c-(a+b+1)x]\,y'-ab\,y=0\;,\quad |x|<1\;.
\smallskip
$$

\item\label{Ex-cp16-6}
Prove the following hypergeometric representations of the classical orthogonal polynomials
of Hermite, Laguerre, Jacobi, and Bessel (with standard normalization):
\begin{snugshade}
$$
\begin{array}{rcl}
H_n(x)&=&\displaystyle (2x)^n\,\pFq{2}{0}{-\frac{n}2,\frac{1-n}2}{\mbox{---}}{-\frac{1}{x^2}}\;, \\ [1.5em]
L_n^{(\alpha)}(x)&=&\displaystyle {{n+\alpha}\choose{n}}\,\pFq{1}{1}{-n}{\alpha+1}{x}\;,  \\ [1.75em]
P_n^{(\alpha,\beta)}(x)&=&
\displaystyle {{n+\alpha}\choose{n}}\,\pFq{2}{1}{-n,n+\alpha+\beta+1}{\alpha+1}{\frac{1-x}2} \\ [1.5em]
&=&\displaystyle {{n+\alpha}\choose{n}}\left(\frac{x+1}{2}\right)^n\,
\pFq{2}{1}{-n,-n-\beta}{\alpha+1}{\frac{x-1}{x+1}}\;, \\ [1.75em]
Y_n^{(\alpha)}(x)&=&\displaystyle\pFq{2}{0}{-n,n+\alpha+1}{\mbox{---}}{-\frac{x}2}\;.
\end{array}\smallskip
$$
\end{snugshade}

\noindent
These formulas hold for every $n\in\mathbb{N}_0$ and $x\in\mathbb{C}$
(with the natural definitions by continuity at the point $x=0$ in the Hermite representation and at the point $x=-1$
in the second representation for the Jacobi polynomials).
\medskip

\item\label{Ex-cp16-6a}
Prove the following hypergeometric representations of the classical discrete orthogonal
polynomials of Charlier and Meixner introduced in exercises {\bf 3} and {\bf 4} of text 4:
\begin{snugshade}
$$
\begin{array}{rcl}
C_n^{(a)}(x)&=& \displaystyle  (-a)^n\;\pFq{2}{0}{-n,-x}{\mbox{---}}{-\frac{1}{a}} \;,  \\ [1.75em]
m_n(x;\beta,c)&=& \displaystyle (x+\beta)_n\;\pFq{2}{1}{-n,-x}{-x-\beta-n+1}{\frac{1}{c}} \;.
\end{array}\smallskip
$$
\end{snugshade}

\noindent
These formulas hold for every $n\in\mathbb{N}_0$ and $x\in\mathbb{C}$.
\end{enumerate}
\medskip
}

\section*{Final remarks}

As we mentioned at the begin of this text, we followed closely chapters 1, 2, and 3 from the book \ref{AndrewsAskeyRoy1999-C16}
by Andrews, Askey, and Roy, %(o qual cont\'em informa\c c\~ao abundante sobre este t\'opico),
with some incursions on the books by Rainville \ref{Rainville1965-C16}, Bailey \ref{Bailey1935-C16}, Whittaker and Watson \ref{WhittakerWatson1927-C16}, and Lebedev \ref{Lebedev1965-C16}, as well as on the so-called {\sl Batman Manuscript Project} \ref{Erdelyi1955-C16} (coordinated by Arthur Erd\'elyi), and on the work \ref{Maroni1994-C16} by Maroni.
The proof of Theorem \ref{GammaF-exist} presented here does not assume any knowledge about infinite products, and it is based on the content appearing at the begin of Chapter XII in the book \ref{WhittakerWatson1927-C16}.
Indeed, assuming the knowledge of some basic facts about infinite products, a more concise proof may be done.
In \ref{AndrewsAskeyRoy1999-C16} only statements (i) and (iii) in Theorem \ref{pFq-conv-thm2} were proved.
Here we presented a detailed proof of (ii), giving the full details of the proof, as well as of the Euler integral representation (Theorem \ref{2F1-Euler-int-rep}).
Exemples 1, 2, and 3 in Section \ref{SumBin} appear in the article \ref{Roy1987-C16} by Ranjan Roy, as well as in the book \ref{AndrewsAskeyRoy1999-C16}.

Exercises {\bf 2}, {\bf 3}, {\bf 4}, {\bf 9}, and {\bf 10} may be found in \ref{AndrewsAskeyRoy1999-C16} (some of them presented here with some minor adjustments, reflecting our style of presentation of full details).
Exercise {\bf 5} is suggested by the need to justify a differentiation under the integral symbol in the proof of Euler's integral representation (usually omitted in the literature).
Exercises {\bf 6}, {\bf 7}, and {\bf 8} may be found in Rainville's book \ref{Rainville1965-C16}.
Exercise {\bf 11} appears in several introductory texts on hypergeometric series.
The ODE which appears in it is in the historical origins of theses series.
The results expressed by exercises {\bf 12} and {\bf 13} are very important and they may be found in several texts mentioned in the bibliography. We point out that the suggestion given for proving the hypergeometric representations in exercise {\bf 12}, based on the explicit formulas for the classical OP deduced in the previous text/chapter, allow us to give very concise proofs for all these formulas (for instance, regarding the hypergeometric representation for the Jacobi OP, compare with Theorem 6.3.3 in \ref{AndrewsAskeyRoy1999-C16}, p.\;295).
\medskip

\section*{Bibliography}
%\medskip

{\small
\begin{enumerate}[label=\emph{\rm [\arabic*]},leftmargin=*]
\item\label{AndrewsAskeyRoy1999-C16} G. Andrews, R. Askey, and R. Roy, {\sl Special Functions},
         Cambridge University Press (1999) [paperback edition: 2000].
\item\label{Bailey1935-C16} W. N. Bailey, {\sl Generalized Hypergeometric Series}, Cambridge University Press (1935).
\item\label{Erdelyi1955-C16} A. Erd\'elyi, Ed., {\sl Higher Transcendental Functions}, vols. 1--3, McGraw-Hill (1955).
\item\label{Lebedev1965-C16} N. N. Lebedev, {\sl Special Functions and their applications}, Prentice-Hall (1965).
[Translated and edited by R. A. Silverman].
\item\label{LutzMattner2001-C16} L. Mattner, {\it Complex differentiation under the integral}, Nieuw Archief voor Wiskunde IV Ser. {\bf 5}/{\bf 2}\,(1) (2001) 32--35.
\item\label{Maroni1994-C16} P. Maroni, {\sl Fonctions eul\'eriennes. Polyn\^omes orthogonaux classiques},
T\'echniques de l'Ing\'enieur, trait\'e G\'en\'eralit\'es (Sciences Fondamentales), A {\bf 154} (1994) 1--30.
\item\label{PetWilfZeilb1996-C16} M. Petkov$\check{{\rm s}}$ek, H. S. Wilf, and D. Zeilberger, {\sl $A=B$}, A\,K Peters, Wellesley (2006), 3th ed.
%, Wellesley, Massachutts (1996). %[http://www.cis.upenn.edu/~wilf/AeqB.html]
\item\label{Rainville1965-C16} E. D. Rainville, {\sl Special Functions}, The Macmillan Company, New York (1960).
\item\label{Roy1987-C16} R. Roy, {\it Binomial identities and hypergeometric identities}, The American Mathematical Monthy {\bf 94} (1987) 36--46.
\item\label{WhittakerWatson1927-C16} E. T. Whittaker and G. N. Watson, {\sl A course of modern analysis}, Cambridge University Press (1963), 4th ed. reprinted.
\end{enumerate}
}

\appendix

\chapter{Topics on locally convex spaces}\label{Appx1}

\pagestyle{myheadings}\markright{Topics on locally convex spaces}
\pagestyle{myheadings}\markleft{J. Petronilho}

In this text we review the most important facts concerning the theory of locally convex spaces (LCS)
needed along the course. %text \ref{OP-foundations}.
%a better understanding of text \ref{OP-foundations}. %All the material
Most of the material presented here is from the book \ref{ReedSimon1972-A1} by M. Reed and B. Simon (specially from Chapter V therein). The subject is also studied in deep detail in the book \ref{Treves1967-A1} by F. Tr\`eves. Other recommended sources of information (containing concise presentations) are the books by M. Al-Gwaiz \ref{Al-Gwaiz1992-A1}, B. Simon \ref{Simon2015-A1}, P. Lax \ref{Lax2002-A1}, and  W. Rudin \ref{Rudin1986-A1}.

\section{Definitions and basic properties}

We denote by $\mathbb{K}$ the scalar field of a given vector space,
being either $\mathbb{K}=\mathbb{R}$ or $\mathbb{C}$.

\begin{snugshade}
\begin{definition}\label{def-semi-norm}
A {\sl seminorm} on a vector space $X$ is a mapping $\,p:X\to[0,\infty)$ obeying the following two conditions:
\begin{enumerate}
\item[{\rm (i)}] $p(x+y)\leq p(x)+p(y)\;,\; \forall x,y\in X\;$;
\item[{\rm (ii)}] $p(\lambda x)=|\lambda|p(x)\;,\; \forall x\in X\;,\;\forall\lambda\in\mathbb{K}$.
\end{enumerate}
A family of seminorms $\{p_\alpha\}_{\alpha\in A}$ is said to {\sl separate points} if
\begin{enumerate}
\item[{\rm (iii)}] $p_\alpha(x)=0\,,\;\forall \alpha\in A\;\;\Rightarrow\;\;x=0$.
\end{enumerate}
\end{definition}
\end{snugshade}

\begin{snugshade}
\begin{definition}\label{def-LCS}
A {\sl locally convex space (LCS)} is a vector space $X$ with a family
$\{p_\alpha\}_{\alpha\in A}$ of seminorms separating points.
The {\sl natural topology} on a LCS is the weakest topology in which
all the seminorms $p_\alpha$ are continuous and in which the operation of addition is continuous.
(Often we will refer to it as the ``$\,\{p_\alpha\}_{\alpha\in A}-$natural topology''.)
%(Often we will refer to it as the ``$p-$natural topology''.)
\end{definition}
\end{snugshade}

\begin{snugshade}
\begin{proposition}\label{LCS-Hausdorff}
The natural topology of a LCS is Hausdorff.
\end{proposition}
\end{snugshade}

%{\it Proof.}\qed
%\medskip

A neighborhood base at $0$ for the $\,\{p_\alpha\}_{\alpha\in A}-$natural topology in a LCS $X$ is given
by the totality of the sets of the form
\begin{snugshade}
$$V\big(0;\epsilon,\{p_{\alpha_1},\ldots,p_{\alpha_N}\}\big):=
\big\{x\in X\;:\; p_{\alpha_i}(x)<\epsilon\;,\;i=1,\ldots,N\big\}\;,$$
\end{snugshade}\noindent
where $\epsilon>0$ and $\{\alpha_1,\ldots,\alpha_n\}$ is a (finite) subset of $A$.
As a consequence, given a sequence $\{x_n\}_{n\geq0}$ in $X$, and being $x\in X$, we deduce

\begin{snugshade}\vspace*{-0.5em}
\begin{equation}\label{xn-conv-LCM}
x_n\to x\;\mbox{\rm in $X$}\quad\Leftrightarrow\quad p_{\alpha}(x_n-x)\to0\;,\;\forall\alpha\in A\;.
\end{equation}
\end{snugshade}

\begin{snugshade}
\begin{definition}\label{def-LCS-equiv-semin}
Two families of seminorms $\{p_\alpha\}_{\alpha\in A}$ and
$\{q_\beta\}_{\beta\in B}$ in a LCS $X$ are called {\sl equivalent}
if they generate the same natural topology in $X$.
\end{definition}
\end{snugshade}

\begin{snugshade}
\begin{proposition}\label{teo-LCS-equiv-semin}
Let $\{p_\alpha\}_{\alpha\in A}$ and $\{q_\beta\}_{\beta\in B}$ be two families of seminorms in a LCS $\,X$.
The following statements are equivalent:
\begin{enumerate}
\item[{\rm (i)}] $\{p_\alpha\}_{\alpha\in A}$ and $\{q_\beta\}_{\beta\in B}$ are equivalent families of seminorms;
\item[{\rm (ii)}] each $p_\alpha$ is continuous in the $\{q_\beta\}_{\beta\in B}-$natural topology, and each $q_\beta$ is continuous in the $\{p_\alpha\}_{\alpha\in A}-$natural topology;
%\item[{\rm (ii)}] each $p_\alpha$ is continuous in the $q-$natural topology, and each $q_\beta$ is continuous in the $p-$natural topology;
\item[{\rm (iii)}] for each $\alpha\in A$, there are $\beta_1,\ldots,\beta_n\in B$ and $C>0$ so that %for all $x\in X$,
$$
p_\alpha(x)\leq C\,\sum_{i=1}^n q_{\beta_i}(x)\; ,\quad\forall x\in X\;;
$$
and for each $\beta\in B$, there are $\alpha_1,\ldots,\alpha_n\in A$ and $K>0$ so that %for all $x\in X$,
$$
q_\beta(x)\leq K\,\sum_{i=1}^n p_{\alpha_i}(x)\; ,\quad\forall x\in X\; .
$$
\end{enumerate}
\end{proposition}
\end{snugshade}

%{\it Proof.}\qed
%\medskip

%\begin{example}
%\rm
%Let $X$ be a vector space and suppose that $\{\ell\,:\,\ell\in\Lambda\}$, with $\Lambda\subseteq X^*$,
%is a nonempty set of linear functionals on $X$ separating points
%(i.e., $\ell\in\Lambda\wedge\ell(x)=0$ for all $x\in X$, implies $x=0$).
%Recall that the $\sigma(X,Y)-$topology is the weaker topology in $X$
%such that all the linear functionals $\ell\in\Lambda$ are continuous.
%With such a topology, $X$ becomes a LCS, where the topology is generated by the family of seminorms
%$\{p_\ell\}_{\ell\in\Lambda}$, being $p_\ell(x):=|\ell(x)|$, $x\in X$.
%Notice, however, that while this topology is given by seminorms, it is never given by norms if $Y$ has
%infinite (algebraic) dimension.
%\end{example}

\section{Fr\'echet spaces}

\begin{snugshade}
\begin{theorem}\label{LCS-metrizable}
Let $X$ be a LCS. The following are equivalent:
\begin{enumerate}
\item[{\rm (i)}] $X$ is metrizable (i.e., the topology in $X$ may be defined by a metric);
\item[{\rm (ii)}] $0$ has a countable neighborhood base;
\item[{\rm (iii)}] the topology in $X$ is generated by some countable family of seminorms.
\end{enumerate}
\end{theorem}
\end{snugshade}

%{\it Proof.}\qed
\medskip

\begin{remark}\label{LCS-metrizable-metric}
\rm
%Under the conditions of Theorem \ref{LCS-metrizable},
If $\{p_k\}_{k\in\mathbb{N}}$ is a countable family of seminorms generating the topology in a LCS $X$,
then the application ${\rm d}:X\times X\to[0,\infty)$ defined by
\begin{snugshade}
$$
{\rm d}(x,y):=\sum_{k=0}^\infty\frac{1}{2^k}\frac{p_k(x-y)}{1+p_k(x-y)}\,,\quad x,y\in X\;,
$$
\end{snugshade}
\noindent
is a metric in $X$ and it generates the same topology in $X$ as the family $\{p_k\}_{k\in\mathbb{N}}$.
\end{remark}

\begin{snugshade}
\begin{definition}\label{def-Frechet-space}
A complete metrizable LCS is called a {\sl Fr\'echet space}.
\end{definition}
\end{snugshade}

\begin{remark}
\rm
Recall that a complete metric space is a metric space in which every Cauchy sequence is convergent
(for some element in that space). %From (\ref{xn-conv-LCM}) we see that
In a metrizable LCS, whose topology is generated by the countable family
of seminorms $\{p_\alpha\}_{\alpha\in\mathbb{N}}$,
a sequence $\{x_n\}_{n\in\mathbb{N}}$ is Cauchy if and only if
\begin{snugshade}\vspace*{-0.5em}
$$%\begin{equation}
\forall \epsilon>0\,,\, \forall \alpha\in\mathbb{N}\, ,\; \exists n_0\in\mathbb{N}\;:\; \forall n,m\in\mathbb{N}\,,\quad
n,m\geq n_0\;\Rightarrow\; p_\alpha(x_n-x_m)<\epsilon\,.
$$%\end{equation}
\end{snugshade}
\end{remark}

\section{The inductive limit topology}

Here we introduce the inductive limit topology in a particularly simple case,
which, however, will be sufficient for our purposes.
Up to some minor modifications (mostly concerning notation), essentially, we pursue following Reed and Simon.

\begin{snugshade}
\begin{definition}\label{def-ind-lim-topology}
Let $X$ be a vector space and $\{X_n\}_{n\in\mathbb{N}}$ a family of subspaces of $X$ such that
$$
X_n\subseteq X_{n+1}\;,\quad X=\bigcup_{n\in\mathbb{N}} X_n\;.
$$
Suppose that each $X_n$ is a LCS and let $i_n:X_n\to X$ $(x\mapsto x)$ be the natural injection from $X_n$ into $X$.
\begin{enumerate}
\item[{\rm (i)}] The {\sl inductive limit topology} (in $X$) of the spaces $X_n$ is the strongest topology in $X$ such that $X$ is a LCS and all the maps $i_n$ are continuous; we write $$X=\mbox{\rm ind\,lim}_n\,X_n\;;$$
\item[{\rm (ii)}] if each $X_{n+1}$ induces in $X_n$ the given topology in $X_n$ (i.e., $X_n$ is a topological subspace of $X_{n+1}$ with the relative topology), the above topology is called the {\sl strict inductive limit topology} of the spaces $X_n$;
\item[{\rm (iii)}] if---in addition to the conditions in {\rm (ii)}---each $X_n$ is a proper closed subspace of $X_{n+1}$, the above topology is called the {\sl hyper strict inductive limit topology} of the spaces $X_n$.
%\item[{\rm (iii)}] if, in addition to the conditions in {\rm (ii)}, each $X_n$ is a proper closed subspace in $X_{n+1}$, the above %topology is called the {\sl hyper strict inductive limit topology} of the spaces $X_n$.
\end{enumerate}
\end{definition}
\end{snugshade}

%\begin{remark}
%When $X$ is a LCS endowed with the inductive limit topology of the LCS spaces $X_n$, we often write
%\vspace*{-0.5em}
%\begin{snugshade}
%$$
%X=\mbox{\rm ind\,lim}_n\,X_n\; .
%$$
%\end{snugshade}
%\end{remark}

\begin{snugshade}
\begin{theorem}\label{TeoLimInd}
Let $X$ be a LCS endowed with the strict inductive limit topology of the LCS $X_n$.
Then the following holds:
\begin{enumerate}
\item[{\rm (i)}] the restriction of the (strict inductive limit) topology on $X$ to each $X_n$ is the given topology on $X_n$;
\item[{\rm (ii)}] the collection of all convex sets $U\subseteq X$ such that $U\cap X_n$ %$[=i_n^{-1}(U)]$
is open in $X_n$ for each $n$ is a neighborhood base at $0$ in $X$;
%\item[{\rm (iii)}] if each $X_n$ is complete, then so is $X$. %(no sentido das sucessões generalizadas, i.e., nets)
%\item[{\rm (iv)}] if $Y$ is a LCS, a linear mapping $T:X\to Y$ is continuous if and only if each of the restrictions $T_n:=T|X_n:X_n\to Y$ is continuous.
\end{enumerate}
\end{theorem}
\end{snugshade}

%{\it Proof.}\qed
%\medskip

\begin{snugshade}
\begin{theorem}\label{Teo-Tcontin-LimInd}
Let $X$ be a LCS with the strict inductive limit topology of the LCS $X_n$, and let
$Y$ be any LCS. Then, a linear mapping $T:X\to Y$ is continuous if and only if
each of the restrictions $T_n:=T|X_n:X_n\to Y$ is continuous.
\end{theorem}
\end{snugshade}

%{\it Proof.}\qed
%\medskip

\begin{snugshade}
\begin{theorem}\label{TeoContinutyLimInd}
Let $X$ be a LCS with the hyper strict inductive limit topology of the LCS $X_n$.
Then the following holds:
\begin{enumerate}
\item[{\rm (i)}] if $\{x_n\}_{n\geq1}$ is a sequence in $X$ and $x\in X$, then
\begin{equation}\nonumber%\label{xn-conv-LimInd}
x_n\to x\;\mbox{\rm in $X$}\quad\Leftrightarrow\quad
\exists k\in\mathbb{N}\;:\; x_n\in X_k\;\mbox{for all $n$}\;\;\wedge\;\;
x_n\to x\;\mbox{\rm in}\; X_k\;;
\end{equation}
\item[{\rm (ii)}] if all the spaces $X_n$ are sequentially complete, then so is $X$;
\item[{\rm (iii)}] $X$ is not a metrizable space.
\end{enumerate}
\end{theorem}
\end{snugshade}

%{\it Proof.}\qed
%\medskip

\begin{remark}
\rm
Those who learned already about LCS (priori to this course) may be a little surprised because here we didn't made any reference to concepts such as ``absorbing'' set, ``balanced'' set, or ``gauge'' (among others). These are indeed very useful tools in the study of LCS---specially for a presentation of the theory of LCS including proofs of all the results stated---, but none for our presentation.%, so we decided to omit any reference to them.
\end{remark}

%\begin{remark}
%\rm
%Those who learned already about LCS (before taking this course) may be a little surprised because of the fact that we didn't spoken here in concepts such as ``absorbing set'' and ``balanced set''. These are indeed concepts of major importance in the framework of the study of LCS. %The reason why such concepts could indeed have been omitted from our presentation
%Such concepts could indeed been omitted from our presentation since we didn't presented the proofs of the results stated (for the proofs they are fundamental).
%\end{remark}

\section{The weak dual topology} %in the dual space of a TVS}

Let $X$ be a vector space (over the field $\mathbb{K}=\mathbb{R}$ or $\mathbb{C}$).
The {\sl algebraic dual} of $X$, denoted by $X^*$,
is the set of all linear functionals ${\bf f}:X\to\mathbb{K}$.
Usually the action of a functional ${\bf f}\in X^*$ over a vector $x\in X$
(i.e., the scalar ${\bf f}(x)$) will be denoted by
\vspace*{-0.5em}\begin{snugshade}$$\langle {\bf f},x\rangle\;.$$\end{snugshade}\noindent
If, besides being a vector space, $X$ is endowed with a compatible topology
(i.e., addition and scalar multiplication are continuous mappings),
$X$ is called a {\sl topological vector space} (TVS).
We denote by $\mathcal{L}(X,Y)$ the set of all linear and continuous operators between two TVS $X$ and $Y$.
%The space of all continuous linear operators between two LCS $X$ and $Y$ will be denote by $\mathcal{L}(X,Y)$.
In particular, the {\sl topological dual} of a TVS $X$ is the set
%$X^\prime:=\mathcal{L}(X,\mathbb{K})$ is called the {\sl topological dual} of the TVS $X$.
\vspace*{-0.5em}\begin{snugshade}
$$X^\prime:=\mathcal{L}(X,\mathbb{K})=\big\{{\bf f}\in X^*\,:\,{\bf f}\;\mbox{\rm is continuous}\big\}\;.$$\end{snugshade}\noindent
%Therefore, $X^\prime$ is the set of all continuous linear functionals ${\bf f}:X\to\mathbb{K}$.
Clearly, $X^\prime\subseteq X^*$.
It is worth mentioning that this inclusion is actually an equality if $X$ is a finite dimensional normed space,
while it is a strict inclusion whenever $X$ is an infinite dimensional normed space
(a fact that can be proved using Zorn's Lemma).
We emphasize, however, that there
are infinite dimensional TVS, $X$, such that the set equality $X'=X^*$ holds.
%An important example, for our purposes, is the space $\mathcal{P}$ of all polynomials
%endowed with an appropriate strict inductive limit topology.

\begin{snugshade}
\begin{definition}\label{wk-dual-topology}
Let $X$ be a TVS. The {\sl weak dual topology} in $X^\prime$ is the
topology in $X^\prime$ generated by the family of seminorms
$\mathcal{S}:=\{s_x|x\in X\}$, where each seminorm $s_x:X^\prime\to[0,+\infty)$ is defined
for each ${\bf f}\in X'$ by
\vspace*{-0.5em}\begin{snugshade}$$s_x({\bf f}):=|\langle {\bf f},x\rangle|\;.$$\end{snugshade}\noindent
\end{definition}
\end{snugshade}

Endowed with the weak dual topology, $X'$ becomes a LCS
(and so it is an Hausdorff space).
Henceforth, according with (\ref{xn-conv-LCM}),
given a sequence $\{{\bf f}_n\}_{n\geq1}$ in $X'$, we have
\vspace*{-0.5em}
\begin{snugshade}
$${\bf f}_n\to\textbf{0}\;\mbox{\rm in $X'$}\quad\mbox{\rm iff}\quad
\langle {\bf f}_n,x\rangle\to0\;,\;\;\forall x\in X\;.$$
\end{snugshade}\noindent
Because of this property, often the name {\sl point convergence topology} is given to the
weak dual topology in $X'$. Another one which we may find in the literature is
{\sl topology of convergence on the finite subsets of $X$}.
This name is due to the fact that the collection of the sets of the form
\vspace*{-0.5em}\begin{snugshade}
$$V_{X'}(\textbf{0};\epsilon,F):=\{{\bf f}\in X'\,:\, s_x({\bf f})<\epsilon\;,\;\forall x\in F\,\big\}\,,$$
\end{snugshade}\noindent
where $\epsilon>0$ and $F$ is a finite subset of $X$, is a neighborhood base at $\textbf{0}\in X'$.
\medskip

\begin{snugshade}
\begin{definition}\label{def-Tdual}
Let $X$ and $Y$ be TVS, and $T\in\mathcal{L}(X,Y)$. %$T:X\to Y$ a linear and continuous operator.
The {\sl dual operator} (or {\sl dual mapping})
of $T$ is the (linear) mapping $$T':Y'\to X'\quad
({\bf g}\in Y'\mapsto T'{\bf g}\in X')\;,$$ where $T'{\bf g}:X\to\mathbb{K}$
is defined by
\vspace*{-0.5em}\begin{snugshade}
$$\langle T'{\bf g},x\rangle:=\langle {\bf g},Tx\rangle\;, \quad x\in X\;.$$
\end{snugshade}\noindent
\end{definition}
\end{snugshade}

%Note that $T'$ is a linear operator, and $T'{\bf g}\in X'$ for each ${\bf g}\in Y'$.

\begin{snugshade}
\begin{theorem}\label{Tdual-contin}
Let $X$ and $Y$ be TVS, and $T\in\mathcal{L}(X,Y)$. % $T:X\to Y$ a linear and continuous operator.
Let $X'$ and $Y'$ be endowed with the weak dual topologies.
Then $T'\in\mathcal{L}(Y',X')$. %the dual operator $T':Y'\to X'$ is also a linear and continuous operator.
\end{theorem}
\end{snugshade}

%{\it Proof.}
%For each $x\in X$ and ${\bf g}\in Y'$, we have
%$$
%s_x(T'{\bf g})=|\langle T'{\bf g},x\rangle|=|\langle {\bf g},Tx\rangle|=s_{Tx}({\bf g})\;.
%$$
%Therefore, for each $x\in X$ and $\epsilon>0$, we may write
%$$
%\begin{array}{rl}
%{\bf T}'^{-1}\big(V_{X'}(\textbf{0};\epsilon,\{s_x\})\big)
%&=\big\{ {\bf g}\in Y'\,:\,T'{\bf g}\in V(\textbf{0};\epsilon,\{s_x\}) \big\}
%=\big\{ {\bf g}\in Y'\,:\,s_{x}(T'{\bf g})<\epsilon \big\} \\ [0.25em]
%&=\big\{ {\bf g}\in Y'\,:\,s_{Tx}({\bf g})<\epsilon \big\}
%=V_{Y'}\big(\textbf{0};\epsilon,\{s_{Tx}\}\big) \;,
%\end{array}
%$$
%hence the desired conclusion follows.
%\qed
\medskip

\section*{Bibliography}
\medskip

{\small
\begin{enumerate}[label=\emph{\rm [\arabic*]},leftmargin=*]
\item\label{Al-Gwaiz1992-A1} M. A. Al-Gwaiz, {\sl Theory of Distributions}, Marcel Dekker, Inc. (1992).
\item\label{Lax2002-A1} Peter D. Lax, {\sl Functional Analysis}, John Wiley $\&$ Sons (2002).
\item\label{ReedSimon1972-A1} M. Reed and B. Simon, {\sl Methods of Modern Mathematical Physics I: Functional Analysis},
Academic Press (1972).
\item\label{Simon2015-A1} B. Simon, {\sl Real Analysis: A Comprehensive Course in Analysis}, Part 1, AMS (2015).
\item\label{Rudin1986-A1} W. Rudin, {\sl Real and Complex Analysis}, McGraw-Hill (1986) [paperback edition: 2004].
\item\label{Treves1967-A1} F. Tr\`eves, {\sl Topological Vector Spaces, Distributions and Kernels}, Academic Press (1967).
\end{enumerate}
}

\end{document}